\documentclass[11pt]{article} 
\usepackage{geometry} 
\usepackage{amsmath,amsthm,amssymb}                
\usepackage{graphicx, color}
\usepackage{float}
\usepackage{amssymb}
\usepackage{epstopdf}
\usepackage{pdfsync}
\usepackage{makeidx}

\newtheorem{definition}{Definition}[section]
\newtheorem{lemma}[definition]{Lemma}

\newtheorem{proposition}[definition]{Proposition}

\newtheorem{remark}[definition]{Remark}

\def\wto{\rightharpoonup}
\def\e{\varepsilon}

\newcommand{\restr}{\ \rule{.4pt}{7pt}\rule{6pt}{.4pt}\ }  
\newcommand{\weak}{\rightharpoonup}

\title{A Game of Scales \\
\small Concentration and oscillations are coming}
\author{Andrea Braides}
\date{Version 1.0}
\makeindex
\begin{document}
\maketitle
%

\centerline{\bf Preface}
\bigskip

\noindent
These notes cover the material of two courses given at SISSA, Trieste, in the academic years 2023-2024 and 2025-2026, together with some material from courses taught in other recent years. The original idea behind these courses emerged from conversations with Andrea Causin and Margherita Solci about the somewhat elusive notion of `multiscale behaviour', which arises in the study of variational problems involving one or more small parameters. Two main mechanisms govern their asymptotic behaviour: concentration and oscillations. Their description is often associated with the emergence of a characteristic `scale', quantifying either the speed of concentration or the period of an optimal oscillation. These mechanisms may originate from different variational principles and may occur separately, coexist, compete, or interact, thus providing a natural interpretation of the term `multiscale'. This perspective inspired the informal title of the courses, which also became the title of these notes.

The natural framework for these analyses is $\Gamma$-convergence. However, these notes are not intended as a general introduction to this now well-established theory. We provide only a brief account of the techniques relevant to our purposes, referring the reader to the author's book chapter {\em A Handbook of $\Gamma$-convergence}. For a broader introduction to $\Gamma$-convergence, and in particular to its topological properties, the standard reference is Dal Maso's monograph {\em An Introduction to $\Gamma$-convergence}.

The courses were highly interactive, and I am grateful to all the students who contributed through their questions, comments, and suggestions: Irene Anello, Fabrizio Cara\-giulo, Simone Carano, Simone Vincini, Giuseppe Cosma Brusca, Davide Donati, Sergio Scalabrino, Chiara Trifone, Edoardo Voglino, Matteo Zanardini, Roman Oleinik, Marco Picerni, Abdilaziz Z.~U.~Komilov, Stefano Mannella, Alec J.~A.~Schiavoni Piazza. Some of them also developed parts of the material presented in these notes, building on research projects proposed during the courses.

\bigskip

Trieste\hfill Andrea Braides

\vfill\eject
\tableofcontents

\section{Introduction} 
In his {\em laudatio} in honour of the newly awarded Fields medalist Artur Avila at the International Congress of Mathematicians in Seoul in 2014, \'Etienne Ghys gave an account of the development of dynamical systems that can equally be applied to the Calculus of Variations, not surprisingly quoting  some of the same reference names. For the mathematicians of the eighteenth century, of the likes of Newton, solving a minimum problem of the form
$$
\min\Big\{\int_a^b f(t, u(t), u'(t))\,dt: u\in X(a,b)\Big\}
$$
would mean explicitly exhibiting a smooth minimizer $u$, or giving a construction that leads to its precise determination. In the nineteenth century the focus shifted to the actual possibility of obtaining a solution, and here we can quote Poincar\'e as Ghys did, together of course with Hilbert and his famous XIX problem on the Calculus of Variations. The twentieth century saw the shift towards the use of weak spaces, with the standpoint of Laurent Schwartz's distributions and even more geometrical objects of Geometric Measure Theory, for which the problem above is just a distant prototype, and the solutions are intended in a more general sense. At that point the focus had turned to determining conditions on the functionals which ensure the existence of solutions, and, more generally, following the so-called Direct Method, the coerciveness and lower semicontinuity of these functionals.
In problems such as the one we have at hand, this translates into questions on the integrand $f$.
In remarkable analogy with Ghys's description, motivated for him by the presence of small perturbation parameters in physical theories, a further step in the development of the Calculus of Variations has been the study of the asymptotic behaviour of families of minimum problems. 
In the language of De Giorgi's Gamma-convergence this is essentially equivalent to the study of the $\Gamma$-limit of families of functionals parameterized by a (small) parameter $\e$. 

In these notes, we interpret $\Gamma$-limits in the light of the properties of $\e$-parameterized families of solutions (or almost solutions) of the corresponding minimum problems. We highlight the possibility of {\em oscillations} and {\em concentration} of the solutions. In the first case, families are increasingly oscillating, and the $\Gamma$-limit describes the averaged effect of such oscillations. In the second case some quantities, mainly gradients or Jacobians, tend to be singular on some lower-dimensional set, and the $\Gamma$-limit describes the nature and energetic cost of the singularity. Our viewpoint is that in these processes the effect of the small parameter is characterized by some {\em scale} that describes, for example, the typical period of oscillation or the size of the region where we have high gradients. The question of how these scales depend on $\e$, and the relation with the general properties of the functionals, will be the main focus of our analysis. 

\smallskip
{\bf A scale for concentration.} To understand what we may mean by the `scale' of a problem, as a first example of $\e$-parameterized minimum problems, we consider one in which $\e$ appears as a {\em singular-perturbation} factor. We consider a so-called {\em double-well potential} $W\colon\mathbb R\to [0,+\infty)$; that is, a sufficiently smooth positive function, with the value $0$ only at two points, which we can assume to be $0$ and $1$. The related perturbed {\em optimal-profile problem} is 
\begin{equation}\label{model-conc}
m_\e=\min\Big\{\int_{\mathbb R}( W(u)+\e^2|u'|^2)\,dt :\lim_{t\to-\infty}u(t)=0,\ \lim_{t\to+\infty}u(t)=1\Big\}.
\end{equation}
Note that the addition of a derivative forbids trivial solutions with $u(t)\in\{0,1\}$ almost everywhere, and this effect, in a sense, weakens as $\e\to 0$. Since the test functions $u$ are continuous, by the conditions at infinity, we deduce that there is at least one point where they are equal to $1/2$. Since the problem is invariant under translations, we can assume that $u(0)=1/2$. We suppose that a minimizer $u_\e$ exists (which can be proved in a suitable space) and is smooth enough so that it satisfies the Euler--Lagrange equation 
$W'(u_\e)-2\e^2u_\e''=0$. Multiplying by $u'_\e$ and integrating, we find that $W(u_\e)-\e^2|u''|^2$ is a constant,
which turns out to be $0$ by the conditions at infinity. We then have
$$
\begin{cases} \e u'_\e=\sqrt{W(u_\e)}\\
u_\e(0)=\frac12,
\end{cases}
$$
which also guarantees uniqueness of the solution.

The change of variables $u(t)= v(t/\e)$ shows that the value $m_\e$ is independent of $\e$.
If $w$ is the minimizer with $w(0)=1/2$ of the problem with $\e=1$; that is, of
$$
m=\min\Big\{\int_{\mathbb R}( W(v)+|v'|^2)\,dt :\lim_{t\to-\infty}v(t)=0,\ \lim_{t\to+\infty}v(t)=1\Big\},
$$
then the function $u_\e(t)= w(t/\e)$ is the corresponding minimizer of the first problem.  
With this remark in mind, we see that the equality $W(w)+|w'|^2=2\sqrt{W(w)}\,w'$ holds, 
from which we can actually compute
$$
m=2\int_{\mathbb R}\sqrt{W(w)}w'dt=2\int_{0}^1\sqrt{W(s)}ds,
$$
after the substitution $s=w(t)$.

The key feature of the $u_\e$ is that they {\em concentrate} at $0$ at {\em scale $\e$}. This statement is to be understood in the sense that, for a given $T>0$, the integral of the energy density  $W(u_\e)+\e^2|u_\e'|^2$ outside $(-\e T,\e T)$ is infinitesimal, while its integral on $(-\e T,\e T)$ tends to $m$, as $T\to+\infty$.

\medskip
{\bf A scale for oscillations.}
The way a `scale' is understood is different in problems with oscillations. As a prototype, we consider the minimum problems
\begin{equation}\label{model-osc}
m_\e=\min\Big\{\int_0^1a\big(\tfrac{t}\e\big)|u'|^2dt: u(0)=0, u(1)=1\Big\}.
\end{equation}
Here $a$ is a periodic function with period $1$, with $0<\alpha\le a(t)\le \beta$. The parameter $\e$ is the period of the scaled coefficient $a(t/\e)$; hence, in this case, it has a geometric meaning. In order to simplify the computations, we suppose that there is an integer number of periods in the interval $(0,1)$; that is, $n=1/\e$ is an integer. Supposing $a$ differentiable, the Euler--Lagrange equation for the minimizer $u_\e$ is 
$$
\big(a(\tfrac{t}\e) u_\e'\big)'=0\hbox{; that is, } \ u_\e'= \frac{C}{a(\frac{t}\e)}.
$$
To determine the constant $C$, using the boundary conditions, we have
$$
1=\int_0^1u'_\e dt=C\int_0^1\frac{1}{a(\frac{t}\e)}dt= \frac{C}n\int_0^n\frac{1}{a(s)}ds\hbox{;\  that is, } \ C=\Big(\int_0^1\frac{1}{a(s)}ds\Big)^{-1}=:\underline a,
$$
the {\em harmonic mean} of $a$,
so that we have 
$$
u_\e(t) = t+\e\,\varphi\big(\tfrac{t}\e\big)\hbox{; where } \varphi(t)=\underline a\int_0^t\Big(\frac{1}{a(s)}-\frac{1}{\underline a}\Big)ds
$$
defines a $1$-periodic function. Moreover, we can compute the minimum value $m_\e=\underline a$, independent of $\e$. This allows us to consider the case $\e=1$ for the characterization of $\varphi$. Rewriting $u(t)=t+\psi(t)$ gives that $\varphi$ is the minimizer of
$$
\underline a
=\min\Big\{\int_0^1a(t)|1+\psi'(t)|^2dt: \psi(0)=\psi(1)=0\Big\};
$$
that is, the minimizer among periodic perturbations of the identity. 

The form of $u_\e$ highlights {\em oscillations}, in the form of a periodic variation around the limit function $u(t)=t$. The scale $\e$ is then interpreted as giving both the period and the scaling of the oscillating term $\e\varphi(\frac{t}\e)$. 

\smallskip
These two examples indicate the two aspects that will be examined in these notes. One is the asymptotic analysis of the functionals that guarantees the convergence of solutions of minimum problems. This requires the use of sufficiently weak notions of convergence of functions that allow for concentration and oscillation effects. The second is the asymptotic description of minimizers. In the terminology of $\Gamma$-convergence, this is described by the computation of a $\Gamma$-limit with respect to these topologies. The object of these lectures will be to improve this description by studying the scales involved in optimal sequences, often highlighting concentration or oscillations of a different nature.

\bigskip
\section{Calculus of Variations in a nutshell}
In this section we introduce very briefly some tools from the Calculus of Variations which will be used in the sequel. The focus will be on the Direct Method, which consists in examining the behaviour of minimizing sequences using  topological arguments. We will not use Euler--Lagrange equations, obtained by examining suitable variations of a minimizer, since in general we will try not to rely on regularity assumptions on the minimizers.
	
\subsection{Calculus of Variations without variations}\label{semicont-sec}  
If a minimum exists for a problem involving a single functional $F$, we expect the nature of a minimizer 
to be revealed by looking at minimizing sequences\index{minimizing sequences}; that is, sequences $\{x_n\}$ such that the values $F(x_n)$ approach the infimum of $F$ on its domain. 
By choosing a notion of convergence on the domain $X$ of the functional, this idea translates into looking at the cluster points of minimizing sequences. Already at this point, the choice of a topology reflects two competing demands: on one side, a weaker topology favours the existence of converging sequences; hence, the existence of these cluster points. On the other hand, a stronger one will provide more information on them. 
Once a minimizing sequence $\{x_n\}$ converging to a point $\overline x\in X$ is found, 
the minimality of $\overline x$ is equivalent to the inequality 
$$F(\overline x)\leq \lim_{n\to+\infty}F(x_n),$$ 
or, equivalently, to the requirement that, for any sequence $\{y_n\}\subset X$ converging to $\overline x$, 
$F$ satisfies 
$$F(\overline x)\leq \liminf\limits_{n\to+\infty}F(y_n),$$
which is the definition of sequential lower semicontinuity\index{lower semicontinuity} of $F$ at $\overline x$.  
This process, ensuring the existence of a minimizer for $F$ in $X$, is a basic version of the {\bf Weierstrass Theorem}\index{Weierstrass Theorem}, which is the core of the so-called {\bf Direct Method}.  

In order to use general functional-analytic methods, it can be useful to view the two conditions 

\smallskip
$\bullet$ existence of a minimizing sequence converging to $\overline x\in X$;  

\smallskip$\bullet$  sequential lower semicontinuity of $F$ at $\overline x$;

\smallskip\noindent
as consequences of properties of the functional $F$ independent of the point $\overline x$.  
The existence of a converging minimizing sequence is ensured by the more general requirement

\smallskip
$\bullet$ (sequential {\bf coerciveness}\index{coerciveness} of $F$) there exists $C\in\mathbb R$ such that 
every sequence $\{y_n\}\subset X$ 

satisfying $F(y_n)\leq C$ admits a converging subsequence. 

\smallskip\noindent
The second requirement is a particular case of the more general property that 

\smallskip
$\bullet$ (sequential {\bf lower semicontinuity} of $F$) for all $x\in X$ we have 
$$F(x)\leq \liminf_{n\to+\infty}F(x_n)$$

\nobreak
for all sequences $\{x_n\}\subset X$ converging to $x$. 
\goodbreak

\smallskip
In these notes, we will always consider topologies for which 
in the relevant sub-level sets $\{x\in X: F(x)\leq C\}$ sequential properties are equivalent to corresponding general topological ones. 
In particular, we will use topologies that are metrizable on such sets. 

\paragraph{Functionals on $L^p$ spaces.}  
As a first example, we examine coerciveness and lower-semicontinuity properties of functionals of the form 
\begin{equation}\label{intLp}
F(u)=\int_{\Omega}f(u)\, dx,
\end{equation}
defined for $u\in L^p(\Omega;\mathbb R^m)$, where $\Omega$ is a bounded open subset of $\mathbb R^d$ and $f\colon\mathbb R^m\to \mathbb R\cup\{+\infty\}$ is a Borel function.  

As for {\bf coerciveness}, if $p>1$, this is usually derived from the {\bf growth condition from below} given by 
$$f(z)\geq c_1|z|^p-c_2,$$
with $c_1>0$ and $c_2\in\mathbb R$. This 
ensures that for any $C\in\mathbb R$ the elements of the set $\{u\in L^p(\Omega): F(u)\leq C\}$ have equibounded $L^p$-norm. Hence, if $F$ is not trivially equal to $+\infty$, then this condition ensures coerciveness with respect to the weak $L^p$ topology. Indeed, 
these sets are weakly pre-compact
in $L^p(\Omega)$, since the weak convergence in bounded sets of $L^p$ is metrizable, and they are not empty for some $C$.    

The {\bf weak lower semicontinuity} (that is, lower semicontinuity with respect to the weak convergence) of these functionals is a classical issue in the Calculus of Variations and is linked to {\bf convexity} properties of $f$. For functionals as in \eqref{intLp}, the weak lower semicontinuity in $L^p(\Omega;\mathbb R^m)$
with $p\ge 1$ is equivalent to the convexity and lower semicontinuity of the integrand $f$. 
The necessity of the lower semicontinuity of $f$ is immediately obtained by looking at sequences of constant functions. The necessity of the convexity of $f$ is derived from the analysis of sequences of oscillating functions, using the so-called 
{\bf Riemann-Lebesgue Lemma}\index{Riemann-Lebesgue Lemma}, describing the weak limit of scaled periodic functions. More precisely,
let $\varphi\in L^p_{\rm loc}(\mathbb R^d;\mathbb R^m)$
be a $T$-periodic function and
$$\overline\varphi=\frac{1}{T^d}\int_{[0,T]^d}\varphi(y)\, dy$$
be its average.
Then, the sequence $u_n(x)=\varphi(nx)$ weakly converges to the constant function $u=\overline\varphi$
(see Lemma~\ref{RL-lemma} in the Appendix to this Chapter for a proof).
We apply this lemma to the $1$-periodic function $\varphi$ defined as follows. We fix $t\in(0,1)$ and choose $z,z^\prime\in\mathbb R^m$. We then define 
\begin{equation}\label{oscifi}
\varphi(x)=\varphi(x_1,\dots, x_d)=\begin{cases}
z & \hbox{\rm if }\ 0\leq x_1<t\\
z^\prime & \hbox{\rm if }\ t\leq x_1<1,  
\end{cases}
\end{equation}
so that $\overline\varphi=tz+(1-t)z^\prime$ and $u_n(x)=\varphi(nx)\weak u(x)=\overline\varphi$ in $L^p(\Omega;\mathbb R^m)$.  
We can apply the Riemann-Lebesgue Lemma also 
to the $1$-periodic function $\psi(x)=f(\varphi(x))$ whose average is 
$\overline\psi=tf(z)+(1-t)f(z^\prime)$. 
From the weak convergence of $\psi(nx)$ to $\overline\psi$, 
in particular we obtain that 
$$\lim\limits_{n\to+\infty} \int_\Omega \psi(nx)\, dx=|\Omega|\overline\psi.$$
Hence, 
from the weak$^*$ lower semicontinuity in $L^\infty$ (which is the case with the strongest convergence)  of the functional $F$, 
we get 
$$|\Omega| f(tz+(1-t)z^\prime)=F(u)\leq \liminf_{n\to+\infty}F(u_n)=|\Omega|  (tf(z)+(1-t)f(z^\prime))$$ 
and the convexity of $f$.  

\smallskip
The sufficiency of the convexity and lower semicontinuity of $f$ for the lower semicontinuity of $F$ with respect to the weak $L^1$-topology (which is in this case the weakest convergence) can be obtained first by considering smooth $f$, and then using an interesting approximation procedure.

We first suppose $f$ to be $C^1$ and Lipschitz. From the convexity inequality 
$f(w)\ge f(z)+\langle\nabla f(z),w-z\rangle$, if $u_n\wto u$ then we have
$$
\liminf_{n\to+\infty} \int_\Omega f(u_n)dx\ge \int_\Omega  f(u)dx+\liminf_{n\to+\infty} \int_\Omega \langle\nabla f(u),u_n-u\rangle dx.
$$
Since $\nabla f(u)\in L^\infty(\Omega;\mathbb R^m)$  is a good test function for the weak convergence, 
the last term is $0$, and we obtain the lower semicontinuity of $F$. 

We then approximate $f$ from below as follows. First, we consider the sequence of convex $n$-Lipschitz functions
$\widehat f_n$ defined by $\widehat f_n(z)=\min\{ f(w)+n|z-w|\}$, and note that $\widehat f_n(z)$ converges increasingly to $f(z)$. Then we regularize the functions $\widehat f_n-\frac1n$ using a convolution argument so that  
the resulting sequence $f_n$ is strictly increasing. This is possible since the term $-\frac1n$ is added so as to separate the graphs of these functions. Note that mollification preserves convexity and Lipschitz continuity, so that
the functionals $F_n$ given by $F_n(u)=\int_\Omega f_n(u)dx$ are weakly lower semicontinuous in $L^1$ as just seen. We finally note that $F_n(u)$ converges increasingly to $F(u)$ by Beppo Levi's Theorem (in this argument we have to note that, up to adding a linear function to $f$, whose corresponding integral gives a weakly continuous perturbation, it is not restrictive to suppose that $f_n\ge 0$). The lower semicontinuity of $F$ is then obtained by noting that the {\bf supremum of a family of lower-semicontinuous functionals} is a lower-semicontinuous functional. 

\smallskip
For future reference, we note that the arguments above can be extended to functionals
\begin{equation}\label{intLp_x}
F(u)=\int_{\Omega}f(x,u)\, dx,
\end{equation}
provided that $f$ is a Borel function (so that the composition with $u$ is measurable).
In this case the necessary and sufficient condition for $F$ to be lower semicontinuous 
is that for almost all $x\in\Omega$ the function $z\mapsto f(x,z)$ is convex and lower 
semicontinuous. \hfill $\Box$

\paragraph{Functionals on Sobolev spaces.}  
We will mainly consider functionals defined on Sobolev spaces. We recall that the Sobolev space
$W^{1,p}(\Omega;\mathbb R^m)$ is defined as the space of functions in $L^p(\Omega;\mathbb R^m)$
such that for each $i\in\{1,\ldots,d\}$ and $j\in\{1,\ldots,m\}$ there exists $v_{ij}\in L^p(\Omega)$ such that
the $j$-th component of $u$ satisfies
$$
\int_\Omega u_j\frac{\partial \varphi}{\partial x_j}dx= -\int_\Omega v_{ij} \varphi \,dx \hbox{ for all }\varphi\in C^\infty_c(\Omega).
$$
The functions $v_{ij}$ are called the weak derivatives, or derivatives in the sense of distributions, of $u$,
and their matrix (or vector if $d=1$ or $m=1$) is the weak gradient of $u$, still denoted by $\nabla u$.
The weak convergence in $W^{1,p}(\Omega;\mathbb R^m)$ is then understood as the weak convergence of 
$u_n\wto u$ and $\nabla u_n\wto \nabla u$ in the appropriate $L^p$-spaces. In this notation we remove $\mathbb R^m$ if $m=1$.

As in the case of functionals on $L^p$-spaces, we fix a bounded domain $\Omega\subset\mathbb R^d$ and consider functionals of the form 
\begin{equation}\label{intW1p}
F(u)=\int_{\Omega}f(\nabla u)\, dx,
\end{equation}
defined for $u\in W^{1,p}(\Omega)$, where $f$ is  a Borel function defined on $\mathbb R^{m\times d}$.  

In the {\bf scalar case}, remarkably, convexity and lower semicontinuity of $f$ are still necessary and sufficient for the lower semicontinuity of $F$. Indeed, while the sufficiency comes from the $L^p$-weak lower semicontinuity, the necessity relies on a {\bf gradient Riemann-Lebesgue Lemma}, stating that if $\psi\in W^{1,p}_{\rm loc}(\mathbb R^d;\mathbb R^m)$ is a function with $\nabla \psi$ $T$-periodic, then the functions $\psi_\e(x)=\e \psi(\frac{x}\e)$ weakly converge in $W^{1,p}(\Omega;\mathbb R^m)$ as $\e\to 0$ to the linear function $\overline\psi(x)= A x$ (or $\overline\psi(x)=\langle A,x\rangle$ in the case $m=1$), where $A=\frac1{T^d} \int_{(0,T)^d}\nabla \psi \,dy$. In the scalar case, the necessity of convexity is then derived in an analogous way as in the $L^p$ case, using $\psi$ with the function $\varphi$ in \eqref{oscifi} as derivative.

\smallskip
As a consequence of the Direct Method, for such an $f$ we immediately obtain solutions for problems of the form
\begin{equation}
\min\Big\{ \int_\Omega f(\nabla u) dx : u-u_0\in W_0^{1,p}(\Omega)\Big\},
\end{equation}
where $W_0^{1,p}(\Omega)$ is the closure of $C^\infty_0$ in $W^{1,p}(\Omega)$, or of the form
\begin{equation}
\min\Big\{ \int_\Omega( f(\nabla u) + gu) dx : u\in W^{1,p}(\Omega)\Big\}
\end{equation}
if $\Omega$ is connected and $\int_\Omega g\,dx=0$. In both cases the coerciveness is obtained by using the Poincar\'e or the Poincar\'e--Wirtinger inequality, as in the second case it is not restrictive to examine only $u$ with zero average.

More in general we can consider functionals of the form 
$$
\int_\Omega( f(\nabla u) + g(x,u) )dx,
$$
provided that $g$ is a Carath\'eodory  function with $|g(x,u)|\le C(1+|u|^q)$ with $q<p^*$ if $p<d$. 
Indeed, in this case the perturbation is a continuous functional with respect to the strong convergence in $L^q(\Omega)$, so that the sum is weakly lower semicontinuous in $W^{1,p}(\Omega)$ by Rellich's theorem.

\smallskip
The determination of necessary and sufficient conditions for weak lower semicontinuity in the vector case $d,m>1$ is more delicate. In particular convexity is not a necessary condition. To understand why the previous argument does not work, we remark that if $z,z'$ are matrices, in general the function $\varphi$ in \eqref{oscifi} is not a derivative of a $W^{1,p}$ function, unless the rank of $z-z'$ is $0$ or $1$, so that the necessity argument cannot be repeated.
\hfill$\Box$

\paragraph{Examples on the effect of lack of lower semicontinuity.} 
Now we look at some prototypical problems, where the lack of lower semicontinuity may translate either into the existence of infinitely many solutions, or of no solutions at all, and we may describe the form of some minimizing sequences.
\smallskip 

\begin{figure}[h!]
	\centerline{\includegraphics [width=2.5in]{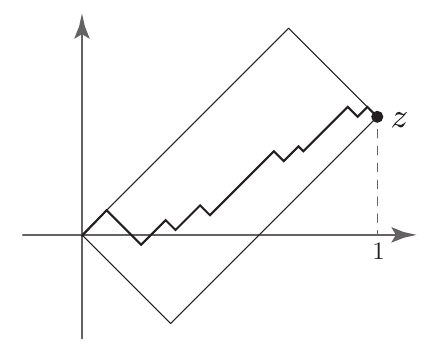}}
	\caption{Minimizers for a double-well problem}\label{Fig1}
\end{figure}
\noindent{\em Infinitely many minimizers.} \ As a first example we consider the minimum problem 
$$\min\Big\{\int_0^1 W(u^\prime)\, dt: u(0)=0, u(1)=z\Big\}$$
where $W$ is a double-well potential as above; that is, $W\colon\mathbb R\to [0,+\infty)$ is a sufficiently smooth positive function vanishing only at two points, which we can suppose to be $-1$ and~$1$. 

We consider the case of a boundary datum $z$ belonging to the interval between the two minimizers; that is, $|z|\leq 1$. 
In this case, we can construct a piecewise-affine test function such that $u^\prime\in\{-1,1\}$ by setting $u(t)=t$ in $(0,\frac{1+z}{2})$ and  $u(t)=1+z-t$ in $(\frac{1+z}{2},1)$; then, the value of the functional is $0$ and $u$ is a minimizer. Similarly, we can also construct a minimizer with $u'=-1$ on $(0, \frac{1-z}{2})$ and $u'=1$ on $(\frac{1-z}{2},1)$. In Fig.~\ref{Fig1} we picture a minimizer which mixes the two slopes $-1$ and $1$. In this way we can construct infinitely many piecewise-affine minimizers. 

\begin{figure}[h!]
	\centerline{\includegraphics [width=5in]{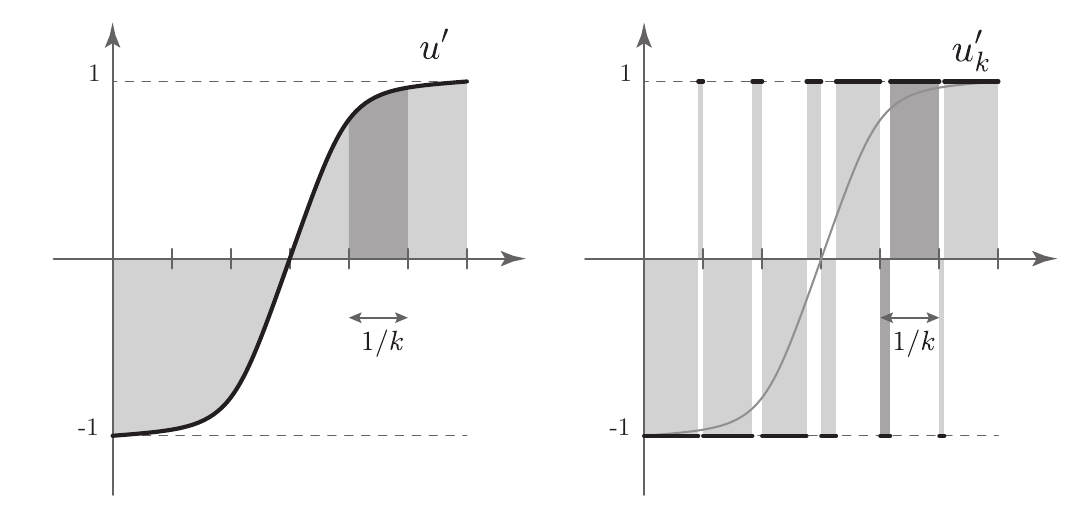}}
	\caption{Construction of minimizers approximating a given $u$}\label{Fig2}
\end{figure}
 We can now consider the problem of characterizing the limits of sequences of minimizers. We note that we can approximate with a sequence of minimizers any function $u$ satisfying the boundary data and the condition $\|u^\prime\|_\infty\leq 1$. A pictorial description of a possible way to explicitly construct such an approximation is given in Fig.~\ref{Fig2}: given such a $u$, in each interval $((i-1)_k, i/k)$ we define the derivative of a function $u_k$ as taking only the values $-1$ and $1$ so that the integral of $u'_k$ on this interval equals that of $u'$. The function $u_k$ is then uniquely defined by the condition $u_k(0)=0$, is a minimizer, and $u_k$ converges weakly in $H^1(0,1)$ and strongly in $L^\infty(0,1)$ to $u$.

 \begin{figure}[h!]
	\centerline{\includegraphics [width=4in]{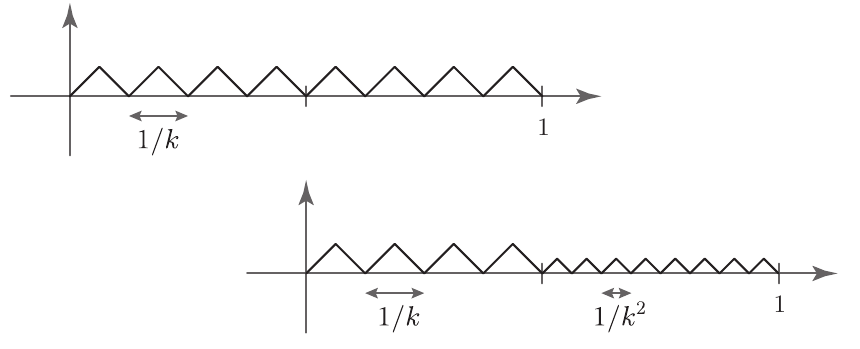}}
	\caption{Oscillations at different scales with $u=0 $.}\label{Fig3}
\end{figure}
We note that, unless $u$ is itself a minimizer, the convergence of $u_k\wto u$ is weak in $H^1(0,1)$. This entails oscillations of the gradient of $u_k$. We note however that these oscillations have no given scale. 
Now we take the boundary condition $z=0$. In Fig.~\ref{Fig3} we take as target function $u=0$, and we picture a sequence $u_k$ constructed as above, oscillating at scale $1/k$, and another sequence of minimizers with oscillations at scale $1/k$ on $(0,1/2)$ and at scale $1/k^2$ on $(1/2,1)$.

 \begin{figure}[h!]
	\centerline{\includegraphics [width=5in]{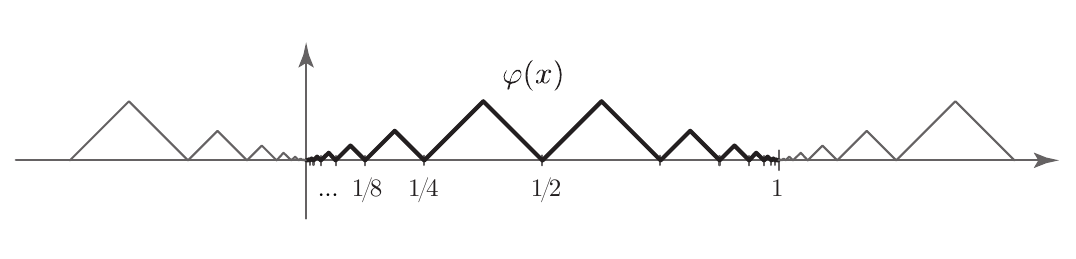}}
	\caption{A function with oscillations at all scales $2^{-j}$ with $j\ge 2$.}\label{Fig4}
\end{figure}
 \begin{figure}[h!]
	\centerline{\includegraphics [width=2.9in]{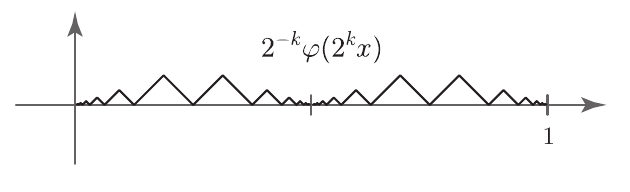}}
	\caption{A test function with oscillations at all scales $2^{-j}$ with $j\ge k+2$.}\label{Fig5}
\end{figure}

We can construct an even more complex example of a sequence of minimizers converging to $0$ and with oscillations at all scales. Their construction is described pictorially in Fig.~\ref{Fig4} and Fig.~\ref{Fig5}. The first one describes the construction of an even one-periodic $\varphi$ with $|\varphi'|=1$ almost everywhere and oscillations
at all scales $2^{-j}$ with $j\ge 2$. The functions $u_k(x)=2^{-k}\varphi(2^{k}x)$ are test functions for the problem with $z=0$ and $u_k$ has oscillations at all scales $2^{-j}$ with $j\ge k+2$.
 
\smallskip 

\noindent{\em Oscillating minimizing sequences without a scale.}\ The second example, due to {\bf Bolza}, is the problem of minimizing the functional 
$$F(u)=\int_0^1 \big(W(u^\prime)+|u|^2\big)\, dt.$$
In this case no boundary conditions are prescribed. The same functions constructed in the previous example and represented in Figs.~\ref{Fig3}--\ref{Fig5} show that
$$\inf\Big\{\int_0^1 \big(W(u^\prime)+|u|^2\big)\, dt\Big\}=0,$$  
and minimizers can develop oscillations at all scales.

Note however that the minimum is not achieved since the minimization of $|u|^2$ would give $u=0$, for which the contribution of $W(u')$ is not zero.

\smallskip 

\noindent{\em Concentrating minimizing sequences without a scale.}\  In this example we see how the problem of minimizing a concave functional will give rise to a phenomenon of concentration and eventually of non-existence.
We minimize
$$F(u)=\int_0^1 \sqrt{|u'|^2+|u'|}\, dt$$ 
over the set of $W^{1,1}$-functions satisfying  $u(0)=0$ and  $u(1)=1$. 
For any such $u$,  
$$F(u)=\int_0^1 \sqrt{|u'|^2+|u'|}\, dt\geq \int_0^1 |u'|\, dt \geq \Big|\int_0^1 u'\, dt\Big|=u(1)-u(0)=1.$$ 
Hence, the infimum of $F(u)$ is greater than $1$.  

We now fix $x\in(0,1)$, and (for all $n$ large enough) we set $u_n(t)=\min\{n(t-x)^+,1 \}$. Since  
$$\int_0^1 \sqrt{|u'_n|^2+|u'_n|}\, dt= \frac{\sqrt{n^2+n}}{n}=1+o(1)_{n\to+\infty},$$
we get 
$$\inf\Big\{\int_0^1 \sqrt{|u'|^2+|u'|}\, dt:u(0)=0, u(1)=1\Big\}=1.$$  
Note that the sequence $\{u_n\}$ tends to the function $u_x=\chi_{(x,1)}$.
In the same way one can find minimizing sequences converging to functions of the form $\sum_{j=1}^N a_j\chi_{(x_j,1)}$ with $a_j>0$ and $\sum_{j=1}^N a_j=1$. The closure (e.g.~in $L^1(0,1)$) of this set is the set of all increasing functions with $u(0)=0$ and $u(1)=1$.

From the inequality above, we see that the equality $F(u)=1$ holds only if $u'=0$ almost everywhere.
This cannot hold if $u(1)\neq u(0)$, so that no minimizer exists.

%

\subsection{Some use of variations}\label{EL-sec}
 Even though we will be careful not to rely on regularity assumptions, sometimes it is convenient to be able to solve a minimum problem by looking at the Euler--Lagrange equation.
 
 Suppose now that $f$ is smooth enough and consider 
$$
F(u)=\int_a^bf(t,u,u^\prime)\, dt.
$$
The Euler--Lagrange equation for this functional is 
$$
-\frac{d}{dt}\frac{\partial f}{\partial \xi}(t,u,u^\prime)+\frac{\partial f}{\partial u}(t,u,u^\prime)=0,
$$
and it is satisfied if $u\in C^2(a,b)$ is a minimizer of $F$ on an interval with respect to functions with its own boundary data. The equation also holds in the sense of distributions under some growth conditions on $f$.

In general, solutions of the Euler--Lagrange equations may not be minimizers. An exception is given by functionals that are convex in the pair $(u,\xi)$.
If the function $(u,\xi)\mapsto f(t,u,\xi)$ is convex, then a function $\overline u$ satisfying the Euler--Lagrange equations is a minimum of the functional $F$ with respect to its own boundary conditions.  Indeed, let $u$ be a test function; that is, $u\in C^2(a,b)$ satisfying 
the same boundary conditions as the solution $\overline u$ of the Euler-Lagrange equation. 
Then, by the convexity condition applied at each point $(\overline u(t), \overline u^\prime(t))$, 
it follows that 
\begin{eqnarray*}
f(t, u(t),u^\prime(t))&\geq&f(t, \overline u(t), \overline u^\prime(t))+
\frac{\partial f}{\partial u}(t, \overline u(t), \overline u^\prime(t))(u(t)-\overline u(t))\\
&&+\frac{\partial f}{\partial \xi}(t, \overline u(t), \overline u^\prime(t))(u^\prime(t)-\overline u^\prime(t)).
\end{eqnarray*}
Here we have used the convexity inequality, and the fact that 
$(\frac{\partial f}{\partial u}, \frac{\partial f}{\partial \xi})$ is the gradient of $f$ with respect to the pair of variables $(u,\xi)$. By integrating the two terms of the inequality in the domain $(a,b)$, and then integrating by parts, since $u-\overline u$ vanishes on the boundary we obtain 
\begin{eqnarray*}
F(u)\ge F(\overline u)+\int_a^b\Big(
\frac{\partial f}{\partial u}(t, \overline u(t), \overline u^\prime(t) -\frac{d}{dt}\frac{\partial f}{\partial \xi}(t, \overline u(t), \overline u^\prime(t)\Big)(u(t)-\overline u(t))\, dt =F(\overline u), 
\end{eqnarray*}
where in the last equality we used the fact that $\overline u$ solves the Euler--Lagrange equation, so that the integrand vanishes and $\overline u$ minimizes the functional in the set of smooth functions satisfying the same boundary conditions.

As an application, we consider a function of the form  
$f(t,u^\prime)=a(t)|u^\prime|^2$, and a solution $\overline u$ of the Euler--Lagrange equation 
$$\big(a(t)u^\prime(t)\big)^\prime=0$$
in an interval $I$ where $a$ is smooth. Then,  
the convexity of $f$ with respect to $u^\prime$ implies that 
$\overline u$ minimizes 
$$F(u)=\int_{I} a(t)|u^\prime|^2\, dt$$
over functions $u$ with the same boundary conditions. 
This can be repeated on all intervals where $a$ is smooth. We choose the simple case when 
$a(t)$ is given by 
$$a(t)=\begin{cases}\alpha & \hbox{\rm if } t\in [0,1/2),\\
\beta & \hbox{\rm if } t\in [1/2,1), 
\end{cases}$$ 
with $\alpha,\beta>0$
and consider the problem 
$$\min\Big\{\int_0^1 a(t)|u^\prime|^2\, dt: u(0)=0, u(1)=1\Big\}.$$ 
The solution $\overline u$, which exists by the Direct Method,  satisfies the Euler--Lagrange equation in both intervals $(0,1/2)$ and $(1/2,1)$, 
hence $u$ is the piecewise-affine function satisfying $\alpha u^\prime(t)=c_\alpha$ in $(0,1/2)$ and $\beta u^\prime(t)=c_\beta$ in $(1/2,1)$ and the imposed boundary conditions. 
The minimum problem for the functional becomes the problem of finding the two optimal constants $c_\alpha$ and $c_\beta$, which can be determined by  
finding the minimizer of 
the function $g(s)=\alpha s^2+\beta (2-s)^2$, so that  
$c_\alpha=c_\beta=\frac{2\alpha\beta}{\alpha+\beta}$, which is also the value of the minimum.  
Note in particular that the solution $\overline u$ 
does not belong to $C^1(0,1)$ if $\alpha\neq\beta$.

  \begin{figure}[h!]
	\centerline{\includegraphics [width=1.8in]{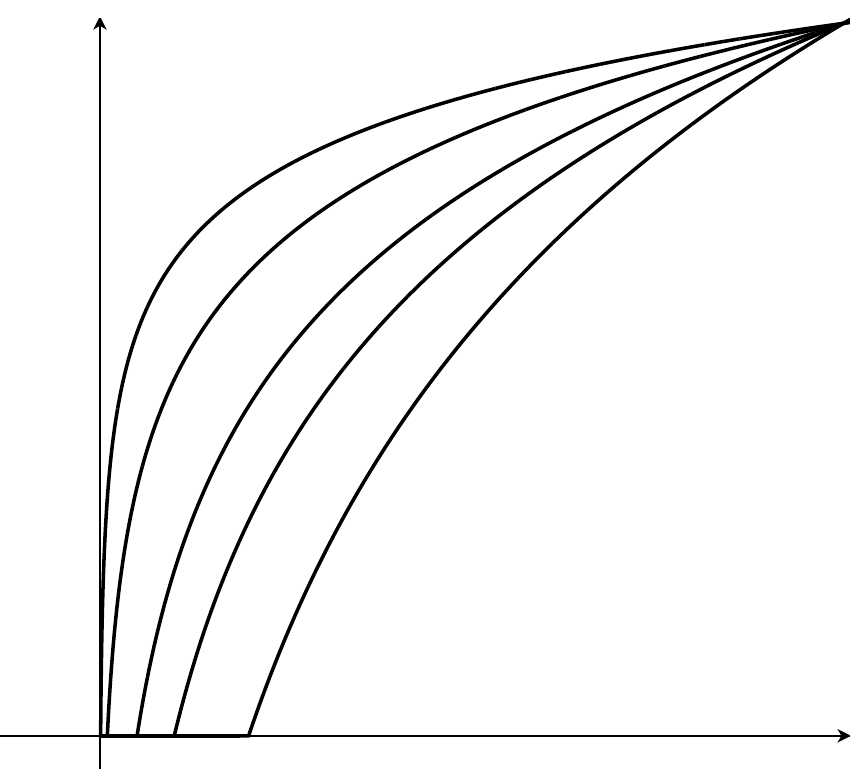}}
	\caption{A minimizing sequence for the Weierstrass example.}\label{Fig8}
\end{figure}
\smallskip
Another example, this time of how the Euler--Lagrange equations can be of use in obtaining minimizing sequences, even if the problem has no solution,
is the so-called {\bf Weierstrass example}. This is the problem of minimizing 
$$F(u)=\int_0^1x|u^\prime|^2\, dx$$ 
over the set of the functions $u\in H^1(0,1)$ satisfying the boundary conditions $u(0)=0$ and $u(1)=1$. 
Since the integrand is smooth, the Euler--Lagrange equation 
$$
(x u')'=0
$$
implies that a solution $u$ of the minimum problem is such that 
$x u^\prime$ is constant; that is, taking into account the condition $u(1)=1$, that
$$
u(x)= 1+C\log x
$$ 
for some $C\in\mathbb R$, which is incompatible with the initial condition $u(0)=0$.
Hence, a smooth minimum does not exist, otherwise it would be characterized by the Euler--Lagrange equation. 
In this case, the Euler--Lagrange equation also helps in proving that the infimum is $0$. To that end, we can test
the minimum problem with $u_C(x)= \max\{1+C\log x, 0\}$ for $C>0$, obtaining $F(u_C)=C$, and the infimum is obtained letting $C\to 0$.
Note that this minimizing sequence tends to the constant $1$ as $C$ tends to $0$ (see Fig.~\ref{Fig8}).

\subsection{Characterization of minimizing sequences: relaxation}\label{Pre-Relax}
With the examples on the effect of lack of lower semicontinuity in Section \ref{semicont-sec} in mind, we now describe the behaviour of minimizing sequences of a non-lower-semicontinuous functional.   
We will see that the limits of minimizing sequences solve a {\bf relaxed minimum problem}. 
In order to describe this latter problem, we start by giving the notion of {\bf lower-semicontinuous envelope} of a functional $F$. If $X$ is a metrizable topological space, and $F\colon X\to \mathbb R\cup\{+\infty\}$, then the lower-semicontinuous envelope of $F$ is the functional $\overline F$ defined on $X$ by  
\begin{eqnarray*}
\overline F(x)&=&\sup\{G(x): \ G\leq F, \ G \ \hbox{\rm lower-semicontinuous}\}.
\end{eqnarray*} 
As we have seen, the supremum of a family of lower-semicontinuous functions is lower semicontinuous, hence the $\sup$ in the definition is in fact a maximum. 
Since $X$ is metrizable, it follows that 
$$\overline F(x)=\inf\Big\{\liminf_{j\to+\infty} F(x_j): \{x_j\}\subset X, x_j\to x\Big\},$$ 
where the convergence $x_j\to x$ is understood as induced by the metric. 
This characterization is simpler to use for the computation of the lower-semicontinuous envelope.  
Note that the infimum is in fact a minimum, since we can use a diagonal argument to construct a sequence achieving the infimum. Then, by considering converging subsequences, we have that 
\begin{equation}\label{car-F}
\overline F(x)=\min\Big\{\lim_{j\to+\infty} F(x_j): \{x_j\}\subset X \ \hbox{\rm such that } x_j\to x \ \hbox{\rm and  } \lim_{j\to+\infty}F(x_j) \hbox{\rm exists} \Big\}.
\end{equation}
To show that a functional $G$ coincides with the lower-semicontinuous envelope of $F$ at a point $x$, we have to prove that $G(x)$ is a lower bound for the infimum; that is, 

\smallskip
$\bullet$ $G(x)\leq \liminf\limits_{j\to +\infty}F(x_j)$ for all sequences $\{x_j\}$ converging to $x$,

\smallskip\noindent
and that it is optimal; that is, 
\begin{enumerate}
\item[$\bullet$] there exists a sequence $\{x_j\}$ converging to $x$ such that 
$G(x)=\lim\limits_{j\to +\infty}F(x_j)$. 
\end{enumerate} 
\paragraph{Cluster points of minimizing sequences.} 
Using the definition of $\overline F$, we can characterize the relaxed minimum problem solved by the limits of minimizing sequences. 
To this end, we note that, for any given $x\in X$, by the characterization of $\overline F(x)$ in \eqref{car-F} there exists $x_j\to x$ such that $F(x_j)\to \overline F(x)$; then, 
$$\overline F(x)=\lim_{j\to+\infty} F(x_j)\geq \inf\{F(y): y\in X\}$$
and we deduce that $\inf\{F(y):y\in X\}=\inf\{\overline F(y):y\in X\}$. 
Now, let $\{\overline x_j\}$ be a minimizing sequence for the problem 
such that $\overline x_j\to \overline x\in X$. We get 
$$\overline F(\overline x)\leq \lim_{j\to+\infty}F(\overline x_j)=\inf\{F(x):x\in X\}=\inf\{\overline F(x): x\in X\}.$$
This argument proves that if there exists a converging minimizing sequence for the problem 
$\inf\{F(x): x\in X\}$, 
then

\smallskip$\bullet$ there exists the minimum of $\{\overline F(x): x\in X\}$,  
and 
$$\min\{\overline F(x): x\in X\}=\inf\{F(x): x\in X\}$$ 

({\bf relaxation} of minimum problems);  

\smallskip
$\bullet$ the minimum points for $\overline F$ are all the cluster points of minimizing sequences for $F$.  

\smallskip\noindent
This is a {\bf Weierstrass theorem for non lower-semicontinuous functionals}.
 
 \smallskip 
Before reconsidering the previous examples of minimum problems for non lower-semi\-continuous functionals, we note the connection between lower semicontinuous envelopes and convexifications. 	The simplest example is the (one-dimensional) lower semicontinuous envelope $\overline F$ of a functional
$$
F(u)=\int_a^b f(u)\,dt
$$
with respect to the weak $L^1$-convergence.

 \begin{figure}[h!]
	\centerline{\includegraphics [width=2.4in]{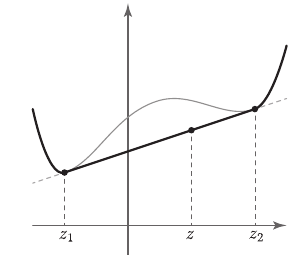}}
	\caption{The convex envelope of a function $f$.}\label{Fig6}
\end{figure}

We define $f^{\ast\ast}$ as the {\bf lower-semicontinuous and convex envelope} of the function $f$ defined as
$$
f^{\ast\ast}(\xi)=\sup\{g(\xi): g\leq f, \ g \ \hbox{\rm lower-semicontinuous and convex}\}.
$$
Then 
$$
\overline F(u)=\int_a^b f^{\ast\ast}(u)\,dt.
$$
The inequality $\overline F(u)\ge\int_a^b f^{\ast\ast}(u)dt$ is due to the lower semicontinuity of the latter.
To understand the converse inequality, we restrict to the case of scalar $u$, and suppose that $f^{\ast\ast}<f$
on an interval $(z_1,z_2)$, with  $f^{\ast\ast}(z_1)=f(z_1)$ and   $f^{\ast\ast}(z_2)=f(z_2)$ (Fig.~\ref{Fig6}). We only consider the case $u(x)=z$ constant with $z=tz_1+(1-t)z_2$. Then we note that $f^{\ast\ast}(z)=tf(z_1)+(1-t)f(z_2)$, and an optimal sequence is obtained by periodically alternating $z_1$ and $z_2$, and using the Riemann--Lebesgue Lemma.

\smallskip
This characterization of the lower semicontinuous envelope can be generalized in many ways.
A general {\bf relaxation} result for functionals {\bf on Sobolev spaces} of the form
$$F(u)=\int_\Omega f(x,\nabla u)\, dx,$$ 
defined for $u\in W^{1,p}(\Omega)$, with respect to the weak convergence is as follows.
Assume that $f$ satisfies a growth condition ensuring that $F$ is coercive in $W^{1,p}(\Omega)$. Then, we can limit the problem of minimizing $F$ to bounded subsets of $W^{1,p}(\Omega)$, where the weak topology is metrizable (and, equivalently, we can use the strong $L^p$-convergence). 
Then, the lower-semicontinuous envelope of $F$ is given by  
$$\overline F(u)=\int_\Omega f^{\ast\ast}(x,\nabla u)\, dx$$
where $f^{\ast\ast}(x,\cdot)$ is the lower-semicontinuous and convex envelope of the function $\xi\mapsto f(x,\xi)$ defined as
$$f^{\ast\ast}(x,\xi)=\sup\{g(\xi): g\leq f, \ g \ \hbox{\rm lower-semicontinuous and convex}\};$$
that is, since convexity and lower-semicontinuity are maintained by taking the supremum, the largest lower-semicontinuous and convex function below $f(x,\cdot)$. 

 \begin{figure}[h!]
	\centerline{\includegraphics [width=1.8in]{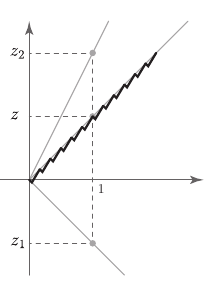}}
	\caption{Optimal sequence for an affine function.}\label{Fig7}
\end{figure}

 In the simplest case of a homogeneous integrand this result is proved in the same way as the corresponding result in $L^p$. In the one-dimensional case an optimal sequence for the affine function $u(x)=zx$ is pictured in Fig.~\ref{Fig7}, with reference to Fig.~\ref{Fig6}. 
We apply this relaxation result to our minimum problems involving the double-well potential $W$. 

\smallskip

\noindent{\em Infinitely many minimizers.}\ We can relax the minimum problem for the integral of $W(u^\prime)$ with prescribed boundary datum $z\in\mathbb R$, obtaining 
$$\inf\Big\{\int_0^1 W(u^\prime)\, dt: u(0)=0, u(1)=z\Big\}=
\min\Big\{\int_0^1 W^{\ast\ast}(u^\prime)\, dt: u(0)=0, u(1)=z\Big\},$$ 
where $W^{\ast\ast}$ is the convex envelope of $W$. Note that 
$W^{\ast\ast}=0$ in the interval $[-1,1]$ between the minimizers.  
We have seen with an explicit construction that all functions $u$ satisfying the given boundary data and $\|u^\prime\|_{\infty}\leq 1$ are limits of minimizing sequences.  If $W$ satisfies a superlinear growth condition at infinity, then  $W^{\ast\ast}$ vanishes exactly in $[-1,1]$ and hence the set of functions satisfying the given boundary data and $\|u^\prime\|_{\infty}\leq 1$ is exactly the set of all limits of minimizing sequences.  

\smallskip 

\noindent{\em Oscillating minimizing sequences without a scale.}\   In Bolza's example, the relaxation gives 
$$\inf\Big\{\int_0^1 (W(u^\prime)+|u|^2)\, dt\Big\}=
\min\Big\{\int_0^1 (W^{\ast\ast}(u^\prime)+|u|^2)\, dt\Big\}=0,$$ 
and the unique minimum point is $u=0$. 
Again, this proves the initial guess that the unique cluster point of minimizing sequences is $0$.

\smallskip

Note that the {\em concentrating minimizing sequences without a scale} example of Section~\ref{semicont-sec} is deliberately not revisited here. There, non-existence is not caused by a failure of convexity of the integrand (which relaxation, in the form of the convex envelope, would repair), but by a loss of coerciveness: along the minimizing sequence the derivatives $u_n^\prime$ do not remain bounded in $L^1$, and concentrate instead, in the limit, as a singular measure. This is a genuinely different phenomenon from the previous two, outside the reach of the relaxation theory developed in this section; it is the prototype, {\em without a scale}, of the concentration phenomena that reappear, now driven by a parameter $\e$, in the model case \eqref{model-conc} of the Introduction and that are studied systematically, with their natural scales, starting from Section~\ref{conc-sec}.

\paragraph{Definitions by relaxation.}
Note that, if we have a problem 
which is well-defined on a suitable space $X$, 
we can extend the definition to a larger set $\tilde X$ 
by relaxation. 

\smallskip

\noindent{\em Sobolev spaces by relaxation.}\ 
As an example of definition by relaxation, we can consider 
the functional 
$$F(u)=\int_{\Omega} \big(|u|^2+|\nabla u|^2\big) \, dx$$
defined in $C^1(\Omega)$, and extend it to $L^2(\Omega)$ by setting $F(u)=+\infty$ in $L^2(\Omega)\setminus C^1(\Omega)$. 
The lower-semicontinuous envelope of $F$, whose domain is the set of all functions in $L^2(\Omega)$ which can be approximated by functions such that $F(u)$ is finite, is the functional $\overline F\colon L^2(\Omega)\to [0,+\infty]$ 
given by 
$$ 
\overline F(u)=\begin{cases}\displaystyle \int_{\Omega} \big(|u|^2+|\nabla u|^2\big) \, dx & \hbox{\rm if } u\in W^{1,2}(\Omega)\\
+\infty & \hbox{\rm if } u\in L^2(\Omega)\setminus W^{1,2}(\Omega).
\end{cases}
$$
This is actually a version of the {\bf $H=W$ Meyers--Serrin Theorem} stating that  $W^{1,2}(\Omega)$ coincides with the space $H^1(\Omega)$; the proof follows from the fact that all functions in the Sobolev space can be approximated by $C^1$ (or $C^\infty$) functions.   
The fact that in the relaxed functional we have the same formula is due to the fact that $F$ is lower-semicontinuous in $C^1(\Omega)$. 

\smallskip

\noindent{\em Sets of finite perimeter by relaxation.}\ 
A measurable set $E\subseteq\Omega\subset\mathbb R^d$ is a set of finite perimeter in $\Omega$ if 
$$\sup\Big\{\int_{\Omega}\chi_E\ \hbox{\rm div}(\varphi)\,dx  : \varphi\in C_c^1(\Omega;\mathbb R^d) \ \hbox{\rm such that } \|\varphi\|_\infty\leq 1 \Big\}<+\infty.$$ 
This quantity is the perimeter of $E$ in $\Omega$, denoted by Per$(E;\Omega)$. 
Note that, if $E$ is smooth, then by integrating by parts we get 
$$\int_{\Omega}\chi_E\ \hbox{\rm div}(\varphi)\,dx=\int_{\Omega\cap\partial E} \langle\varphi, \nu\rangle\, d\Sigma,$$ 
where $\nu$ is the exterior normal to $\partial E$ and $d\Sigma$ denotes integration with respect to the $d-1$-dimensional surface measure. 
Then, the perimeter coincides with the $d-1$-dimensional surface area of the boundary of $E$ in $\Omega$. 

If we define the functional 
$$ 
F(E)=\begin{cases} 
\hbox{\rm Per}(E;\Omega) & \hbox{\rm if } E \ \hbox{\rm smooth} \ \hbox{\rm (or if $E$ is a polytope)}\\
+\infty & \hbox{\rm otherwise,}
\end{cases}
$$
and consider the convergence $E_j\to E$ given by the convergence to $0$ of the measure of the symmetric differences $E\Delta E_j$, or by the convergence $\chi_{E_j}\to \chi_E$ in $L^1_{\rm loc}(\Omega)$, then the lower-semicontinuous envelope of $F$ is given by 
$$ 
\overline F(E)=\begin{cases} 
\hbox{\rm Per}(E;\Omega) & \hbox{\rm if } E \ \hbox{\rm is a set of finite perimeter in $\Omega$}\\
+\infty & \hbox{\rm otherwise.} 
\end{cases}
$$
Note that also in this case the functional is lower-semicontinuous in the space of smooth sets with respect to the same convergence, since it is the supremum of a family of lower-semicontinuous functions. 
This relaxation result states that every set of finite perimeter $E$ can be approximated by smooth sets $E_h$ (or by polytopes) in such a way that the perimeters of $E_h$ converge to the perimeter of $E$.

\subsection{Limits of minimizers and limit functionals: $\Gamma$-convergence}\label{Gconv-sec}
We go back to the two model cases \eqref{model-conc} and \eqref{model-osc} described in the Introduction, where we have seen concentration and oscillation, respectively, of the solutions of minimum problems depending on a parameter $\e$. We now want to extend the analysis to a general setting, including cases when the minimum problems have no solutions. 

Given a family of minimum problems $m_\e$ depending on a parameter $\e$, the idea is to introduce suitable functionals $F_\e$ such that $\inf F_\e=m_\e$, and analyze the limit behaviour of {\em minimizing sequences} \index{minimizing sequences}$u_\e$, such that $$F_\e(u_\e)=m_\e+o(1)_{\e\to 0},$$  
showing that the cluster points of such sequences are solutions of a minimum problem depending on a ``limit'' functional independent of $\e$. In the same spirit as the definition of the relaxation, the definition of such functional will be given in topological terms.    
We recall that $\overline F$ is the lower-semicontinuous envelope of $F$ at $x$ if two conditions are satisfied; that is, for all $x_j\to x$ the lower estimate $\liminf_{j}F(x_j)\geq \overline F(x)$ holds, and there exists a sequence $x_j\to 0$ such that $\lim_{j}F(x_j)=\overline F(x)$. 
We now generalize this definition to a sequence of functionals $F_j$ in such a way that a Weierstrass theorem still holds. Let $X$ be a metrizable space, and let $F_j\colon X\to\mathbb R\cup \{+\infty\}$ be a sequence of functionals depending on a parameter $j$. Then, a functional $F$ is the {\bf  $\Gamma$-limit of $F_j$} at $x\in X$, denoted by $F(x)=\Gamma$-$\lim\limits_{j\to+\infty}F_j(x)$, if 
\begin{enumerate}
\item[$\bullet$] {(\em $\liminf$ inequality)} $\liminf\limits_{j\to+\infty}F_j(x_j)\geq F(x)$ for all sequences $x_j\to x$; 
\item[$\bullet\bullet$] {(\em $\limsup$ inequality)} there exists a sequence $\overline x_j\to x$ such that $\limsup\limits_{j\to+\infty}F_j(\overline x_j)\leq F(x)$. 
\end{enumerate}
Note that, once the $\liminf$ inequality holds, 
the $\limsup$ inequality is equivalent to 
\begin{enumerate} 
\item[$\bullet\bullet$] {(\em existence of a recovery sequence)} there exists $\overline x_j\to x$ such that $\lim\limits_{j\to+\infty}F_j(\overline x_j)\leq F(x)$. 
\end{enumerate} 
The following equivalent formulation of the $\limsup$ inequality can in general be useful in the computation of $\Gamma$-limits: 
\begin{enumerate} 
\item[$\bullet\bullet$] {(\em approximate $\limsup$ inequality)} for all $\eta>0$, there exists a sequence $\overline x_j\to x$ such that $\limsup\limits_{j\to+\infty}F_j(\overline x_j)\leq F(x)+\eta$.  
\end{enumerate}  

We will use the corresponding definition for a family of functionals $F_\e$ depending on a parameter $\e>0$, meaning that, for all sequences $\e_j\to 0$, the functionals $F_j=F_{\e_j}$ $\Gamma$-converge to the same $F$. Moreover, we mention that $\Gamma$-convergence enjoys the {\bf Urysohn property}; that is, the $\Gamma$-convergence of $F_j$
is equivalent to the property that
for all sequences $\{j_k\}$ we can extract a subsequence of $\{F_{j_k}\}$ which $\Gamma$-converges to $F$ (see Proposition~\ref{Urysohn-prop} in the Appendix to this Chapter for a proof).

With the definition of $\Gamma$-convergence for a sequence of functionals, we can go back to the problem of the limit behaviour of a minimizing sequence for a family of functionals $F_j$. 
Assume then that there exists a converging minimizing sequence for the sequence $\{F_j\}$; that is, 
a sequence $\{\overline x_j\}$ converging to some $\overline x\in X$ and such that $F_j(\overline x_j)=\inf F_j+o(1)_{j\to+\infty}$, and suppose that $F$ is the $\Gamma$-limit of $F_j$.  
Then, by the $\liminf$ inequality, 
$$\liminf\limits_{j\to+\infty}\inf\{F_j(x): x\in X\}=\liminf\limits_{j\to+\infty}F_j(\overline x_j)\geq F(\overline x).$$ 
For all $x\in X$, there exists a sequence $x_j\to x$ such that  
$$F(x)\geq \limsup_{j\to+\infty}F_j(x_j)\geq \limsup_{j\to+\infty} \inf\{F_j(y): y\in X\};$$ 
taking the infimum over $x\in X$, by the previous inequality we get 
$$\inf\{F(x):x\in X\}=F(\overline x).$$ 
We have then shown that, if $F=\Gamma$-$\lim_j F_j$ and there exists a converging minimizing sequence for $F_j$, then: 

\smallskip
$\bullet$ there exists the minimum $\min\{F(x): x\in X\}$, and 
$$\min\{F(x): x\in X\}=\lim_{j\to+\infty}\inf\{F_j(x):x\in X\}$$ 

(convergence of minima); 

\smallskip
$\bullet$ the minimum points for $F$ are all the cluster points of minimizing sequences for $F_j$ 

(convergence of minimizing sequences);      

\smallskip\noindent
that is, a {\bf Weierstrass Theorem along a $\Gamma$-converging sequence} holds. 

\bigskip

A simple example of convergence that can be stated in terms of $\Gamma$-convergence is given by problems involving {\em finite-difference approximations}. The simplest case, but still bearing the characteristics of the general case, is the one-dimensional discrete problem, where the parameters are functions $u:\{0,\ldots, n\}\to \mathbb R$, $i
\mapsto u_i$, and the energy is the quadratic energy
$$
F_n(u)=\sum_{i=1}^n  n(u_i-u_{i-1})^2	.
$$
We can reinterpret such a $u$ as an element of the subspace $\mathcal P_n$ of the piecewise-affine functions defined on $[0,1]$ and with values $u(\frac in)=u_i$.
We then have
$$
F_n(u)=\int_0^1|u'|^2dt \qquad \hbox{ if } u\in \mathcal P_n,
$$
so that the $\Gamma$-limit of $F_n$ is simply $F(u)=\int_0^1|u'|^2dt$ with domain $H^1(0,1)$. The liminf inequality is a consequence of the lower semicontinuity of the norm, while the limsup inequality follows from the strong density of functions in $\mathcal P_n$. Note that similarly we can consider energies 
$$
G(u)=\sum_{i=1}^n \frac1n g_iu_i  \qquad \hbox{ if } u\in \mathcal P_n,
$$
for example if $g$ is continuous on $[0,1]$ and $g_i=g(\frac in)$. If $u_n\to u$ in $L^2(0,1)$ then  $G(u_n)\to G(u)$, so that $F_n+G_n$ $\Gamma$-converges to $F+G$, and problems of the form
$$
\min\Big\{\int_0^1 (|u'|^2-gu)dt: u(0)=\overline u_0, u(1)=\overline u_1\Big\}
$$
are approximated by discrete problems of the form 
$$
\min\bigg\{\sum_{i=1}^n  n(u_i-u_{i-1})^2-\sum_{i=1}^n \frac1n g_iu_i : u=\overline u_0, u_n=\overline u_1\bigg\}.
$$
The same holds in any dimension, also for arbitrary reference sets $\Omega$, once we consider a local convergence, and we identify discrete functions with interpolations on suitable triangulations.
The only thing to remember is that the notation is a little different, taking, in place of $F_n$, the functionals
$$
F_\e(u)=\sum_{i,j\in\mathbb Z^d, |i-j|=1}\e^{d-1} (u_i-u_j)^2
$$
and the functions $\{u_i\}$ are identified with interpolations of the functions defined on $\e\mathbb Z^d$ by $u(\e i)=u_i$.
For finite-difference schemes, $\e$ is a {\bf scale of discretization}, in the sense that the gradients of the related functions $u$ are constant on sets (elements of the triangulation) at the scale $\e$.

These approximations can be coupled with relaxation, giving a {\bf relaxation with a scale}. For example, if $f$ is not convex then we can use finite-difference approximation to characterize the lower-semicontinuous envelope of $\int f(u)dt$, obtaining (in the case of $f$ continuous for simplicity) a discrete formula for the convex envelope of $f$
$$
f^{\ast\ast}(z)=\lim_{n\to+\infty} \min\bigg\{\frac 1n\sum_{i=1}^n f(z_i): \frac 1n\sum_{i=1}^n z_i=z\bigg\}.
$$

\bigskip
In the rest of these notes, our analyses will in general be translated into 

\smallskip
$\bullet$  computation of the $\Gamma$-limit of some $F_\e$; 

\smallskip
$\bullet$ the relation of the form of the recovery sequences to oscillation or concentration phenomena 

at some scale, described in terms of the parameter $\e$.

\subsection{Appendix: relevant tools of Chapter 2}\label{App-Ch2}

\paragraph{References.} The Direct Method and the characterization of weak lower semicontinuity through convexity go back to Tonelli \cite{Tonelli}; modern self-contained treatments, including the vectorial (quasiconvex) case only touched upon here, are given in the books by Dacorogna \cite{Dacorogna}, Giusti \cite{Giusti}, and Fonseca and Leoni \cite{FonsecaLeoni}. Relaxation and lower-semicontinuous envelopes are treated systematically in Ekeland and Temam \cite{ET} and, in the integral-functional setting used here, in the monograph by Buttazzo \cite{Buttazzo89}. The notion of $\Gamma$-convergence is due to De Giorgi and Franzoni \cite{DGF75}; the standard references are the books by Dal Maso \cite{DM93} and by the author of these notes \cite{Bra02}, to which we refer for all the properties of $\Gamma$-convergence used without proof below and in the rest of these notes (compactness with respect to $\Gamma$-convergence, stability under continuous perturbations, metrizability of $\Gamma$-convergence on separable spaces, etc.). The $H=W$ theorem is due to Meyers and Serrin \cite{MS64}. For the final example on sets of finite perimeter (whose systematic theory will only be used starting from Chapter~\ref{conc-sec}, see the Appendix to that Chapter) we refer to Ambrosio, Fusco and Pallara \cite{AFP} or to Maggi \cite{Maggi}. The Weierstrass and Bolza examples are classical; both, together with several variants, can be found for instance in the book by Buttazzo, Giaquinta and Hildebrandt \cite{BGH}, specifically devoted to one-dimensional problems as in most of this chapter; for the second, see also Bolza's own lecture notes \cite{Bolza}, although we have not been able to trace the precise original source of the specific form used here, and would be grateful for a more precise pointer.

\paragraph{A remark on suprema of lower-semicontinuous functionals.} We have twice used the elementary fact that {\em the pointwise supremum $F=\sup_{i\in I} F_i$ of an arbitrary family of lower-semicontinuous functionals $F_i\colon X\to\mathbb R\cup\{+\infty\}$ is itself lower semicontinuous}. This follows at once from the characterization of lower semicontinuity through sublevel sets: for every $c\in\mathbb R$,
$$
\{x\in X: F(x)\le c\}=\bigcap_{i\in I}\{x\in X: F_i(x)\le c\},
$$
an intersection of closed sets (since each $F_i$ is lower semicontinuous), hence closed; this is exactly the statement that $F$ is lower semicontinuous.

\paragraph{Proof of the Riemann--Lebesgue-type Lemma.} We restate and prove the lemma used in Section~\ref{semicont-sec} to derive the necessity of convexity for weak lower semicontinuity.

\begin{lemma}\label{RL-lemma} Let $1\le p<\infty$, let $\varphi\in L^p_{\rm loc}(\mathbb R^d;\mathbb R^m)$ be $1$-periodic (i.e.\ $\varphi(y+k)=\varphi(y)$ for a.e.\ $y$ and all $k\in\mathbb Z^d$), and let $\overline\varphi=\int_{(0,1)^d}\varphi(y)\,dy$ be its average. Then, for every bounded open set $\Omega\subset\mathbb R^d$, the functions $u_n(x)=\varphi(nx)$ converge weakly in $L^p(\Omega;\mathbb R^m)$ (weakly$^\ast$ in $L^\infty(\Omega;\mathbb R^m)$ if $p=\infty$) to the constant function $\overline\varphi$ as $n\to+\infty$, $n\in\mathbb N$.
\end{lemma}

\begin{proof}
Since $\varphi$ is $1$-periodic, for every $k\in\mathbb Z^d$ and every $g\in L^1_{\rm loc}(\mathbb R^d;\mathbb R^\ell)$ built from $\varphi$ (e.g.\ $g=\varphi$ or $g=|\varphi|^p$) we trivially have $\int_{(0,1)^d+k}g(y)\,dy=\int_{(0,1)^d}g(y)\,dy$, since the change of variables $y\mapsto y+k$ leaves both the domain of integration and, by periodicity, the integrand unchanged. We use this fact to compute the limit, as $n\to+\infty$, of $n^{-d}\int_{nB}g(y)\,dy$ for a box $B=\prod_{i=1}^d(a_i,b_i)\subset\mathbb R^d$ and $g$ as above: the box $nB$ is covered, up to a boundary layer of $O(n^{d-1})$ unit cells (as $n\to+\infty$), by the $N_n\sim n^d|B|$ unit cells $(0,1)^d+k$, $k\in\mathbb Z^d$, entirely contained in $nB$; on each such cell the integral of $g$ equals $\int_{(0,1)^d}g\,dy$ by the observation above, while the boundary cells contribute a total of at most $O(n^{d-1})\int_{(0,1)^d}|g|\,dy$. Hence
$$
\int_{nB}g(y)\,dy=N_n\int_{(0,1)^d}g(y)\,dy+O(n^{d-1}),\qquad\hbox{so that}\qquad n^{-d}\int_{nB}g(y)\,dy\ \longrightarrow\ |B|\int_{(0,1)^d}g(y)\,dy
$$
as $n\to+\infty$.

Applying this with $g=|\varphi|^p$ and $B$ any box containing $\overline\Omega$, and using the change of variables $y=nx$, we get
$$
\int_\Omega |u_n(x)|^p\,dx\le \int_B|\varphi(nx)|^p\,dx=n^{-d}\int_{nB}|\varphi(y)|^p\,dy\ \longrightarrow\ |B|\int_{(0,1)^d}|\varphi|^p\,dy<+\infty,
$$
so that $(u_n)$ is bounded in $L^p(\Omega;\mathbb R^m)$ (trivially so if $p=\infty$, since then $\|u_n\|_\infty=\|\varphi\|_\infty$ for all $n$).

Applying the same computation with $g=\varphi$ and $B$ an arbitrary box contained in $\Omega$, we get, again by the change of variables $y=nx$,
$$
\int_B u_n(x)\,dx=n^{-d}\int_{nB}\varphi(y)\,dy\ \longrightarrow\ |B|\,\overline\varphi,
$$
so that $u_n\wto\overline\varphi$ when tested against the indicator function of any box $B\subset\Omega$, and hence, by linearity, against any simple function on $\Omega$.

Since simple functions are dense in $L^{p'}(\Omega;\mathbb R^m)$ for $1\le p'<\infty$ (respectively, in $L^1(\Omega;\mathbb R^m)$, for the weak$^\ast$ statement with $p=\infty$), and $(u_n)$ is bounded in $L^p(\Omega;\mathbb R^m)$, a standard density argument (approximate a general test function $\psi$ by a simple function $\psi_\eta$ with $\|\psi-\psi_\eta\|_{p'}<\eta$, and use the uniform bound on $u_n$ to control the error) upgrades the convergence tested on simple functions to convergence tested on all of $L^{p'}(\Omega;\mathbb R^m)$ (resp.\ $L^1(\Omega;\mathbb R^m)$). This is exactly the statement $u_n\wto\overline\varphi$ in $L^p(\Omega;\mathbb R^m)$ (resp.\ weakly$^\ast$ in $L^\infty$).
\end{proof}

\begin{remark}\rm The same proof, with obvious changes (covering $B/\e$ instead of $nB$ by unit cells), shows that if $\varphi$ is $1$-periodic then $\varphi_\e(x):=\varphi(x/\e)$ converges weakly to $\overline\varphi$ as $\e\to0^+$, $\e$ ranging over all of $(0,+\infty)$ and not only over the reciprocals of integers; this is the form of the lemma used, e.g., in the discussion of oscillating coefficients in Section~\ref{hom-sec} and elsewhere.
\end{remark}

\paragraph{The Urysohn property of $\Gamma$-convergence.} We prove the equivalence stated in Section~\ref{Gconv-sec}. We recall that, for a fixed $x\in X$, the $\Gamma$-$\liminf$ and $\Gamma$-$\limsup$ of a sequence of functionals $G_j\colon X\to\mathbb R\cup\{+\infty\}$ at $x$ are defined by
$$
\Gamma\hbox{-}\liminf_j G_j(x)=\inf\Big\{\liminf_j G_j(x_j): x_j\to x\Big\},\qquad \Gamma\hbox{-}\limsup_j G_j(x)=\inf\Big\{\limsup_j G_j(x_j): x_j\to x\Big\},
$$
so that $G_j\to G$ at $x$ in the sense of $\Gamma$-convergence if and only if these two quantities coincide and equal $G(x)$. We will use two elementary facts.

\smallskip
$\bullet$ ({\em monotonicity along subsequences}) If $(j_k)$ is a subsequence of $(j)$, then, since for any fixed sequence $(y_j)$ the liminf (resp.\ limsup) along a subsequence of indices is always $\ge$ (resp.\ $\le$) the liminf (resp.\ limsup) along the full sequence, taking the infimum over $y_j\to x$ on both sides gives
$$
\Gamma\hbox{-}\liminf_j G_j(x)\ \le\ \Gamma\hbox{-}\liminf_k G_{j_k}(x), \qquad \Gamma\hbox{-}\limsup_k G_{j_k}(x)\ \le\ \Gamma\hbox{-}\limsup_j G_j(x).
$$

\smallskip
$\bullet$ ({\em attainment}) Since $X$ is metrizable, the infima defining $\Gamma\hbox{-}\liminf_j G_j(x)$ and $\Gamma\hbox{-}\limsup_j G_j(x)$ are always attained, by a standard diagonal argument: for a countable, decreasing (in $h$) local base $\{B_h(x)\}_{h\ge1}$ of neighbourhoods of $x$, choose $x_j^h\to x$ with $\liminf_j G_j(x_j^h)$ within $1/h$ of $\Gamma\hbox{-}\liminf_jG_j(x)$ (truncating if the latter is $-\infty$), and let $N_h\uparrow+\infty$ be such that $x_j^h\in B_h(x)$ for $j\ge N_h$; the diagonal sequence $\overline x_j:=x_j^h$ for $N_h\le j<N_{h+1}$ converges to $x$ and realizes the infimum. The same construction applies to $\Gamma\hbox{-}\limsup$.

\begin{proposition}[Urysohn property]\label{Urysohn-prop} Let $X$ be metrizable, let $(F_j)$ be a sequence of functionals on $X$, and let $F\colon X\to\mathbb R\cup\{+\infty\}$. Then $F_j$ $\Gamma$-converges to $F$ if and only if every subsequence $(F_{j_k})$ of $(F_j)$ admits a further subsequence $\Gamma$-converging to $F$.
\end{proposition}

\begin{proof}
If $F_j$ $\Gamma$-converges to $F$ then, by the monotonicity fact above applied to any subsequence $(j_k)$ (with the full sequence playing the role of the ``subsequence'' of itself), we get
$$
F(x)=\Gamma\hbox{-}\liminf_j F_j(x)\le \Gamma\hbox{-}\liminf_k F_{j_k}(x)\le \Gamma\hbox{-}\limsup_k F_{j_k}(x)\le \Gamma\hbox{-}\limsup_j F_j(x)=F(x),
$$
so that $(F_{j_k})$ itself (with no need to pass to a further subsequence) $\Gamma$-converges to $F$.

Conversely, suppose every subsequence of $(F_j)$ has a further subsequence $\Gamma$-converging to $F$, and fix $x\in X$. We show $\Gamma\hbox{-}\liminf_jF_j(x)=F(x)$ (the proof that $\Gamma\hbox{-}\limsup_jF_j(x)=F(x)$ is entirely symmetric, using recovery sequences in place of the sequence $\overline x_j$ below). Suppose, for contradiction, that $\Gamma\hbox{-}\liminf_jF_j(x)\neq F(x)$; since by the monotonicity fact (applied to the trivial subsequence $(j_k)=(j)$, using the hypothesis that some subsequence of the full sequence $\Gamma$-converges to $F$) we already know $\Gamma\hbox{-}\liminf_jF_j(x)\le F(x)$, we must have $\Gamma\hbox{-}\liminf_jF_j(x)<F(x)$. By the attainment fact, there exists $\overline x_j\to x$ realizing $\liminf_jF_j(\overline x_j)=\Gamma\hbox{-}\liminf_jF_j(x)<F(x)$; passing to a subsequence $(j_k)$ along which $F_{j_k}(\overline x_{j_k})$ converges to this same value (possible since it is the liminf of the full sequence $F_j(\overline x_j)$), the hypothesis provides a further subsequence $(j_{k_i})$ along which $F_{j_{k_i}}$ $\Gamma$-converges to $F$; the $\Gamma\hbox{-}\liminf$ inequality for this subsequence, applied to $\overline x_{j_{k_i}}\to x$, gives
$$
\Gamma\hbox{-}\liminf_jF_j(x)=\lim_i F_{j_{k_i}}(\overline x_{j_{k_i}})=\liminf_i F_{j_{k_i}}(\overline x_{j_{k_i}})\ge F(x),
$$
contradicting $\Gamma\hbox{-}\liminf_jF_j(x)<F(x)$.
\end{proof}

\paragraph{A remark on the $H=W$ theorem.} The inclusion $\overline F(u)\le\int_\Omega(|u|^2+|\nabla u|^2)dx$ for $u\in W^{1,2}(\Omega)$ used in the main text amounts exactly to the density of $C^1(\Omega)\cap W^{1,2}(\Omega)$ (in fact of $C^\infty(\Omega)\cap W^{1,2}(\Omega)$) in $W^{1,2}(\Omega)$ with respect to the $W^{1,2}$-norm; this holds on an arbitrary open set $\Omega$ (no regularity of $\partial\Omega$ is required) and is the content of the Meyers--Serrin theorem \cite{MS64}, proved by a clever partition-of-unity argument that mollifies $u$ at a scale depending on the distance to $\partial\Omega$, so that the local mollifications patch together into a single globally defined smooth approximation.

\vfill\eject

\section{Concentration on interfaces}\label{conc-sec}
Let $W\colon\mathbb R\to \mathbb R$ be a non-convex function, with superlinear growth, and consider 
the problem  
$$\min\Big\{\int_{\Omega}W(u)\, dx: \int_\Omega u\, dx=c|\Omega|\Big\},$$
where $\Omega$ is a bounded domain in $\mathbb R^d$ and $u\colon\Omega\to\mathbb R$. 
We can compare this problem with the convex-minimization problem 
$$\min\Big\{\int_{\Omega}W^{\ast\ast}(u)\, dx: \int_\Omega u\, dx=c|\Omega|\Big\},$$ 
whose minimizer is trivially given by the constant function $u=c$, which is then also a solution of the initial non-convex problem if $W(c)=W^{\ast\ast}(c)$. Otherwise, there exists an affine function $r(z)$ such that $W^{\ast\ast}(z)=r(z)$ in a neighbourhood of $c$. Noting that, if we replace $W(u)$ with $W(u)-r(u)$, minimizers do not change, it is not restrictive to suppose that  $W\geq 0$ and that $c$ lies between two consecutive zeroes, and, up to a change of variables, that these two zeroes are $0$ and $1$.  With these observations in mind, it follows immediately that all functions of the form $u=\chi_E$ with $|E|=c|\Omega|$ solve the minimum problem.   

In order to rule out minimizers with ``large boundaries'', a possible strategy is the one described in the previous chapter; that is, the choice of those solutions that are limits of minimizers of suitable approximating problems depending on a parameter $\e$. We then introduce a family of {\bf singularly perturbed} energies, which are referred to as {\bf Cahn--Hilliard free energies} in the theory of phase separation, and the corresponding 
minimum problems  
$$\min\Big\{\int_{\Omega}(W(u)+\e^2|\nabla u|^2)\, dx: \int_\Omega u\, dx=c|\Omega|\Big\}.$$ 
This choice is motivated by the Van der Waals conjecture that the sets corresponding to limits of such minimizers satisfy a {\bf minimal interface criterion}.   

We can consider the Cahn--Hilliard functionals, or, more precisely, the functionals  
$$
F_\e^0(u)=\begin{cases}
\displaystyle\int_{\Omega}(W(u)+\e^2|\nabla u|^2)\, dx & \hbox{\rm if } u\in H^1(\Omega)\\
+\infty & \hbox{\rm otherwise in } L^1(\Omega) 
\end{cases}
$$
in terms of $\Gamma$-convergence, neglecting for the moment the integral constraint. 
If $u$ is the characteristic function of a smooth set $E$, and $u_\e=\phi_\e\ast u$ with $\phi_\e$ a standard mollifier supported in $B_\e$, 
we can check that the corresponding $\Gamma$-limit with respect to the $L^1$ convergence is $0$. By approximation, it holds for any set $E$, so that $\Gamma$-convergence does not give information about the limits of the minimizers of the approximating functionals. 

In order to make this approach work, we have to find a meaningful {\bf scaling of the energies}; indeed, multiplying by a constant $\e$-dependent scaling factor does not affect the minimizers while it allows for a different form of the $\Gamma$-limit, hence restricting the limit minimum problem. 
We then consider the functionals $F^\alpha_\e$ defined by setting $F^\alpha_\e(u)=\e^{-\alpha}F^0_\e(u)$. Since for $\alpha<1$ the same computation performed for $\alpha=0$ gives a vanishing $\Gamma$-limit, we can conjecture that the meaningful scaling is obtained by dividing by $\e$.   
These arguments lead to the analysis in terms of $\Gamma$-convergence of the functionals 
\begin{equation}\label{MMf}
F_\e(u)=\begin{cases}
\displaystyle\int_{\Omega}\Big(\frac{1}{\e}W(u)+\e |\nabla u|^2\Big)\, dx & \hbox{\rm if } u\in H^1(\Omega)\\
+\infty & \hbox{\rm otherwise in } L^1(\Omega) 
\end{cases}
\end{equation}
as a model of concentration on interfaces.

\subsection{Minimizers as scaled phase transitions}  
The key argument in the computation of the $\Gamma$-limit of 
\eqref{MMf} is based on the analysis of optimal sequences close to a limit interface, showing that 
the transition essentially occurs on a layer of size of order $\e$ and follows a
scaled one-dimensional profile in the direction orthogonal to the interface. 
The $\Gamma$-convergence result when the dimension $d$ is equal to $1$ 
consists in the characterization of this optimal profile, while the general case is obtained by using 
the structure of the sets of finite perimeter and the one-dimensional result 
applied to one-dimensional sections. 

\smallskip 

In the computations, we assume that $W$ is a continuous function, and that $W(z)=0$ if and only if $z\in\{0,1\}$.

\paragraph{Compactness and $\Gamma$-convergence in one dimension.} 
It is not restrictive to fix $\Omega=(0,1)$. We prove that all sequences 
 $\{u_\e\}$ such that 
$F_\e(u_\e)\leq C<+\infty$ admit a subsequence strongly converging in $L^1$ to some $u$, and that such $u$ is the characteristic function of a finite union of intervals. 

Let then $\{u_\e\}$ be a sequence with equibounded energy. Since we reason up to subsequences, we do not change the notation and still use $u_\e$. 

We first note that the values of $u_\e$ are in measure close to the set $\{0,1\}$. Indeed, setting $S=\sup_\e F_\e(u_\e)$, 
from the bound $\int_{0}^1 W(u)\, dx\leq S\e$, we deduce that for each $\eta>0$ there exists a $c_\eta>0$ such that 
\begin{equation}\label{MMmeas}
|\{x: |u_\e(x)-1|>\eta\}\cap \{x: |u_\e(x)|>\eta\}|\leq c_\eta \e. 
\end{equation}
In particular, $c_\eta$ can be taken equal to $S (\min\{W(s): |s-1|\geq \eta \ \hbox{\rm and } |s|\geq \eta\})^{-1}$.
Next, we note that the norms $\|u_\e\|_\infty$ are equibounded. This can be proven by contradiction. Indeed, suppose that $\max u_\e=u_\e(x_\e)$ tends to $+\infty$. By \eqref{MMmeas}, there exists $y_\e$ such that $u_\e(y_\e)=2$. We may suppose $y_\e<x_\e$ and $u_\e(x)\geq 2$ in $[y_\e,x_\e]$, which implies $|x_\e-y_\e|\leq c_1\e$ by \eqref{MMmeas}.  
Then 
$$F_\e(u_\e)\geq \int_{y_\e}^{x_\e}\e (u^\prime_\e)^2 \, dx \geq   
\e \frac{(u_\e(x_\e)-u_\e(y_\e))^2}{|x_\e-y_\e|}\geq  \frac{(\|u_\e\|_\infty-2)^2}{c_1}.$$
This gives the contradiction. 

Now, we examine transitions from close to $0$ to close to $1$. Each such transition will cross a region where $W$ is bounded from below by a positive constant, say where $u_\e\in ({1\over 4},{3\over 4})$. Consider a {\em transition interval} $I=(x^-,x^+)$ such that $\{u_\e(x^-), u_\e(x^+)\}=\{{1\over 4},{3\over 4}\}$ and $u(x)\in ({1\over 4},{3\over 4})$ in $(x^-,x^+)$. Let now $\kappa= \min \{W(s): s\in ({1\over 4},{3\over 4})\}$ and $\delta=x^+-x^-$. We estimate  
\begin{equation}\label{MMinf14}\int_{x^-}^{x^+}\Big(\frac{W(u_\e)}{\e}+\e (u^\prime_\e)^2\Big) \, dx \geq
\kappa\frac{\delta}{\e}+\frac{\e}{4\delta}
\geq \sqrt{\kappa},\end{equation}
where the last inequality is obtained by optimizing in $\delta$. 
Hence, the number of such transition intervals is equibounded independently of $\e$. 

\begin{figure}[h!]
	\centerline{\includegraphics [width=3.4in]{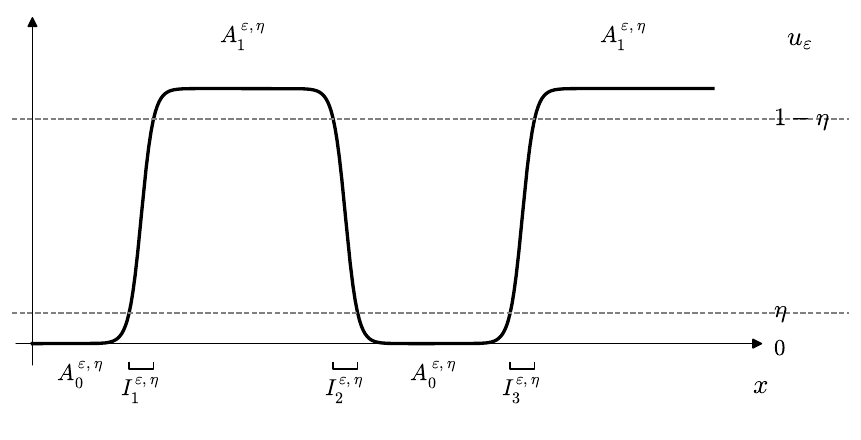}}
	\caption{The decomposition of $u_\e$ into transition intervals $I^{\e,\eta}_k$ and the complementary sets $A^{\e,\eta}_0$, $A^{\e,\eta}_1$.}\label{Fig13}
\end{figure}

Now, we fix $\eta\in(0,\frac{1}{4})$ and we estimate the energy of each transition between $\eta$ and $1-\eta$. 
With fixed $\e$, if now $I=(x^-,x^+)$ is a {\em transition interval between $\eta$ and $1-\eta$}; that is, such that $\{u_\e(x^-), u_\e(x^+)\}=\{\eta,1-\eta\}$ and $u(x)\in (\eta,1-\eta)$ in $(x^-,x^+)$,
we can apply inequality \eqref{MMinf14} to some subinterval of $(x^-,x^+)$, deducing that the number of such transition intervals is equibounded independently of $\e$ and $\eta$. 
This argument allows us to describe some structure of $u_\e$ by subdividing the domain in 
a finite number of (disjoint) transition intervals $I^{\e,\eta}_k=(x^{\e,\eta,-}_k, x^{\e,\eta,+}_k)$ for $k\in\{1,\dots, N\}$,  with $N$ independent (up to subsequences) of the parameters, and the complement. 
For each interval of the form $(x^{\e,\eta,+}_k,  x^{\e,\eta,-}_{k+1})$, the values of the function $u_\e$ are either strictly less than $1-\eta$ or strictly larger than $\eta$ at any point, and the same holds for the boundary intervals $(0,x^{\e,\eta,-}_1)$ and $(x^{\e,\eta,+}_N,1)$. Hence, the complement of the transition set can be again decomposed in a finite union of intervals $A^{\e,\eta}_0$ where $u_\e<1-\eta$ and a finite union of intervals $A^{\e,\eta}_1$ where $u_\e>\eta$ (see Fig.~\ref{Fig13}).

\begin{figure}[h!]
	\centerline{\includegraphics [width=5.6in]{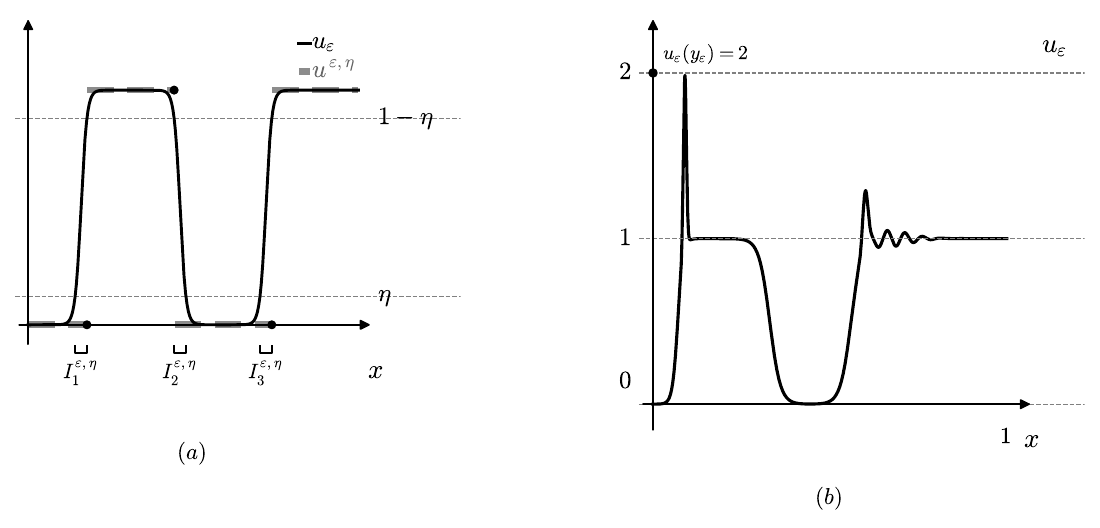}}
	\caption{(a) $u_\e$ (thin) together with the piecewise-constant $u^{\e,\eta}$ (thick, dashed) obtained from it: an asymmetric $\{0,1\}$-valued projection of $u_\e$, agreeing with it up to $\eta$ outside the transition intervals $I_k^{\e,\eta}$, but set to $0$ throughout each $I_k^{\e,\eta}$ regardless of the transition's direction. (b) A non-monotone competitor, overshooting above the level $2$ used in the equiboundedness argument, and oscillating before settling near a plateau: the same construction applies unchanged.}\label{Fig14}
\end{figure}

We can then define a sequence of piecewise-constant functions with $N$ discontinuities by setting
$$u^{\e,\eta}(x)=
\begin{cases}
1 & \hbox{\rm if } x\in A^{\e,\eta}_1\\
0 & \hbox{\rm if } x\in (0,1)\setminus A^{\e,\eta}_1.
\end{cases}
$$
Fig.~\ref{Fig14}(a) shows $u_\e$ together with the resulting $u^{\e,\eta}$, for the same function pictured in Fig.~\ref{Fig13}.

Since by \eqref{MMmeas} the measure of
 the union of the transition intervals is bounded by $c_\eta \e$, the choice of the value of $u^{\e,\eta}$ in this set does not affect the convergence of the sequence. Up to subsequences, by letting $\e\to 0$ and then $\eta\to 0$ we obtain a limit function $u$ valued in $\{0,1\}$ with at most $N$ discontinuities. Moreover, 
 $$\int_{0}^1|u^{\e,\eta}(x)-u{_\e}(x)|\, dx\leq \eta +\|u_\e\|_\infty c_\eta\e$$ 
 by \eqref{MMmeas}, so that by the equiboundedness of $\|u_\e\|_\infty$ we also obtain the convergence 
 $u_\e\to u$ in $L^1$, proving the claimed coerciveness. 
 
 \smallskip 
 
 To give a meaningful lower bound for the $\Gamma$-limit, we note that we can estimate the energy integrated on each transition interval between $\eta$ and $1-\eta$ by minimizing over the length of the interval, 
 after setting $u_\e(x)=v_\e(\frac{x}{\e})$, as 
 \begin{eqnarray*}
&&\int_{x^-}^{x^+}\Big(\frac{W(u_\e)}{\e}+\e (u^\prime_\e)^2\Big) \, dx
=\int_{\frac{x^-}{\e}}^{\frac{x^+}{\e}}
(W(v_\e)+(v^\prime_\e)^2) \, dt\\
\nonumber &\geq&   \inf_T
 \min\Big\{
 \int_{-T}^{T}(W(v)+|v'|^2)\, dt:
 v(-T)=\eta, v(T)=1-\eta\Big\}=c_W+o(1)    
  \end{eqnarray*}
as $\eta\to 0$, where
\begin{equation}\label{defcW}
c_W=\inf\Big\{ \int_{-\infty}^{+\infty}({W(v)}+|v'|^2)\, dt:
 \lim_{t\to-\infty}v(t)=0, \lim_{t\to+\infty}v(t)= 1\Big\}.
 \end{equation}
 This minimum problem is referred to as an {\bf optimal-profile problem}, describing the least-energy transition from $0$ to $1$.
Note that we have used the fact that $W(0)=W(1)=0$ to pass from conditions of the form $v(-T)=0$, $v(T)=1$ to asymptotic conditions at $\pm\infty$.
If $S(u)$ denotes the set of discontinuity points of the limit function $u$, since $\#S(u)\le N$, we then obtain the estimate 
 \begin{equation}\label{MMinf2} 
 \Gamma\hbox{\rm-}\liminf_{\e\to 0} F_\e(u)\geq  c_W\, \# S(u). 
 \end{equation}

In order to show that this estimate is sharp, we can use the optimal-profile problem in order to construct recovery sequences. If $u=\chi_{(x_0,+\infty)}$ with $x_0\in (0,1)$, for given $\delta>0$ we choose 
an admissible test function for the minimum problem in \eqref{defcW}
realizing the value $c_W$ up to an error $\delta$, and define 
$$u_\e(x)= v\big({x-x_0\over\e}\big).$$ 
Then, $u_\e\to u$ and $\limsup\limits_{\e\to 0} F_\e(u_\e)\le c_W+\delta$. If $S(u)=\{x_j\}_j$, the construction can be extended by locally using functions of the form $v\big({x-x_j\over\e}\big)$ and $v\big({x_j-x\over\e}\big)$.

We can then conclude that 
\begin{equation}\label{MMsup1} 
 \Gamma\hbox{\rm-}\lim_{\e\to 0} F_\e(u)=  
 \begin{cases}
 c_W\, \# S(u) &\hbox{\rm if } u \ \hbox{\rm piecewise constant and } u\in\{0,1\} \ \hbox{\rm a.e.}\\
 +\infty & \hbox{\rm otherwise in } L^1. 
 \end{cases}
 \end{equation}

\smallskip 

The construction of the recovery sequences, obtained by scaling optimal transitions, 
shows that the parameter $\e$ has a geometric interpretation as the typical size of the intervals where the energy of minimizers for problems involving $F_\e$ concentrates. 
In this sense we say that the energy is {\bf at scale $\e$}. 

%

\smallskip 

The coefficient $c_W$ can be expressed directly as 
\begin{equation}\label{cW}
c_W=2\int_0^1\sqrt{W(s)}\, ds.
\end{equation}
This can be seen as a consequence of the so-called {\bf Modica--Mortola trick} based on the simple inequality $W(v)+(v^\prime)^2\geq 2\sqrt{W(v)}|v^\prime|$, which gives the lower bound 
$$c_W\geq 2 \inf\Big\{ \int_{-\infty}^{+\infty}\sqrt{W(v)}|v'|\, dt:
 \lim_{t\to-\infty}v(t)=0, \lim_{t\to+\infty}v(t)= 1\Big\}=2\int_0^1\sqrt{W(s)}\, ds,$$ 
 achieved on increasing functions $v$, using the change of variable $s=v(t)$. 
The converse inequality is obtained by using as a test function a solution of
$$\begin{cases}
v^\prime(t)=\sqrt{W(v(t))}\\
v(0)=\frac{1}{2},
\end{cases}$$
for which $W(v)+(v^{\prime})^2=2\sqrt{W(v)}|v^\prime|$ (see Fig.~\ref{Fig12}; the existence of such a $v$, with the correct behaviour as $t\to\pm\infty$, is checked in the Appendix to this Chapter).

\begin{figure}[h!]
	\centerline{\includegraphics [width=4.6in]{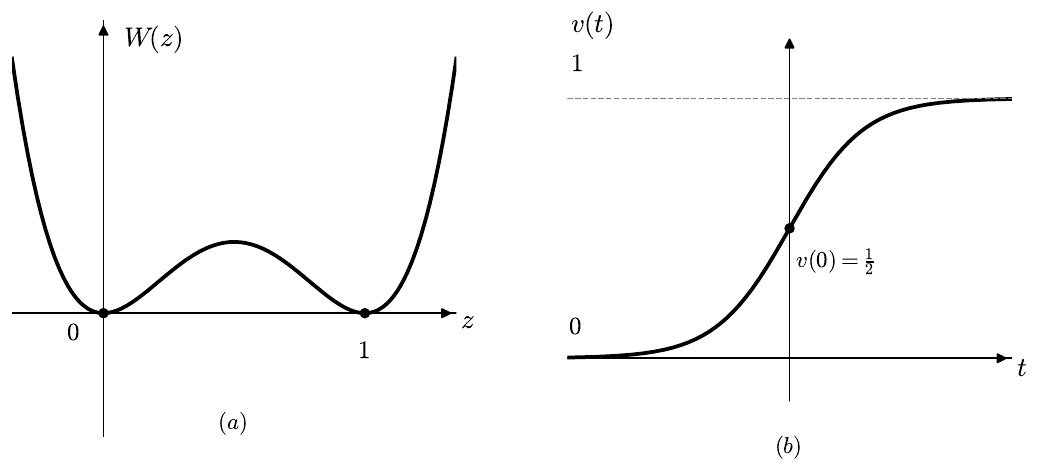}}
	\caption{(a) A double-well potential $W$, with wells at $0$ and $1$; (b) the corresponding optimal profile $v$.}\label{Fig12}
\end{figure}

\paragraph{Compactness and $\Gamma$-convergence in higher dimension.}  
We note that the precompactness of a sequence $\{u_\e\}$ with equibounded energy can be directly obtained by using the Modica--Mortola trick, which gives 
$$F_\e(u_\e)\geq 2\int_{\Omega}\sqrt{W(u_\e)}|\nabla u_\e|\, dx.$$
Setting $\Phi(s)=2\int_{0}^{s}\sqrt{W(\tau)}\, d\tau$, 
and $v_\e=\Phi(u_\e)$, we then have the boundedness of $v_\e$ in $W^{1,1}(\Omega)$. 
Hence, up to subsequences, $v_\e$ converge to some $v$ in $L^1$, up to translations. Since, again by \eqref{MMmeas}, $u_\e$ converge to $\{0,1\}$ in measure, we also deduce that $v_\e\to \{\Phi(0),\Phi(1)\}$ in measure.   
Since $\Phi(0)=0$, the limit of $v_\e$ is given by $\Phi(1)\chi_E$ with $E$ a subset of $\Omega$.   We then deduce that $u_\e$ converge in $L^1$ to $\chi_E$ since $u_\e=\Phi^{-1}(v_\e)$. 
Finally, we note that the boundedness of $v_\e$ in $W^{1,1}$ implies that $E$ is a set of finite perimeter.
\smallskip

\begin{figure}[h!]
	\centerline{\includegraphics [width=2.6in]{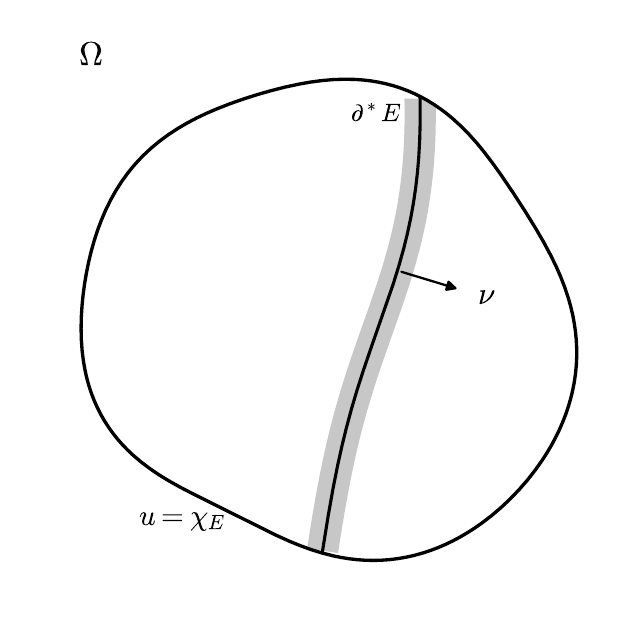}}
	\caption{The limit set $E\subset\Omega$ (so that $u_\e\to\chi_E$), with the $\e$-neighbourhood of its reduced boundary $\partial^\ast E$ and the outer normal $\nu$ used in the blow-up argument below.}\label{Fig15}
\end{figure}

We estimate the $\Gamma$-$\liminf$ by using the {\bf Fonseca--M\"uller blow-up method}, which essentially consists in interpreting the energy densities as densities of a sequence of weakly converging measures, and estimating the derivative of the limit measure with respect to the relevant fixed measure.

In this case, given a sequence $\{u_\e\}$ with equibounded energy and converging in $L^1$ to a set $E$ of finite perimeter, the relevant fixed measure is $\mathcal H^{d-1}\restr \partial^\ast E$, and the approximating sequence of measures is given by  
$$\mu_\e(A)=\int_{A}\Big(\frac{W(u_\e)}{\e}+\e|\nabla u_\e|^2\Big)\, dx.$$  
We can suppose that $\mu_\e$ weakly$^\ast$ converge to a bounded measure $\mu$. 
For $\mathcal H^{d-1}$-almost all $x_0\in \Omega\cap \partial E$, we have that
$$
    \frac{d\mu}{d\mathcal H^{d-1}\restr \partial^\ast E}(x_0)=\lim_{\rho\to 0}\frac{\mu(Q^\nu_\rho(x_0))}{\rho^{d-1}}=\lim_{\rho\to 0}\lim_{\e\to 0}\frac{\mu_\e(Q^\nu_\rho(x_0))}{\rho^{d-1}},
$$
where in the last equality we discard a countable set of $\rho$. By a diagonal argument and by the structure of sets of finite perimeter we may choose $\rho=\rho_\e>\!\!>\e$ such that
\begin{equation}
    \frac{d\mu}{d\mathcal H^{d-1}\restr \partial^\ast E}(x_0)=\lim_{\e\to 0}\frac{\mu_\e(Q^\nu_\rho(x_0))}{\rho^{d-1}}, 
\end{equation}
and $u_\e({y-x_0\over\rho})$ converge in $L^1(Q^\nu)$ to the characteristic function of $S^\nu$. We now use a Fubini argument. To that end, for each $y\in Q^\nu\cap\Pi^\nu$ we denote
$v^y_\e(t)=u_\e(x_0+y+t\rho\nu) $ defined for $t\in (-{1\over 2},{1\over 2})$ and note that $v^y_\e$ converge to $\chi_{(0,1/2)}$ in $L^1(-{1\over 2},{1\over 2})$ for almost all $y$. We then have
$$
{1\over \rho^{d-1}}\int_{Q^\nu_\rho(x_0)}\Big(\frac{W(u_\e)}{\e}+\e|\nabla u_\e|^2\Big)\, dx
\ge 
\int_{Q^\nu\cap\Pi^\nu}\int_{-1/2}^{1/2}
\Big(\frac{W(u^y_\e)}{\e/\rho}+{\e\over\rho}|(u^y_\e)'|^2\Big)\, dt\, d\mathcal H^{d-1}(y).
$$
By the one-dimensional result we have
$$\liminf_{\e\to0}\int_{-1/2}^{1/2}
\Big(\frac{W(u^y_\e)}{\e/\rho}+{\e\over\rho}|(u^y_\e)'|^2\Big)\, dt\ge c_W;$$
hence, by Fatou's Lemma, we obtain
$$\liminf_{\e\to 0}
{1\over \rho^{d-1}}\int_{Q^\nu_\rho(x_0)}\Big(\frac{W(u_\e)}{\e}+\e|\nabla u_\e|^2\Big)\, dx\ge c_W H^{d-1}(Q^\nu\cap\Pi^\nu)= c_W,
$$
and conclude the blow-up argument since
$$\liminf_{\e\to 0}F_\e(u_\e)=\liminf_{\e\to 0}\mu_\e(\Omega)\ge \mu(\Omega)\ge \int_{\Omega\cap \partial^\ast E} \frac{d\mu}{d\mathcal H^{d-1}\restr \partial^\ast E} d\mathcal H^{d-1}$$
$$\ge \int_{\Omega\cap \partial^\ast E} \frac{d\mu}{d\mathcal H^{d-1}\restr \partial^\ast E}c_W d\mathcal H^{d-1}
= c_W{\rm Per} (E,\Omega).
$$
This inequality is obtained by estimating from below the contribution of $d$-dimensional transition energies by the corresponding one-dimensional transition energies on sections in the direction of the normal to $E$. We now see that this argument is sharp, and allows us to construct recovery sequences. To that end, for fixed $T$, consider a minimizer $v_T$ of 
$$
\min\Big\{
 \int_{-T}^{T}(W(v)+|v'|^2)\, dt:
 v(-T)=0, v(T)=1\Big\},
 $$
 extended to $0$ if $t<-T$ and to $1$ if $t>T$.
 If $P$ is a polyhedron, we then consider the functions
$$u^T_\e(x)=v_T\Bigl(\frac{{\rm d}_P (x,\mathbb R^d\setminus P)}\e \Bigr),
$$
where $d_P$ is the signed distance from $\partial P$ defined by
$${\rm d}_P(x)=\begin{cases}
{\rm dist} (x,\mathbb R^d\setminus P) & \hbox{ if } x\in P\\
- {\rm dist} (x, P) & \hbox{ if } x\in \mathbb R^d\setminus P.
    \end{cases}
$$
Then $u^T_\e$ converge to $\chi_P$, and we have
$$\limsup_{\e\to0} F_\e(u^T_\e)\le c_W {\rm Per }(P,\Omega) +o(1)$$
as $T\to+\infty$, proving the upper bound if $P$ is a polyhedron. By the density of polyhedral sets in the family of sets of finite perimeter, the argument is complete.

Summarizing, in the higher-dimensional case we have seen that optimal sequences produce {\bf transition layers} at {\bf scale} $\e$  tending to an interface, that follow a {\bf one-dimensional profile} in the direction orthogonal to the interface itself.

\subsection{Internal and external boundary layers} 
Singular perturbations can also be considered for problems with boundary conditions; for example of the form
\begin{equation}\label{pert-bc-e}
\min\Big\{\int_0^1W(u)dt+\e^2\int_0^1|u'|^2dt: u(0)=z_0\Big\}.
\end{equation}
Consider a continuous non-negative potential $W\colon\mathbb R\to [0,+\infty)$, with a single minimum point  $\overline z$ such that $W(\overline z)=0$.
Then the minimizers of \eqref{pert-bc-e} converge in measure to the constant function $\overline u(x)=\overline z$.
As a consequence, the domain of the $\Gamma$-limit $F$ of 
\begin{equation}
F_\e(u)=\frac1\e \int_0^1W(u)dt+\e\int_0^1|u'|^2dt \qquad  u(0)=z_0
\end{equation}
consists only of the function $\overline u$. Arguing as for the phase-transition energies, we deduce that the value of $F$ on this function is given by the one-sided optimal-transition problem
\begin{equation}\label{pert-bc-of}
m(z_0)=\inf\Big\{\int_0^{+\infty}(W(v)+ |v'|^2)dt: v(0)=z_0, v(+\infty)=\overline z\Big\},
\end{equation}
and minimizing sequences {\bf concentrate on the boundary} at scale $\e$, their profile following a scaled {\bf boundary layer}; that is, $u_\e(x)= v(x/\e)$ with $v$ a minimizer of \eqref{pert-bc-of}. Note that this is a kind of double-well problem, where the well at $\overline z$ is due to minimization, while the other $z_0$ is imposed by
the boundary condition. 

The same argument can be used for two-sided boundary conditions, for which the functionals are
\begin{equation}\label{feuzo}
F_\e(u)=\frac1\e \int_0^1W(u)dt+\e\int_0^1|u'|^2dt \qquad  u(0)=z_0, u(1)=z_1,
\end{equation}
and the value of their $\Gamma$-limit at $\overline u$ is $m(z_0)+m(z_1)$.

This argument can be extended, for example, to the case of problems of the form  
\begin{equation}\label{pert-bc-e-}
\min\Big\{\int_0^1W(u)dt+\e^2\int_0^1|u'|^2dt: u(0)=z_0, \int_0^1 u\,dt=\overline z\Big\},
\end{equation}
where now $\overline z$ is not necessarily a minimum. If $W$ is strictly convex and $C^1$, then we consider
$$
F_\e(u)=\frac1\e \int_0^1\big(W(u)-W(\overline z)-W'(\overline z)(u-\overline z)\big)dt+\e\int_0^1|u'|^2dt 
$$
with the constraints $ u(0)=z_0$ and $ \int_0^1 u\,dt=\overline z$. Note that considering $\overline W(u)=W(u)-W(\overline z)-W'(\overline z)(u-\overline z)$ in place of $W(u)$ only adds a constant to the integrals in \eqref{pert-bc-e-}, and $\overline W$ has the single minimizer $\overline z$.
Hence, again the domain of the $\Gamma$-limit $F$ of $F_\e$ is the only function $\overline u$, and its value there 
is 
\begin{equation}\label{pert-bc-of-bar}
\overline m(z_0)=\inf\Big\{\int_0^{+\infty}(W(v)-W(\overline z)-W'(\overline z)(v-\overline z)+ |v'|^2)dt: v(0)=z_0, v(+\infty)=\overline z\Big\}.
\end{equation}

\begin{figure}[h!]
	\centerline{\includegraphics [width=5.2in]{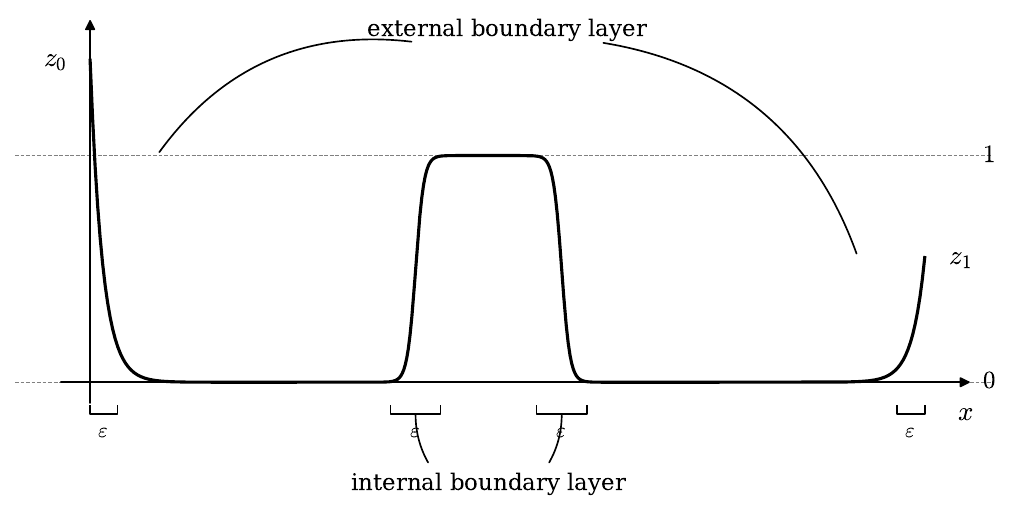}}
	\caption{A minimizer $u_\e$ with two internal jumps and boundary values $z_0,z_1$ not equal to $0$ or $1$. Near $x=0$ and $x=1$, external boundary layers of width $O(\e)$ connect $z_0$ and $z_1$ to the nearest well; the two internal jumps are realized, as in Figs.~\ref{Fig13}--\ref{Fig14}, by optimal-profile transitions between $0$ and $1$.}\label{Fig16}
\end{figure}

In the case of a double-well potential $W$ with minimum $0$ attained at $0$ and $1$, we can couple this analysis with the analysis of phase transitions. In this case the functionals \eqref{feuzo} $\Gamma$-converge to the sharp-transition functional
$$
F(u)=c_W \#S(u)+ m(z_0, u(0^+))+ m(z_1,u(1^-)),
$$
where $c_W$ is given by \eqref{defcW}, and $m(z,w)$ is defined for $w\in\{0,1\}$ and $z\in\mathbb R$ by
\begin{equation}\label{pert-bc-of-two}
m(z,w)=\inf\Big\{\int_0^{+\infty}(W(v)+ |v'|^2)dt: v(0)=z, v(+\infty)=w\Big\},
\end{equation}
and represents a boundary layer between the boundary value $z$ and the minimum state $w$.
Note that the minimization of $F$ is non-local since the behaviour at $0$ may influence the behaviour at $1$,
and minimizers exhibit both {\bf internal and external boundary layers} at scale $\e$, as pictured in Fig.~\ref{Fig16}.

\bigskip
Finally, we briefly remark that one can also consider a problem such as \eqref{pert-bc-e} in the general $d$-dimensional case. Suppose that $\Omega$ is a bounded $C^2$-domain and $\phi:\partial\Omega\to\mathbb R$ is sufficiently smooth (e.g., the trace of a $C^1$-function), and as above $W\colon\mathbb R\to [0,+\infty)$, with a single minimum point  $\overline z$ such that $W(\overline z)=0$. Then the $\Gamma$-limit $F$ of 
\begin{equation}
F_\e(u)=\frac1\e \int_\Omega W(u)dx+\e\int_\Omega|\nabla u|^2dx 
\end{equation}
with constraint $u=\phi$ on $\partial\Omega$ is again finite only on the function $\overline u$, where
$$
F(u)=\int_{\partial\Omega} m(\phi(x))d{\mathcal H}^{d-1}(x),
$$ 
where $m$ is given by \eqref{pert-bc-of}. In the case of a double-well potential, the additional term on the boundary depends on the integral of $m(\phi(x), {\rm Tr}u(x))$, where Tr$u$ is the trace of $u$ on $\partial\Omega$.

\subsection{Competing effects: concentration and uniformity}\label{sec-competing} 
We now examine the behaviour of minimum problems 
\begin{equation}\label{pert-bc-e-m}
m_\e(\lambda):=\min\Big\{\int_0^1W(u)dt+\e^2\int_0^1|u'|^2dt:  \int_0^1 u\,dt=\lambda\Big\},
\end{equation}
with respect to the parameter $\lambda$ in the case of a double-well potential $W$ with minimum $0$ achieved at $0$ and $1$, and such that $W^{\ast\ast}(z)=W(z)$ for $z\not\in(0,1)$. 
The $\Gamma$-limit analysis gives

\smallskip
$\bullet$ for $\lambda\not\in(0,1)$ we have $m_\e(\lambda)= W(\lambda)$;

$\bullet$ for $\lambda\in(0,1)$ we have  $m_\e(\lambda)=\e c_W+o(\e)$.

\begin{figure}[h!]
	\centerline{\includegraphics [width=3.6in]{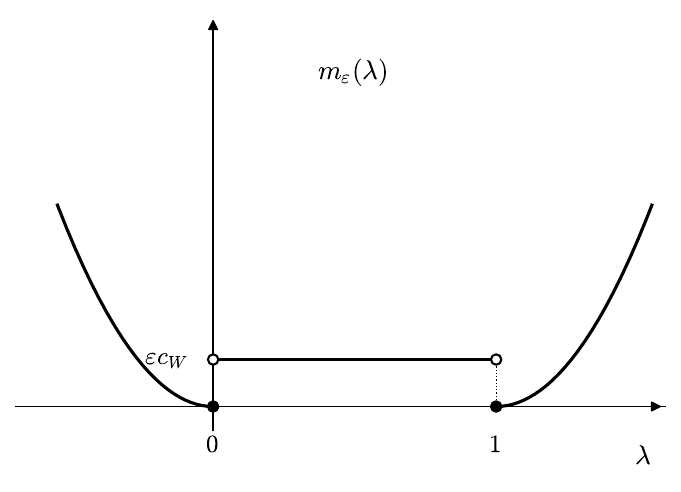}}
	\caption{The graph of $m_\e(\lambda)$ suggested by the two bullets above, for $\e>0$ fixed: it jumps from $0$, attained at $\lambda=0,1$, up to $\e c_W$ as soon as $\lambda$ enters $(0,1)$, and back down to $0$ at $\lambda=1$.}\label{Fig17}
\end{figure}

\smallskip
\noindent For $\e>0$ fixed, the function $m_\e$ is discontinuous at $0$ and $1$ (Fig.~\ref{Fig17}). This is in contrast with the definition of $m_\e$ in \eqref{pert-bc-e-m}.

\begin{figure}[h!]
	\centerline{\includegraphics [width=4.6in]{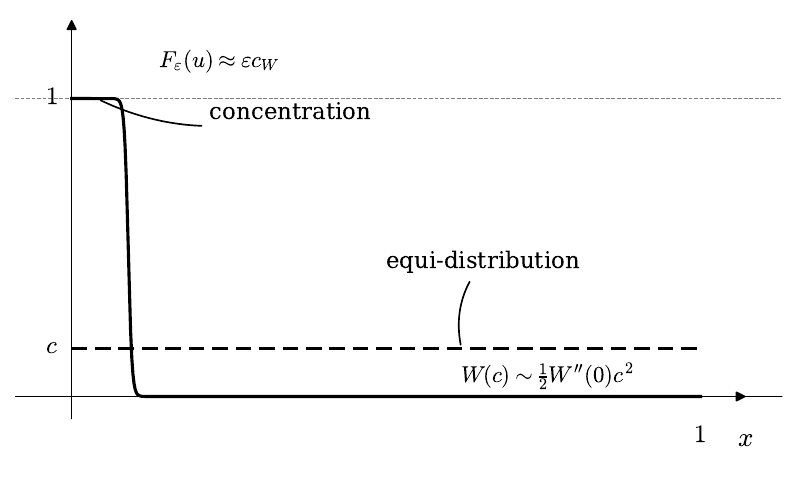}}
	\caption{The two competing constructions realizing a small mean value $\lambda=c$: concentration (a single internal transition, of energy $\approx\e c_W$ regardless of $c$) and equi-distribution ($u\equiv c$, of energy $W(c)\sim\frac12 W''(0)c^2$). Whichever is smaller determines $m_\e(\lambda)$ near $\lambda=0$, and symmetrically near $\lambda=1$.}\label{Fig18}
\end{figure}

Actually, the sharp-discontinuity energy must be compared with the energy of a uniform minimizer, which is $W(\lambda)$: for $\lambda$ close to $0$, say, concentrating the mass $\lambda$ on a single internal transition costs $\approx\e c_W$, independently of $\lambda$, while spreading it uniformly, with $u\equiv\lambda$, costs $W(\lambda)$, which is small precisely when $\lambda$ is small (Fig.~\ref{Fig18}).

\begin{figure}[h!]
	\centerline{\includegraphics [width=3.6in]{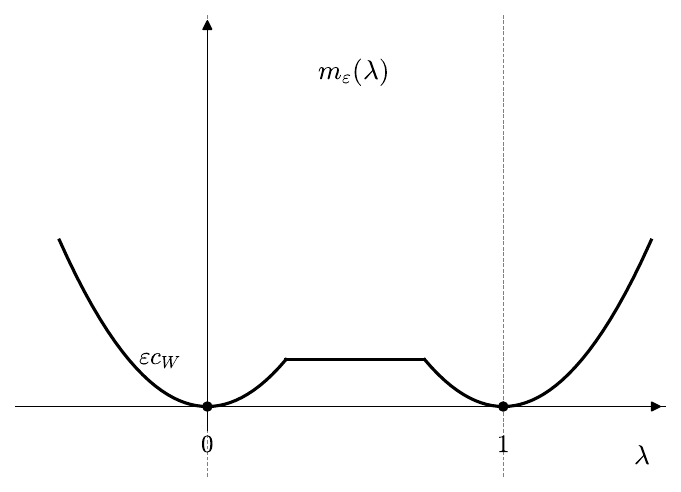}}
	\caption{The continuous graph of $m_\e(\lambda)$ given by \eqref{uniform}: equi-distribution wins close to $0$ and $1$, concentration (the plateau at $\e c_W$) wins away from them.}\label{Fig19}
\end{figure}

Hence, restricting our analysis to $\lambda\in [0,1]$ for simplicity, a better description of $m_\e$ is
$$
m_\e(\lambda)=\min\{W(\lambda), \e c_W\} +o(\e),
$$
or, if $W$ is twice differentiable at $0$ and $1$,
\begin{equation}\label{uniform}
m_\e(\lambda)=\min\Big\{\frac12W''(0)\lambda^2, \frac12W''(1)(1-\lambda)^2,\e c_W\Big\}+ o(\e),
\end{equation}
whose graph, continuous this time, is pictured in Fig.~\ref{Fig19}.

This description is sharp, in the sense that the functionals depending on two parameters $\e$ and $\lambda$
$$
F_{\e,\lambda}(u)= \frac1\e\int_0^1W(u)dt+\e\int_0^1|u'|^2dt \hbox{ with the constraint }  \int_0^1 u\,dt=\lambda
$$
$\Gamma$-converge as $(\e,\lambda)\to (0,0)$ to the functional finite only for $u=0$, on which its value is equal to $\min\Big\{\frac12W''(0)\lambda^2,\e c_W\Big\}$, and analogously for $(\e,\lambda)\to (0,1)$, while if $(\e,\lambda)\to (0,\lambda_0)$ with $\lambda_0\in(0,1)$ it converges to the usual sharp phase-transition functional. Then equality \eqref{uniform} can be obtained by a Weierstrass theorem along the sequence $F_{\e,\lambda}$.


\subsection{Appendix: relevant tools of Chapter 3}\label{App-Ch3}

\paragraph{References.} The variational, $\Gamma$-convergence formulation of the Van der Waals--Cahn--Hilliard theory of phase transitions is due to Modica and Mortola \cite{MM77}, who first proved the one-dimensional case of the $\Gamma$-convergence result \eqref{MMsup1}; the higher-dimensional case, together with the identification of the sharp-interface limit with the minimal-interface criterion conjectured by Van der Waals, is due to Modica \cite{Mod87}. The blow-up method used above to prove the $\Gamma$-liminf inequality in dimension $d>1$ is due to Fonseca and M\"uller \cite{FM92}; the general, pedagogically-oriented account on which the summary below is modeled is given by Braides, Maslennikov and Sigalotti \cite{BMS08}. For the systematic theory of sets of finite perimeter, only briefly touched upon in Section~\ref{Pre-Relax}, we refer again to Ambrosio, Fusco and Pallara \cite{AFP} or to Maggi \cite{Maggi}.

\paragraph{Sets of finite perimeter: the structure theorem.} Recall that a measurable set $E\subset\mathbb R^d$ has finite perimeter in $\Omega$ if the distributional gradient $D\chi_E$ of its characteristic function is a (vector-valued) Radon measure on $\Omega$ with finite total variation, this total variation being by definition $\mathrm{Per}(E;\Omega)$. De Giorgi's structure theorem (see e.g.\ \cite[Theorem~3.59]{AFP} or \cite[Chapter~15]{Maggi}) asserts that
$$
D\chi_E=\nu\,\mathcal H^{d-1}\restr\partial^\ast E,
$$
where $\partial^\ast E$ is the {\bf reduced boundary} of $E$; that is, the set of points $x$ at which the blow-ups $(E-x)/\rho$ converge, as $\rho\to0^+$, to the half-space $\{y: y\cdot\nu(x)<0\}$ for some unit vector $\nu(x)$ (the generalized outer normal). In particular,
$$
\mathrm{Per}(E;\Omega)=\mathcal H^{d-1}(\Omega\cap\partial^\ast E).
$$
The reduced boundary $\partial^\ast E$ coincides with the classical boundary $\partial E$ when $E$ is smooth, but for a general set of finite perimeter it is a proper, $\mathcal H^{d-1}$-measurable, subset of $\partial E$. It is exactly on $\Omega\cap\partial^\ast E$ that the blow-up argument used above applies: by definition of reduced boundary, $\mathcal H^{d-1}$-almost every $x_0\in\partial^\ast E$ is a point at which $E$ has an asymptotically flat (half-space) blow-up, which is what allows the $d$-dimensional transition energy to be estimated, on sections through $x_0$ in the direction $\nu(x_0)$, by the one-dimensional optimal-profile problem.

\paragraph{Existence of the optimal profile.} In the computation of $c_W$ above it is implicitly used that the infimum in \eqref{defcW} is attained, by the increasing solution of $v^\prime=\sqrt{W(v)}$, $v(0)=\frac12$, pictured in Fig.~\ref{Fig12}(b). We check that this is always the case, for an arbitrary continuous $W\colon\mathbb R\to[0,+\infty)$ with $W(z)=0$ if and only if $z\in\{0,1\}$ (no further regularity of $W$ is needed).

\begin{proposition}\label{opt-profile-exist}
The infimum $c_W$ in \eqref{defcW} is a minimum.
\end{proposition}
\begin{proof}
Since $W$ is continuous and strictly positive on $(0,1)$, the function
$$
g(z)=\int_{1/2}^z\frac{ds}{\sqrt{W(s)}}, \qquad z\in(0,1),
$$
is well defined, of class $C^1$, and strictly increasing, hence a bijection of $(0,1)$ onto an interval $(t^-,t^+)$ with $-\infty\le t^-<0<t^+\le+\infty$. Let $v_0=g^{-1}\colon(t^-,t^+)\to(0,1)$ be its inverse; then $v_0$ is strictly increasing, $v_0(0)=\frac12$, and, differentiating the identity $g(v_0(t))=t$,
$$
v_0^\prime(t)=\sqrt{W(v_0(t))}\qquad\hbox{\rm for every } t\in(t^-,t^+).
$$
Define $v\colon\mathbb R\to[0,1]$ by $v=v_0$ on $(t^-,t^+)$, and $v(t)=0$ for $t\le t^-$, $v(t)=1$ for $t\ge t^+$ (only relevant if $t^->-\infty$, respectively $t^+<+\infty$). Then $v$ is continuous, non-decreasing, locally absolutely continuous on $\mathbb R$, with $v(-\infty)=0$ and $v(+\infty)=1$; that is, $v$ is an admissible test function for \eqref{defcW}. Since $v^\prime=0$ a.e.\ on $\{v\in\{0,1\}\}$, where also $W(v)=0$, using $W(v)+(v^\prime)^2=2\sqrt{W(v)}\,v^\prime$ on $(t^-,t^+)$ (an equality, since there $v^\prime=\sqrt{W(v)}$) and the change of variable $s=v(t)$, we get
$$
\int_{-\infty}^{+\infty}(W(v)+|v^\prime|^2)\,dt=\int_{t^-}^{t^+}2\sqrt{W(v)}\,v^\prime\,dt=2\int_0^1\sqrt{W(s)}\,ds,
$$
which equals $c_W$ by \eqref{cW}. Hence $v$ is a minimizer.
\end{proof}

\begin{remark}\rm
If $\int_0^{1/2}(W(s))^{-1/2}ds=\int_{1/2}^1(W(s))^{-1/2}ds=+\infty$ --- as is the case, for instance, whenever $W$ is $C^2$ with $W^{\prime\prime}(0)>0$ and $W^{\prime\prime}(1)>0$, i.e.\ $W$ vanishes to second order at the two wells --- then $t^\pm=\pm\infty$, and the minimizer $v$ constructed above coincides with $v_0$ on all of $\mathbb R$, exactly as pictured in Fig.~\ref{Fig12}(b). If instead $W$ vanishes to first order at (say) $0$ --- a case compatible with the mere continuity assumed here, e.g.\ $W(z)\sim cz$ as $z\to0^+$ --- then $t^->-\infty$: the increasing solution of the ODE reaches the value $0$ in finite (reversed) time, and the minimizer constructed above is the one obtained by extension with the constant value $0$. In both cases $c_W$ is attained, with no change in its value.
\end{remark}

\paragraph{The blow-up method.} We summarize, in the abstract form of Fonseca and M\"uller \cite{FM92} (see also \cite{BMS08} for a general exposition, on which the scheme below is modeled), the blow-up technique used above to prove the $\Gamma$-liminf inequality in dimension $d>1$. The method applies whenever functionals $F_j$ can be written as set functions, $F_j(u)=\mu_j(\Omega)$, for set functions $\mu_j$ defined on the open subsets of $\Omega$ (here, $F_j=F_\e$ and $\mu_\e(A)=\int_A(\e^{-1}W(u_\e)+\e|\nabla u_\e|^2)\, dx$), and consists of the following steps.
\begin{enumerate}
\item[\em Step 1.] Extract a subsequence along which the $\mu_j$ weakly$^\ast$ converge, as (equibounded) measures, to some finite measure $\mu$ on $\Omega$; fix a reference positive measure $\lambda$, suited to the target function $u_0=\lim_ju_j$ (here, $\lambda=\mathcal H^{d-1}\restr(\Omega\cap\partial^\ast E)$), and write the Radon--Nikod\'ym decomposition $\mu=\frac{d\mu}{d\lambda}\lambda+\mu^s$.
\item[\em Step 2.] At a $\lambda$-a.e.\ point $x_0$ (here, $\mathcal H^{d-1}$-a.e.\ $x_0\in\Omega\cap\partial^\ast E$), use the differentiation theorem for measures to express $\frac{d\mu}{d\lambda}(x_0)$ as a limit of ratios $\mu(D_\rho(x_0))/\lambda(D_\rho(x_0))$ over a suitable family of shrinking sets $D_\rho(x_0)$ (here, cubes $Q^\nu_\rho(x_0)$ with a face orthogonal to $\nu=\nu(x_0)$).
\item[\em Step 3 (blow-up).] By a diagonal argument, choose $\rho=\rho_j\to0$, with $\rho_j\gg\e_j$, such that $\frac{d\mu}{d\lambda}(x_0)=\lim_j\mu_j(D_{\rho_j}(x_0))/\lambda(D_{\rho_j}(x_0))$. The rescaled functions converge to a model profile $v_0$ determined only by $u_0$ and $x_0$ (here, $v_0=\chi_{S^\nu}$, by the very definition of reduced boundary).
\item[\em Step 4 (local estimate).] Minimize the rescaled energy among competitors agreeing with $v_0$ on $\partial D$, reducing the local ratio to a minimum problem, independent of $j$ in the limit (a homogenization formula in general; here, by a Fubini argument on the one-dimensional sections in the direction $\nu$ and the one-dimensional result of this section, the relevant minimum problem gives exactly $c_W$).
\item[\em Step 5 (global estimate).] Integrate the local lower bound over $x_0$ with respect to $\lambda$, to obtain $\mu(\Omega)\ge\int\frac{d\mu}{d\lambda}\,d\lambda$, and conclude since $\liminf_jF_j(u_j)\ge\mu(\Omega)$ (here, $\mu(\Omega)\ge c_W\,\mathrm{Per}(E;\Omega)$).
\end{enumerate}
This is precisely the scheme followed, with the choices indicated above, in the proof of the $\Gamma$-liminf inequality in the higher-dimensional case.


\section{Concentration on singularities}  
The analysis of Van der Waals interactions with Dirichlet data shows that 
imposed boundary conditions 
may result in the appearance of an external boundary layer at the scale $\e$. This then affects the limit through a contribution concentrated on the boundary, which is, in that case, a $(d-1)$-dimensional set.  
Starting from this observation, in a sense we can revert the argument and impose boundary conditions on $\e$-dependent sets that allow for concentration effects on different types of sets, not necessarily limited to the boundary. 

\begin{figure}[h!]
	\centerline{\includegraphics [width=6.2in]{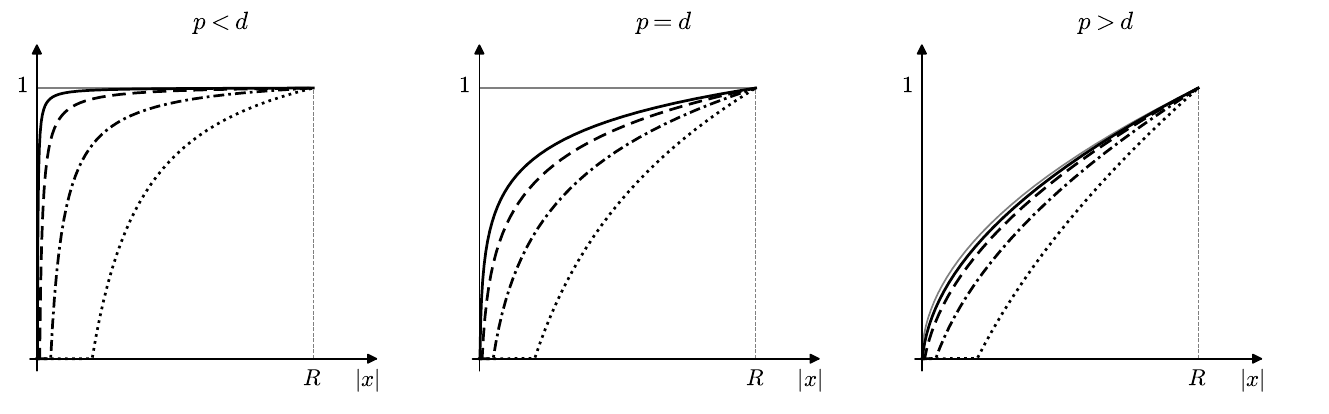}}
	\caption{Radial profiles of the minimizers $u_\e$, with $z=1$, extended by $0$ in $B_\e$, for $\e=0.2,\,0.05,\,0.01,\,0.002$ (dotted, dash-dotted, dashed and solid line, respectively) and $d=3,\,p=2$ ($p<d$), $d=p=2$ (formula \eqref{uelog}), $d=2,\,p=3$ ($p>d$). The thin gray line is the pointwise limit: the constant $1$ if $p\le d$, and $|x/R|^{\frac{p-d}{p-1}}$ if $p>d$.}\label{Fig20}
\end{figure}

To explain how the choice of these $\e$-dependent sets is reflected in the scale of the corresponding minima, and then in the limit, we look at a prototypical minimum problem, whose asymptotic behaviour can be explicitly computed.   
We consider $p>1$ and the minimum problem in $W^{1,p}(B_R\setminus \overline{B_\e})$ given by 
$$\min\Big\{\int_{B_R\setminus \overline{B_\e}}|\nabla u|^p\, dx: 
u=z \ \hbox{\rm on } \partial B_R \ \hbox{\rm and }   u=0 \ \hbox{\rm on } \partial B_\e \Big\},$$ 
where $B_r\subset \mathbb R^d$ is the open ball of radius $r$ centered at $0$, $z\in\mathbb R$ and $p>1$. 
The minimizer $u_\e$, computed from the Euler-Lagrange equation, is given by 
$$
u_\e(x)=z 
\displaystyle\frac{|x|^{\frac{p-d}{p-1}}-\e^{\frac{p-d}{p-1}}}{R^{\frac{p-d}{p-1}}-\e^{\frac{p-d}{p-1}}}
$$
if $p\neq d$, and by 
\begin{equation}\label{uelog}
u_\e(x)=z \displaystyle\frac{\log |x|-\log\e}{\log R-\log\e}
\end{equation}
in the case $p=d$. 
Note that, if $p\leq d$, then $u_\e$ tends to the constant function $z$; if $p>d$, on the contrary, the limit of $u_\e$ is the function $z\big|\frac{x}{R}\big|^{\frac{p-d}{p-1}}$ (see Fig.~\ref{Fig20}, where $z=1$).

The corresponding minimum values are 
$$m_\e=
\begin{cases} 
\omega_d |z|^p\Big(\frac{p-1}{p-d} \big(R^{\frac{p-d}{p-1}}-\e^{\frac{p-d}{p-1}}\big)\Big)^{1-p}
& \hbox{\rm if } p\neq d\\
\omega_d |z|^p (\log R-\log\e)^{1-d}& \hbox{\rm if } p=d,\end{cases}
$$   
and therefore they tend to $0$ if $p\leq d$. 

In order to explain the singular behaviour for $p\leq d$, for each fixed $p$ we can embed all minimum problems in the common space $W^{1,p}(B_R)$, independent of $\e$, by extending all test functions to $0$ in $B_\e$. 
Since the condition $u=0$ is imposed on a set shrinking to a point, this corresponds to considering as $\e$ tends to $0$ a minimum problem with the condition $u(0)=0$, which can be neglected if $p\leq d$. 
The precise computation of the minima $m_\e$ 
allows one to focus on the effects of the perturbation; that is, on the relevant $\e$-dependent terms in the approximation.  
We note that the minimum is positively homogeneous of degree $p$ in $z$, so that, after the change of variables $x= \e y$, 
it equals  
\begin{equation}\label{minpd}
|z|^p \e^{d-p} \ \min\Big\{\int_{B_{R/\e}\setminus \overline{B_1}}|\nabla u|^p\, dy: 
u=1 \ \hbox{\rm on } \partial B_{R/\e} \ \hbox{\rm and }   u=0 \ \hbox{\rm on } \partial B_1 \Big\},
\end{equation}
where $u\in W^{1,p}(B_{R/\e}\setminus \overline{B_1})$. This last minimum is the so-called {\bf $p$-capacity}\index{capacity} of $B_1$ in $B_{R/\e}$. 

\smallskip
If $p<d$, these minima decreasingly converge to the strictly positive value 
$$\min\Big\{\int_{\mathbb R^d}|\nabla u|^p\, dy: u\in W^{1,p}(\mathbb R^d) \ \hbox{\rm with } 
u=1 \ \hbox{\rm on } B_1 \Big\}=:\hbox{\rm Cap}_p(B_1;\mathbb R^d),$$ 
the $p$-capacity of $B_1$ in $\mathbb R^d$.   
Then, the energy effect of the perturbation of the boundary is at a scale $\e^{d-p}$, and it is concentrated on the point $x=0$, in the sense that it does not depend on $R$.  The minimizers of the problem $m_\e$ are approximately described as $u_\e(x)= 
v({x\over\e})$, where $v$ is the solution of the minimum problem giving $\hbox{\rm Cap}_p(B_1;\mathbb R^d)$.

If $p=d$, the problems are scale-invariant; that is, the values are the same after the change of variables, and the information is carried by the capacitary term, depending on $\log\big(\frac{R}{\e}\big)$. This highlights that, contrary to the case $p<d$, we may have a dependence on $R$ even if $R>\!>\e$.
In this sense, the effect of $\e$ is still singular, but it is not concentrated at a point. 

For our purposes, this example is prototypical in two directions. 
\begin{enumerate}
\item[$\bullet$] Choice of the boundaries where Dirichlet conditions are imposed, and of their scales, to recover asymptotic effects of order $1$, competing with the integral term. 
\item[$\bullet$] Analysis of scale-invariant minimum problems at a scale where the arising singularities drive the behaviour of minimizers. In this case, a vectorial setting allows us to substitute Dirichlet boundary conditions on shrinking sets with compatibility requirements of topological nature. 
\end{enumerate}

\subsection{Strange terms}\label{strangeterms} 
As suggested by the described example, we generalize the $\e$-dependent domain $B_R\setminus B_\e$ by considering a general set $\Omega_\e$, on part of whose boundary, corresponding to $\partial B_\e$, we can impose 
homogeneous Dirichlet boundary conditions. 
In order to give a heuristic argument of how we can choose these boundaries to produce a limit effect of order $1$, we can consider them as composed of a union of disjoint $(d-1)$-dimensional spheres with a common radius $\e$. 
By suitably increasing the number of these sets as $\e\to 0$, the overall limit energy concentrating on the centers of the spheres can be chosen to be of order $1$. 

\medskip\noindent{\bf Periodic arrays of boundaries.} 
As a first example, we consider a fixed bounded open set $\Omega\subset\mathbb R^d$ and define $\Omega_\e$ 
as the set obtained by removing from $\Omega$ an array of closed balls of radius $\e$ 
centered at the points in the cubic lattice $\delta \mathbb Z^d\cap\Omega$. The Dirichlet boundary condition $u=0$ is imposed on the union of the boundaries of these balls. 
Since the number of such balls is proportional to $\delta^{-d}$, by recalling the prototypical example (hence, in particular, with $p<d$) we expect a contribution of order $1$ if $\e^{d-p}\delta^{-d}$ is of order $1$. This suggests choosing the scale $\delta$ such that 
$\delta\sim \e^{\frac{d-p}{d}}$. 

With the choice $\delta_\e=\e^{\frac{d-p}{d}}$, we then consider the functionals $F_\e\colon W^{1,p}(\Omega)\to [0,+\infty]$ defined by 
$$F_\e(u)=\int_{\Omega}|\nabla u|^p\, dx \quad \hbox{\rm if }\ u=0 \ \hbox{\rm in } \bigcup_{k\in\mathbb Z^d}B_\e(\delta_\e k)$$ 
and $F_\e(u)=+\infty$ otherwise in $W^{1,p}(\Omega)$. 
For simplicity, we do not impose boundary conditions on the remaining part of $\partial\Omega$ since it may intersect
the boundary of the balls, thus  
they may be in contrast with the given homogeneous Dirichlet conditions. This choice makes the minimum problems for $F_\e$ trivial, so, to have a relevant application of the study of the asymptotic behaviour of $F_\e$, we can for instance think of minimum problems of the form 
$$\min\Big\{\int_{\Omega}\big(|\nabla u(x)|^p-g(x)u(x)\big)\, dx: u=0 \ \hbox{\rm in } \bigcup_{k\in\mathbb Z^d}B_\e(\delta_\e k) \Big\},$$ 
where $g\in L^{p^\prime}(\Omega)$. 
Note that the addition of the last term is compatible with the $\Gamma$-convergence of $F_\e$ since it is 
continuous with respect to the strong $L^p$-convergence. 

To heuristically understand the form of the $\Gamma$-limit of $F_\e$ (with respect to the strong $L^p$-convergence), 
we first assume that the sequence $\{u_\e\}$ converges to a continuous $u$, and suppose that $u_\e=u(k\delta_\e)$ on the boundary of each scaled ball with center $\delta_\e k$ and radius $T \delta_\e$ for a small enough fixed $T$. By applying the previous considerations for the case of one single ball, we obtain a global contribution of 
$$\hbox{\rm Cap}_p(B_1;\mathbb R^d) \sum_{\delta_\e k\in\delta_\e \mathbb Z^d\cap \Omega}  \delta_\e^d |u(k\delta_\e)|^p=\hbox{\rm Cap}_p(B_1;\mathbb R^d) \int_{\Omega}|u|^p\, dx +o(1).$$ 
Noting that the constant extensions of $u_\e$ to each ball $B_\e(k\delta_\e)$ with the corresponding constant value still converge to the same $u$, by the lower semicontinuity of the gradient part (which contributes outside the perforations), we can conjecture an expression of the form 
$$F(u)=\int_\Omega |\nabla u|^p\,dx + \hbox{\rm Cap}_p(B_1;\mathbb R^d) \int_{\Omega}|u|^p\, dx. $$
for the $\Gamma$-limit of $F_\e$ on the whole space $W^{1,p}(\Omega)$. 
The additional lower order term appearing in the $\Gamma$-limit has been poetically called a {\em strange term coming from nowhere}, or more philologically {\em un terme \'etrange venu d'ailleur}.  

\begin{figure}[h!]
	\centerline{\includegraphics [width=6.2in]{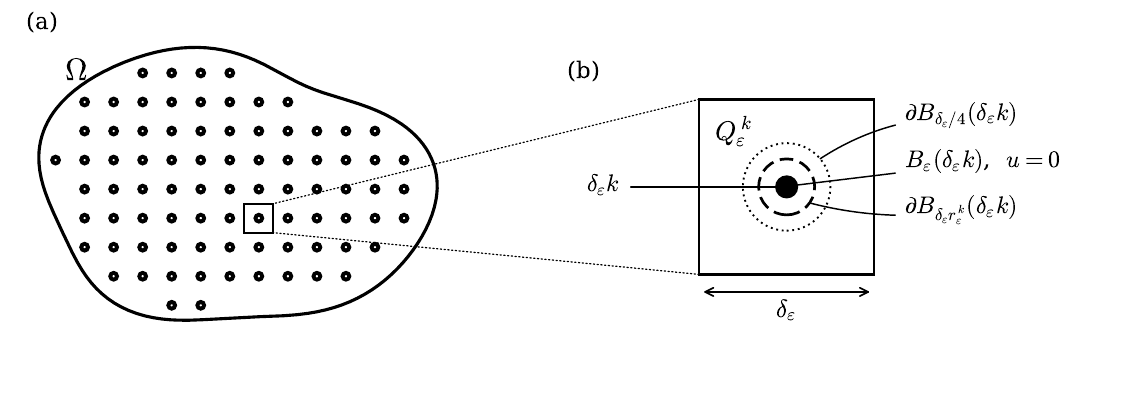}}
	\caption{The geometry of the periodic array of perforations for $d=2$. (a) The balls $B_\e(\delta_\e k)$, where $u=0$, are centred at the points of the lattice $\delta_\e\mathbf Z^2\cap\Omega$; the radius $\e$ is much smaller than the period $\delta_\e$. (b) A single cell $Q_\e^k$ of side $\delta_\e$: the ball $B_\e(\delta_\e k)$ (black), the sphere $\partial B_{\delta_\e r_\e^k}(\delta_\e k)$ with $r_\e^k\in(0,\frac14)$ (dashed), on which $u_\e$ is compared with a constant in the lower bound, and the sphere $\partial B_{\delta_\e/4}(\delta_\e k)$ (dotted).}\label{Fig21}
\end{figure}

We prove the lower estimate for the $\Gamma$-limit of $F_\e$ by making this heuristic argument rigorous; the geometry is pictured in Fig.~\ref{Fig21}. 
If $u_\e\weak u$, we can suppose that 
\begin{enumerate} 
\item[{\rm (i)}] for any $k\in\mathbb Z^d$ there exist $r_\e^k\in(0,\frac{1}{4})$ and a constant $u_{\e,k}$ such that 
$u_{\e}=u_{\e,k}$ on $\partial B_{\delta_\e r_\e^k}(\delta_\e k)$ in the sense of traces;  
\item[{\rm (ii)}] $\lim\limits_{\e\to 0}\delta_\e^d \displaystyle\sum_{k\in\mathbb Z^d}|u_{\e,k}|^p=\displaystyle\int_{\Omega}|u|^p\, dx$.  
\end{enumerate} 
This follows by a {\bf changing of boundary data} argument applied on each suitable sphere 
$\partial B_{\delta_\e r_\e^k}(\delta_\e k)$, taking into account that the distance of $u_\e$ from a constant can be controlled by a Poincar\'e-type inequality on sets at a scale $\delta_\e$. The details of this argument are given in the Appendix to this chapter, since they are of independent interest and will often be used in asymptotic problems involving Dirichlet conditions. 
Setting $P_\e=\Omega\setminus\bigcup_{k\in\mathbb Z^d} B_{\delta_\e r_\e^k}(\delta_\e k)$, 
the lower estimate can be written as 
\begin{eqnarray*}\liminf_{\e\to 0}F_\e(u_\e)&\geq&\liminf_{\e\to 0}\int_{P_\e}|\nabla u_\e|^p\, dx+\liminf_{\e\to 0}\sum_{Q_\e^k\subset\Omega}\int_{B_{\delta_\e r_\e^k}(\delta_\e k)}|\nabla u_\e|^p\, dx\\
&\geq&\liminf_{\e\to 0}\int_{P_\e}|\nabla u_\e|^p\, dx
+\liminf_{\e\to 0}\sum_{Q_\e^k\subset\Omega}\min\Big\{\int_{B_{\delta_\e r_\e^k}(\delta_\e k)} |\nabla v|^p\, dx: 
v=0 \ \hbox{\rm in } B_{\e}(\delta_\e k)\Big\}\\
&\geq&\liminf_{\e\to 0}\int_{\Omega}|\nabla w_\e|^p\, dx+\liminf_{\e\to 0}\sum_{Q_\e^k\subset\Omega}
\delta_\e^d |u_{\e,k}|^p \hbox{\rm Cap}_p(B_1; \mathbb R^d), 
\end{eqnarray*}
where $w_\e$ equals $u_\e$ on $P_\e$ and the constant $u_{\e,k}$ in each ball $B_{\delta_\e r_\e^k}(\delta_\e k)$. Note that $w_\e\weak u$ in $W^{1,p}(\Omega)$; hence, we can pass to the limit by using the lower semicontinuity of the first integral term, and by condition (ii), concluding the proof of the lower bound. 

To prove the upper estimate, we fix $u\in C^1(\Omega)$ and $\eta>0$. 
Let $T>0$ and $\varphi\in W^{1,p}(B_T)$ be a test function such that 
$$\int_{B_T}|\nabla \varphi|^p\, dx\leq  (1+\eta) \hbox{\rm Cap}_p(B_1;\mathbb R^d)$$ 
with $\varphi=1$ on $\partial B_T$ and $\varphi=0$ in $B_1$. 
We define $\varphi_\e$ by setting  
$$\varphi_\e(x)=\varphi\Big(\frac{x-\delta_\e k}{\e}\Big) \quad \hbox{\rm in } \ B_{\e T}(\delta_\e k),$$
and $\varphi_\e=1$ outside the union of such balls. 
\begin{figure}[h!]
	\centerline{\includegraphics [width=6.2in]{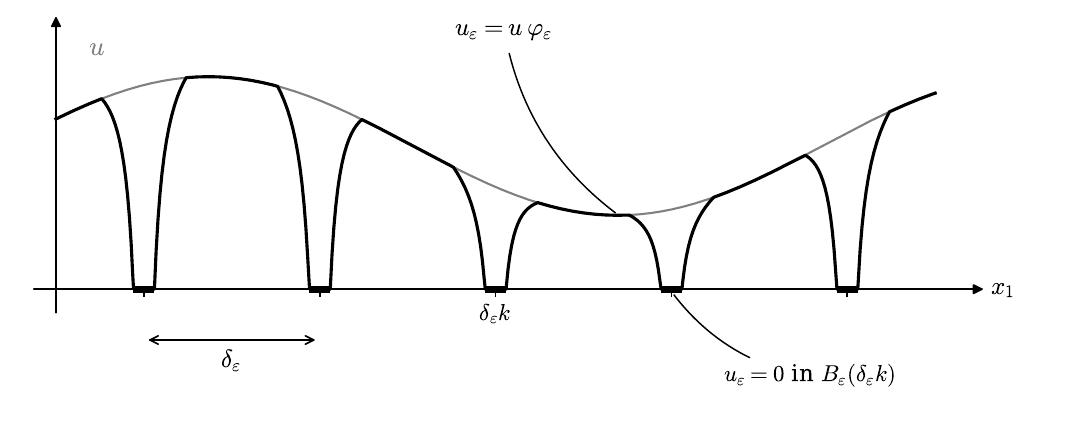}}
	\caption{A section, along a line through the centres of a row of balls, of the recovery sequence $u_\e=u\varphi_\e$ for a non-constant function $u$ (gray): $u_\e$ coincides with $u$ away from the balls, vanishes in each $B_\e(\delta_\e k)$, and in between follows the capacitary potential $\varphi_\e$ rescaled at scale $\e$ around each centre. The ratio $\e/\delta_\e$ is exaggerated.}\label{Fig23}
\end{figure}
Choosing $u_\e=u \varphi_\e$ (see Fig.~\ref{Fig23}), 
on each ball we can estimate  
\begin{eqnarray*}
\int_{B_{\e T}(\delta_\e k)}|\nabla u_\e|^p\, dx&\leq& 
\int_{B_{\e T}(\delta_\e k)}|\nabla u \varphi_\e+u\nabla\varphi_\e|^p\, dx\\
&\leq&\frac{C}{\eta}\int_{B_{\e T}(\delta_\e k)}|\nabla u|^p\, dx+ (1+\eta)\int_{B_{\e T}(\delta_\e k)}|u \nabla \varphi_\e|^p\, dx\\
&\leq&\frac{C}{\eta}\e^d T^d\|\nabla u\|^p_\infty+  (1+\eta) |u(\delta_\e k)+\omega(\e T)|^p \int_{B_{\e T}(\delta_\e k)}|\nabla \varphi_\e|^p\, dx\\
&\leq&\frac{C}{\eta}\e^d T^d\|\nabla u\|^p_\infty+  (1+\eta)^2 |u(\delta_\e k)+\omega(\e T)|^p \e^{d-p} \hbox{\rm Cap}_p(B_1;\mathbb R^d), 
\end{eqnarray*}
where $\omega$ is the modulus of continuity of $u$. 
By summing up over $k$, we then obtain 
\begin{eqnarray*}
F_\e(u_\e)&\leq&\int_{\Omega}|\nabla u|^p\, dx+  
\frac{C}{\eta}\Big(\frac{\e}{\delta_\e}\Big)^d T^d  \|\nabla u\|^p_\infty\\
&&+  
(1+\eta)^2 \Big(\int_{\Omega}|u|^p\, dx +o(1)_{\e\to 0}\Big)\hbox{\rm Cap}_p(B_1;\mathbb R^d). 
\end{eqnarray*}
Since $\e=o(\delta_\e)$, the upper estimate follows by the arbitrariness of $\eta>0$. 
We extend the result to a general $u\in W^{1,p}(\Omega)$ by density. The recovery sequence highlights concentration at scale $\e$, {\bf averaged at scale} $\e^{(d-p)/d}$.

\paragraph{Concentration ``at all $\e$-scales''.}
We now go back to the case $p=d$, where the behaviour is singular. In particular, the power of $\e$ disappears in \eqref{minpd}. Plugging the solution $u_\e$ described in \eqref{uelog} in \eqref{minpd}, we then obtain
\begin{equation}\label{mindd}
\min\Big\{\int_{B_{R/\e}\setminus \overline{B_1}}|\nabla u|^d\, dy: 
u=1 \ \hbox{\rm on } \partial B_{R/\e} \ \hbox{\rm and }   u=0 \ \hbox{\rm on } \partial B_1 \Big\}=
\omega_d |\log(R/\e)|^{1-d},
\end{equation}
where $\omega_d=\mathcal H^{d-1}(S^{d-1})$ is the $d-1$-dimensional measure of the unit sphere in $\mathbb R^d$ and results from integration in polar coordinates.

We can repeat the arguments in the previous case $p<d$ looking for a period $\delta>\!>\e$ for a periodic arrangement of balls of radius $\e$ where to impose $0$ boundary conditions. Using the computation \eqref{minpd} with $R=\delta/4$, so that the contribution of the balls can be decoupled, the resulting relation is $\delta^d= |\log (\delta/4\e)|^{1-d}$. Neglecting the asymptotically negligible contributions, we can choose
$$
\delta=|\log \e|^{\frac{1-d}d}.
$$
We do not repeat the computation of the $\Gamma$-limit, which gives an analogous result with the ``strange term'' $\int_\Omega|u|^ddx$. As far as scales are concerned, we note that recovery sequences for a function $u$ are constructed of the form $u\varphi_\e$, with $\varphi_\e\colon \mathbb R^d\to [0,1]$ $\delta$-periodic, and 
\begin{equation}\label{fep}
\varphi_\e(x)=\displaystyle\frac{\log |x|-\log\e}{\log c\delta-\log\e}\wedge 1
\end{equation}
on $[0,1]^d$, where $c>0$ is a small constant. The relevant scaling property of such a function is that for all $0<\alpha<\beta<1$ we have
$$
|\log\e|^{d-1}\int^{\e^\alpha}_{\e^\beta}|\nabla \varphi_\e|^2dx \sim \omega_d(\beta-\alpha);
$$
that is, the energy at scales between $\e^\beta$ and $\e^{\alpha}$ is a given portion of the total energy concentrating at the perforation. Hence, we may say that the energy {\bf concentrates at all scales $\e^\alpha$} around the perforation, {\bf modulated at scale} $|\log \e|^{\frac{1-d}d}$.

\subsection{Topological singularities} 
In order to introduce a different way in which singularities may arise, and the relevant corresponding scales,  
we return to problems of the form 
\begin{equation}\label{minW}
\min\Big\{\int_{\Omega}\Big(W(u)+\e^2|\nabla u|^2\Big)\, dx: u=\varphi \ \hbox{\rm on } \partial \Omega\Big\},
\end{equation} 
where $u\in H^1(\Omega)$.
If $W$ has a unique minimum point, the minimizers of these problems converge to a constant function possibly developing external boundary layers. 
In this chapter, a similar behaviour is that of minimizers of singular capacitary problems, which again converge to a constant function, given by an imposed  Dirichlet datum, except for a boundary layer concentrated at a singularity.   

We now consider a case in which the boundary datum $\varphi$ pointwise minimizes $W$ but forces a singularity in $\Omega$ similar to that in the capacitary problem. For such a phenomenon to be possible, we have to consider a vectorial setting; that is, in the simplest case, $\Omega$ is a bounded and simply connected domain in $\mathbb R^2$ and $u\in H^1(\Omega;\mathbb R^2)$, with the function $W\colon\mathbb R^2\to\mathbb R$ attaining its minimum on the set $S^1$. We fix, for simplicity, 
$$W(u)=(|u|^2-1)^2.$$ 
We choose a smooth function $\varphi\colon\partial\Omega\to S^1$ as a boundary datum, and consider the minimum problem in \eqref{minW}. 
To briefly explain the topological nature of the problem, we first assume that a minimizing sequence exists such that it converges to a continuous $u$ such that $W(u)=0$ pointwise and the boundary condition is satisfied. Then, $u\in C^0(\overline\Omega; S^1)$, which is possible only if the degree of $u$ (that is, of the datum $\varphi$) on the boundary of $\Omega$ is $0$, since $\Omega$ is simply connected.  We recall that, loosely speaking, the degree is the number of times $u$ winds around $S^1$. Writing $\varphi$ in complex notation as $\varphi(x)=e^{i\theta(x)}$, the degree is given by 
$$\frac{1}{2\pi}\int_{\partial \Omega}\frac{\partial\theta}{\partial\tau}(x)\, d\tau,$$ 
where $\tau$ parameterizes $\partial\Omega$. 
The degree on a closed curve $\gamma\subset\overline\Omega$  
can in general be defined in a similar way for functions belonging to $H^1(\Omega; S^1)$; if $\Omega$ is not simply connected, the degree may be a non-zero integer.

\paragraph{Ginzburg-Landau functionals\index{Ginzburg-Landau functionals}.}  
We consider the energy 
$$
F_\e(u)=\int_\Omega ((|u|^2-1)^2+\e^2|\nabla u|^2)dx
$$
defined on $u\in H^1(\Omega;\mathbb R^2)$. The pointwise limit of $F_\e$ as $\e\to 0$ is equal to 
$\int_\Omega (|u|^2-1)^2dx$, which is minimized by all functions in $H^1(\Omega; S^1)$.

The similarity with the van der Waals model may suggest scaling the energies by some $\e^\alpha$ with $\alpha<2$. We can check that the $\Gamma$-limit (for example in the sense of the $L^2$-convergence) of such scaled energies is trivial. To that end, choose $z_1,z_2\in S^1$ and $u(x)=u(x_1)= z_1+(z_2-z_1)\chi_{(0,+\infty)}(x_1)$; that is, a function with an interface in the hyperplane $\{x_1=0\}$, which we suppose intersects $\Omega$.
If in complex notation $z_1=e^{i\theta_1}$ and $z_2=e^{i\theta_2}$, for $\delta>0$ we define $$u_\delta(x)= 
{\rm exp}\big(i(\theta_1+ x_1(\theta_2-\theta_1)/\delta\big)\qquad \hbox{ for } x_1\in (0,\delta),
$$
and $u_\delta(x)=u(x)$ otherwise. Since $u(x)\in S^1$ for all $x$,
we then have
$$
\frac1{\e^\alpha}F_\e(u_\delta)= \frac{\e^{2-\alpha}}\delta |\theta_2-\theta_1|^2.
$$
Choosing $\delta=\delta_\e>\!>\e^{2-\alpha}$, we get that the $\Gamma$-limit on such $u$ is trivial. The same argument shows that it is zero on all piecewise-constant functions in $L^\infty(\Omega; S^1)$, and then by density on the whole $L^\infty(\Omega; S^1)$. On the other hand, dividing by $\e^2$ trivially gives that the functionals are equi-coercive with respect to the weak convergence in $H^1$ and strong convergence in $L^2$, and their $\Gamma$-limit has domain $H^1(\Omega; S^1)$, where it is simply $\int_\Omega |\nabla u|^2dx$. Hence, also this scaling has a trivial effect.

\begin{figure}[h!]
	\centerline{\includegraphics [width=6.2in]{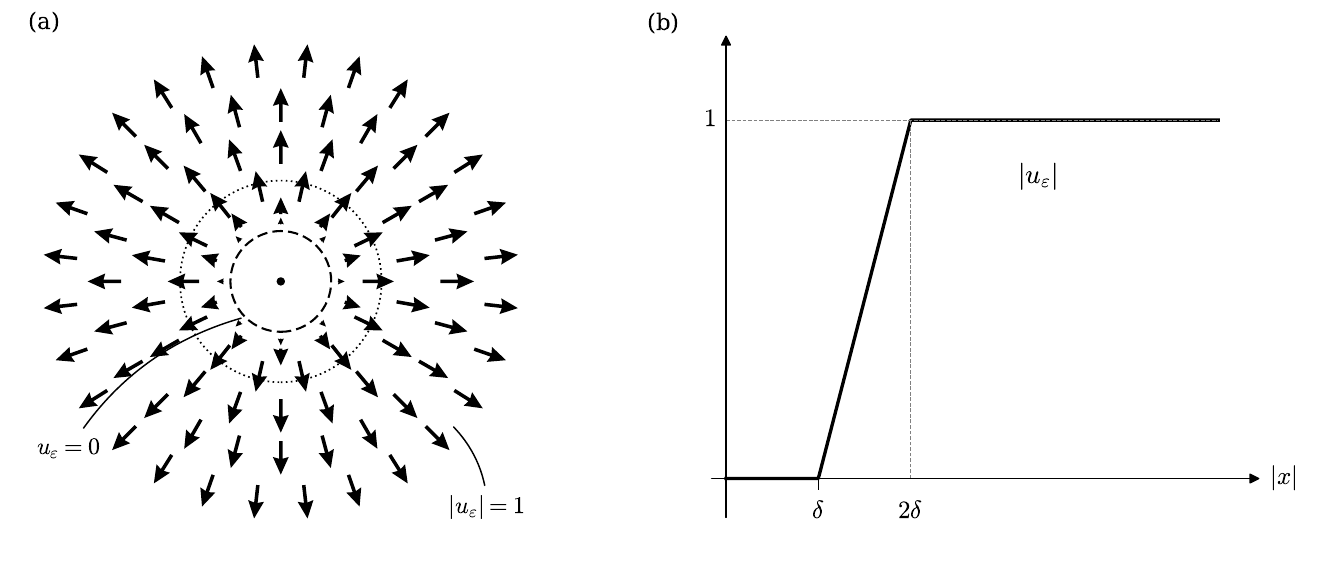}}
	\caption{A vortex: the function $u_\e$ in \eqref{uep}. (a) The vector field $u_\e$ in the plane: it vanishes in $B_\delta$ (dashed circle), and has modulus $1$ and points radially outwards, $u_\e=x/|x|$, outside $B_{2\delta}$ (dotted circle). (b) The modulus $|u_\e|$ as a function of $|x|$.}\label{Fig22}
\end{figure}

In order to highlight a non-trivial scaling we consider $u(x)=x/|x|$ as a target function, and a sequence $u_\e\to u$ 
in $H^1_{\rm loc}(B_r\setminus\{0\}; \mathbb R^2)$, such that there exists a ball of radius $2\delta$ with
$|u_\e|=u(x)$ in $B_r\setminus B_{2\delta}$ and $u_\e=0$ on $B_{\delta}$. Such a $u_\e$ can be constructed explicitly as
\begin{equation}\label{uep}
u_\e(x)= \Big(\Big(\frac{|x|}{\delta}-1\Big)^+\wedge 1\Big)\frac{x}{|x|},
\end{equation}
as pictured in Fig.~\ref{Fig22}.

Then we estimate
$$
F_\e(u_\e)\sim\int_{B_{2\delta}}(|u_\e|^2-1)^2dx+ \e^2\int_{B_{2\delta}}\frac1{\delta^2}dx+ \e^2\int_{B_1\setminus B_{\delta}}|\nabla u|^2dx.
$$
Noting that $|\nabla u|^2=1/|x|^2$, we see that this term is of order
$\delta^2+ \e^2+ \e^2|\log\delta|$, which is minimized when $\delta$ is of order $\e$.

This scaling gives the Ginzburg-Landau functionals at the {\bf vortex scale},
$$
F^{GL}_\e(u)=\frac{1}{\e^2|\log\e|}\int_{\Omega}(|u|^2-1)^2\, dx+
\frac{1}{|\log\e|}\int_{\Omega}|\nabla u|^2\, dx,
$$ 
whose $\Gamma$-limit must be understood in an appropriate sense.
Without entering into details, this limit is computed with respect to a convergence 
$u_\e\to \mu$, where $\mu$ is a measure of the form
$$
\mu=\sum_{j=1}^N d_j \delta_{x_j} \hbox{ (vortex measure),}
$$  
and $d_j\in\mathbb Z$ can be understood as the degree of a limit function at $x_j$. This 
convergence may be defined as the convergence of the Jacobian determinants 
of $u_\e$ in the dual space of Lipschitz functions ({\em flat sense}). The limit 
is then the {\em vortex energy}
$$
F^{GL}(u)=\pi \sum_{j=1}^N |d_j|.
$$

As far as scales are concerned, we remark that the recovery sequences, constructed as in 
\eqref{uep}, can be compared with those for perforated domains in the critical scaling 
\eqref{fep}, showing again {\bf concentration at all scales $\e^\alpha$} around the vortices.

\subsection{Appendix: relevant tools of Chapter 4}\label{App-Ch4}

\paragraph{References.} The ``strange term coming from nowhere'' of Section~\ref{strangeterms} is due to Cioranescu and Murat \cite{CM82}; for the general theory of $p$-capacity we refer to the monograph of Heinonen, Kilpel\"ainen and Martio \cite{HKM93}. The study of Ginzburg-Landau functionals at the vortex scale originates with the book of Bethuel, Brezis and H\'elein \cite{BBH94}; the sharp lower bounds from which the vortex energy is derived are due to Sandier \cite{Sandier98} and Jerrard \cite{Jerrard99}, while the full $\Gamma$-convergence result stated without proof above, with convergence of the Jacobians to the vortex measure, is due to Alberti, Baldo and Orlandi \cite{ABO05}.

\paragraph{Changing of boundary data.}
In the analysis of asymptotic problems a recurring issue is that of matching given boundary conditions. This occurs when Dirichlet data are assigned on a portion of the boundary of the domain, and also when test functions for asymptotic formulas have to satisfy a prescribed condition on some auxiliary set, such as the spheres $\partial B_{\delta_\e r_\e^k}(\delta_\e k)$ used in the capacitary formula of Section~\ref{strangeterms}.
While in phase-transition problems imposing a boundary datum results in a nontrivial additional term (an external boundary layer, as discussed in the previous chapter), for functionals of $p$-growth in Sobolev spaces the same constraint is {\em compatible} with $\Gamma$-convergence: recovery sequences can always be modified to satisfy a prescribed boundary datum at an asymptotically vanishing extra cost. We prove this fact by a quick slicing argument, due to De Giorgi: the near-boundary region is divided into $N$ disjoint zones, and, by the pigeonhole principle, at least one of them carries an energy no larger than $1/N$-th of the total; interpolating between the given sequence and the target datum only across that (thin, low-energy) zone changes the energy by an amount that vanishes as $N\to\infty$.

\begin{lemma}\label{slicing-lemma}
Let $p>1$, let $A$ be a bounded open subset of $\mathbb R^d$, $S\subset\partial A$ a hypersurface, and $\rho>0$ such that the collar $A_\rho=\{x\in A:\mathrm{dist}(x,S)<\rho\}$ is well defined. Let $v\in W^{1,p}(A_\rho)$ and let $\{u_\e\}\subset W^{1,p}(A_\rho)$ be such that $u_\e\to v$ in $L^p(A_\rho)$ and
$$
M:=\sup_\e\int_{A_\rho}|\nabla u_\e|^p\,dx+\int_{A_\rho}|\nabla v|^p\,dx<+\infty.
$$
Then there exists a sequence $\{w_\e\}\subset W^{1,p}(A_\rho)$ such that $w_\e=v$ on $S$ in the sense of traces, $w_\e\to v$ in $L^p(A_\rho)$, and
$$
\limsup_{\e\to0}\int_{A_\rho}|\nabla w_\e|^p\,dx\leq\limsup_{\e\to0}\int_{A_\rho}|\nabla u_\e|^p\,dx.
$$
\end{lemma}

\begin{proof}
Fix $N\in\mathbb N$ and partition $A_\rho$ into $N$ consecutive thin sub-collars $S_1,\dots,S_N$ of width $\rho/N$, labelled so that $S_1$ is adjacent to $S$. For every $\e$, since
$$
\sum_{j=1}^N\int_{S_j}(|\nabla u_\e|^p+|\nabla v|^p)\,dx\leq M,
$$
the pigeonhole principle gives an index $j_\e\in\{1,\dots,N\}$ with
$$
\int_{S_{j_\e}}(|\nabla u_\e|^p+|\nabla v|^p)\,dx\leq\frac MN.
$$
Let $\eta_\e$ be the Lipschitz cut-off equal to $1$ on $S_1\cup\dots\cup S_{j_\e-1}$, equal to $0$ on $S_{j_\e+1}\cup\dots\cup S_N$, and affine in $\mathrm{dist}(\cdot,S)$ on $S_{j_\e}$ (so that $|\nabla\eta_\e|\leq N/\rho$ there), and set $w_\e^N=\eta_\e v+(1-\eta_\e)u_\e$. Then $w_\e^N=v$ on $S$ and $w_\e^N=u_\e$ on $S_{j_\e+1}\cup\dots\cup S_N$, while on $S_{j_\e}$
$$
\int_{S_{j_\e}}|\nabla w_\e^N|^p\,dx\leq C\int_{S_{j_\e}}(|\nabla u_\e|^p+|\nabla v|^p)\,dx+C\Big(\frac N\rho\Big)^p\int_{A_\rho}|u_\e-v|^p\,dx\leq\frac{CM}N+C\Big(\frac N\rho\Big)^p\int_{A_\rho}|u_\e-v|^p\,dx.
$$
Since $u_\e\to v$ in $L^p(A_\rho)$, for $N$ fixed the last term vanishes as $\e\to0$, and $w_\e^N\to v$ in $L^p(A_\rho)$ (being equal to $v$ on part of the collar and to $u_\e\to v$ on the rest); hence
$$
\limsup_{\e\to0}\int_{A_\rho}|\nabla w_\e^N|^p\,dx\leq\limsup_{\e\to0}\int_{A_\rho}|\nabla u_\e|^p\,dx+\frac{CM}N.
$$
Choosing, by a standard diagonal argument, $N=N_\e\to\infty$ slowly enough as $\e\to0$ and setting $w_\e=w_\e^{N_\e}$, all the stated properties hold for the single sequence $\{w_\e\}$.
\end{proof}

\paragraph{Dirichlet conditions on periodic perforations.} We now complete the lower-bound proof of Section~\ref{strangeterms}, where $\delta_\e=\e^{\frac{d-p}{d}}$ is the period of the lattice of balls. Finding, around each ball $B_\e(\delta_\e k)$, a sphere on which $u_\e$ is close to a constant requires, in addition to Lemma~\ref{slicing-lemma}, a control of the oscillation of $u_\e$: this is provided by the Poincar\'e inequality, applied to a dyadic family of annuli around $\delta_\e k$, each of radius comparable to its own thickness.

\begin{proposition}\label{perforated-prop}
Let $\Omega\subset\mathbb R^d$ be open and bounded, let $p\in(1,d)$, and let $\{u_\e\}\subset W^{1,p}(\Omega)$ be such that $u_\e=0$ in $\bigcup_{k\in\mathbb Z^d}B_\e(\delta_\e k)$, $u_\e\weak u$ in $W^{1,p}(\Omega)$, and $F_\e(u_\e)=\int_\Omega|\nabla u_\e|^p\,dx\leq C$. Then, for every $\eta>0$, there exists $\{v_\e^\eta\}\subset W^{1,p}(\Omega)$ such that
\begin{enumerate}
\item[{\rm (i)}] $v_\e^\eta\weak u$ in $W^{1,p}(\Omega)$;
\item[{\rm (ii)}] $F_\e(v_\e^\eta)\leq(1+\eta)F_\e(u_\e)$;
\item[{\rm (iii)}] for every $k\in\mathbb Z^d$ there exist $r_\e^k\in(0,\frac14)$ and a constant $v_{\e,k}^\eta$ such that $v_\e^\eta=v_{\e,k}^\eta$ on $\partial B_{\delta_\e r_\e^k}(\delta_\e k)$ in the sense of traces;
\item[{\rm (iv)}] $\displaystyle\lim_{\e\to0}\delta_\e^d\sum_{k\in\mathbb Z^d}|v_{\e,k}^\eta|^p=\int_\Omega|u|^p\,dx$.
\end{enumerate}
In particular, replacing $u_\e$ by $v_\e^\eta$ and letting $\eta\to0$ along a diagonal sequence gives properties {\rm(i)--(ii)} assumed in Section~\ref{strangeterms}.
\end{proposition}

\begin{proof}
Fix $k\in\mathbb Z^d$ and consider, around $\delta_\e k$, the dyadic family of circular coronas
$$
C_j(\delta_\e k)=B_{2^{1-j}\delta_\e/4}(\delta_\e k)\setminus B_{2^{-j}\delta_\e/4}(\delta_\e k),\qquad j=1,\dots,N,
$$
each with radius of the same order as its own thickness, where $N\ge1/\eta$ is a fixed integer depending only on $\eta$ (we take {\em only $N$ coronas}). Since the $C_j(\delta_\e k)$ are pairwise disjoint subsets of $B_{\delta_\e/4}(\delta_\e k)$,
$$
\sum_{j=1}^N\int_{C_j(\delta_\e k)}|\nabla u_\e|^p\,dx\leq\int_{B_{\delta_\e/4}(\delta_\e k)}|\nabla u_\e|^p\,dx=:M_\e^k,
$$
so that, by the averaging argument of De Giorgi (the most convenient corona carries at most $\frac1N$ of the total), there exists an index $j_\e^k\in\{1,\dots,N\}$ such that $\int_{C_{j_\e^k}(\delta_\e k)}|\nabla u_\e|^p\,dx\leq\frac1NM_\e^k\leq\eta M_\e^k$. The Poincar\'e inequality on $C_{j_\e^k}(\delta_\e k)$, rescaled to the fixed reference corona $B_2\setminus B_1$, controls the oscillation of $u_\e$ on the spheres crossing this corona by
$$
\mathrm{osc}\,u_\e\leq C\big(2^{-j_\e^k}\delta_\e\big)^{1-\frac dp}\Big(\int_{C_{j_\e^k}(\delta_\e k)}|\nabla u_\e|^p\,dx\Big)^{1/p}\leq C_\eta\,\delta_\e^{1-\frac dp}(M_\e^k)^{1/p},
$$
where $C_\eta=C\,2^{N(\frac dp-1)}$ depends only on $\eta$ (here $p<d$ and $j_\e^k\le N$), and in particular not on $\e$ and $k$.
Setting $r_\e^k=2^{-j_\e^k-2}\in(0,\tfrac18]\subset(0,\tfrac14)$, let $v_{\e,k}^\eta$ be the average of $u_\e$ on $\partial B_{\delta_\e r_\e^k}(\delta_\e k)$: the estimate above shows that $u_\e$ is, up to an $\eta$-controlled, $L^p$-small error, already equal to this constant on the sphere. Applying Lemma~\ref{slicing-lemma} on the corona $C_{j_\e^k}(\delta_\e k)$, with $S=\partial B_{\delta_\e r_\e^k}(\delta_\e k)$ and target value $v\equiv v_{\e,k}^\eta$, replaces $u_\e$ by a function equal to $v_{\e,k}^\eta$ exactly on $S$, at an extra cost of the order of $\eta M_\e^k$, since the corona carrying the transition was chosen, by the averaging argument above, to hold only a fraction $\eta$ of the local energy $M_\e^k$. Performing this replacement independently on each (pairwise disjoint) ball $B_{\delta_\e/4}(\delta_\e k)$, and leaving $u_\e$ unchanged elsewhere, defines $v_\e^\eta$, with
$$
F_\e(v_\e^\eta)\leq F_\e(u_\e)+C\eta\sum_{k}M_\e^k\leq(1+C\eta)F_\e(u_\e),
$$
which gives (ii) up to renaming $\eta$. Properties (i) and (iii) hold by construction, and (iv) follows from $u_\e\to u$ in $L^p(\Omega)$ together with the smallness of the oscillation of $u_\e$ on the spheres $\partial B_{\delta_\e r_\e^k}(\delta_\e k)$: by the estimate above, $\delta_\e^d\sum_k\mathrm{osc}^p\le C_\eta\,\delta_\e^{p}\sum_kM_\e^k\le C_\eta\,\delta_\e^pF_\e(u_\e)\to0$, which identifies the constants $v_{\e,k}^\eta$, in the limit, with the local averages of $u$.
\end{proof}


\section{Concentration on free discontinuities} 
We now consider a class of problems where the main variable is the set of discontinuity points $S(u)$ of some function $u$. In order to detect such a set it is convenient to use some small scale of approximation, in which the set is substituted by some other variable or corresponds to a feature of the approximating sequence.

\subsection{A {\em phase-field} approach to free-discontinuity problems}
We consider a simple one-dimensional approximation problem, with the target functional
\begin{equation}\label{MS=0}
F(u)=\begin{cases}\displaystyle\int_{-1}^1|u'|^2dt &\hbox{ if } u\in H^1(-1,1)\\
\displaystyle\int_{-1}^1|u'|^2dt+1 &\hbox{ if $u\in H^1((-1,0)\cup(0,1))$ and $u(0^-)\neq u(0^+)$};\end{cases}
\end{equation}
that is, we want to ``detect'' whether $u$, elsewhere $H^1$, has a jump or not at $0$, by approximation with elliptic functionals. The function $u'$ in the second formula is considered as an $L^2$ function defined piecewise.
Note that $F$ is not a convex functional; hence, it cannot be approximated by functionals of the form 
$\int_{-1}^1a_\e(x)|u'|^2dt$, which are convex. 
The idea to overcome this difficulty is to introduce a new variable $v$ which concentrates at $0$ and is of ``phase-transition type'', for which the convexity restriction does not hold. Such an approach is usually called a {\em phase-field approach} to the treatment of discontinuity. 
The functionals we introduce are of the form
$$
F_\e(u,v)=\int_{-1}^1v^2|u'|^2dt+ \int_{-1}^1\Big( \frac{(1-v)^2}{4\e}+\e|v'|^2\Big)dt,
$$
defined on pairs of functions $u,v\in H^1(-1,1)$ such that $v(t)=1$ if $|t|>\sqrt\e$. Note that this functional is separately convex (even quadratic) in both its variables.

Consider a sequence $u_\e,v_\e$ with $F_\e(u_\e,v_\e)\le S<+\infty$ and $u_\e\to u$ in $L^2(-1,1)$. Then, up to subsequences, we have 
$u_\e\wto u$ in $ H^1((-1,0)\cup(0,1))$.
Note that 
$$
 \liminf_{\e\to 0}\int_{(-1,1)\setminus (-\sqrt\e,\sqrt\e)}|u'_\e|^2dt
\ge \liminf_{\e\to 0}\int_{(-1,1)\setminus (-\delta,\delta)}|u'_\e|^2dt\ge \int_{(-1,1)\setminus (-\delta,\delta)}|u'|^2dt
$$
for all $\delta>0$, so that the inequality
\begin{eqnarray}\label{MS-d1}
\liminf_{\e\to 0}\int_{(-1,1)\setminus (-\sqrt\e,\sqrt\e)}|u'_\e|^2dt\ge \int_{(-1,1)}|u'|^2dt
\end{eqnarray}
always holds by the arbitrariness of $\delta$.

 If $v_\e\ge c>0$ on $(-1,1)$ then, up to subsequences, $u_\e\wto u$ in $H^1(-1,1)$, and the lower inequality is achieved. Conversely, if, up to subsequences, there exist $t_\e$ such that $v_\e(t_\e)\to0$, 
then we have 
$$
\int_{(-\sqrt \e,\sqrt \e)} \Big( \frac{(1-v_\e)^2}{4\e}+\e|v_\e'|^2\Big)dt
\ge \int_{(-\sqrt \e,\sqrt \e)}(1-v_\e)|v_\e'|dt
$$
$$
\ge \bigg(\Big|\int_{(-\sqrt \e,t_\e)}(1-v_\e)v_\e'dt\Big|+ \Big|\int_{(t_\e,\sqrt \e)}(1-v_\e)v_\e'dt\Big|\bigg)
= 2\int_{v_\e(t_\e)}^1(1-s)ds= 1+o(1).
$$
In the terminology that we have used in the previous chapters, these contributions are given by two boundary layers connecting points where $v_\e$ is almost $0$ and $1$, respectively. The presence of the first point is ensured by the necessity of having a discontinuity in the limit. The second is ensured by the vanishing of the contribution of the ``single-well potential'' $(1-v)^2$. Gathering the estimate just obtained and \eqref{MS-d1}, we then have, in this second case,
$$
\liminf_{\e\to 0} F_\e(u_\e,v_\e)\ge \int_{(-1,1)}|u'|^2dt +1.
$$ 
Since in the first case $u\in H^1(-1,1)$ and the liminf is at least $\int_{(-1,1)}|u'|^2dt$, in both cases we have $\liminf_{\e\to 0} F_\e(u_\e,v_\e)\ge F(u)$.

\begin{figure}[h!]
	\centerline{\includegraphics [width=6.2in]{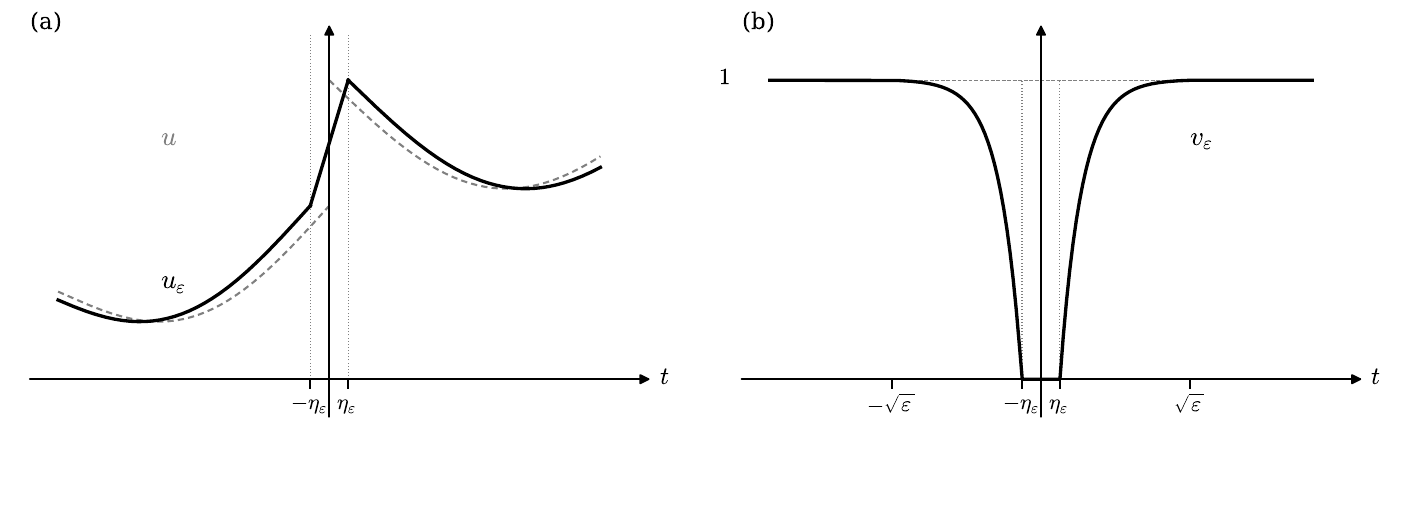}}
	\caption{The recovery sequence: (a) a function $u$ with a jump at $0$ (gray, dashed) and the corresponding $u_\e$, which is a translation of $u$ outside $(-\eta_\e,\eta_\e)$ and is affine in this interval; (b) the function $v_\e$, which vanishes in $(-\eta_\e,\eta_\e)$, is equal to $1$ for $|t|\ge\sqrt\e$, and connects the two values in two boundary layers at scale $\e$ at $\pm\eta_\e$. The scales ($\eta_\e\ll\e\ll\sqrt\e$) are not respected.}\label{Fig24}
\end{figure}

We now fix $u\in H^1((-1,0)\cup(0,1))$ (the case $u\in H^1(-1,1)$ being trivial with $u_\e=u$ and $v_\e=v=1$), and a parameter $\eta_\e>0$ with $\eta_\e=o(\e)$ (for instance, $\eta_\e=\e^2$). We define the recovery sequences as follows
$$
u_\e(t)= \begin{cases}u(t-\eta_\e) & \hbox{ if } t>\eta_\e\\
u(t+\eta_\e) & \hbox{ if } t<-\eta_\e\\ \displaystyle
u(0^-)+\frac{u(0^+)-u(0^-)}{2\eta_\e}(t+\eta_\e)& \hbox{ if } |t|\le\eta_\e,
\end{cases}
$$
and
$$
v_\e(t)= \begin{cases}1 & \hbox{ if } |t|>\sqrt\e
\\
0 &  \hbox{ if } |t|<\eta_\e
\\ \displaystyle
\frac{1-e^{-\frac1{2\e}(|t|-\eta_\e)}}{1-e^{-\frac1{2\e}(\sqrt\e-\eta_\e)}} &  \hbox{ if } |t|\in[\eta_\e,\sqrt\e],
\end{cases}
$$
as pictured in Fig.~\ref{Fig24}.

The function $u_\e$ is the translation of $u$ by $\eta_\e$ on both sides of $0$, with an affine interpolation on $(-\eta_\e,\eta_\e)$, so that $u_\e\to u$ in $L^2(-1,1)$. Since $v_\e=0$ on $(-\eta_\e,\eta_\e)$, the only region where $u_\e$ has a large derivative does not contribute to the first term of $F_\e$, and
$$
\int_{(-1,1)}v_\e^2|u_\e'|^2dt\le\int_{(-1,1)\setminus(-\eta_\e,\eta_\e)}|u_\e'|^2dt\le \int_{(-1,0)\cup(0,1)}|u'|^2dt.
$$
Up to a multiplicative constant converging to $1$, the function $v_\e$ is defined, on each of the two intervals $(\eta_\e,\sqrt\e)$ and $(-\sqrt\e,-\eta_\e)$, by a version, scaled by $\e$, of the function $v(s)=1-e^{-s/2}$, solution of $v'=\frac12(1-v)$ with initial condition $v(0)=0$; that is, $v_\e(t)=v\big(\frac{|t|-\eta_\e}{\e}\big)/N_\e$ with $N_\e=1-e^{-\frac1{2\e}(\sqrt\e-\eta_\e)}\to1$. Setting $S_\e=\frac{\sqrt\e-\eta_\e}{\e}\to+\infty$ and changing variables, we obtain
\begin{eqnarray*}
\int_{(-\sqrt \e,\sqrt \e)\setminus(-\eta_\e,\eta_\e)} \Big( \frac{(1-v_\e)^2}{4\e}+\e|v_\e'|^2\Big)dt
&=&\frac{2}{4N_\e^2}\int_0^{S_\e}\Big[\big(e^{-\frac s2}-e^{-\frac{S_\e}2}\big)^2+e^{-s}\Big]ds
\\&\longrightarrow&\frac12\int_0^{+\infty}2e^{-s}ds=1,
\end{eqnarray*}
since $N_\e\to1$ and $e^{-S_\e/2}\to0$. The region $(-\eta_\e,\eta_\e)$, where $v_\e=0$, contributes
$$
\int_{(-\eta_\e,\eta_\e)}\frac{(1-v_\e)^2}{4\e}dt=\frac{\eta_\e}{2\e}\to0,
$$
which is the reason for choosing $\eta_\e=o(\e)$, while $v_\e=1$ on $\{|t|>\sqrt\e\}$ gives no contribution. Hence $\limsup_{\e\to0}F_\e(u_\e,v_\e)\le\int_{(-1,1)}|u'|^2dt+1$.
Thus, we have used a phase-field variable at the {\bf scale} $\e$, in the same way as for the phase-transition case.
Note that formally we have shown the $\Gamma$-convergence of $F_\e$ to a different functional than $F$; namely, with a little abuse of notation
$$
F(u,v)=\begin{cases}F(u) &\hbox{ if  $v=1$}\\
+\infty &\hbox{ otherwise. }\end{cases}
$$
Clearly, the two functionals are equivalent.

\smallskip
Functional \eqref{MS=0} is a particular case of a more general functional, in which the location of the jump points is an unknown, the (one-dimensional version of the) {\em Mumford--Shah functional} \cite{MS89,DGA88,AFP}. In dimension one its domain consists of piecewise-$H^1$ functions; that is, functions $u$ which are $H^1$ on the complement of a finite set of jump points $S(u)$. The functional on an interval $(a,b)$ can then be written
\begin{equation}\label{MS=ab}
F^{MS}(u)=\alpha\int_{(a,b)}|u'|^2dt+\beta \#(S(u)).
\end{equation}
 This functional is lower semicontinuous with respect to the strong $L^2$ convergence, and it is coercive on sets of functions bounded in $L^\infty(a,b)$ (see the Appendix to this chapter).
 Note that for $\alpha=\beta=1$, $(a,b)=(-1,1)$ and on functions $u\in H^1((-1,1)\setminus \{0\})$, the functional is the same as \eqref{MS=0}. The same proof as in that particular case shows that an approximation of $F^{MS}$ is given by the {\em Ambrosio--Tortorelli functionals} \cite{AT90,AT92}
$$
F^{AT}_\e(u,v)=\alpha\int_{(a,b)}v^2|u'|^2dt+ \beta\int_{(a,b)}\Big( \frac{(1-v)^2}{4\e}+\e|v'|^2\Big)dt,
$$
defined on pairs of functions $u,v\in H^1(a,b)$, and with respect to the $L^2$ convergence for both. The function $u'$ is considered as an $L^2$ function defined piecewise. Again, the limit is formally a functional depending on $v$, but it is finite, and equal to the Mumford--Shah functional, only if $v$ is identically $1$. The phase-field variable concentrates at {\bf scale} $\e$ on the set of discontinuity points $S(u)$ of $u$.

\paragraph{Alternative Approximations}
Due to its importance in Visual Reconstruction and, in its higher-dimensional versions, in Fracture Mechanics, several different approximations of the Mumford--Shah functionals have been proposed in the literature (see also \cite{Bra98}). Here we list a few of them in the simpler one-dimensional case \eqref{MS=ab}.

\bigskip
\noindent{\bf Discrete Approximation (Blake--Zisserman/Chambolle)}.
We consider the discrete energies, introduced by Blake and Zisserman \cite{BZ87} and studied by Chambolle \cite{Cha95},
\begin{equation}\label{BZC}
F_\e(u)=\sum_{i=1}^{N_\e} \min\bigg\{\frac{\alpha(u_i-u_{i-1})^2}\e,\beta\bigg\}=\sum_{i=1}^{N_\e} \e\min\bigg\{\alpha\Big(\frac{u_i-u_{i-1}}{\e}\Big)^2,\frac\beta\e\bigg\},
\end{equation}
where $\{t_0,\ldots,t_{N_\e}\}$ is a parameterization of $[a,b]\cap\e\mathbb Z$. If we identify each such function $u$ with the function $v$ given by

\smallskip
$\bullet$ its affine interpolation of the values $u_{i-1}$ and $u_i$ on $[t_{i-1}, t_i]$ if $|u_i-u_{i-1}|\le\sqrt{\beta\e/\alpha}$;

\smallskip
$\bullet$ the constant value $u_{i-1}$ on $[t_{i-1}, t_i)$ if $|u_i-u_{i-1}|>\sqrt{\beta\e/\alpha}$;

\smallskip
$\bullet$ the value $u_0$ in $(a,t_0)$ and $u_{N_\e}$ in $(t_{N_\e},b)$.

\smallskip\noindent
\begin{figure}[h!]
	\centerline{\includegraphics [width=6.2in]{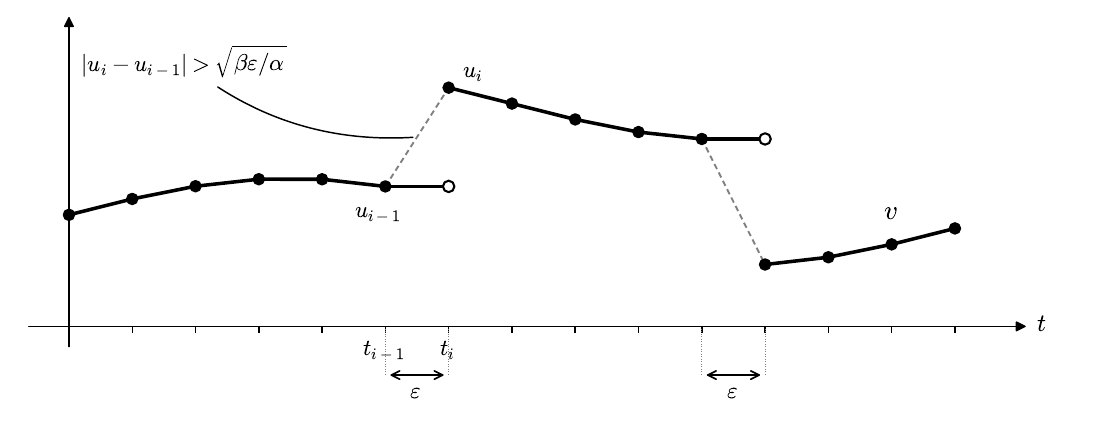}}
	\caption{The identification of a discrete function $u=(u_0,\dots,u_{N_\e})$ (dots) with the piecewise-affine function $v$ (solid line): $v$ is the affine interpolation of the values $u_{i-1}$ and $u_i$ on $[t_{i-1},t_i]$ if $|u_i-u_{i-1}|\le\sqrt{\beta\e/\alpha}$; otherwise, $v=u_{i-1}$ on $[t_{i-1},t_i)$ and $v$ has a jump at $t_i$ (the dashed segment is the affine interpolation, which is not used). The intervals of length $\e$ where this happens form the set $S_\e$.}\label{Fig25}
\end{figure}

With this identification (Fig.~\ref{Fig25}), we have $F_\e(u)=F^{MS}(v)$, and the liminf inequality follows since $u_\e\to u$ implies $v_\e\to u$. The upper bound is simply obtained with the recovery sequence $u_\e=u$; that is, $u_\e$ is the discretization of the target function $u$. Hence, the set $S(u)$ is approximated by the sets $S_\e$ given by the union of the intervals $[t^\e_{i-1}, t^\e_i]$ {\bf of length $\e$} where $|u_i-u_{i-1}|>\sqrt{\beta\e/\alpha}$.

\bigskip
\noindent{\bf Finite-difference Approximation (De Giorgi/Gobbino)}. This is a continuum approximation of the previous one, conjectured by De Giorgi and proved by Gobbino \cite{Gob98}. For simplicity of notation we consider the case of functionals defined on the whole line. We fix a symmetric mollifier $\varphi:\mathbb R\to [0,+\infty)$ with support in $[-1,1]$, and define
$$
F_\e(u)=\frac1{\e^2}\int_{\mathbb R\times \mathbb R} \varphi\Big(\frac{x-y}\e\Big)\min\bigg\{\frac{|u(x)-u(y)|^2}\e,1\bigg\}dx\,dy.
$$
This functional can be rewritten as
\begin{eqnarray*}
F_\e(u)&=&\int_{\mathbb R}\varphi(\xi) \frac1{\e}\int_{\mathbb R}\min\bigg\{\frac{|u(x+\e\xi)-u(x)|^2}\e,1\bigg\}dx\,d\xi
\\
&=&\int_{\mathbb R}\varphi(\xi) \frac1{\e}\int_0^{\e|\xi|}\sum_{j\in\mathbb Z}\min\bigg\{\frac{|u(\e\xi j+t)-u(\e\xi(j-1)+t)|^2}\e,1\bigg\}dt\,d\xi
\\
&\ge&\int_{\mathbb R}\varphi(\xi)  |\xi|^2 \sum_{j\in\mathbb Z}\min\bigg\{\frac{|u(\e\xi j+t_\e)-u(\e\xi(j-1)+t_\e)|^2}{\e|\xi|},\frac1{|\xi|}\bigg\}d\xi
\end{eqnarray*}
for some $t_\e\in [0,\e|\xi|]$. By the previous discrete approximation, if $u_\e\to u$, we have
$$
\liminf_{\e\to 0} \sum_{j\in\mathbb Z}\min\bigg\{\frac{|u(\e\xi j+t_\e)-u(\e\xi(j-1)+t_\e)|^2}{\e|\xi|},\frac1{|\xi|}\bigg\}
\ge \int_{\mathbb R} |u'|^2dt+\frac1{|\xi|}\#(S(u)).
$$
Hence, by applying Fatou's Lemma, we obtain the liminf inequality with 
$$
\alpha= \int_{\mathbb R}\varphi(\xi)  |\xi|^2d\xi,\qquad \beta=\int_{\mathbb R}\varphi(\xi)  |\xi|d\xi.
$$

The recovery sequence is again obtained by taking $u_\e=u$. At fixed $\xi$, each jump point of $u$ corresponds to an interval of length $\e|\xi|$; hence, the {\bf scale is $\e$ modulated by $\xi$}.

\bigskip
\noindent{\bf Non-local Approximation (Braides--Dal Maso)}. In this approximation \cite{BDM97} the gradient in a non-convex energy can be replaced by some averaged quantity which is lower semicontinuous for the weak convergence. A simple such energy is
$$
F_\e(u)=\int_{\mathbb R} \min\bigg\{\frac1{\e}\int_{x-\frac\e2}^{x+\frac\e2}|u'|^2dt,\frac1\e\bigg\}dx,
$$
again integrated on the whole real line for simplicity of notation. The domain can be taken to be $H^1_{\rm loc}$ functions or simply $L^2_{\rm loc}$ functions, with the convention that the integrand takes the value $1/\e$ if $u\not\in H^1(x-\frac\e2,x+\frac\e2)$. 

By Jensen's inequality and a discretization argument we have that there exists $x_\e\in[0,\e]$ such that, setting
$u_j=u(j\e+x_\e)$, we obtain
$$
F_\e(u)\ge\int_{\mathbb R} \min\Big\{\Big(\frac{u(x+\frac\e{2})-u(x-\frac\e{2})}{\e}\Big)^2,\frac1\e\Big\}dx
\ge\sum_{j\in\mathbb Z} \min\Big\{\frac{(u_j-u_{j-1})^2}\e,1\Big\};
$$
that is, we have a lower bound by the discrete approximation, obtaining the Mumford--Shah functional with $\alpha=\beta=1$. Again, a recovery sequence for a piecewise $H^1$-function is obtained with the trivial sequence $u_\e=u$. In the computation of $F_\e(u_\e)$, note that if $t\in S(u)$ then $u\not\in H^1(x-\frac\e2,x+\frac\e2)$ if $x\in(t-\frac\e2,t+\frac\e2)$, on which interval we then have the value $1/\e$ giving the constant $1$. This shows a concentration of the discontinuity at {\bf scale} $\e$.

\subsection{Interaction of concentration and free discontinuities}
The Mumford--Shah functional is a particular {\em free-discontinuity functional}. General (spatially homogeneous one-dimensional) free-discontinuity functionals take the form 
\begin{equation}\label{fdf}
F(u)=\int_{(a,b)} f(u')dt +\sum_{t\in S(u)} \vartheta (u(t^+)-u(t^-)).
\end{equation}
If $\vartheta\ge c>0$ and $f(z)/|z|\to+\infty$ as $|z|\to+\infty$, then $F$ is well defined (and coercive) on piecewise-$W^{1,1}$ functions. In order to remove the first of these two conditions, we have to enlarge the domain of $F$ to one-dimensional {\em special functions of bounded variation} $SBV(a,b)$; that is, $u\in L^1(a,b)$ whose distributional derivative can be written as 
$$
Du=u' d{\mathcal L^1}+ \sum_{t\in S(u)}(u(t^+)-u(t^-))\delta_t,
$$
where $t^\pm$ denote the approximate limits at $t$, and $S(u)=\{t: u(t^+)\neq u(t^-)\}$ is the set of discontinuity points of $u$. Then necessary and sufficient conditions for $F$ to be lower-semicontinuous on $SBV(a,b)$ with respect to the $L^1$ convergence, and coercive on sets of functions bounded in $L^\infty(a,b)$, are that

$\bullet$  $f$ is lower semicontinuous and convex, and $\lim\limits_{z\to \infty} \frac{f(z)}{|z|}=+\infty$;

$\bullet$ $\vartheta$ is  lower semicontinuous and {\em subadditive}, and  $\lim\limits_{z\to 0} \frac{\vartheta(z)}{|z|}=+\infty$.

\smallskip
The subadditivity condition on $\vartheta$, that is,
$$
\vartheta(w_1+w_2)\le \vartheta(w_1)+\vartheta(w_2),
$$
guarantees that it is not convenient to split jumps. The superlinearity of $\vartheta$ at $0$ implies that it is not possible to approximate gradients by jumps, just as the superlinearity of $f$ at $\infty$ implies that the converse is not possible.

\paragraph{An approximation by singular perturbation.}
Following \cite{ABG98}, we consider a continuum version of the functional depending on a {\em truncated quadratic potential} as in \eqref{BZC}, with a singular perturbation of second order
\begin{equation}
F_\e(u)=\int_{(a,b)} \min\Big\{ |u'|^2,\frac1\e\Big\}dt +\e^3\int_{(a,b)} |u''|^2dt
\end{equation}
defined on $H^2(a,b)$. The reason for the scaling $\e^3$ will be clear in the following, as we explain how these functionals approximate a free-discontinuity functional.
Let $u\in H^2$ and let $[a',b']\subset (a,b)$ be a maximal interval where $|u'|\ge 1/{\sqrt\e}$. Then $u'(a')=u'(b')$. For our arguments it is not restrictive to suppose that $u'(a')=1/{\sqrt\e}$. We then consider the function $v(t)= u(a'+t)- u(a')-\frac1{\sqrt\e}t$ defined on $[0,T]$ with $T=b'-a'$. We have 
$$
\int_{(a',b')} \min\Big\{ |u'|^2,\frac1\e\Big\}dt +\e^3\int_{(a',b')} |u''|^2dt=
\frac T\e+\e^3\int_0^T|v''|^2dt.
$$
If we set $z=v(T)$, and change variables by setting $s=t/\e\sqrt{z} $,  $v(t)= zw(t/\e\sqrt{z} )$, $S=T/\e\sqrt{z}$,
we have 
$$
\frac T\e+\e^3\int_0^T|v''|^2dt =\sqrt{z} \Big(S+\int_0^S|w''|^2ds\Big).
$$
Noting that $w(0)=0$, $w(S)=1$ and $w'(0)=w'(S)=0$, setting 
\begin{equation}\label{c2}
C_2=\inf_{S>0}\min\bigg\{S+\int_0^S|w''|^2ds: w(0)=0, w(S)=1, w'(0)=w'(S)=0\bigg\},
\end{equation}
we have 
$$
\int_{(a',b')} \min\Big\{ |u'|^2,\frac1\e\Big\}dt +\e^3\int_{(a',b')} |u''|^2dt\ge C_2\sqrt{\Big|u(b')-u(a')-
\frac{b'-a'}{\sqrt\e}\Big|}
$$
If we define the function $\overline u$ as equal to $u$ outside $[a',b']$ and as
$$
\overline u(t)= u(a')+u'(a')(t-a') \hbox{ on } [a',b'),
$$
then $\overline u(b'^+)-\overline u(b'^-)= u(b')-u(a')-\frac{b'-a'}{\sqrt\e}$.
Hence, if we repeat this construction on all such intervals $[a',b']$, then we construct a function $\overline u$ such that
$$
F_\e(u)\ge \int_{(a,b)} |\overline u'|^2dt + C_2\sum_{t\in S(\overline u)} \sqrt{|\overline u(t^+)-\overline u(t^-)|}.
$$
In this construction, if also $[a,b']$ is a maximal interval where $|u'|\ge 1/{\sqrt\e}$ then we simply set $\overline u(t)=u(b')+u'(b')(t-b')$ on this interval, and analogously on intervals $[a',b]$. 

If $u_\e\to u$ in $L^1(a,b)$ then also the corresponding $\overline u_\e$ constructed as above converge to $u$, so that $u\in SBV(a,b)$, and 
$$
\liminf_{\e\to 0}F_\e(u_\e)\ge  \int_{(a,b)} |u'|^2dt + C_2\sum_{t\in S(u)} \sqrt{|u(t^+)-u(t^-)|}=:F_2(u).
$$

\begin{figure}[h!]
	\centerline{\includegraphics [width=6.2in]{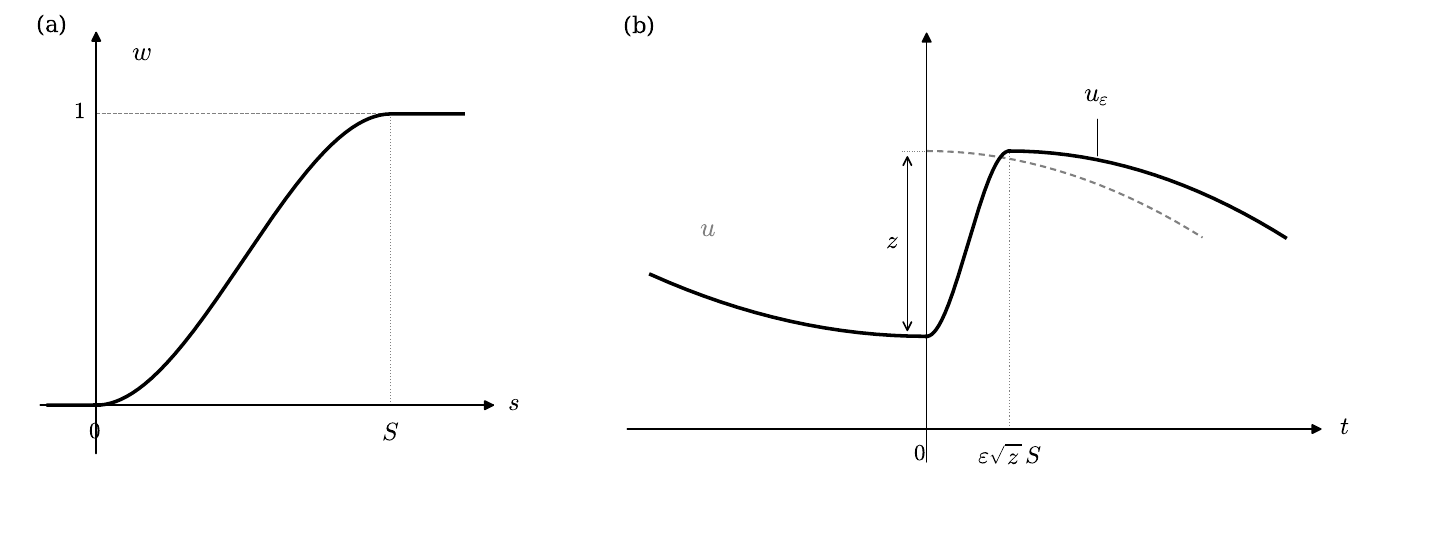}}
	\caption{The recovery sequence in \eqref{ABG-ue}. (a) The profile $w$ in the definition of $C_2$: $w(0)=0$, $w(S)=1$ and $w'(0)=w'(S)=0$ (extended by $0$ and $1$ outside $[0,S]$). (b) A function $u$ with $u'(0^\pm)=0$ and a jump of size $z$ at $0$ (gray, dashed), and the corresponding $u_\e$: it coincides with $u$ for $t\le0$, follows the graph of $w$, rescaled by $z$ in height and by $\e\sqrt z$ in length, on $[0,\e\sqrt zS]$, and is a translation of $u$ afterwards. The width $\e\sqrt zS$ of the transition is not to scale.}\label{Fig26}
\end{figure}

The construction of a recovery sequence describes the scaling properties of this approximation. It suffices to consider the case of $(a,b)=(-1,1)$ and $S(u)=\{0\}$. It is not restrictive to suppose that $u\in H^2((-1,1)\setminus\{0\})$, that $u'(0^\pm)=0$, and $z:=u(0^+)-u(0^-)>0$. Let $S,w$ be a minimizing pair for the definition of $C_2$ (see Fig.~\ref{Fig26}(a)). Then we set
\begin{equation}\label{ABG-ue}
u_\e(t)=\begin{cases} u(t) & \hbox{ if  $t\le 0$}\\ u(0^-)+
z w\Big(\frac{t}{\e\sqrt z}\Big) &\hbox{ if  $0\le t\le \e\sqrt z S$}\\
u(t- \e\sqrt z S) & \hbox{ if  $t\ge \e\sqrt z S$,}
\end{cases}
\end{equation}
and note that for $\e$ small enough
$$
F_\e(u_\e)\le \int_{(a,b)} (|u'|^2+\e^3|u''|^2)dt + \sqrt{|u(0^+)-u(0^-)|}\Big(S+\int_0^S |w''|^2ds\Big),
$$
which shows the upper bound with $F_2(u)$. The construction in \eqref{ABG-ue} shows the concentration on discontinuities at {\bf scale $\e$ with factor} $\sqrt z$.

\paragraph{Perturbations of higher order.}
The same analysis can be carried out for singular perturbations of order $k\ge2$,
$$
F_\e(u)=\int_{(a,b)} \min\Big\{ |u'|^2,\frac1\e\Big\}dt +\e^{2k-1}\int_{(a,b)} |u^{(k)}|^2dt,
$$
defined on $H^k(a,b)$, where the power $2k-1$ is chosen so that a transition giving a jump of size $z$ has an energy of order one. Repeating the scaling argument above with $v(t)=zw(t/\lambda)$ and $\lambda=\e|z|^{1/k}$, the energy of such a transition on an interval of length $\lambda S$ equals $|z|^{1/k}\big(S+\int_0^S|w^{(k)}|^2ds\big)$. The jump energy density is then $m_k|z|^{1/k}$, with
$$
\begin{array}{c}
m_k=\min_{S>0}\min\Big\{S+\int_0^S|w^{(k)}|^2ds:\ w(0)=0,\ w(S)=1,\\ w^{(\ell)}(0)=w^{(\ell)}(S)=0,\ \ell=1,\dots,k-1\Big\},
\end{array}
$$
so that $m_2=C_2$ and the power $\frac12$ is replaced by $\frac1k$; the concentration occurs at scale $\e|z|^{1/k}$. For $k=2$ the only boundary condition needed for the profile is on the first derivative, which is given by the definition of the interval where $|u'|\ge1/\sqrt\e$. For $k>2$, instead, nothing is known a priori on the derivatives $u^{(\ell)}$ of intermediate order $\ell=2,\dots,k-1$ at the endpoints of such an interval, and the homogeneous conditions $w^{(\ell)}(0)=w^{(\ell)}(S)=0$ in the definition of $m_k$ can only be obtained by changing the boundary data of $u_\e$ on slightly larger intervals at a negligible cost. This requires refined interpolation estimates, developed by Solci \cite{Sol25}, who also treats the higher-dimensional case by a slicing argument and obtains the Mumford--Shah functional as a limit of these approximations as $k\to+\infty$; see also \cite{BF26} for the analogous result for Perona--Malik-type energies.

\subsection{Interaction of relaxation and phase-field discontinuities}
\begin{figure}[h!]
	\centerline{\includegraphics [width=4.2in]{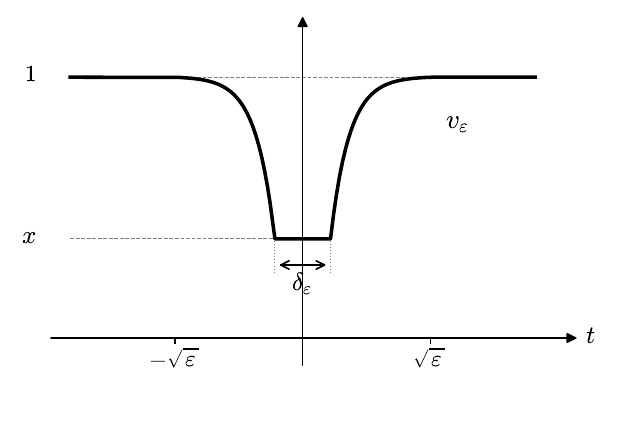}}
	\caption{The phase-field variable $v_\e$ optimizing the jump energy: it is equal to the constant $x$ on an interval of length $\delta_\e$, where the gradient of $u_\e$ concentrates, and connects $x$ to $1$ on both sides by two boundary layers at scale $\e$. The scales are not respected.}\label{Fig27}
\end{figure}

We now examine the jump energies due to an interaction between a phase-field approach as in the Ambrosio--Tortorelli functionals and the scale-free relaxation due to a bulk energy of linear growth \cite{ABS99}.
We consider the functionals 
$$
F_\e(u,v)=\int_{(a,b)}\psi(v)|u'|dt+ \int_{(a,b)}\Big( \frac{(1-v)^2}{4\e}+\e|v'|^2\Big)dt,
$$
where $\psi:[0,1]\to [0,1]$ is an increasing lower-semicontinuous function with $\psi(0)=0$, $\psi(1)=1$ and $\psi>0$ on $(0,1]$; in particular, we can consider $\psi(v)=v^2$ as in the Ambrosio--Tortorelli case.  The main difference here is the linear dependence on $u'$ in the first variable, so that it allows for gradients concentrating to a jump discontinuity of size $z$. If we take $v=x$ constant on the sets where the gradient tends to a jump, we have a jump energy of $\psi(x)|z|$. On the other hand, if the variable $v$ has a transition from $1$ to $x$ on an interval $[t_1,t_2]$, the second term in $F_\e$ gives an optimal contribution 
$$
\int_{t_1}^{t_2} |1-v||v'|dt = \frac12 (1-x)^2.
$$
We can then optimize this process for $v$ by concatenating a transition layer between $1$ and $x$, then a transition giving a jump maintaining $v=x$ constant on an interval of length $\delta$ (see Fig.~\ref{Fig27}), and then a symmetric transition layer between $x$ and $1$. The resulting value is then 
$$
\psi(x) |z|+ (1-x)^2+ \frac\delta{4\e} (1-x)^2.
$$

By the arbitrariness of $\delta$ and optimizing in $x$ we obtain the jump energy
$$
\vartheta(z)=\min\{ \psi(x) |z|+ (1-x)^2\}.
$$
The optimal form of the phase-field variable $v$ highlights a concentration at {\bf scale $\e$ optimized on the well depth}.

In this model, we have the freedom to choose $\psi$, which influences the resulting function $\vartheta$. For example,

$\bullet$ if $\psi(x)=x^2$, then $\displaystyle\vartheta (z)=\frac{|z|}{1+|z|}$;

$\bullet$ if $\psi=1-\chi_{\{0\}}$, then $\vartheta (z)=\min\{|z|,1\}$.

\begin{figure}[h!]
	\centerline{\includegraphics [width=5.2in]{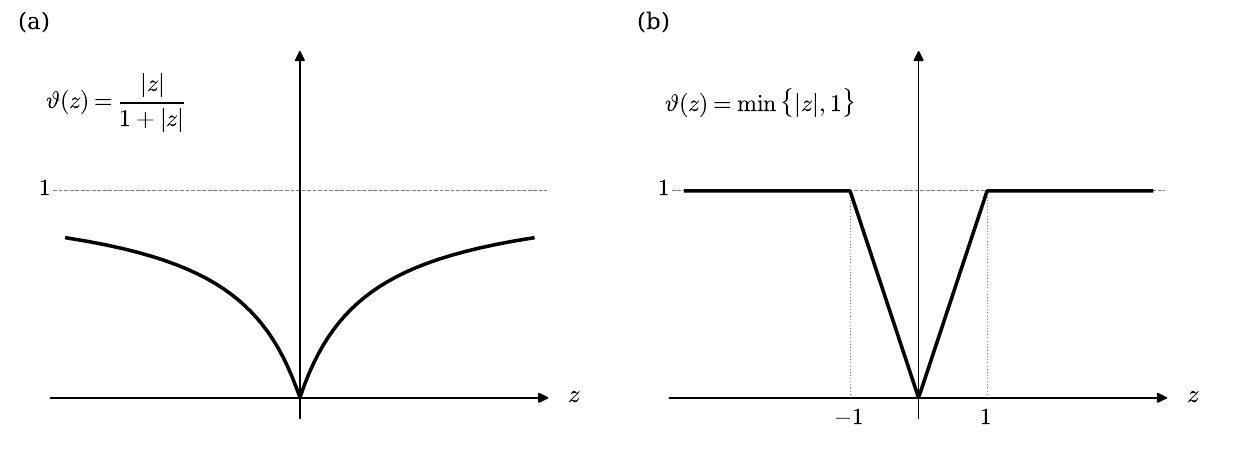}}
	\caption{The jump energy density $\vartheta$ in the two examples: (a) $\psi(x)=x^2$, for which $\vartheta(z)=\frac{|z|}{1+|z|}$; (b) $\psi=1-\chi_{\{0\}}$, for which $\vartheta(z)=\min\{|z|,1\}$. In both cases $\vartheta$ has slope $1$ at $0$ and is bounded by $1$.}\label{Fig28}
\end{figure}

The two corresponding functions $\vartheta$ are pictured in Fig.~\ref{Fig28}. Note that in both cases $\vartheta(z)\sim|z|$ as $z\to0$: the functionals $F_\e$ have linear growth in $u'$, and their $\Gamma$-limit is defined on the space of functions of bounded variation (more generally, of generalized functions of bounded variation, in order to allow for unbounded jumps) and is given by $\int|u'|dt+|D^cu|+\sum\vartheta(|u^+-u^-|)$, where $|D^cu|$ is the total variation of the Cantor part of the derivative of $u$, with the same coefficient $1$ as the absolutely continuous part \cite{ABS99}. The jump energy $\vartheta$ is concave, bounded, and with slope $1$ at $0$; this is not in contradiction with the superlinear behaviour at $0$ required in the previous section, since concave densities with infinite slope at $0$ are recovered from these by a further limit procedure, see the Appendix to this chapter.

\subsection{Appendix: relevant tools of Chapter 5}\label{App-Ch5}

\paragraph{References.} The Mumford--Shah functional was introduced in image segmentation by Mumford and Shah \cite{MS89}; its weak formulation in the space $SBV$ of special functions of bounded variation, and the proof of existence of minimizers for it, are due to De Giorgi and Ambrosio \cite{DGA88} (after the compactness and lower semicontinuity results of Ambrosio \cite{Amb89}); a complete treatment of $SBV$ and of free-discontinuity problems is given in the book by Ambrosio, Fusco and Pallara \cite{AFP}. The approximation by the phase-field functionals of Section~5.1 is due to Ambrosio and Tortorelli \cite{AT90,AT92}; the discrete approximation is based on the {\em weak membrane} model of Blake and Zisserman \cite{BZ87} and was studied by Chambolle \cite{Cha95}; the finite-difference approximation was conjectured by De Giorgi and proved by Gobbino \cite{Gob98}, and the non-local one is due to Braides and Dal Maso \cite{BDM97}. A systematic account of these and other approximations is given in the monograph \cite{Bra98}. The lower-semicontinuity and compactness theory for the functionals \eqref{fdf} is due to Ambrosio \cite{Amb89} and to Bouchitt\'e, Buttazzo and the author \cite{BBB95}, see also \cite{Bra95}. The approximation of such functionals by singular perturbations of the type of Section~5.2 is due to Alicandro, Braides and Gelli \cite{ABG98}, and its $d$-dimensional version to Alicandro and Gelli \cite{AG00}; the extension to perturbations of order $k>2$ is due to Solci \cite{Sol25}. The interaction of relaxation and phase-field discontinuities of Section~5.3 is studied by Alicandro, Braides and Shah \cite{ABS99}.

\paragraph{Special functions of bounded variation in dimension one.}
Let $(a,b)$ be a bounded interval. A function $u\in BV(a,b)$ belongs to $SBV(a,b)$ if its derivative is the sum of an absolutely continuous part and of a jump part,
$$
Du=u'\,dt+\sum_{t\in S(u)}(u^+(t)-u^-(t))\delta_t,
$$
that is, if there is no Cantor part; here $u'\in L^1(a,b)$ is the pointwise derivative and $S(u)$ is the (at most countable) set of jump points of $u$. If $\#S(u)<+\infty$ and $u'\in L^2$, then $u$ is a function which is $H^1$ on each of the finitely many intervals into which $S(u)$ divides $(a,b)$; the one-dimensional Mumford--Shah functional is
$$
F^{MS}(u)=\alpha\int_{(a,b)}|u'|^2dt+\beta\,\#S(u),
$$
and is $+\infty$ outside this class. We prove its compactness and lower semicontinuity.

\begin{proposition}\label{MS-comp}
Let $u_j\in SBV(a,b)$ be such that $\|u_j\|_{L^\infty}\le M$ and $F^{MS}(u_j)\le C$. Then, up to subsequences, $u_j$ converges in $L^1(a,b)$ to some $u$ with $F^{MS}(u)\le\liminf_jF^{MS}(u_j)$.
\end{proposition}

\begin{proof}
Up to passing to a subsequence we may assume that $\liminf_jF^{MS}(u_j)$ is a limit, that $\#S(u_j)=m$ for all $j$ (with $m\le C/\beta$), and that the ordered jump points $t^1_j<\dots<t^m_j$ converge to points $t^1\le\dots\le t^m$ in $[a,b]$. Denote by $T$ the set of the distinct points $t^i\in(a,b)$, so that $\#T\le m$. Let $I$ be a compact interval contained in $(a,b)\setminus T$; for $j$ large $I$ does not contain any $t^i_j$, so that $u_j\in H^1(I)$ with $\|u_j\|_{L^\infty}\le M$ and $\|u_j'\|_{L^2(I)}^2\le C/\alpha$. Hence $\{u_j\}$ is bounded in $H^1(I)$, and is therefore precompact in the uniform topology of $I$, with weak limits in $H^1(I)$. By a diagonal argument on an exhaustion of $(a,b)\setminus T$ by compact intervals we obtain a subsequence converging locally uniformly in $(a,b)\setminus T$ to a function $u$ with $u\in H^1_{\rm loc}((a,b)\setminus T)$, $|u|\le M$. By dominated convergence $u_j\to u$ in $L^1(a,b)$ (the set $T$ is negligible). By the weak lower semicontinuity of the Dirichlet integral, for every such $I$
$$
\int_I|u'|^2dt\le\liminf_{j}\int_I|u_j'|^2dt\le\liminf_j\int_{(a,b)\setminus S(u_j)}|u_j'|^2dt,
$$
and taking the supremum over $I$ we get $\int_{(a,b)\setminus T}|u'|^2dt\le\liminf_j\int_{(a,b)\setminus S(u_j)}|u_j'|^2dt$. Finally, $u$ can only be discontinuous at points of $T$ (it is continuous elsewhere), so that $S(u)\subset T$ and $\#S(u)\le\#T\le m=\#S(u_j)$. Summing the two inequalities gives the claim.
\end{proof}

The $L^\infty$ bound cannot be dropped: the functions $u_j=j\chi_{(a,a+1/j)}$ have bounded $L^1$ norm and $F^{MS}(u_j)=\beta$, but converge to a Dirac delta. In the general case of functionals \eqref{fdf} the same strategy applies, but the lower semicontinuity of the bulk and surface part requires, besides the convexity of $f$ and the subadditivity of $\vartheta$, an argument of {\em slicing/blow-up} type which we do not reproduce here, see \cite{Amb89,BBB95,Bra95}. Note that the assumption $\vartheta(z)/|z|\to+\infty$ as $z\to0$ guarantees that the Cantor part of the derivative is not created in the limit, which is the reason why the natural space is $SBV$.

\paragraph{The constant $C_2$.}
We compute the constant $C_2$ in \eqref{c2}. Fix $S>0$; the minimum problem
$$
\min\Big\{\int_0^S|w''|^2dt:\ w(0)=0,\ w(S)=1,\ w'(0)=w'(S)=0\Big\}
$$
is strictly convex, and its Euler--Lagrange equation is $w''''=0$, so that the minimizer is the cubic polynomial
$$
w(t)=3\Big(\frac tS\Big)^2-2\Big(\frac tS\Big)^3.
$$
Since $w''(t)=\frac6{S^2}\big(1-\frac{2t}S\big)$, we have
$$
\int_0^S|w''|^2dt=\frac{36}{S^4}\int_0^S\Big(1-\frac{2t}S\Big)^2dt=\frac{12}{S^3}.
$$
We are left with the minimization of $S+12S^{-3}$ over $S>0$: the derivative vanishes if and only if $S^4=36$, that is $S=\sqrt6$, and
$$
C_2=\sqrt6+\frac{12}{6\sqrt6}=\sqrt6+\frac{\sqrt6}3=\frac{4\sqrt6}3\approx3.27.
$$

\paragraph{Scaling for perturbations of order $k$.}
Let $F_\e^k$ be the functional with the perturbation $\e^{2k-1}\int|u^{(k)}|^2$. Consider a transition $v(t)=zw(t/\lambda)$ on an interval of length $\lambda S$, with $w(0)=0$, $w(S)=1$. On this interval the term $\min\{|v'|^2,1/\e\}$ equals $1/\e$ (if $|v'|^2\ge1/\e$), which gives the contribution $\lambda S/\e$, while the singular perturbation contributes
$$
\e^{2k-1}\int_0^{\lambda S}|v^{(k)}|^2dt=\e^{2k-1}z^2\lambda^{1-2k}\int_0^S|w^{(k)}|^2ds.
$$
The two contributions are of the same order if $\lambda/\e=\e^{2k-1}z^2\lambda^{1-2k}$, that is $\lambda^{2k}=\e^{2k}z^2$, i.e. $\lambda=\e|z|^{1/k}$. For this choice both of them are multiplied by $|z|^{1/k}$ and the total energy is
$$
|z|^{1/k}\Big(S+\int_0^S|w^{(k)}|^2ds\Big),
$$
which is the expression defining $m_k|z|^{1/k}$ in Section~5.2. For $k=2$ we recover the scale $\e\sqrt{|z|}$ and the density $C_2\sqrt{|z|}$.

\paragraph{The jump energy density of Section~5.3.}
The density $\vartheta$ is obtained by minimizing, for a jump of size $z$, the cost of the possible values of the phase variable on the jump layer,
$$
\vartheta(z)=\min_{x\in[0,1]}\big\{\psi(x)|z|+(1-x)^2\big\}.
$$
For $\psi(x)=x^2$ the minimum of $x\mapsto x^2|z|+(1-x)^2$ is attained when $2x|z|=2(1-x)$, that is at $x=1/(1+|z|)$, and
$$
\vartheta(z)=\frac{|z|+|z|^2}{(1+|z|)^2}=\frac{|z|}{1+|z|}.
$$
For $\psi=1-\chi_{\{0\}}$ the value $x=0$ gives $1$, while for $x\in(0,1]$ we get $|z|+(1-x)^2$ whose minimum is $|z|$, attained at $x=1$; hence $\vartheta(z)=\min\{|z|,1\}$. In general, if $\psi\ge0$ is lower semicontinuous with $\psi(1)=1$:
\begin{itemize}
\item $\vartheta$ is the infimum of the family of functions $|z|\mapsto\psi(x)|z|+(1-x)^2$, which are affine and nondecreasing; hence $\vartheta$ is concave and nondecreasing in $|z|$, and $\vartheta(0)=0$ (take $x=1$);
\item choosing $x=1$ we have $\vartheta(z)\le|z|$. If $x_z$ is a minimizer, then $(1-x_z)^2\le\vartheta(z)\le|z|$ so that $x_z\to1$ as $z\to0$; by lower semicontinuity $\liminf\psi(x_z)\ge\psi(1)=1$ and thus $\vartheta(z)/|z|\to1$ as $z\to0$;
\item concavity and $\vartheta(0)=0$ give $\vartheta(\lambda r)\ge\lambda\vartheta(r)$ for $\lambda\in[0,1]$, hence $\vartheta(r)/r$ is nonincreasing, and for $r,s>0$
$$
\vartheta(r+s)=\frac{r}{r+s}\vartheta(r+s)+\frac{s}{r+s}\vartheta(r+s)\le\vartheta(r)+\vartheta(s),
$$
that is, $\vartheta$ is subadditive.
\end{itemize}
In particular the slope of $\vartheta$ at $0$ is $1$, so that $\vartheta(z)/|z|\not\to+\infty$ as $z\to0$: the superlinearity condition required in Section~5.2 is not satisfied, and the energies of Section~5.3 must be considered on the larger space $BV$ (or $GBV$), where a Cantor part of the derivative is allowed.

Let now $\psi:[0,1]\to[0,1]$ be increasing and lower semicontinuous, with $\psi(0)=0$, $\psi(1)=1$ and $\psi>0$ on $(0,1]$ (which is the situation of Section~5.3, and implies the properties above). The theorem of Alicandro, Braides and Shah \cite{ABS99} states that the functionals
$$
F_\e(u,v)=\int_{(a,b)}\Big(\psi(v)|u'|+\frac1\e W(v)+\e|v'|^2\Big)dt,
$$
with $W$ continuous, $W\ge0$ and $W(x)=0$ only for $x=1$, defined on $H^1$ functions with $0\le v\le1$, $\Gamma$-converge with respect to the $L^1\times L^1$ convergence to the functional
$$
F(u,v)=\int_{(a,b)}|u'|dt+|D^cu|(a,b)+\sum_{t\in S(u)}\vartheta(|u(t^+)-u(t^-)|)
$$
if $u\in GBV(a,b)$ and $v=1$ almost everywhere, and $+\infty$ otherwise. Here $GBV$ is the space of functions whose truncations are in $BV$, $|D^cu|$ is the supremum over the truncations of the total variation of the Cantor part, and
$$
\vartheta(z)=\min_{x\in[0,1]}\big\{\psi(x)z+2\big(c_W-c_W(x)\big)\big\},\qquad c_W(x)=2\int_0^x\sqrt{W(s)}\,ds,
$$
For $W(s)=\frac14(1-s)^2$, as in Section~5.3, we have $c_W(x)=x-\frac12x^2$ and $c_W=c_W(1)=\frac12$, so that $2(c_W-c_W(x))=(1-x)^2$, which gives the formula for $\vartheta$ used there: the term $c_W-c_W(x)=\frac12(1-x)^2$ is the cost of one transition layer from $1$ to $x$, computed by the Modica--Mortola inequality $\frac1\e W(v)+\e|v'|^2\ge2\sqrt{W(v)}|v'|$. The properties of $\vartheta$ listed above are those of \cite{ABS99}: $\vartheta$ is increasing, subadditive, Lipschitz with constant $1$, $\vartheta(z)\le z$, $\vartheta(z)/z\to1$ as $z\to0$ and $\vartheta(z)\to2c_W$ as $z\to+\infty$, where $c_W=c_W(1)=2\int_0^1\sqrt{W(s)}\,ds$ is the constant of the previous chapters.

The same paper shows how to obtain, by a double limit, the free-discontinuity functionals of Section~5.2. If the energy $\psi(v)|u'|$ is replaced by $\psi(v)f(|u'|)$, with $f$ convex, increasing and with $f(t)/t\to1$ as $t\to+\infty$, the limit is the same with $\int f(|u'|)dt$ in place of $\int|u'|dt$. Moreover, every concave $\vartheta$ which is the minimum of a finite number of affine functions $a_iz+b_i$ can be obtained by choosing $\psi$ piecewise constant (and multiplying the phase-field terms by a constant). Taking finer and finer piecewise-affine approximations of a concave $\vartheta$ with $\vartheta(z)/z\to+\infty$ as $z\to0$, and of a convex superlinear $\varphi$, one obtains a family of functionals $G^j_\e$ with
$$
\lim_{j\to+\infty}\Big(\Gamma\hbox{-}\lim_{\e\to0}G^j_\e\Big)(u,v)=\int_{(a,b)}\varphi(|u'|)dt+\sum_{t\in S(u)}\vartheta(|u(t^+)-u(t^-)|)
$$
on $GSBV(a,b)$, $v=1$, so that every functional of the form \eqref{fdf} is a double limit of phase-field functionals of the type of Section~5.3.

\section{Periodicity by concentration}

In this section we explore some examples in which lower-order terms may interact with concentrating variables or free discontinuities, generating periodic patterns. The oscillations and the concentration effect will occur at different scales showing a {\bf multi-scale behaviour}.

\subsection{A perturbation of the Bolza example}
We consider a double-well potential $W$ with wells at $-1$ and $1$, and study the minimization of the functional
$$
\int_0^1(W(u')+|u|^2+\e^2|u''|^2)dt.
$$
In the case $\e=0$ this is Bolza's example, in which no minimizer exists and minimizing sequences oscillate at no precise scale. The introduction of $\e>0$ makes the problem coercive (and lower semicontinuous) in $H^2$, and minimizers do exist. We will now argue that they also have a precise scale of oscillations (this problem has been studied in detail by M\"uller \cite{Mu93} and by Alberti and M\"uller \cite{AM01}).

\begin{figure}[h!]
	\centerline{\includegraphics [width=2.5in]{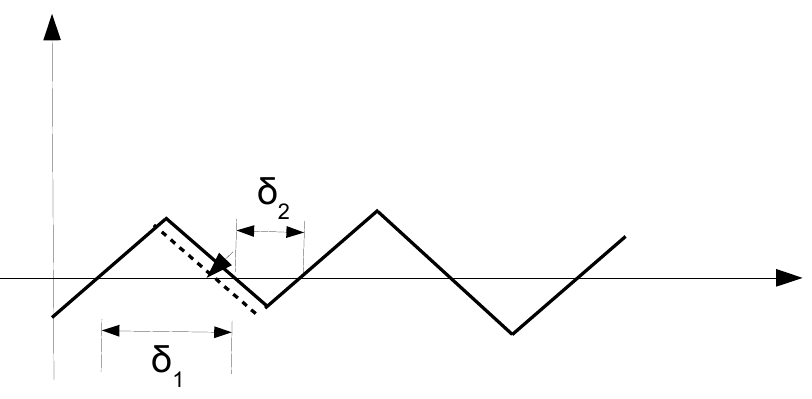}}
	\caption{A test function for the perturbed Bolza problem}\label{Fig9}
\end{figure}

In order to simplify the problem, we write the term $\int(W(u')+\e^2|u''|^2)dt$ as $\e\cdot\frac1\e\int(W(u')+\e^2|u''|^2)dt$ and substitute the second factor with its limit $c_W\#(S(u'))$ as $\e\to 0$, given by the Modica--Mortola theory applied to $u'$ (see the Appendix to this chapter), with $c_W=2\int_{-1}^1\sqrt{W(s)}ds$. This is not justified analytically a priori, but gives us a more tractable family of functionals
$$
\int_0^1|u|^2dt +\e c_W\#(S(u'))\hbox{ with } |u'|=1 \hbox{ almost everywhere}.
$$
A test function for this functional is pictured in Fig.~\ref{Fig9}. The picture highlights that if we have oscillations at two different periods $\delta_1$ and $\delta_2$, then moving the two neighbouring angles lowers the energy. We can then limit ourselves to functions with $2\delta$-periodic oscillations, and, taking into account that the number of discontinuities of $u'$ is approximately $1/\delta$, for such a function the energy is (up to small errors due to boundary approximations) equal to 
$$
\frac{\delta^2}{12}+c_W \frac{\e}\delta,
$$
which is minimized for 
$
\delta^3=6c_W \e
$; that is, the period is approximately 
$$
2\delta= (48c_W \e)^{1/3}.
$$
This simple approximate analysis suggests a {\bf periodic scale triggered by concentration}, and that minimizers of the perturbed Bolza example show a {\bf concentration effect on a scale $\e$ modulated periodically on a scale $\e^{1/3}$}.
Note that even in this simplified case this superposition of effects is not easily described in terms of a $\Gamma$-limit.

\subsection{Staircasing}
We want to describe the scaling properties of the solutions of the minimum problem 
\begin{equation}\label{stairs}
\min\Big\{ \int_0^1(\min\{|u'|^2,1\}+ |u-at|^2 + \e^4|u''|^2)dt\Big\}.
\end{equation}
In order to present a simplified approximate version of this problem, we examine
the asymptotic behaviour of 
$$
F_\e(u)= \frac1\e\int_0^1\min\{|u'|^2,1\}dt+\e^3\int_0^1|u''|^2dt.
$$
Since $F_\e\ge F^\alpha_\e$ for every $\alpha>0$ and $\e<1/\alpha$, by comparison with the functionals 
$$
F^\alpha_\e(u)= \int_0^1\min\Big\{\alpha|u'|^2,\frac1\e\Big\}dt+\e^3\int_0^1|u''|^2dt,
$$
we obtain a lower bound with 
$$
F^\alpha(u)= \alpha\int_0^1|u'|^2dt+C_2\sum_{t\in S(u)}\sqrt{|u(t^+)-u(t^-)|}\qquad u\in SBV(0,1),
$$
with $C_2$ independent of $\alpha$ and in particular obtained by formula \eqref{c2} with $\alpha=1$. By taking the supremum over all $\alpha>0$ we argue that the $\Gamma$-limit is given by
$$
F(u)= C_2\sum_{t\in S(u)}\sqrt{|u(t^+)-u(t^-)|}\quad u\in SBV(0,1) \hbox{ and } u'=0 \hbox{ almost everywhere},
$$
recovery sequences for piecewise-constant target $u$ being constructed exactly as in the case $\alpha=1$. 

We can then provide a simplification of the problem in \eqref{stairs} as 
\begin{eqnarray}\label{stairs-2}\nonumber
&&\min\Big\{ \e C_2\sum_{t\in S(u)}\sqrt{|u(t^+)-u(t^-)|}+ \int_0^1|u-at|^2 dt\Big\}\\
&=&a^2\min\Big\{ \e a^{-3/2}C_2\sum_{t\in S(u)}\sqrt{|u(t^+)-u(t^-)|}+ \int_0^1|u-t|^2 dt\Big\}.
\end{eqnarray}

\begin{figure}[h!]
	\centerline{\includegraphics [width=2.5in]{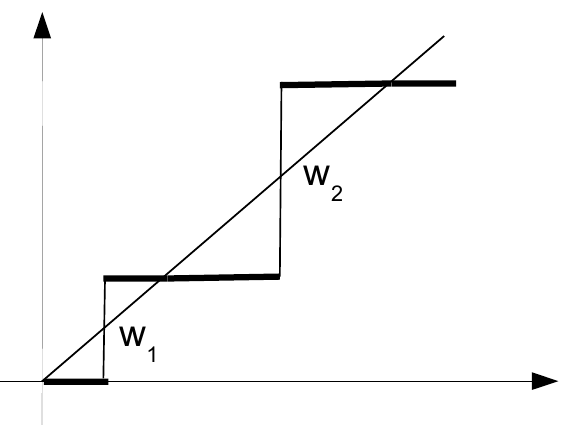}}
	\caption{A test function of `staircase' type}\label{Fig10}
\end{figure}

Assuming that $u(0)\sim 0$ and $u(1)\sim 1$ (neglecting boundary effects), noting that $u$ is increasing, the minimum is determined by its jumps of size $\{w_i:i\in\{1,\ldots,N\}\}$, each of which gives a term due to the integral part equal to $w_i^3/12$. The corresponding energy is then given by
$$
a^2\sum_{i=1}^N\Big(\e \frac{C_2}{a^{3/2}} \sqrt{w_i}+\frac{w_i^3}{12}\Big), 
$$
with $\sum_iw_i=1$. The function $f(z)= K \sqrt{z}+\frac{z^3}{12}$, with $K=\e C_2a^{-3/2}$, defined for $z\ge 0$, is concave close to $0$ and then becomes convex, the inflection point depending on $\e$. If $\e$ is small, the optimal $w_i$ are all equal to a value $\delta$ and $N=1/\delta$: the optimal periodic distribution of jump points is obtained by minimizing the energy per unit length $f(\delta)/\delta$, that is
$$
\frac{K}{\sqrt{\delta}}+\frac{\delta^2}{12}, 
$$
obtaining $3\e C_2a^{-3/2}=\delta^{5/2}$ (geometrically, $\delta$ is the point of tangency of the graph of $f$ with the line through the origin); that is, 
$$
\delta=(3 C_2)^{2/5}a^{-3/5}\e^{2/5}.
$$
(If instead $\e$ is so large that $\delta>1$, no jump is convenient and the minimizer is the constant $u=a/2$; this is outside the regime $\e\to0$ we are interested in.) 
Recalling that a jump of size $z$ is a concentration at scale $\e\sqrt z$, and that here the jumps of $u$ have size $a\delta\sim\e^{2/5}$, we conclude that this analysis suggests, for the shape of minimizers of problem \eqref{stairs}, concentration on jump points at scale $\e^{6/5}$, periodically distributed with a period at scale $\e^{2/5}$. With respect to the previous example, this one shows a dependence of the period on the variable $a$: the period is proportional to $a^{-3/5}$: jump points are denser (the concentration occurs at a {\bf higher rate}) for large values of $a$. A rigorous quantitative analysis of the staircasing phenomenon for minimizers of this type of second-order regularizations of the Perona--Malik functional is due to Gobbino and Picenni \cite{GP23,GP25}; see also \cite{BF08} for the $\Gamma$-limit of related functionals, and \cite{BF26} for perturbations of higher order.

\subsection{Periodic micro-fracture selected by the strain}
In this last example the periodic pattern is generated by the interaction of an elastic variable $u$ with a second variable $v$ having free discontinuities; the jumps of $v$ are the concentration phenomenon, and their spacing is selected through a cell problem which depends on the macroscopic strain, rather than being given by a power of $\e$. The functionals are in the spirit of the Mumford--Shah functionals of the previous chapter, but with a small coefficient of the jump part. 

The functional we consider is defined on pairs $u,v\colon(0,1)\to \mathbb R$ with $u\in H^1(0,1)$ and $v\in SBV(0,1)$. More precisely, $v$ is a piecewise-$H^1$ function, and the energy is 
$$
F_\e(u,v)=\int_0^1 |u'|^2dt+ \int_0^1 |v'|^2dt+\e\#(S(v))+ \frac1{\e^2} \int_0^1 |u-v|^2dt.
$$
From a standpoint similar to that of thin films, this can be considered as the energy of a thin object with only two layers, one with a quadratic dependence on the variable $u$ (a linearly elastic medium) and the second one a Mumford--Shah functional depending on the variable $v$ with a small coefficient of the jump part. The two layers interact through a quadratic term in the `vertical' difference quotient $\frac{u(t)-v(t)}\e$, which can be seen as a gradient term. The choice of the scaling of the coefficient of the discontinuity part is such that the discontinuity part is of the same order as the interaction part on optimizing sequences.

The coerciveness properties of the functional $F_\e$ derive from its dependence on $u$. Indeed, sequences $u_\e$ such that $F_\e(u_\e,v_\e)$ is equibounded for some $v_\e$ have equi-bounded gradients in $L^2(0,1)$ so that they are precompact, up to addition of constants. If $u_\e\to u$ then the last term ensures that also $v_\e\to u$ in $L^2(0,1)$, even though for $v_\e$ we do not directly have sufficient information on their gradients.

Since the domain is $H^1(0,1)$ and compactness is only determined by the behaviour with respect to $u$, we may minimize out the dependence on $v$ and consider equivalently
$$
F_\e(u)=\int_0^1 |u'|^2dt+ \min_v\Big\{\int_0^1 |v'|^2dt+\e\#(S(v))+ \frac1{\e^2} \int_0^1 |u-v|^2dt\}.
$$

For fixed $u$, the inner minimum over $v$ is itself a scalar Mumford--Shah-type problem with an elastic forcing term $\frac1{\e^2}|u-v|^2$, penalizing jumps at a rate $\e$ per jump. The optimal strategy is to place the jumps of $v$ \emph{periodically}: if $u$ is close to affine with slope $z$ over a window of unit order, it does not pay to space the jumps irregularly, and the optimal period, rescaled by $\e$, converges to some $S=S(z)>0$ (possibly $+\infty$, meaning no jumps at all). Scaling the cell $[0,\e S]$ back to unit size, the cost per unit length of one period of this periodic pattern is, to leading order,
$$
\psi(S,z)=\frac1S\min\Big\{\int_0^S\big((u')^2+(v')^2+(u-v)^2\big)\,dt\ :\ u,v\in H^1(0,S),\ u(0)=0,\ u(S)=Sz\Big\}
$$
(no boundary condition is imposed on $v$, which is free to relax at the ends of the cell), plus a fracture cost $1/S$ for the one jump of $v$ per period. Optimizing over the period itself we get
$$
F_\e(u,v)\ \ \Gamma\hbox{-converges to}\ \ F_0(u)=\int_0^1 g(u')\,dt,\qquad g(z)=\inf_{S>0}\Big\{\psi(S,z)+\frac1S\Big\}.
$$
This is a particular case of a general homogenization phenomenon for energies coupling an elastic variable to a brittle one through a vanishing-scale interaction, studied by Braides, Causin and Solci \cite{BCS19}:
{\bf elastic  micro-oscillation and micro-fracture compete through a single auxiliary length scale} $S$, optimized at each macroscopic strain $z$.

\begin{figure}[h!]
	\centerline{\includegraphics [width=4.8in]{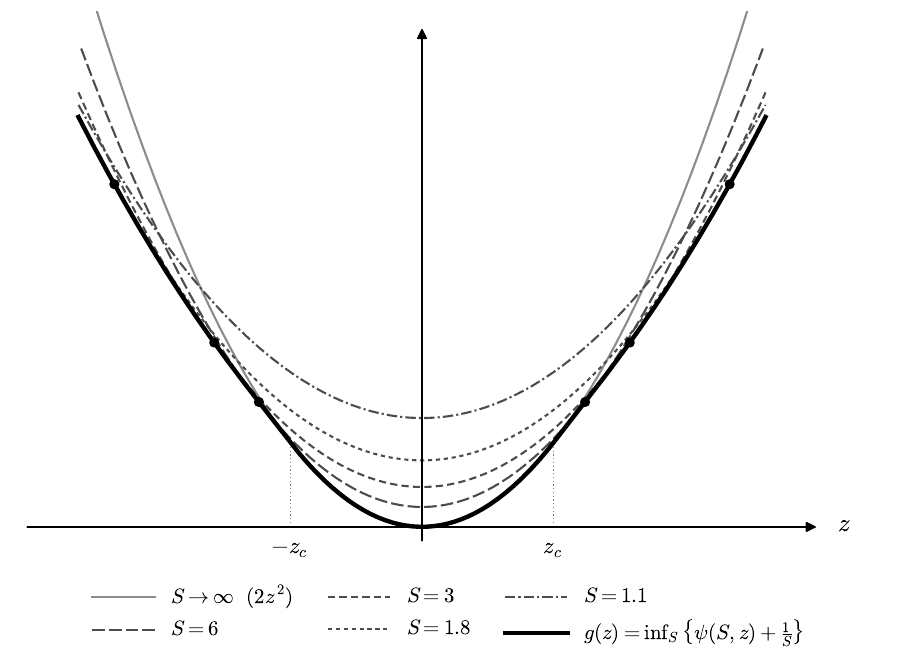}}
	\caption{The function $g$ as the lower envelope of the functions $z\mapsto\psi(S,z)+\frac1S$ for some values of $S$ (the thin lines) and for $S\to+\infty$ (the parabola $2z^2$). The dots mark the points of tangency, where $S=S(z)$ is the optimal period. For $|z|\le z_c$ the envelope coincides with $2z^2$; for larger $|z|$ the optimal period decreases with $|z|$.}\label{Fig29}
\end{figure}

Note that, by the quadratic dependence on $u$ of the cell problem, $\psi(S,z)=\kappa(S)z^2$ for a function $\kappa$ of $S$ only (given explicitly below), so that $g$ is the lower envelope of the family of parabolas $z\mapsto\kappa(S)z^2+\frac1S$, $S>0$, and of the limit parabola $2z^2$ for $S\to+\infty$; see Fig.~\ref{Fig29}.

Here the cell problem can be solved explicitly (see the Appendix to this chapter). The Euler--Lagrange system for $\psi(S,z)$ is $u''=u-v$, $v''=v-u$ (together with the natural condition $v'(0)=v'(S)=0$), so that $u+v$ is affine while $w=u-v$ solves $w''=2w$; imposing the boundary and natural conditions, and observing that at a solution the energy reduces (after an integration by parts using the equations themselves) to the boundary term $Sz\,u'(S)$, gives the explicit formula
$$
\psi(S,z)=\frac{2Sz^2}{S+\sqrt2\,\tanh\big(\frac{\sqrt2}2\,S\big)}.
$$
One can then check that $g$ is even and strictly convex, and that:
\begin{itemize}
\item[$\bullet$] for $|z|\le z_c:=2^{-3/4}$, the infimum defining $g(z)$ is approached only as $S\to\infty$: it never pays to fracture, and $g(z)=2z^2$ is simply the energy of the two layers elastically glued together;
\item[$\bullet$] for $|z|>z_c$, the infimum is attained at a finite, unique optimal period $S(z)$, strictly decreasing in $|z|$; as $z\to\infty$,
$$
S(z)\sim \sqrt[3]{6}\,|z|^{-2/3},\qquad\qquad g(z)\sim z^2+\frac{\sqrt[3]{36}}{4}\,|z|^{2/3}.
$$
\end{itemize}
Below the critical strain $z_c$, then, the coupled two-layer system behaves macroscopically as a single ordinary elastic medium; above it, the brittle layer spontaneously fragments into finer and finer periodic micro-cracks as $|z|\to\infty$, with the crack density $1/S(z)$ growing like $|z|^{2/3}$ and the resulting excess energy only like $|z|^{2/3}$ rather than $|z|^2$. In this example a genuine geometric period adjusts itself to the imposed macroscopic strain. Compared with the two previous examples, the period is not given by a power of $\e$ with an unknown constant, but is exactly $\e S(z)$, with $S$ determined by the cell problem and depending on the macroscopic strain $z$; the decrease $S(z)\sim|z|^{-2/3}$ of the period for large strains is analogous to the dependence of the period on the variable $a$ in the staircasing problem.

\subsection{Appendix: relevant tools of Chapter 6}\label{App-Ch6}

\paragraph{References.} The perturbation of the Bolza example by a second-order term, and the selection of periodic minimizers at scale $\e^{1/3}$, are due to M\"uller \cite{Mu93}, who proves that minimizers are periodic at this scale in a suitable sense; the Young-measure description of the limit micro-patterns, which captures the sawtooth structure and its period, is due to Alberti and M\"uller \cite{AM01}. The one-dimensional Modica--Mortola estimate used below is from \cite{MM77}, and the original example of Bolza can be found in \cite{Bolza,BGH}. The staircasing phenomenon for second-order regularizations of the Perona--Malik functional $\int\log(1+|u'|^2)$, of which \eqref{stairs} is a model case with a truncated quadratic energy, has been analyzed in detail by Gobbino and Picenni \cite{GP23}, who describe the minimizers at the scale of the steps, and by the same authors at finer scales in \cite{GP25}; see also Bellettini and Fusco \cite{BF08} for the $\Gamma$-limit of the corresponding functionals and Braides and Fonseca \cite{BF26} for perturbations of order $k\ge2$. The homogenization of the elastic-brittle model of the last example is due to Braides, Causin and Solci \cite{BCS19}.

\paragraph{The transition cost of the derivative.}
Let $W:\mathbb R\to[0,+\infty)$ be continuous with $W(s)=0$ if and only if $s=\pm1$, and set $c_W=2\int_{-1}^1\sqrt{W(s)}\,ds$.

\begin{lemma}\label{MM-trans}
Let $\e>0$ and let $v\in H^1(\alpha,\beta)$. Then
$$
\int_\alpha^\beta\big(W(v)+\e^2|v'|^2\big)dt\ge2\e\Big|\int_{v(\alpha)}^{v(\beta)}\sqrt{W(s)}\,ds\Big|.
$$
In particular, if $v(\alpha)\le-1+\eta$ and $v(\beta)\ge1-\eta$, then the left-hand side is not less than $\e c_W-2\e\int_{[-1,-1+\eta]\cup[1-\eta,1]}\sqrt W$, which is $\e c_W+o(\e)$ as $\eta\to0$.
\end{lemma}

\begin{proof}
By Young's inequality $a^2+b^2\ge2ab$ applied with $a=\sqrt{W(v)}$ and $b=\e|v'|$ we have $W(v)+\e^2|v'|^2\ge2\e\sqrt{W(v)}|v'|$, and by the change of variables $s=v(t)$ (valid for absolutely continuous $v$) $\int_\alpha^\beta\sqrt{W(v)}|v'|\,dt\ge\big|\int_{v(\alpha)}^{v(\beta)}\sqrt{W(s)}\,ds\big|$.
\end{proof}

Equality is (approximately) achieved by the heteroclinic profile solving $\e v'=\sqrt{W(v)}$, as in the one-dimensional Modica--Mortola theorem of Chapter 3. Applying this to $v=u'$ on each interval where $u'$ goes from $-1$ to $1$ or conversely, we get that the energy $\int(W(u')+\e^2|u''|^2)dt$ is at least $\e c_W$ times the number of such transitions (up to $o(\e)$), and that $\frac1\e\int(W(u')+\e^2|u''|^2)dt$ $\Gamma$-converges to $c_W\#S(u')$ on the functions with $u'\in\{\pm1\}$ almost everywhere. This justifies the reduction made in the text.

\paragraph{An optimal-partition lemma.} Both examples rely on the same elementary argument.

\begin{lemma}\label{partition}
Let $g:(0,+\infty)\to[0,+\infty)$ be a function, and set $h_*=\inf_{\delta>0}g(\delta)/\delta$. Then for every finite family $\delta_1,\dots,\delta_N>0$ with $\sum_i\delta_i=1$ we have $\sum_ig(\delta_i)\ge h_*$, and equality holds (when $1/\delta_*$ is an integer) for $N=1/\delta_*$ and $\delta_i=\delta_*$, if $\delta_*$ is a minimizer of $g(\delta)/\delta$.
\end{lemma}

\begin{proof}
Since $g(\delta_i)\ge h_*\delta_i$ for every $i$, summing we get $\sum_ig(\delta_i)\ge h_*\sum_i\delta_i=h_*$. If $\delta_*$ is a minimizer and $N=1/\delta_*$ then $\sum_ig(\delta_*)=N g(\delta_*)=g(\delta_*)/\delta_*=h_*$.
\end{proof}

Geometrically, $h_*$ is the slope of the line through the origin which is tangent from below to the graph of $g$.

\paragraph{The period of the perturbed Bolza example.}
Let $u$ be such that $|u'|=1$ almost everywhere, with $N-1$ corners (points of $S(u')$), at which the derivative changes sign, dividing $(0,1)$ into $N$ intervals of lengths $\delta_1,\dots,\delta_N$ ($\sum_i\delta_i=1$). On the $i$-th interval $u$ is affine with slope $\pm1$, so that $u(t)=c_i\pm(t-t_i)$, with $t_i$ the midpoint, and
$$
\int_{I_i}u^2dt=\int_{-\delta_i/2}^{\delta_i/2}(c_i+s)^2ds=\delta_ic_i^2+\frac{\delta_i^3}{12}\ge\frac{\delta_i^3}{12},
$$
with equality if and only if $c_i=0$, that is, $u$ vanishes at the midpoint of the interval. Hence
$$
\int_0^1u^2dt+\e c_W\#S(u')\ge\sum_{i=1}^N\Big(\frac{\delta_i^3}{12}+\e c_W\Big)-\e c_W\ge\min_{\delta>0}\Big(\frac{\delta^2}{12}+\frac{\e c_W}{\delta}\Big)-\e c_W
$$
by Lemma~\ref{partition} applied to $g(\delta)=\frac{\delta^3}{12}+\e c_W$. The function $\frac{\delta^2}{12}+\frac{\e c_W}\delta$ is minimized for $\frac\delta6=\frac{\e c_W}{\delta^2}$, that is $\delta_*^3=6\e c_W$, and the minimum is $\frac{\delta_*^2}{12}+\frac{\e c_W}{\delta_*}=\frac{\delta_*^2}4$. This bound is attained, up to boundary terms, by the zigzag function with all corners at distance $\delta_*$ and with zeros at all midpoints of the intervals (which is a $2\delta_*$-periodic sawtooth function, as in Fig.~\ref{Fig9}). The corresponding period is $2\delta_*=(48c_W\e)^{1/3}$, and the minimal energy is $\frac14(6c_W\e)^{2/3}+o(\e^{2/3})$, of order $\e^{2/3}$.

\paragraph{The $\Gamma$-limit of $F_\e$.}

\begin{proposition}\label{Gamma-stairs}
The functionals $F_\e(u)=\frac1\e\int_0^1\min\{|u'|^2,1\}dt+\e^3\int_0^1|u''|^2dt$, defined on $H^2(0,1)$, $\Gamma$-converge, with respect to the $L^1(0,1)$ convergence, to
$$
F(u)=C_2\sum_{t\in S(u)}\sqrt{|u(t^+)-u(t^-)|}
$$
if $u\in SBV(0,1)$ with $u'=0$ almost everywhere, and $+\infty$ otherwise, where $C_2$ is the constant in \eqref{c2}.
\end{proposition}

\begin{proof}
{\em Lower bound.} Let $\alpha>0$ and $\e<1/\alpha$. Then $\frac1\e\min\{|\xi|^2,1\}=\min\{\frac{|\xi|^2}\e,\frac1\e\}\ge\min\{\alpha|\xi|^2,\frac1\e\}$ for every $\xi$, so that $F_\e\ge F^\alpha_\e$. The functionals $F^\alpha_\e$ differ from those of Section~5.2 only by the factor $\alpha$ in front of the quadratic part. This factor modifies the threshold $|u'|\ge1/\sqrt{\alpha\e}$ defining the intervals on which the gradient is large; on the endpoints of such an interval the derivative is of order $\e^{-1/2}$, which is negligible with respect to the order $\e^{-1}\sqrt z$ of the derivative of a transition giving a jump of size $z$ (at scale $\e\sqrt z$). Hence the argument of Section~5.2 (see \cite{ABG98}) applies without changes and gives the lower bound $\Gamma\hbox{-}\liminf_{\e\to0}F^\alpha_\e\ge F^\alpha$, with $F^\alpha(u)=\alpha\int|u'|^2dt+C_2\sum_{t\in S(u)}\sqrt{|u(t^+)-u(t^-)|}$ and the same constant $C_2$. Taking the supremum over $\alpha>0$ gives $\Gamma\hbox{-}\liminf F_\e\ge F$, with $F=+\infty$ if $u'\not\equiv0$.

{\em Upper bound.} Let $u$ be piecewise constant with finitely many jumps. We use the recovery sequence \eqref{ABG-ue} in the neighbourhood of each jump $t_i$, with $u_\e=u$ elsewhere. In the transition layer of $t_i$, which has length $\e\sqrt{z_i}S$, we have $\frac1\e\min\{|u_\e'|^2,1\}\le\frac1\e$, while this term vanishes outside the layers, and the energy of the layer is at most $\sqrt{z_i}\big(S+\int_0^S|w''|^2ds\big)=C_2\sqrt{z_i}$ for the optimal pair $(S,w)$. If $u$ has infinitely many jumps, we conclude by the lower semicontinuity of the $\Gamma$-limsup and by the density of piecewise constant functions with finitely many jumps in the energy, since $\sum_i\sqrt{z_i}<+\infty$.
\end{proof}

\paragraph{The staircase.}
Let $K>0$ and $f(z)=K\sqrt z+\frac{z^3}{12}$, so that the energy of the staircase problem \eqref{stairs-2} is, after the rescaling $u=a\tilde u$,
$$
a^2\sum_{i=1}^Nf(w_i),\qquad K=\e C_2a^{-3/2},\ \ \sum_iw_i=1,
$$
neglecting boundary effects. The minimum of $f(z)/z=Kz^{-1/2}+\frac{z^2}{12}$ is attained at the only point $\delta_*$ at which $-\frac K2\delta_*^{-3/2}+\frac{\delta_*}6=0$, that is $\delta_*^{5/2}=3K$, and equals $h_*=\frac{K}{\sqrt{\delta_*}}+\frac{\delta_*^2}{12}=\frac{5}{12}\delta_*^2$. By Lemma~\ref{partition} the minimal energy is $a^2h_*=\frac5{12}a^2(3K)^{4/5}=\frac5{12}(3C_2a\e)^{4/5}$, and it is attained by a staircase with equal steps of length $\delta_*=(3C_2)^{2/5}a^{-3/5}\e^{2/5}$. In particular the minimal energy is of order $(a\e)^{4/5}$ and depends on $a$ and $\e$ only through the product $a\e$. If $\delta_*>1$ the infimum is attained by the constant $u=a/2$, with energy $a^2/12$. The transition layers have size $\e\sqrt{a\delta_*}$, which scales as $\e^{6/5}$.

\paragraph{Plateaus with a positive slope.}
In the simplified model we have supposed that $u$ is constant between jumps. This is justified to leading order. Consider a step of $\tilde u$ of length $\delta$ with slope $\sigma$ instead of constant: the best affine function on it gives the energy
$$
a^2\delta\Big(\frac{(1-\sigma)^2\delta^2}{12}+\sigma^2\Big),
$$
where the first term is the tracking error and the second one is the gradient energy (here $|u'|=a\sigma<1$, so that the truncation plays no role). Minimizing in $\sigma$ with $m=\delta^2/12$ we get $\sigma=m/(1+m)$ and the value $a^2\delta\frac m{1+m}$, which differs from the value $a^2\delta m$ for $\sigma=0$ by a relative factor $1+O(\delta^2)=1+O(\e^{4/5})$. The slope is therefore of order $\delta^2\sim\e^{4/5}$, and neglecting it does not affect the scaling laws above.

\paragraph{The cell problem of the elastic-brittle model.}
We first check the scaling of the cell energy. If $u$ is close to an affine function of slope $z$, with the jumps of $v$ at distance $\e S$, we write on one period $u(t)=\e\tilde u(t/\e)$ and $v(t)=\e\tilde v(t/\e)$, with $\tilde u,\tilde v$ defined on $(0,S)$. Then
$$
\int|u'|^2dt=\e\int|\tilde u'|^2ds,\quad \int|v'|^2dt=\e\int|\tilde v'|^2ds,\quad \frac1{\e^2}\int|u-v|^2dt=\e\int|\tilde u-\tilde v|^2ds,
$$
so that the energy of one period is $\e\big(\int_0^S(|\tilde u'|^2+|\tilde v'|^2+|\tilde u-\tilde v|^2)ds+1\big)$, where the last term accounts for the jump of $v$ at the end of the period. Dividing by the length $\e S$ of the period gives the energy per unit length $\psi(S,z)+\frac1S$, with $\psi$ as in the text (the boundary condition $\tilde u(S)-\tilde u(0)=Sz$ expresses that $u$ has slope $z$).

\paragraph{Computation of $\psi(S,z)$.}
The Euler--Lagrange equations of the minimum problem defining $\psi(S,z)$ are $u''=u-v$, $v''=v-u$, with $u(0)=0$, $u(S)=Sz$ and the natural conditions $v'(0)=v'(S)=0$. The sum $p=u+v$ is affine, $p(t)=p_0+ct$, and $w=u-v$ solves $w''=2w$. Since $v'=\frac12(p'-w')$, the natural conditions give $w'(0)=w'(S)=c$. Setting $k=\sqrt2$ we obtain
$$
w(t)=\frac ck\Big(\sinh kt-\tanh\frac{kS}2\cosh kt\Big),\qquad w(0)+w(S)=0.
$$
The boundary conditions on $u=\frac12(p+w)$ then read $p_0+w(0)=0$ and $p_0+cS+w(S)=2Sz$, which give $cS+\frac{2c}k\tanh\frac{kS}2=2Sz$, that is
$$
c=\frac{2Sz}{S+\sqrt2\tanh\big(\frac{\sqrt2}2S\big)}.
$$
Multiplying the equations by $u$ and $v$, respectively, integrating by parts and adding, we obtain that the minimal value of $\int_0^S(|u'|^2+|v'|^2+(u-v)^2)dt$ equals the boundary term $u(S)u'(S)=Sz\,u'(S)$, and $u'(S)=\frac12(c+w'(S))=c$. Hence $\psi(S,z)=zc$, which is the formula in the text.

\paragraph{The critical strain and the large-strain asymptotics.}
Set $\tau=\tanh\big(\frac{\sqrt2}2S\big)\in(0,1)$. A direct computation gives
$$
\psi(S,z)+\frac1S-2z^2=\frac{S+\sqrt2\,\tau-2\sqrt2\,z^2\tau S}{S\,(S+\sqrt2\,\tau)}.
$$
If $2\sqrt2\,z^2\le1$, that is $|z|\le z_c=2^{-3/4}$, the numerator is not less than $S+\sqrt2\,\tau-\tau S>0$, so that the quantity is positive for every $S$, while it tends to $0$ as $S\to+\infty$: the infimum $g(z)=2z^2$ is not attained. If $|z|>z_c$ the numerator is negative for $S$ large, since $\tau\to1$ and $1-2\sqrt2z^2<0$, so that $g(z)<2z^2$ and the infimum is attained at a finite $S(z)$. For $|z|$ large the optimal $S$ is small, and since $\tanh x=x-\frac{x^3}3+o(x^3)$,
$$
\psi(S,z)=\frac{z^2}{1-S^2/12+o(S^2)}=z^2\Big(1+\frac{S^2}{12}+o(S^2)\Big).
$$
We then minimize $z^2+\frac{z^2S^2}{12}+\frac1S$: the derivative vanishes for $S^3=6/z^2$, that is $S(z)\sim\sqrt[3]6\,|z|^{-2/3}$, and at this value $\frac{z^2S^2}{12}=\frac1{2S}$, so that
$$
g(z)\sim z^2+\frac3{2S(z)}=z^2+\frac{3}{2\sqrt[3]6}|z|^{2/3}=z^2+\frac{\sqrt[3]{36}}4|z|^{2/3}.
$$
(The convexity of $g$, and the uniqueness of $S(z)$ and its monotonicity in $|z|$, which we have stated in the text, can be checked numerically from the formula for $\psi$, but are not proved here.)

\section{Continuum Homogenization}\label{hom-sec} 
Solutions of minimum problems depending on a small parameter $\e$ may develop 
increasingly fast oscillations as $\e$ goes to $0$. The effect of such oscillations tends to be homogeneously distributed, in a sense at the opposite end to concentration, which can be seen as a localized effect.  
In the simplest case, oscillations can be forced by spatial inhomogeneity of the energy densities at scale $\e$, which in particular can be modeled by scaling a periodic integrand. 
A basic example, already presented in the Introduction, is that of functionals of the form 
\begin{equation}\label{hom-base}
F_\e(u)=\int_{0}^1 a\Big(\frac{t}{\e}\Big)|u^\prime|^2\, dt 
\end{equation}
defined for $u\in H^1(0,1)$, 
where $a$ is a $1$-periodic function such that $0<\alpha\leq a(t)\leq \beta$. In terms of $\Gamma$-convergence, the limit is then given by 
$$F_0(u)=\underline a \int_0^1 |u^\prime|^2\, dt, \quad \hbox{\rm where }\ \  \frac{1}{\underline a}=\int_0^1\frac{1}{a(s)}\, ds.$$ 
This limit is homogeneous in the sense that the integrand is translation invariant with respect to the spatial variable $t$, and the appearance of the coefficient $\underline a$ is a consequence of the optimization of oscillations at scale $\e$. 
The asymptotic analysis of energies with increasing spatial inhomogeneity 
as $\e$ tends to $0$ is commonly called {\em homogenization} (see the Appendix to this chapter for references). 
Usually, the integrands of the limit energies are characterized through suitable {\em homogenization formulas}, entailing 
oscillations at scale $\e$ of the minimizers.  

\begin{figure}[h!]
	\centerline{\includegraphics [width=5.0in]{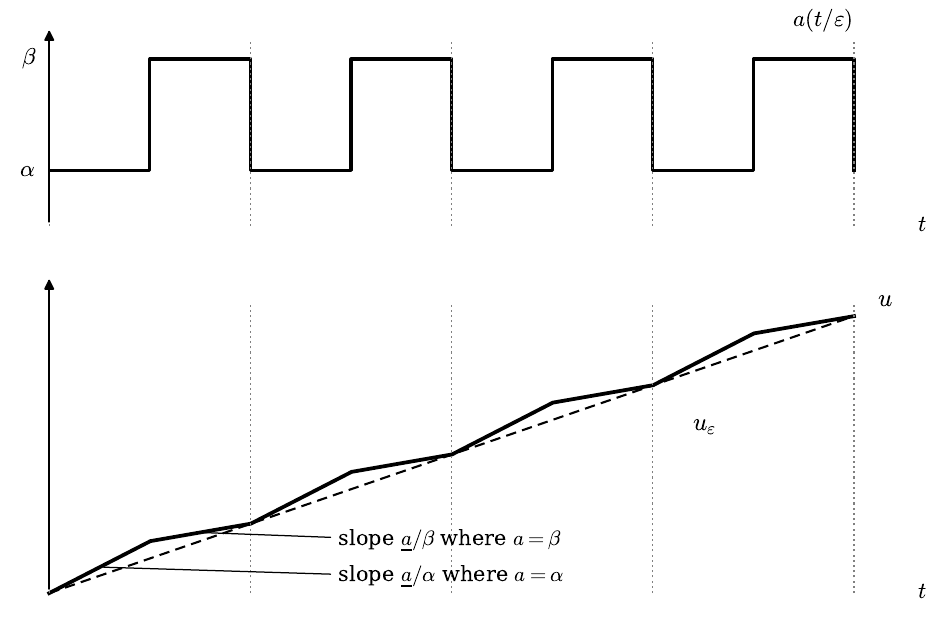}}
	\caption{The one-dimensional example \eqref{hom-base}: the coefficient $a(t/\e)$ (top) and the optimal $u_\e$ for the target $u(t)=t$ (bottom). The slope $u_\e'=\underline a/a(t/\e)$ is larger where $a$ is smaller, and $u_\e$ coincides with $u$ at the end-points of the periodicity cells, so that $F_\e(u_\e)\to\underline a\int_0^1|u'|^2dt$.}\label{Fig35}
\end{figure}

\subsection{Convexity effects: scales reduce to periods}\label{sec-L2hom}
In one dimension, the analysis of energies $F_\e$ as in \eqref{hom-base} 
is equivalent to dealing with energies 
$$
F_\e(u)=\int_{0}^1 a\Big(\frac{t}{\e}\Big)u^2\, dt 
$$
with respect to the weak $L^2$-convergence. The treatment of energies on $L^p$ modeled on this example allows us to highlight the asymptotic effect of periodic oscillations in a sufficiently general setting, without having to cope with other effects, such as those given by the structure of gradients in higher dimensions.

\paragraph{Homogenization formula in $L^p$.} 
Let then $p>1$, and let $f\colon\mathbb R^d\times\mathbb R^m\to [0,+\infty)$ be a Borel function such that 

\smallskip 
\hspace{5mm} (i) for all $z$, the map $x\mapsto f(x,z)$ is $1$-periodic; 

\hspace{5mm} (ii) $f(x,z)\geq c|z|^p-\frac{1}{c}$ for some $c>0$. 
\smallskip 

\noindent Let $\Omega$ be a bounded domain in $\mathbb R^d$. 
For every $\e>0$ we define the functional $F_\e\colon L^p(\Omega;\mathbb R^m)\to [0,+\infty]$ given by 
\begin{equation}\label{Fehom1}
F_\e(u)=\int_{\Omega}f\big(\frac{x}{\e},u\big)\, dx.
\end{equation}
Up to the relaxation argument recalled in Section \ref{Pre-Relax}, we can suppose that the map $z\mapsto f(x,z)$ is convex for almost all $x$. 

A lower bound for the asymptotic behaviour of $F_\e$ can be given by a coarse graining at the level of the period. This is done by comparing the energy of a weakly converging sequence $\{u_\e\}$ with another energy obtained by integrating out oscillations at scale $\e$, depending only on the averages of $u_\e$ on periodicity cubes. 
More precisely, let $u_\e$ weakly converge to $u\in L^p(\Omega;\mathbb R^m)$.
For all $k\in\mathbb Z^d$, let $Q^\e_k$ be the scaled (periodicity) cube defined by $Q^\e_k=\e k+\e [0,1]^d$ and $I^\e=\{k\in\mathbb Z^d: Q^\e_k\subset\Omega\}$. 
Then 
$$F_\e(u_\e)\geq\sum_{k\in I^\e}\int_{Q_k^\e}f\Big(\frac{x}{\e},u_\e(x)\Big)\,dx=
\sum_{k\in I^\e}\e^d\int_{[0,1]^d}f(y,u_\e(\e k+\e y))\,dy.
$$
If we define 
\begin{equation}\label{fhom1}f_{\rm hom}(z)=\min\Big\{\int_{[0,1]^d}f(y,v(y))\,dy: \int_{[0,1]^d} v(y)\, dy=z\Big\},
\end{equation}
then we obtain the inequality 
\begin{equation}\label{zek}F_\e(u_\e)\geq 
\sum_{k\in I^\e}\e^d f_{\rm hom}(z^\e_k), \quad \hbox{\rm where }\ \ z^\e_k=\frac{1}{\e^d}\int_{Q_k^\e} u_\e(x)\, dx, 
\end{equation}
which we can rewrite as 
\begin{equation}\label{coarse}F_\e(u_\e)\geq 
\int_{\Omega}f_{\rm hom}(\overline u_\e)\, dx-f_{\rm hom}(0)\Big|\Omega\setminus \bigcup_{k\in I^\e}Q^\e_k\Big|, \quad \hbox{\rm where }\ \ \overline u_\e(x)=\begin{cases}
z^\e_k & \hbox{\rm if } x\in Q^\e_k, k\in I^\e\\
0 & \hbox{\rm if } x\in \Omega\setminus \bigcup_{k\in I^\e}Q^\e_k.
\end{cases}
\end{equation}
Noting that $f_{\rm hom}$ is convex, that the measure of $\Omega\setminus \bigcup_{k\in I^\e}Q^\e_k$ tends to $0$, and that $\overline u_\e$ weakly converge to the same limit $u$ as $u_\e$, it follows that 
$$\liminf_{\e\to 0}F_\e(u_\e)\geq \liminf_{\e\to 0} \int_{\Omega}f_{\rm hom}(\overline u_\e)\, dx\geq \int_{\Omega}f_{\rm hom}(u)\, dx.$$
\begin{figure}[h!]
	\centerline{\includegraphics [width=5.0in]{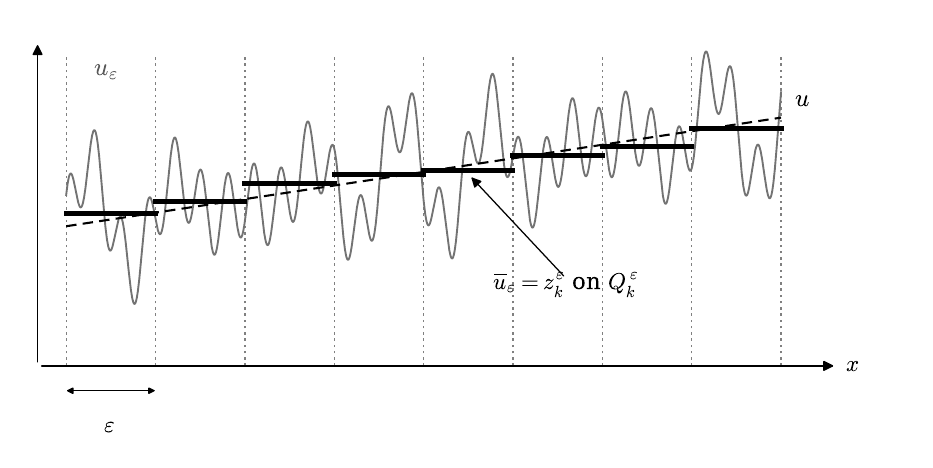}}
	\caption{Coarse graining at the level of the period: a function $u_\e$ oscillating below the scale $\e$ (thin, grey) and the piecewise-constant function $\overline u_\e$ of its averages $z^\e_k$ on the periodicity cubes $Q^\e_k$ (thick). Both converge weakly to the same limit $u$ (dashed).}\label{Fig34}
\end{figure}

To check that this bound is sharp, we first treat the case of a constant $u$ equal to $z$. To that end, let $v$ be a minimizer of the problem in \eqref{fhom1}, and set 
$u_\e(x)=v(\frac{x}{\e})$. By the Riemann-Lebesgue Lemma, $u_\e$ weakly converge to $u$, and 
$F_\e(u_\e)$ tend to $|\Omega|f_{\rm hom}(z)$. If $u$ is piecewise constant, this argument 
can be repeated on each set where $u$ is constant (there is no constraint of continuity between different sets, since the functions are only in $L^p$). For a general $u$ with $F_{\rm hom}(u)<+\infty$ we take the piecewise-constant functions $u_k$ obtained by averaging $u$ on the cubes of a grid of side $2^{-k}$: then $u_k\to u$ strongly in $L^p$, and by Jensen's inequality for the convex function $f_{\rm hom}$ we have $\int_\Omega f_{\rm hom}(u_k)\,dx\le\int_\Omega f_{\rm hom}(u)\,dx$, so that the lower semicontinuity of the $\Gamma$-limsup gives $\Gamma\hbox{-}\limsup_{\e\to0}F_\e(u)\le\int_\Omega f_{\rm hom}(u)\,dx$. Note that no upper bound on $f$ is needed here: $f_{\rm hom}$ may take the value $+\infty$, and in this case so does the $\Gamma$-limit. This shows that 
\begin{equation}\label{Fhom1}
F_{\rm hom}(u)=\int_{\Omega}f_{\rm hom}(u)\, dx
\end{equation}
is the {\em homogenized functional} of the family $\{F_\e\}$. 
The form of the recovery sequence shows that the 
scale of the oscillations is $\e$; that is, we have {\bf oscillations at the scale of the period}, as expected. 
Moreover, the recovery sequences are characterized by the single {\em cell-problem formula} \eqref{fhom1}, reducing the analysis to periodic test functions with the same period as that of the spatial oscillation of the corresponding energy density.   

As an example, if $f(x,z)=a(x)|z|^2$ we can show that  
$$f_{\rm hom}(z)=\underline a z^2, \quad \hbox{\rm where }\ \  
\underline a =\Big(\int_{[0,1]^d} \frac{1}{a(y)}\, dy\Big)^{-1}.$$ 
Up to reasoning componentwise, it is sufficient to deal with the case $m=1$. For all test functions $v$ in \eqref{fhom1},   
$$z^2=\Big(\int_{[0,1]^d}\frac{1}{\sqrt{a(y)}}\sqrt{a(y)} v(y)\, dy\Big)^2
\leq \int_{[0,1]^d}\frac{1}{a(y)}\, dy\  
\int_{[0,1]^d}a(y)v^2(y)\, dy,$$ 
which gives $f_{\rm hom}(z)\geq \underline a z^2$.  
The converse inequality is obtained by choosing $v(y)=z \underline a (a(y))^{-1}$ as test function.

\paragraph{Two-scale convergence in $L^p$.} 
The possibility of characterizing the limit behaviour of functionals by only minimizing over (locally) periodic functions allows one to use arguments which are specifically designed to treat periodic problems, such as formal asymptotic expansions. In this case, recovery sequences are supposed to take the form 
$$u_\e(x)=u(x)+\e u_1(x,\frac{x}{\e})+\e^2 u_2(x,\frac{x}{\e})+\dots,$$ with the functions $u_1,u_2,\dots$ $1$-periodic in the second variable. This approach has been formalized by introducing {\em two-scale convergence} defined as the convergence of a sequence $\{u_\e\}\subset L^p(\Omega)$ to a function $u\in L^{p}(\Omega\times [0,1]^d)$ requiring that 
$$\int_{\Omega}u_\e(x)\phi\Big(x,\frac{x}{\e}\Big)\, dx\longrightarrow \int_{\Omega\times[0,1]^d}u(x,y)\phi(x,y)\,dx dy$$ 
for all admissible test functions $\phi\in L^{p^\prime}(\Omega;C_\#([0,1]^d))$, where $C_\#([0,1]^d)$ denotes the continuous $1$-periodic functions \cite{Ngu89,All92}. This convergence can be interpreted (see the Appendix to this chapter and \cite{CDG08}) as a weak convergence in $L^{p}(\Omega\times [0,1]^d)$ of the functions $\mathcal U_\e(u_\e)$, where, using the notation introduced in  \eqref{zek},
$$
\mathcal U_\e(u_\e)(x,y)= \begin{cases} z^\e_k + v^\e_k(y)& \hbox{\rm if } x\in Q^\e_k, k\in I^\e\\
0 & \hbox{\rm if } x\in \Omega\setminus \bigcup_{k\in I^\e}Q^\e_k,
\end{cases}\qquad
v^\e_k(y)= u_\e(\e k+\e y)-z^\e_k.
$$
The operator $\mathcal U_\e: L^{p}(\Omega)\to L^{p}(\Omega\times [0,1]^d)$ is called the {\em unfolding operator}. It formalizes the coarse-graining operation introduced in \eqref{coarse}, retaining the information about oscillations of period $\e$ as a dependence on $y$. Note that the unfolding operator is tailored to deal with functions exactly of period $\e$ (as opposed to scale $\e$). Indeed, if $u_\e(x)= v({x\over\e})$ then the two-scale limit is $v(y)$ if it is $1$-periodic, but is the average of $v$ if its period is irrational. 

With respect to this convergence, the $\Gamma$-limit of the functionals 
$F_\e$ in \eqref{Fehom1} reduces to an $\e$-independent functional given by 
$$F_{\rm hom}^{2s}(u)=\int_{\Omega\times[0,1]^d}f(y,u(x,y))\, dx dy.$$
The use of two-scale convergence and $F_{\rm hom}^{2s}$ allows one to obtain information about microscopical oscillations of minimizing sequences directly from the dependence on the variable $y$ of minimizers $u$.

This functional can be compared with the functional $F_{\rm hom}$ in \eqref{Fhom1}; that is, the $\Gamma$-limit with respect to the weak convergence in $L^p$. 
Indeed, let $u_0\in L^{p}(\Omega)$, and $u\in L^p(\Omega\times [0,1]^d)$ be any function such that 
$\int_{[0,1]^d}u(x,y)\, dy=u_0(x)$. 
Then we take $\{u_\e\}$ to be a recovery sequence for $F^{2s}_{\rm hom}(u)$. By testing over $\phi$ depending only on $x$, it follows that $u_\e$ weakly converge in $L^p(\Omega)$ to $u_0$, so that by the lower estimate for $F_{\rm hom}(u_0)$, we obtain that 
$$F_{\rm hom}^{2s}(u)=\lim_{\e\to 0}F_\e(u_\e)\geq F_{\rm hom}(u_0),$$ 
and the inequality 
$$\inf\Big\{F_{\rm hom}^{2s}(v): u_0(x)=\int_{[0,1]^d}v(x,y)\, dy\Big\}\geq F_{\rm hom}(u_0).$$
For the converse inequality, let $\{u_\e\}$ be a recovery sequence for $F_{\rm hom}(u_0)$. Up to subsequences, we can suppose that $u_\e$ also two-scale converge to some $w\in L^p(\Omega\times [0,1]^d)$ satisfying the integral condition. Hence, 
$$F_{\rm hom}(u_0)=\lim_{\e\to 0}F_\e(u_\e)\geq 
F_{\rm hom}^{2s}(w)\geq \inf\Big\{F_{\rm hom}^{2s}(v): u_0(x)=\int_{[0,1]^d}v(x,y)\, dy\Big\}.$$ 
This shows that the functional $F_{\rm hom}$ in \eqref{Fhom1} can be 
recovered by integrating out the variable $y$ in $F_{\rm hom}^{2s}$. 


\subsection{Effects at scale $\e$ over all frequencies}\label{sec-W1phom} 
We now turn to homogenization in a space where functions satisfy more structure constraints than $L^p$-functions. As a result, even considering functionals with an $\e$-periodic dependence in the space variable, the limit will not be characterized by only using $\e$-periodic test functions. This is in a sense a {\em nonlocal effect} felt at scale $\e$ and not limited to the period $\e$. 

\paragraph{Homogenization formula in $W^{1,p}$.} 
Let $p>1$ and let $f\colon\mathbb R^d\times \mathbb R^{m\times d}\to [0,+\infty)$ be a Borel function such that 

\smallskip 
\hspace{5mm} (i) for all $\xi$, the map $x\mapsto f(x,\xi)$ is $1$-periodic; 

\hspace{5mm} (ii) $f(x,\xi)\geq c|\xi|^p-\frac{1}{c}$ for some $c>0$; 

\hspace{5mm} (iii) $f(x,\xi)\leq c^\prime(|\xi|^p+1)$ for some $c^\prime>0$. 
\smallskip 

\noindent Let $\Omega$ be a bounded domain in $\mathbb R^d$. 
For every $\e>0$, let $F_\e\colon W^{1,p}(\Omega;\mathbb R^m)\to [0,+\infty)$ be given by 
\begin{equation}\label{Fehom2}
F_\e(u)=\int_{\Omega}f\big(\frac{x}{\e},\nabla u\big)\, dx.
\end{equation}
Note that the gradient is an $m\times d$ matrix. As a consequence of the structure of $W^{1,p}$, convexity is not a necessary condition for the lower semicontinuity of $F_\e$ if $m>1$ and $d>1$, and hence, contrary to the $L^p$ case, it cannot be assumed without loss of generality.  

In order to obtain a lower bound producing a homogenization formula, we use the blow-up technique. Let then $u_\e$ weakly converge to $u\in W^{1,p}$ with $\sup_{\e>0}F_\e(u_\e)=M<+\infty$. In the following, extraction of subsequences will be understood, even if retaining the subscript $\e$. 
We consider the sequence of bounded measures 
$$\mu_\e=f\big(\frac{x}{\e},\nabla u_\e\big)\, d\mathcal L^d$$
which we may suppose weakly$^\ast$ converge to $\mu$. 
We now estimate the absolutely continuous part of $\mu$ with respect to the Lebesgue measure. For almost all $x_0\in\Omega$, there exists the approximate gradient $\nabla u(x_0)$, and the density of $\mu$ with respect to $\mathcal L^d$ at the point $x_0$ is given by 
$$\frac{d\mu}{d\mathcal L^d}(x_0)=\lim_{\rho\to 0}\frac{\mu(Q_\rho(x_0))}{\rho^d}=\lim_{\rho\to 0}\lim_{\e\to 0}\frac{\mu_\e(Q_\rho(x_0))}{\rho^d},$$
where we neglect a countable set of $\rho$. We then choose  $\rho=\rho_\e\to 0$ with $\rho_\e>\!\!>\e$ such that 
$$\frac{d\mu}{d\mathcal L^d}(x_0)=\lim_{\e\to 0}\frac{\mu_\e(Q_{\rho_\e}(x_0))}{\rho_\e^d}=\lim_{\e\to 0}\frac{1}{\rho_\e^d}\int_{Q_{\rho_\e}(x_0)}f\Big(\frac{x}{\e},\nabla u_\e(x)\Big)\, dx.$$ 
We may assume that the functions $v_\e(y)=\frac{u_\e(x_0+\rho_\e y)-a_\e}{\rho_\e}$, where $a_\e$ is the average of $u_\e$ on $Q_{\rho_\e}(x_0)$, strongly converge to the affine function $v(y)=\nabla u(x_0) y$ in $L^p(Q_1;\mathbb R^m)$. This convergence and both growth conditions on $f$ ensure that we may indeed assume that $u_\e(x)$ equals the value $\nabla u(x_0) x$ at the boundary of $Q_{\rho_\e}(x_0)$ up to a fixed translation. Hence,  
\begin{eqnarray*}
&&\hspace{-1cm}\frac{1}{\rho_\e^d}\int_{Q_{\rho_\e}(x_0)}f\Big(\frac{x}{\e},\nabla u_\e(x)\Big)\, dx\\
&&\geq \frac{1}{\rho_\e^d} \inf\Big\{
\int_{Q_{\rho_\e}(x_0)}f\Big(\frac{x}{\e},\nabla u(x_0)+\nabla w(x)\Big)\, dx: w\in W^{1,p}_0(Q_{\rho_\e}(x_0);\mathbb R^m)\Big\}. 
\end{eqnarray*}
After the change of variable $y=\frac{x}{\e}$, setting $T_\e=\frac{\rho_\e}{\e}$ we have that 
\begin{eqnarray*}
\frac{d\mu}{d\mathcal L^d}(x_0)&\geq& \liminf_{\e\to 0}\frac{1}{T_\e^d}\inf\Big\{   \int_{Q_{T_\e}(\frac{x_0}{\e})}f(y,\nabla u(x_0)+\nabla w(y))\, dy: w\in W^{1,p}_0(Q_{T_\e}(\frac{x_0}{\e});\mathbb R^m)      \Big\}\\
&\geq& \liminf_{\e\to 0}\frac{1}{S_\e^d}\inf\Big\{   \int_{Q_{S_\e}(k_\e)}f(y,\nabla u(x_0)+\nabla w(y))\, dy: w\in W^{1,p}_0(Q_{S_\e}(k_\e);\mathbb R^m)      \Big\}, 
\end{eqnarray*}
where $k_\e\in\mathbb Z^d$ is such that $|\frac{x_0}{\e}-k_\e|_\infty\le\frac12$ and $S_\e=T_\e+1$, so that $Q_{T_\e}(\frac{x_0}{\e})\subset Q_{S_\e}(k_\e)$ (the cubes are centered at the indicated points). In the last inequality we have extended $w$ by $0$ on $Q_{S_\e}(k_\e)\setminus Q_{T_\e}(\frac{x_0}{\e})$, which by the upper bound (iii) on $f$ adds an error of order $T_\e^{d-1}$, negligible after division by $T_\e^d$. Note that in the last formula we can substitute $k_\e\in\mathbb Z^d$ with $0$ by the $1$-periodicity of $f$. 

\begin{figure}[h!]
	\centerline{\includegraphics [width=5.4in]{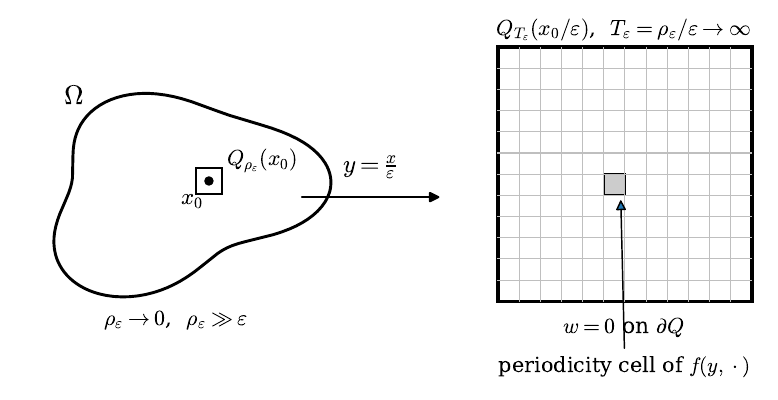}}
	\caption{The blow-up argument: at a point $x_0$ we zoom in on a cube $Q_{\rho_\e}(x_0)$ with $\e\ll\rho_\e\ll1$, on which $u_\e$ is close to the affine function $\nabla u(x_0)(x-x_0)$. After the change of variable $y=x/\e$ we obtain a cube of side $T_\e=\rho_\e/\e\to\infty$, made of many periodicity cells of $f(y,\cdot)$, with affine boundary conditions: this gives the asymptotic homogenization formula.}\label{Fig36}
\end{figure}

For fixed $\xi\in\mathbb R^{m\times d}$, we set 
$$g(S)=\frac{1}{S^d}\inf\Big\{   \int_{Q_{S}}f(y,\xi+\nabla w(y))\, dy: w\in W^{1,p}_0(Q_{S};\mathbb R^m)      \Big\}.$$ 
By a {\em subadditivity argument}, we then have $g(S)\leq g(T)+r(S,T)$ for all $S>T$, where 
$$\lim_{T\to+\infty}\limsup_{S\to+\infty}r(S,T)=0,$$ 
which guarantees the existence of the limit (see the Appendix to this chapter) in the {\em asymptotic homogenization formula}, due to Marcellini, Braides and M\"uller \cite{Mar78,Bra85,Mu87},
$$f_{\rm hom}(\xi)=\lim_{S\to+\infty}
\frac{1}{S^d}\inf\Big\{   \int_{Q_{S}}f(y,\xi+\nabla w(y))\, dy: w\in W^{1,p}_0(Q_{S};\mathbb R^m)      \Big\}.
$$ 
We then have the inequality $\frac{d\mu}{d\mathcal L^d}(x_0)\ge f_{\rm hom}(\nabla u(x_0))$ for almost all $x_0\in\Omega$, which can be integrated and gives the lower bound
$$
\liminf_{\e\to 0} F_\e(u_\e)\ge \mu(\Omega)\ge \int_\Omega \frac{d\mu}{d\mathcal L^d}(x_0)\,dx\ge \int_\Omega f_{\rm hom}(\nabla u(x))\,dx.
$$

This formula suggests a construction for approximate recovery sequences in analogy with the $L^p$ case. For a target linear function $u(x)=\xi x$, with fixed $\eta>0$ we can choose $S\in\mathbb N$ and $w\in W^{1,p}_0(Q_S;\mathbb R^m)$ such that 
$$\frac{1}{S^d}  \int_{Q_{S}}f(y,\xi+\nabla w(y))\, dy\le f_{\rm hom}(\xi)+\eta. $$
We then extend $w$ $S$-periodically and define the sequence $u_\e(x)= \xi x+\e w\big({x\over\e}\big)$ converging to $u$. By the Riemann-Lebesgue Lemma we then have 
$\lim_{\e\to 0} F_\e(u_\e)\le |\Omega|(f_{\rm hom}(\xi)+\eta)$.
Since we can modify the construction of $u_\e$ in such a way that $u_\e=u$ on $\partial\Omega$, and the functionals are invariant by addition of constants, the construction can be repeated for piecewise-affine $u$, and then extended to arbitrary $u$ by density, obtaining that the $\Gamma$-limit is 
$$F_{\rm hom} (u)= \int_\Omega f_{\rm hom}(\nabla u(x))\,dx.
$$

\smallskip\noindent
{\it A ``periodic version'' of the homogenization formula.} We can give an alternate formula for $f_{\rm hom}$ as
$$f_{\rm hom}(\xi)=\lim_{S\to+\infty}
\frac{1}{S^d}\inf\Big\{   \int_{Q_{S}}f(y,\xi+\nabla w(y))\, dy: w\in W^{1,p}_{\rm loc}(\mathbb R^d;\mathbb R^m), w \hbox{ S-periodic}  \Big\},
$$ 
where we can even restrict to $S\in\mathbb N$.
Indeed, optimizing over $S$-periodic functions provides a lower bound for the asymptotic homogenization formula, while, on the other hand, the formula provides an upper bound, as the construction of recovery sequences can be carried out using minimizing $S$-periodic functions.
It is interesting to note that this formula (in this version, or in the one with zero boundary conditions) can be used also in the case when $f=f(\xi)$ is a non-(quasi)convex function, giving a characterization of the {\em(quasi)convex envelope} $Qf$ of $f$ as an {\bf optimization process over all scales of periodic oscillations}. This observation is helpful in that it suggests that integrands of $\Gamma$-limits can sometimes be characterized by an asymptotic formula over periodic functions even though not deriving from a periodic setting.

\medskip
The asymptotic homogenization formula shows that the effect of the periodicity is at {\bf scale $\e$ but not necessarily confined to the period $\e$}. A counterexample by Stefan M\"uller \cite{Mu87} shows that the analysis indeed cannot be reduced to considering only $\e$-periodic functions. This can be explained by recalling that a class of integrands giving lower-semicontinuous functionals are {\em polyconvex} functions; that is, convex functions of the minors. Hence, in particular if $m=d=2$ we have convex functions of $\det\nabla u= D_1u_1D_2u_2-D_1u_2D_2u_1$, and we may have functions $u$ whose determinant is $1$-periodic such that $\nabla u$ is periodic with a larger period.

\paragraph{Cell-problem formula in the convex case.} 
In order to deduce some properties in the case of a convex function, we consider, for a general $f$,  
the simpler cell-problem formula on $1$-periodic functions given by  
$$f^\#(\xi)=\inf\Big\{   \int_{Q_{1}}f(y,\xi+\nabla v(y))\, dy: v\in W^{1,p}_{\rm loc}(\mathbb R^d;\mathbb R^m), v \hbox{ 1-periodic}   \Big\}. 
$$ 
If $v$ is a 1-periodic test function, and $\phi_S(x)=$dist$(x,\partial Q_S)\wedge 1$, we can use 
$\phi_S v$ as a test function for $g(S)$ and deduce that $g(S)\le f^\#(\xi)+ o(1)$ as $S\to+\infty$. 
Hence, 
$$f_{\rm hom}(\xi)\leq f^\#(\xi).$$

If the function $f$ is convex in $\xi$, then we can prove the equality.  Indeed, 
if $S\in\mathbb N$ and $w$ is a test function for $g(S)$ extended by $S$-periodicity, we can use 
as test function for $f^\#(\xi)$ the $1$-periodic function 
$$v(x)={1\over S^d}
\sum_{k\in \{1,\ldots, S\}^d} w(x+k), 
$$
obtained as a convex combination of all integer translations of $w$. 
Using convexity and 
the $1$-periodicity of $v$, and the $S$-periodicity of both $v$ and $w$, by the convexity of $f$ we have 
\begin{eqnarray*}
&&\hspace{-5mm}\int_{Q_1} f(y, \xi+\nabla v)dy=
{1\over S^d} \int_{Q_S} f(y, \xi+\nabla v)dy\\
&&\hspace{5mm}=
{1\over S^d} \int_{Q_S} f\Big(y, {1\over S^d}
\sum_{k\in \{1,\ldots, S\}^d} (\xi +\nabla w(y+k))\Big)dy\\
&&\hspace{5mm}\le {1\over S^d} {1\over S^d} \sum_{k\in \{1,\ldots, S\}^d}\int_{Q_S} f(y, 
 \xi +\nabla w(y+k))dy = 
 {1\over S^d} \int_{Q_S} f(y,
 \xi +\nabla w(y))dy, 
\end{eqnarray*}
which shows that $f^\#(\xi)\le g(S)$ and then $f^\#(\xi)=f_{\rm hom}(\xi)$.  
Hence, also in the Sobolev case, in a sense convexity trivializes the oscillations to those of $\e$-periodic functions, so that {\bf the scale $\e$ reduces to the period $\e$}. On the other hand, having a single problem instead of an asymptotic formula makes it more handy for computations and stability issues.

\paragraph{Non-one dimensional dependence on the microstructure.} 
The dependence of $f$ on the spatial variable is often referred to as its {\em microstructure}. In the one-dimensional case $d=1$ and 
$f(t,z)=a(t)z^2$, we have seen that the limit coefficient is the harmonic mean of $a$, so that the effect of the microstructure only depends on this average quantity. In the $L^p$ case, we have noted that this extends to higher dimension. In the Sobolev case, this is not true anymore, and in particular {\bf the scale is not one-dimensional}; that is, we cannot recover $f_{\rm hom}(\xi)$ by looking at one-dimensional sections in the direction of $\xi$. This can be seen through an example, for which we can use the equivalent cell-problem formula giving $f^\#(\xi)$. We fix a one-dimensional coefficient $a$ and consider the function $f(x_1,x_2,\xi)= a(x_1)|\xi|^2$ for $d=2$ and $m=1$. The formula for the homogenized function is
$$f^\#(\xi)=\inf\Big\{\int_{Q_{1}}a(y_1)|\xi+\nabla v(y)|^2\, dy: v\in W^{1,p}_{\rm loc}(\mathbb R^2), v \hbox{ 1-periodic}   \Big\}. 
$$ 
This formula defines a quadratic form, so that there exists a $2\times 2$ matrix $A_{\rm hom}$ such that $f^\#(\xi)=\langle A_{\rm hom}\xi,\xi\rangle $. By the independence of the coefficient from $x_2$, we deduce that $f^\#(\xi_1,-\xi_2)=f^\#(\xi_1,\xi_2)$, which implies that $ A_{\rm hom}$ is diagonal. In order to determine the coefficients $A_{11}$ and $A_{22}$ on the diagonal we note that minimizers for $f^\#(e_1)$ and $f^\#(e_2)$ are both one-dimensional. In the first case the minimizer is $v(x_1,x_2)= v_1(x_1)$, where $v_1$ is the minimizer of the one dimensional case related to $a(x)|\xi_1|^2$. Hence $A_{11}=\underline a$. In the second case the minimizer is $v(x_1,x_2)= 0$ since sections of the coefficients in the direction $x_2$ are constant for all fixed $x_1$ and $A_{22}=\int_0^1a(y)\,dy$. 

Now, we analyze the functionals corresponding to one-dimensional sections of $a$ in each fixed direction $\xi$. If $\xi_1\neq 0$, for every $\xi^\prime$ with $\langle\xi,\xi^\prime\rangle=0$ we consider the $1$-dimensional homogenization problem related to 
$$F^{\xi^\prime}_\e(u)=\int_{0}^1a\Big({\langle\xi^\prime+t\xi, e_1\rangle\over\e}\Big)|u^\prime|^2\, dt.$$
Since the map $t\mapsto a(\langle\xi^\prime+t\xi, e_1\rangle)$ equals $a$ up to an affine change of variables, its harmonic mean equals $\underline a$. Hence all homogenized functionals 
$F_{\rm hom}^{\xi^\prime}$ equal the functional $\overline F$ given by 
$$\overline F(u)=\underline a \int_0^1 |u^\prime|^2\, dt,$$
which only gives the trivial lower estimate 
$f^\#(\xi)\geq \underline a |\xi|^2$, which is not sharp if $\xi_2\neq 0$, as we have shown.

Another classical example, from which we can deduce an even more strictly non-one dimensional behaviour of the homogenized functional in the Sobolev case, is that of the two-dimensional periodic checkerboard, with the values of the function $a$ alternating between $\alpha$ and $\beta$ (both strictly positive) on the black and white squares, respectively. If $f(x,\xi)=a(x)|\xi|^2$, it can be shown (see the Appendix to this chapter) that $f^\#(\xi)=\sqrt{\alpha\beta}|\xi|^2$, while the analysis on one-dimensional sections still gives the coefficient $\underline a={2\alpha\beta\over \alpha+\beta}<\sqrt{\alpha\beta}$.

\begin{figure}[h!]
	\centerline{\includegraphics [width=5.4in]{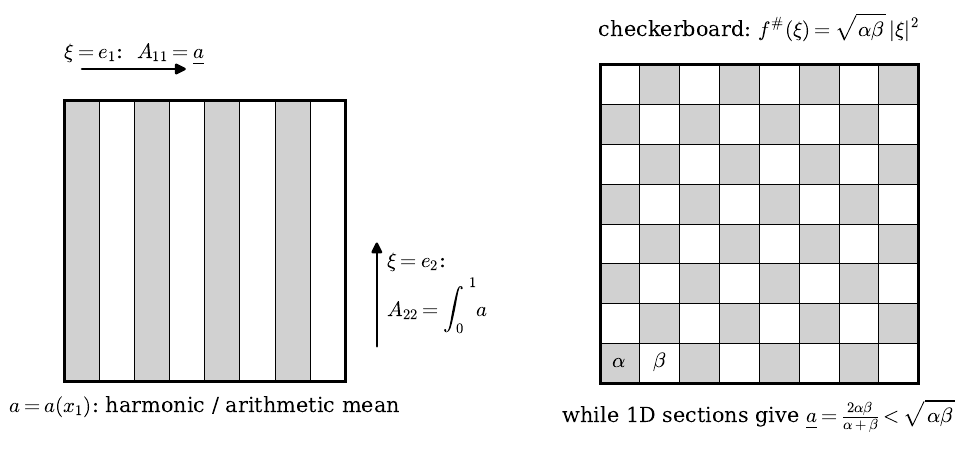}}
	\caption{Non-one-dimensional microstructure. Left: for a laminate $a=a(x_1)$ the homogenized matrix is diagonal, with the harmonic mean $\underline a$ for $\xi=e_1$ (oscillations in the direction of $\xi$) and the arithmetic mean for $\xi=e_2$ (no oscillations). Right: the checkerboard, whose homogenized coefficient $\sqrt{\alpha\beta}$ is larger than the harmonic mean given by one-dimensional sections.}\label{Fig37}
\end{figure}

\subsection{Other continuum homogenization processes}  
The general principle of homogenization can be applied to a broad class of energies, with an oscillating dependence on some variable, whose effect can be integrated out in the limit, obtaining a ``simpler'' homogeneous energy describing, in a sense, a macroscopic behaviour. The main characteristic required to make a homogenization approach work is that the functionals have to be defined on spaces allowing oscillations, endowed with a structure that may capture average behaviours.   From our standpoint, we choose to highlight the general procedure of recovering a characteristic scale by the analysis of the limit energy density. In this perspective, we now examine the behaviour of families of energies defined on different spaces, where the effect of weak convergences must be suitably understood. 

\paragraph{Homogenization of perimeter functionals.} 
We consider oscillating energies defined on sets of finite perimeter inside a fixed domain $\Omega\subset\mathbb R^d$, given by
\begin{equation}\label{homper}
F_\e(E)=\int_{\Omega\cap\partial^\ast E}\varphi\Big(\frac{x}{\e}, \nu\Big)\, d\mathcal H^{d-1},
\end{equation}
with $\varphi\colon \mathbb R^d\times S^{d-1}\to [\alpha,\beta]$ a Borel function, even in the second variable, and $\alpha,\beta>0$ \cite{AB90}. Note that, if $\varphi$ is constant, then $F_\e$ is a multiple of the perimeter, so that we can refer to this family of energies as {\bf oscillating perimeters}.
In order to envisage the form of a homogenization formula, we can use the structure of sets of finite perimeter. 
To that end, we let $E_\e$ (with equibounded perimeters) converge to $E$ and focus on points $x_0\in\partial^\ast E$. Using the blow-up technique, we let 
$$\mu_\e=\varphi\Big(\frac{x}{\e},\nu_\e\Big)\, \mathcal H^{d-1}\restr \partial^\ast E_\e,$$
where $\nu_\e$ is the inner normal to $E_\e$, and suppose, without loss of generality, that $\mu_\e$ weakly$^\ast$ converge to a bounded measure $\mu$. 
For $\mathcal H^{d-1}$-almost all $x_0\in\partial^\ast E$, 
\begin{eqnarray*}
\frac{d\mu}{d\mathcal H^{d-1}\restr\partial^\ast E}(x_0)=\lim_{\e\to 0}\frac{1}{\rho_\e^{d-1}}\int_{Q^{\nu}_{\rho_\e}(x_0)\cap\partial^\ast E_\e}\varphi\Big(\frac{x}{\e},\nu_\e\Big)\, d\mathcal H^{d-1}
\end{eqnarray*}
where $\nu=\nu(x_0)$ is the inner normal to $\partial^\ast E$ at $x_0$ and $Q^{\nu}_{T}(y)$ is any cube with side length $T$, centered at $y$ and with one face orthogonal to $\nu$. 
The infinitesimal parameter $\rho_\e$ is chosen, after a diagonal argument, such that $\rho_\e>\!\!>\e$ and also $\frac{1}{\rho_\e}(E_\e-x_0)$ converge to the half space $H^\nu=\{x: \langle x,\nu \rangle >0\}$. 

As usual, the second step in the application of the blow-up method is the variation of boundary data, which in this case implies that we can assume that 
$E_\e=x_0+H^\nu$ in a neighbourhood of the boundary of the cube $Q^\nu_{\rho_\e}(x_0)$. 
By minimizing over all sets satisfying this boundary condition, by a change of variable we get the lower inequality 
\begin{eqnarray*}
&&\hspace{-1cm}\frac{d\mu}{d\mathcal H^{d-1}\restr\partial^\ast E}(x_0)\\
&&\geq\liminf_{\e\to 0}
\frac{1}{T_\e^{d-1}}\inf\Big\{ \int_{Q^{\nu}_{T_\e}(\frac{x_0}{\e})\cap\partial^\ast A}\varphi(y,\nu_A)\, d\mathcal H^{d-1} : A=\frac{x_0}{\e}+H^\nu \ \hbox{\rm on } \partial Q^{\nu}_{T_\e}\big(\frac{x_0}{\e}\big) \Big\},
\end{eqnarray*}
where $T_\e=\frac{\rho_\e}{\e}$. As in the Sobolev case, by the periodicity it can be shown that we may suppose without loss of generality that $x_0=0$.
In order to give a lower bound independent of subsequences, we note that for all $\nu$ we can define the {\em asymptotic homogenization surface formula} 
\begin{eqnarray}\label{homformper}
\varphi_{\rm hom}(\nu)=\lim_{T\to +\infty}
\frac{1}{T^{d-1}}\inf\Big\{ \int_{Q^{\nu}_{T}\cap\partial^\ast A}\varphi(y,\nu_A)\, d\mathcal H^{d-1} : A=H^\nu \ \hbox{\rm on } \partial Q^{\nu}_{T}\Big\},\end{eqnarray}
where $Q_T^\nu=Q_T^\nu(0)$. 
The existence of the limit is again proved by a subadditive argument, which we outline below. The validity of this formula allows us to give a lower bound independent of subsequences
$$
\liminf_{\e\to0} F_\e(E_\e)\ge \mu(\Omega)\ge
\int_{\Omega\cap\partial^* E}\frac{d\mu}{d\mathcal H^{d-1}\restr\partial^\ast E} d\mathcal H^{d-1}\ge \int_{\Omega\cap\partial^* E} \varphi_{\rm hom}(\nu(x))d\mathcal H^{d-1}.
$$

In order to show the upper bound, it is sufficient to treat in detail the case $E=H^\nu\cap \Omega$. This construction also describes the subadditive argument leading to the proof of the existence of the limit in 
\eqref{homformper}.
We fix $T>0$ and consider a(n almost) minimizer $A_T$ of the problem in \eqref{homformper}.
Let $\mathbb Z^{d-1}_\nu$ be a lattice congruent to  $\mathbb Z^{d-1}$ in the plane $\Pi^\nu$. For each $n\in\mathbb Z^{d-1}$ we choose $k_n\in \mathbb Z^d$ at distance not larger than ${1\over 2}\sqrt d$ from $(T+1)\mathbb Z^{d-1}_\nu$ and define the set $E_T$ as equal to $A_T+k_n$ inside $Q_T+k_n$ and equal to $H^\nu$ otherwise.
The sequence $E^T_\e= \e E_T$ converges to $H^\nu$ and 
$$
\limsup_{\e\to 0} F_\e(E^T_\e)\le {1\over T^{d-1}}\int_{Q^{\nu}_{T}\cap\partial^\ast A_T}\varphi(y,\nu_{A_T})\, d\mathcal H^{d-1} + r(T),$$
where the remainder term $r(T)$ is due to the extra contribution on the boundary of the cubes $Q_T+k_n$, where we have a mismatch with the boundary condition in a neighbourhood of $\Pi^\nu$ of size at most $\sqrt d$, and to the contribution outside the union of the cubes, where the interface coincides with $\Pi^\nu$. In both cases we obtain $r(T)=O({1\over T})$, so that, by the arbitrariness of $T$ we have 
$$\Gamma\hbox{-}\limsup_{\e\to 0}F_\e(E)\le \varphi_{\rm hom}(\nu)\,{\mathcal H}^{d-1}(\Pi^\nu\cap \Omega)=\int_{\Omega\cap\partial^* E} \varphi_{\rm hom}(\nu(x))d\mathcal H^{d-1}.
$$
This construction highlights that recovery sequences develop {\bf oscillations at scale $\e$ in the direction orthogonal to $\nu$}, confined to a neighbourhood of the boundary of scale $\e$. Moreover, in general such oscillations are not periodic but in a sense {\bf almost periodic}, unless $\nu$ is a {\em rational direction}; that is, there exists $t\neq 0$ such that $t\nu\in\mathbb Z^d$.
The upper bound can be extended to polyhedral sets by using the construction above close to each face of the polyhedron, and then to all sets of finite perimeter by density, thus showing that
$$\Gamma\hbox{-}\limsup_{\e\to 0}F_\e(E)=\int_{\Omega\cap\partial^* E} \varphi_{\rm hom}(\nu(x))d\mathcal H^{d-1}.
$$

\bigskip As an example we can take the simple energies in the plane ($d=2$)
$$
F_\e(E)=\int_{\Omega\cap\partial^*E} a\Bigl({x\over\e}\Bigr)d{\mathcal H}^{1},
$$
with {\em isotropic} (that is, not depending on $\nu$) energy densities, in the particular case when $a$ is the one-periodic function 
$$
a(x_1,x_2)=\begin{cases} \alpha &\hbox{ if } x_1\in\mathbb Z \hbox { or } x_2\in\mathbb Z\\
\beta &\hbox{ otherwise,} \end{cases}
$$
and $\beta>\sqrt2 \alpha$. In this case the optimal sequences $E_\e$ are such that $a(x_1,x_2)=\alpha$ on $\partial^*E_\e$ and the homogenized energy density is $\varphi_{\rm hom} (\nu)= \alpha(|\nu_1|+|\nu_2|)$, which is homogeneous but not isotropic (for the study of the minimizers of such oscillating perimeters, see \cite{CdlL05,CT09}).

\begin{figure}[h!]
	\centerline{\includegraphics [width=5.6in]{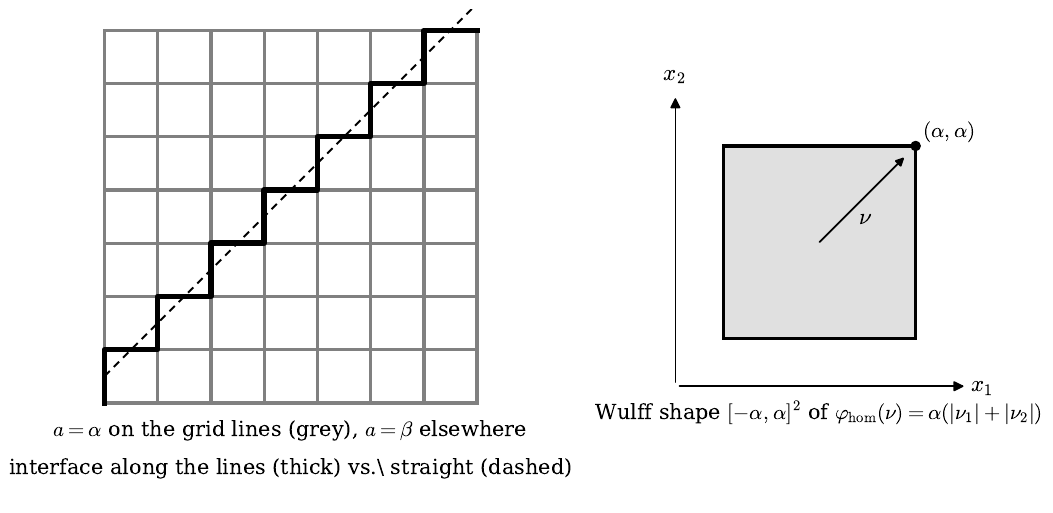}}
	\caption{Oscillating perimeters. Left: the density is $\alpha$ on the lines of the grid and $\beta>\sqrt2\alpha$ elsewhere; a straight interface at $45^\circ$ (dashed) pays $\beta$ almost everywhere, while an interface running along the lines (thick, a staircase at the scale $\e$) pays only $\alpha$ times its length, which gives $\varphi_{\rm hom}(\nu)=\alpha(|\nu_1|+|\nu_2|)$. Right: the corresponding Wulff shape, the square $[-\alpha,\alpha]^2$.}\label{Fig38}
\end{figure}

\goodbreak

\paragraph{Homogenization of free-discontinuity problems.} We now consider an interesting class of energies that mix bulk and surface parts; namely, functionals of the form
$$
F_\e(u)=\int_\Omega f\Big({x\over\e'},\nabla u\Big)\,dx +\int_{\Omega\cap S(u)} \varphi\Big({x\over\e''},\nu_u\Big)d{\mathcal H}^{d-1},
$$
defined on functions in $SBV(\Omega)$, where $\e',\e''\to 0$ independently. For simplicity, we consider the scalar case, and surface energy densities not depending on the traces of $u$ on either side of $S(u)$, and symmetric in $\nu$. The functions $f$ and $\varphi$ satisfy the same hypotheses as for the bulk and perimeter homogenizations just considered. 

The main issue in the treatment of such energies is an {\bf independence of scale effect}; that is, that in a sense the bulk and surface parts can be considered independently. We briefly explain this phenomenon leaving out the technical details.

\smallskip
It is interesting to note that the blow-up technique may be applied also in this case, since its ingredients are differentiation of measures and characterization of pointwise quantities by scaling, which is the case also for SBV-functions. We then take a sequence $u^\e$ of equibounded energy. Upon truncation we can suppose that $u^\e$ is uniformly bounded and converges to a function $u\in SBV(\Omega)$ in $L^1(\Omega)$, and $\e'=\e'(\e)$, $\e''=\e''(\e)$. We define the measures
$$
\mu_\e =f\Big({x\over\e'},\nabla u_\e\Big){\mathcal L}^d + \varphi\Big({x\over\e''},\nu_{u_\e}\Big) {\mathcal H}^{d-1}\restr S(u_\e),
$$
which we may assume weakly$^*$ converge to some finite measure $\mu$, whose absolutely continuous part with respect to the target measure $\lambda= {\mathcal L}^d + {\mathcal H}^{d-1}\restr S(u)$ we want to represent. Note that $\lambda$ is the sum of two mutually singular measures. 

We concentrate on the bulk part first, and note that for almost all $x_0\in\Omega$ we may suppose that the functions $v_\rho(y)={u(x_0+\rho y)- u(x_0)\over \rho}$ converge in $L^1_{\rm loc} (\mathbb R^d)$ to $\langle \nabla u(x_0),y\rangle$.
We can then extract a sequence $\rho=\rho_\e>>\e$ such that still $w_\e(y)={u_\e(x_0+\rho y)- u_\e(x_0)\over \rho}$ converge in $L^1_{\rm loc} (\mathbb R^d)$ to $\langle \nabla u(x_0),y\rangle$ and 
$$
{d\mu\over d\lambda}(x_0)={d\mu\over d{\mathcal L}^d}(x_0)=\lim_{\e\to 0}{1\over \rho^d}
\Big(\int_{Q_\rho(x_0)}  f\Big({x\over\e'},\nabla u_\e\Big)\,dx +\int_{Q_\rho(x_0)\cap S(u_\e)} \varphi\Big({x\over\e''},\nu_{u_\e}\Big)d{\mathcal H}^{d-1}\Big)
$$
$$
\ge \liminf_{\e\to 0}\int_{Q_1}  f\Big({x_0+\rho y\over\e'},\nabla w_\e(y)\Big)\,dy.
$$

Since the $(d-1)$-Hausdorff measure is homogeneous of degree ${d-1}$ by homothety, and by the inequality $\varphi\ge \alpha>0$, we note that ${\mathcal H}^{d-1}(Q_1\cap S(w_\e)) \le  C\rho$ by the finiteness of $
{d\mu\over d\lambda}(x_0)$. This implies that we may find a sequence $\widetilde w_\e\in W^{1,p}(Q_1)$ converging to $\langle \nabla u(x_0),y\rangle$, and such that
$$
\int_{Q_1} f\Big({x_0+\rho y\over\e'},\nabla \widetilde w_\e(y)\Big)\,dy\le  \int_{Q_1}  f\Big({x_0+\rho y\over\e'},\nabla w_\e(y)\Big)\,dy+o(1)
$$
as $\e\to 0$. We do not prove this construction, which expresses the fact that we may substitute a sequence in SBV with vanishing jump set with a sequence with no jump set at all, and hence in $W^{1,p}$. Arguing as in the proof of the homogenization theorem in $W^{1,p}$ we can suppose that $x_0=0$ by the periodicity of $f$, and then obtain that  
$$
{d\mu\over d\lambda}(x_0)\ge \liminf_{\e\to 0}\int_{Q_1}  f\Big({y\over\e'/\rho},\nabla \widetilde  w_\e(y)\Big)\,dy\ge f_{\rm hom}(\nabla u(x_0))
$$
by the bulk homogenization theorem applied with the small parameter $\e'/\rho$ in place of $\e'$.

The surface part again uses the different scaling properties of the Lebesgue and Hausdorff parts. 
Now, for almost all $x_0\in\Omega\cap S(u)$ we may suppose that the functions $v_\rho(y)=u(x_0+\rho y)- u(x_0)$ converge in $L^1_{\rm loc} (\mathbb R^d)$ to $u^{\nu}$ defined by
$$
u^{\nu}(y)=\begin{cases} u^+(x_0) & \hbox{ if } \langle y, \nu\rangle>0\\
u^-(x_0) & \hbox{ otherwise, }
\end{cases}
$$
where $\nu=\nu_u(x_0)$.
We can then extract a sequence $\rho=\rho_\e>\!>\e$ such that still $w_\e(y)=u_\e(x_0+\rho y)- u_\e(x_0)$ converge in $L^1_{\rm loc} (\mathbb R^d)$ to $u^{\nu}$ and 
\begin{eqnarray*}
{d\mu\over d\lambda}(x_0)&=&{d\mu\over d{\mathcal H}^{d-1}\restr S(u)}(x_0)
\\&=&\lim_{\e\to 0}{1\over \rho^{d-1}}
\Big(\int_{Q^\nu_\rho(x_0)}  f\Big({x\over\e'},\nabla u_\e\Big)\,dx +\int_{Q^\nu_\rho(x_0)\cap S(u_\e)} \varphi\Big({x\over\e''},\nu_{u_\e}\Big)d{\mathcal H}^{d-1}\Big)
\\
&\ge& \liminf_{\e\to 0}\int_{Q^\nu_1\cap S(w_\e)}  \varphi\Big({x_0+\rho y\over\e''},\nu_{w_\e}\Big)\,d{\mathcal H}^{d-1}.
\end{eqnarray*}
Thanks to the $d$-homogeneity of the Lebesgue measure and the inequality $f\ge c|z|^p-{1\over c}$ with $c>0$, we note that $\int_{Q^\nu_1}|\nabla w_\e|^pdx \le  C\rho^{p-1}$ by the finiteness of ${d\mu\over d\lambda}(x_0)$. This implies that we may find a sequence $\widetilde w_\e\in SBV(Q_1)$ converging to $u^{\nu}$, such that $\nabla \widetilde w_\e=0$ almost everywhere, and 
$$
\mathcal H^{d-1}((Q^\nu_1\cap S(\widetilde w_\e))\setminus S(w_\e))\le C\rho.
$$
as $\e\to 0$. We do not prove this construction, which is an application of the coarea formula, which expresses the fact that we may substitute a sequence in SBV with vanishing gradient with a sequence with no gradient at all, that is, a piecewise-constant function. 
If we set 
$$
E_\e=\{ x\in\Omega: |\widetilde w_\e(x)- u^+(x_0)|\le  |\widetilde w_\e(x)- u^-(x_0)|\}
$$ 
then this is a sequence of sets of finite perimeter converging to $H^{\nu}$ and 
$$
\int_{Q^\nu_1\cap \partial E_\e}  \varphi\Big({x_0+\rho y\over\e''},\nu_{E_\e}\Big)\,d{\mathcal H}^{d-1}\le
\int_{Q^\nu_1\cap S(w_\e)}  \varphi\Big({x_0+\rho y\over\e''},\nu_{w_\e}\Big)\,d{\mathcal H}^{d-1}.
$$
From this convergence and the homogenization theorem for perimeters we deduce
$$
{d\mu\over d\lambda}(x_0)\ge \liminf_{\e\to 0}\int_{Q_1\cap \partial E_\e}  \varphi\Big({x_0+\rho y\over\e''},\nu_{E_\e}\Big)\,d{\mathcal H}^{d-1}\ge \varphi_{\rm hom}(\nu).
$$

Summing up this inequality and the one for the bulk part, we then have the lower bound
\begin{eqnarray*}
\liminf_{\e\to 0} F_\e(u_\e)\ge \int_\Omega {d\mu\over d\lambda}d\lambda
\ge \int_\Omega f_{\rm hom}(\nabla u)\,dx
+ \int_{\Omega\cap S(u)}\varphi_{\rm hom}(\nu_u)
\,d{\mathcal H}^{d-1}.
\end{eqnarray*}

The separation of scales is reflected in the construction of recovery sequences, which we only hint at: it is sufficient to consider the case when $S(u)$ is a $(d-1)$-dimensional polytope. This can be considered as a subset of some $\Pi^\nu$, and hence as a boundary, for which an almost-periodic almost-optimal construction can be performed as in the case of perimeters.
After this modification a recovery sequence as in the case of a Sobolev function can be constructed in the complement of this small deformation of $S(u)$. We do not include the details. The conclusion is that, as $\e',\e''\to0$, the functionals $F_\e$ $\Gamma$-converge (with respect to the $L^1$ convergence) to 
$$
F_{\rm hom}(u)=\int_\Omega f_{\rm hom}(\nabla u)\,dx+\int_{\Omega\cap S(u)}\varphi_{\rm hom}(\nu_u)\,d\mathcal H^{d-1},
$$
where $f_{\rm hom}$ and $\varphi_{\rm hom}$ are given by the homogenization formulas for the bulk and surface terms separately, independently of the relative speed of $\e'$ and $\e''$ \cite{BDV96}.

\subsection{Appendix: relevant tools of Chapter 7}\label{App-Ch7}

\paragraph{References.} The homogenization of integral functionals by $\Gamma$-convergence starts from the quadratic case, treated by De Giorgi and Spagnolo \cite{DGS73} in the framework of $G$-convergence of elliptic operators; the convex periodic case is due to Marcellini \cite{Mar78}, and the non-convex one, with the asymptotic homogenization formula, to Braides \cite{Bra85} and M\"uller \cite{Mu87}, the latter also proving that the cell problem cannot be restricted to functions with period $1$. Complete accounts are given in the books by Bensoussan, Lions and Papanicolaou \cite{BLP78}, Jikov, Kozlov and Oleinik \cite{JKO94}, and Braides and Defranceschi \cite{BD98}, the last one being specifically devoted to integral functionals and following the point of view of these notes. Two-scale convergence is due to Nguetseng \cite{Ngu89} and Allaire \cite{All92}, and the unfolding operator, which gives a convenient reformulation, to Cioranescu, Damlamian and Griso \cite{CDG08}. The homogenization of oscillating perimeters is due to Ambrosio and Braides \cite{AB90}; for the structure of minimizers of such energies see Caffarelli and de la Llave \cite{CdlL05} and Chambolle and Thouroude \cite{CT09}. The independence of scales for free-discontinuity energies is studied in \cite{BDV96}. The formula for the checkerboard is classical and is due to Keller \cite{Kel64} and Dykhne \cite{Dyk71}; see also \cite{MS85}.

\paragraph{Oscillating functions.} The following lemma is used in all the constructions of recovery sequences.

\begin{lemma}[Riemann--Lebesgue]\label{RL}
Let $g\in L^1_{\rm loc}(\mathbb R^d)$ be $1$-periodic and set $\langle g\rangle=\int_{[0,1]^d}g\,dy$. Then:

(a) $\int_Ag(x/\e)\,dx\to|A|\langle g\rangle$ for every bounded set $A$ with $|\partial A|=0$;

(b) if $g\in L^p_{\rm loc}$ with $1<p\le+\infty$, then $g(\cdot/\e)\rightharpoonup\langle g\rangle$ weakly in $L^p_{\rm loc}$ (weakly$^*$ if $p=+\infty$).
\end{lemma}

\begin{proof}
(a) First let $A=\prod_i(a_i,b_i)$ be a box. The number of cubes $Q^\e_k$ contained in $A$ is $|A|\e^{-d}+O(\e^{1-d})$, each of them contributes $\e^d\langle g\rangle$ to the integral, and the remaining cubes, which meet $\partial A$ and are $O(\e^{1-d})$ in number, contribute at most $\e^d\int_{[0,1]^d}|g|\,dy$ each. Hence $\int_Ag(x/\e)dx=|A|\langle g\rangle+O(\e)$. For a general $A$ with $|\partial A|=0$ we approximate $A$ from inside and from outside by finite unions of boxes, using $|g|$ to control the difference.

(b) By (a) applied to $|g|^p$ the sequence $g(\cdot/\e)$ is bounded in $L^p(K)$ for every compact $K$. By (a) again $\int g(x/\e)\psi\,dx\to\langle g\rangle\int\psi\,dx$ for all $\psi$ which are finite linear combinations of characteristic functions of boxes; these are dense in $L^{p'}(K)$ if $p>1$ (in $L^1$ if $p=+\infty$), which, together with the boundedness, gives the claim.
\end{proof}

In particular, if $v$ is $1$-periodic and $f(y,v(y))\in L^1_{\rm loc}$, then $\int_\Omega f(x/\e,v(x/\e))dx\to|\Omega|\int_{[0,1]^d}f(y,v(y))dy$, which is the computation used for the energy of the recovery sequences of Section~7.1.

\paragraph{The subadditive argument.}
Let $f$ satisfy (i)--(iii) of Section~7.2, fix $\xi\in\mathbb R^{m\times d}$ and set, for $S>0$ and $Q_S=(0,S)^d$,
$$
I(S)=\inf\Big\{\int_{Q_S}f(y,\xi+\nabla w)\,dy:\ w\in W^{1,p}_0(Q_S;\mathbb R^m)\Big\},\qquad g(S)=\frac{I(S)}{S^d}.
$$

\begin{proposition}\label{subadd}
The limit $\lim_{S\to+\infty}g(S)$ exists.
\end{proposition}

\begin{proof}
Let $S,T\in\mathbb N$ with $S\ge T$, let $q=\lfloor S/T\rfloor$, and let $\eta>0$ and $w_T$ be a test function with $\int_{Q_T}f(y,\xi+\nabla w_T)dy\le I(T)+\eta$. We define $w$ on $Q_S$ as $w_T(y-Tk)$ on each cube $Tk+Q_T$, $k\in\{0,\dots,q-1\}^d$, and as $0$ on the remaining part of $Q_S$. Since $T$ is an integer and $f$ is $1$-periodic in the first variable, each of the $q^d$ cubes gives the same contribution $\int_{Q_T}f(y,\xi+\nabla w_T)dy$, while on the remaining set, of measure $S^d-q^dT^d$, we use the upper bound (iii), obtaining
$$
I(S)\le q^d\big(I(T)+\eta\big)+C(|\xi|^p+1)(S^d-q^dT^d),
$$
and, dividing by $S^d$ and using $qT\ge S-T$,
$$
g(S)\le \Big(\frac{qT}S\Big)^d(g(T)+\eta T^{-d})+C(|\xi|^p+1)\Big(1-\Big(1-\frac TS\Big)^d\Big).
$$
Letting $\eta\to0$ and then $S\to+\infty$ with $T$ fixed we get $\limsup_{S\to+\infty}g(S)\le g(T)$ for all $T\in\mathbb N$ (here $S$ ranges in the integers). Taking the $\liminf$ in $T$ we obtain that $\lim_{S\in\mathbb N}g(S)$ exists (it is finite, since $0\le g(S)\le C(|\xi|^p+1)$ by (iii) with $w=0$). If $S$ is real, then, extending by $0$ test functions on $Q_{\lfloor S\rfloor}\subset Q_S\subset Q_{\lceil S\rceil}$, we get $|I(S)-I(\lfloor S\rfloor)|\le C(|\xi|^p+1)(\lceil S\rceil^d-\lfloor S\rfloor^d)=O(S^{d-1})$ and the same for $\lceil S\rceil$, so that $g(S)-g(\lfloor S\rfloor)\to0$ and the limit exists along the reals.
\end{proof}

The same argument shows that the limit does not change if the cubes $Q_S$ are replaced by cubes with centre in $\mathbb Z^d$ (or by any family of cubes with $O(1)$ distance from them, as in Section~7.2): by periodicity such cubes can be translated to the origin up to a translation of the integer vertices, and the two corresponding quantities differ by $O(S^{d-1})$ by the same extension argument. The proof for the perimeter formula \eqref{homformper} is the same, with cubes $Q^\nu_T$ translated by the lattice points close to the hyperplane $\Pi^\nu$, as in the construction of the upper bound.

\paragraph{Changing the boundary data for sets of finite perimeter.}
For Sobolev functions with growth of order $p$ the lemma on the change of boundary data was proved in the Appendix to Chapter 4 (Lemma~\ref{slicing-lemma}). We prove here the analogous statement for sets, used in the lower bound for perimeters.

\begin{lemma}\label{slicing-perim}
Let $Q$ be an open cube of side $1$ and centre $0$, let $H$ be a half space with $0\in\partial H$, let $\varphi\colon\mathbb R^d\times S^{d-1}\to[\alpha,\beta]$ be a Borel function, and let $E_j\subset Q$ be sets of finite perimeter with $E_j\to H\cap Q$ in $L^1(Q)$ and $\sup_j\mathrm{Per}(E_j,Q)<+\infty$. Then for every $\delta\in(0,1)$ there exist sets $E_j'$ of finite perimeter, with $E_j'=H$ in a neighbourhood of $\partial Q$, such that
$$
\limsup_{j\to+\infty}\int_{Q\cap\partial^*E_j'}\varphi(x,\nu_{E_j'})\,d\mathcal H^{d-1}\le\liminf_{j\to+\infty}\int_{Q\cap\partial^*E_j}\varphi(x,\nu_{E_j})\,d\mathcal H^{d-1}+C\beta\delta,
$$
with $C$ depending only on $d$.
\end{lemma}

\begin{proof}
Let $Q_t$ be the concentric cube of side $t$, $t\in(1-\delta,1)$, and let $\psi(x)=\|x\|_\infty$. By the coarea formula applied to the function $\psi$, whose gradient has modulus $1$ almost everywhere,
$$
\int_{1-\delta}^1\mathcal H^{d-1}\big(\partial Q_{t}\cap(E_j\triangle H)\big)\,dt\le|E_j\triangle H|\to0.
$$
Hence, for fixed $\delta$, we can choose $t_j\in(1-\delta,1)$ (not belonging to the exceptional set of $t$ for which the slice is not well defined) with $\mathcal H^{d-1}(\partial Q_{t_j}\cap(E_j\triangle H))\le|E_j\triangle H|/\delta\to0$. Set $E_j'=(E_j\cap Q_{t_j})\cup(H\setminus Q_{t_j})$. Its boundary in $Q$ consists of the part of $\partial^*E_j$ inside $Q_{t_j}$, of the part of $\partial^*H$ in $Q\setminus Q_{t_j}$, whose $\mathcal H^{d-1}$ measure is at most $C\delta$, and of the part of $\partial Q_{t_j}$ where $E_j$ and $H$ differ, whose measure tends to $0$. Using $\varphi\le\beta$ on the last two parts we get the claim.
\end{proof}

By a diagonal argument on $\delta\to0$ the energy of the modified sets is not larger than that of the original ones up to an infinitesimal error. Applying this to the rescaled sets $\frac{1}{\rho_\e}(E_\e-x_0)$ gives the change of boundary data used in the lower bound for perimeters.

\paragraph{Two-scale convergence and unfolding.} Let $\Omega_\e=\bigcup_{k\in I^\e}Q^\e_k$, $Y=[0,1]^d$, and let $\mathcal U_\e$ be the unfolding operator of Section~7.1, $\mathcal U_\e(u)(x,y)=u(\e k+\e y)$ for $x\in Q^\e_k$, $y\in Y$.

\begin{proposition}\label{unfold}
Let $1<p<+\infty$ and let $(u_\e)$ be bounded in $L^p(\Omega)$. Then, up to subsequences, $\mathcal U_\e(u_\e)$ converges weakly in $L^p(\Omega\times Y)$ to some $u$, and $u_\e$ two-scale converges to $u$. Conversely, if $u_\e$ two-scale converges to $u\in L^p(\Omega\times Y)$ then $\mathcal U_\e(u_\e)\rightharpoonup u$ weakly in $L^p(\Omega\times Y)$.
\end{proposition}

\begin{proof}
By the change of variables $x=\e k+\e y$ in each cube we have
$$
\int_{\Omega_\e\times Y}|\mathcal U_\e(u_\e)|^p\,dx\,dy=\sum_{k\in I^\e}\e^d\int_Y|u_\e(\e k+\e y)|^pdy=\int_{\Omega_\e}|u_\e|^pdx,
$$
so that $\mathcal U_\e(u_\e)$ is bounded in $L^p(\Omega\times Y)$ and has weakly convergent subsequences. Let now $\phi\in C_c(\Omega;C_\#(Y))$ (continuous in $x$ with compact support, periodic in $y$). Since $x/\e=k+y$ at the point $x=\e k+\e y$ and $\phi$ is periodic in $y$,
$$
\int_\Omega u_\e(x)\phi\Big(x,\frac x\e\Big)dx=\int_{\Omega_\e\times Y}\mathcal U_\e(u_\e)(x,y)\,\phi(\e k_x+\e y,y)\,dx\,dy+o(1),
$$
where $k_x$ is the index of the cube containing $x$, and the $o(1)$ takes into account the contribution of $\Omega\setminus\Omega_\e$, whose measure tends to $0$ (and $\phi$ has compact support). Since $|\e k_x+\e y-x|\le\e\sqrt d$ the uniform continuity of $\phi$ in $x$ and the boundedness of $\mathcal U_\e(u_\e)$ allow us to replace $\phi(\e k_x+\e y,y)$ by $\phi(x,y)$ with an error that tends to $0$. Hence
$$
\int_\Omega u_\e(x)\phi\Big(x,\frac x\e\Big)dx=\int_{\Omega\times Y}\mathcal U_\e(u_\e)\,\phi\,dx\,dy+o(1)
$$
for all such $\phi$. If $\mathcal U_\e(u_\e)\rightharpoonup u$ the right-hand side converges to $\int_{\Omega\times Y}u\phi$, which is the two-scale convergence for these test functions, and by density and boundedness for all admissible test functions in $L^{p'}(\Omega;C_\#(Y))$. Conversely, two-scale convergence gives the convergence of the integrals of $\mathcal U_\e(u_\e)\phi$ for a class of $\phi$ which is dense in $L^{p'}(\Omega\times Y)$; since $\mathcal U_\e(u_\e)$ is bounded, this gives the weak convergence to $u$.
\end{proof}

If $u_\e(x)=v(x/\e)$ with $v$ $1$-periodic, then $\mathcal U_\e(u_\e)(x,y)=v(y)$ for all $\e$ and the two-scale limit is $v$. If instead $v(x)=e^{2\pi i\omega x}$ with $\omega\in\mathbb R\setminus\mathbb Q$ ($d=1$), then $\mathcal U_\e(u_\e)(x,y)=e^{2\pi i\omega y}e^{2\pi i\omega\lfloor x/\e\rfloor}$, and for every interval $(a,b)$ the sum $\sum_{a/\e<k<b/\e}e^{2\pi i\omega k}$ is bounded by $2/|1-e^{2\pi i\omega}|$, so that $\int_a^be^{2\pi i\omega\lfloor x/\e\rfloor}dx\to0$ and the two-scale limit is $0$, the mean value of $v$.

\paragraph{The checkerboard.} Let $a(y)=\alpha$ if $\lfloor2y_1\rfloor+\lfloor2y_2\rfloor$ is even and $a(y)=\beta$ otherwise ($y\in\mathbb R^2$), where $\alpha,\beta>0$; this is a $1$-periodic checkerboard of squares of side $\frac12$. Let $f^\#(\xi)=\min\{\int_{Q_1}a|\xi+\nabla v|^2dy:\ v\ 1\hbox{-periodic}\}$, which is a quadratic form $\langle A\xi,\xi\rangle$ with $A$ symmetric.

\begin{proposition}\label{checker}
$A=\sqrt{\alpha\beta}\,I$.
\end{proposition}

\begin{proof}
{\em Step 1: $A$ is a multiple of the identity.} The function $a$ is invariant (up to negligible sets) under the maps $y\mapsto(y_2,y_1)$ and $y\mapsto(\frac12-y_1,y_2)$. If $M$ is the orthogonal matrix of one of these maps, then composing the test functions with the map shows that $f^\#(M\xi)=f^\#(\xi)$, that is $M^TAM=A$. The reflection gives $A_{12}=0$, and the exchange of the variables gives $A_{11}=A_{22}$. Hence $A=\lambda(\alpha,\beta)I$.

{\em Step 2: scaling and exchange.} Clearly $\lambda(t\alpha,t\beta)=t\lambda(\alpha,\beta)$ for $t>0$. Moreover, the translation $y\mapsto y+(\frac12,0)$ exchanges the roles of $\alpha$ and $\beta$, and the cell formula is translation invariant, so that $\lambda(\alpha,\beta)=\lambda(\beta,\alpha)$.

{\em Step 3: duality.} Let $a$ be a general periodic coefficient with $0<c\le a\le C$ and let $A_a$ be the corresponding homogenized matrix. For $\xi\in\mathbb R^2$ let $v$ be the periodic minimizer in the cell problem, and set $E=\xi+\nabla v$, $J=aE$. Then $\mathrm{div}\,J=0$ and the mean of $J$ is $A_a\xi$. Let $R(x_1,x_2)=(-x_2,x_1)$. Since $\mathrm{div}\,J=\mathrm{curl}(RJ)$ the periodic field $RJ$ is curl-free, hence $RJ=RA_a\xi+\nabla w$ with $w$ periodic. On the other hand $RE$ is divergence-free, with mean $R\xi$, and $RE=RJ/a$. Setting $b=1/a$ we then have $\mathrm{div}(b\,(RA_a\xi+\nabla w))=0$, which means that $w$ is the minimizer of the cell problem for $b$ with datum $RA_a\xi$, and that the mean of $b(RA_a\xi+\nabla w)$, which is $R\xi$, equals $A_bRA_a\xi$. Hence $A_bRA_a=R$ and $A_{1/a}=RA_a^{-1}R^{-1}$.

{\em Step 4.} Applying Step 3 to $a$ (with values $\alpha,\beta$) we have $1/a$ with values $1/\alpha,1/\beta$ and $\lambda(1/\alpha,1/\beta)=1/\lambda(\alpha,\beta)$. By Step 2, $\lambda(1/\alpha,1/\beta)=\frac1{\alpha\beta}\lambda(\beta,\alpha)=\frac1{\alpha\beta}\lambda(\alpha,\beta)$, so that $\lambda(\alpha,\beta)^2=\alpha\beta$.
\end{proof}

The coefficient $\sqrt{\alpha\beta}$ is the geometric mean, which is larger than the harmonic mean $\frac{2\alpha\beta}{\alpha+\beta}$ given by the one-dimensional sections (as it must be), and smaller than the arithmetic mean $\frac{\alpha+\beta}2$ (the energy of the constant test function $v=0$).

\section{Concentration and homogenization of domains}
We now examine two instances where the small parameter $\e$ describes some ``extreme'' geometric features of the domain of the functionals. In the first case we treat {\em thin films}; that is, domains with one small dimension $\e$, while in the second one we treat periodic {\em perforated domains} (with Neumann boundary conditions); that is, energies whose integrand is $0$ on a set obtained by scaling by $\e$ a given periodic Lipschitz domain in $\mathbb R^d$.

\subsection{Thin objects and dimension reduction}\label{sec-thin}
A subset of $\mathbb R^d$ of the form
$$
\Omega_\e=\omega\times (0,\e),
$$
where $\omega\subset\mathbb R^{d-1}$ is a Lipschitz bounded open set,
is called a {\em thin film}. Similarly one can consider other {\em thin objects} of the form
$\omega_{d-k}\times \e\, \omega_{k}$, where $\omega_k$ and $\omega_{d-k}$ are Lipschitz bounded open sets in $\mathbb R^k$ and $\mathbb R^{d-k}$, respectively.

\begin{figure}[h!]
	\centerline{\includegraphics [width=5.8in]{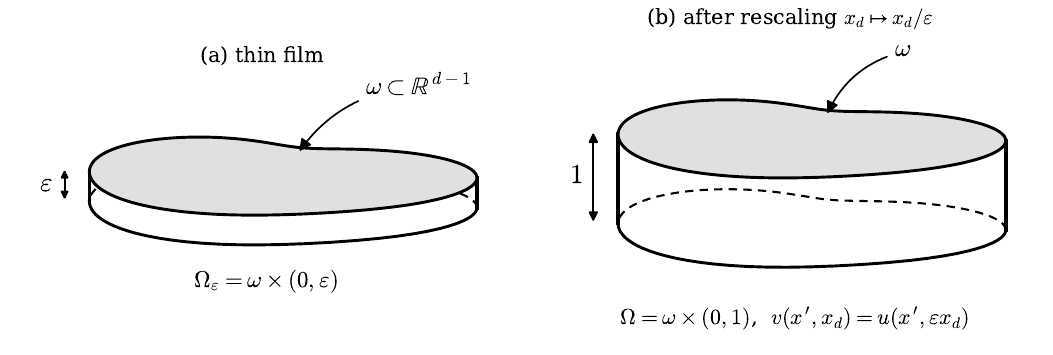}}
	\caption{(a) A thin film $\Omega_\e=\omega\times(0,\e)$, with $\omega\subset\mathbb R^{d-1}$ a bounded Lipschitz open set. (b) The same set after the rescaling $x_d\mapsto x_d/\e$ of the thin direction, which transforms functions $u$ on $\Omega_\e$ into functions $v(x',x_d)=u(x',\e x_d)$ on the fixed domain $\Omega=\omega\times(0,1)$.}\label{Fig47}
\end{figure}

Functionals defined on $\Omega_\e$ can be seen as approximations of corresponding functionals defined on $\omega$ as $\e\to 0$. In order to make this dimension-reduction process rigorous, we have to introduce a notion of {\em dimension-reduction convergence} of functions defined on $\Omega_\e$ to functions defined on $\omega$, which makes the functionals we are interested in equi-coercive. To that end we consider the simple case of the Dirichlet integral
$$
\int_{\Omega_\e} |\nabla u|^2dx,
$$
defined on $H^1(\Omega_\e)$, and, in order to work on a single set, scale the domain to $\Omega=\omega\times(0,1)$.
We now introduce the notation $x'=(x_1,\ldots,x_{d-1})$ and write $x=(x',x_d)$. Similarly, quantities or operations with a prime will refer to their $d-1$-dimensional version; e.g.~$\nabla'$ or $dx'$.
We can then write
$$
\int_{\Omega_\e} |\nabla u|^2dx=\e\int_{\Omega} |\nabla' v|^2dx+\frac1\e \int_{\Omega} \Big|\frac{\partial v}{\partial x_d}\Big|^2dx,
$$
where
$
v(x)= u(x',\e x_d).
$
This shows that, if we consider the scaled energies
$$
F_\e(u)=\frac1\e\int_{\Omega_\e} |\nabla u|^2dx,
$$
and $F_\e(u_\e)\le S<+\infty$, then, upon defining
$
v_\e(x)= u_\e(x',\e x_d)
$
as above, we obtain
$$
\int_{\Omega} |\nabla' v_\e|^2dx+\frac1{\e^2} \int_{\Omega} \Big|\frac{\partial v_\e}{\partial x_d}\Big|^2dx
\le S<+\infty,
$$
from which we first obtain that $\int_{\Omega} |\nabla v_\e|^2dx\le S<+\infty$ for $\e\le 1$, so that, up to addition of constants and extraction of subsequences, the sequence $v_\e$ converges weakly to some $v\in H^1(\Omega)$, and
$$
\int_{\Omega} \Big|\frac{\partial v}{\partial x_d}\Big|^2dx\le \liminf_{\e\to 0}\int_{\Omega} \Big|\frac{\partial v_\e}{\partial x_d}\Big|^2dx\le  \liminf_{\e\to 0}
\e^2 S=0,
$$
so that there exists $u\in H^1(\omega)$ such that $v(x)=u(x')$. This function $u$ is the {\em dimensionally reduced limit} of $u_\e$. More precisely, we say that $u_\e$ converges to $u$ {\em in the sense of dimension reduction} if $v_\e\to u$ in $L^p(\Omega)$, where $u$ is seen as a function on $\Omega$ independent of $x_d$ (here $p=2$; the same argument works for every $p>1$). Note that the $\Gamma$-limit of $F_\e$ with respect to this
convergence is the ``dimensionally reduced Dirichlet integral'' $\int_\omega|\nabla'u|^2dx'$.

\begin{figure}[h!]
	\centerline{\includegraphics [width=4.8in]{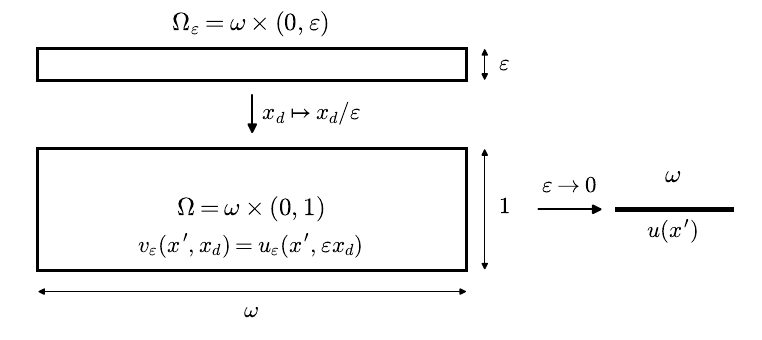}}
	\caption{A thin film $\Omega_\e=\omega\times(0,\e)$ is rescaled in the thin direction to $\Omega=\omega\times(0,1)$, where $v_\e(x',x_d)=u_\e(x',\e x_d)$; as $\e\to0$ the functions $v_\e$ converge to a function $u(x')$ on $\omega$ (dimension-reduction convergence).}\label{Fig30}
\end{figure}

\smallskip
We now examine the effect of this dimension-reduction process on more general energies  of the form
$$
F_\e(u)=\frac1\e\int_{\Omega_\e} f(\nabla u)dx,
$$
where $f$ is a continuous function of $p$-growth with $p>1$; that is, $c|\xi|^p-1\le f(\xi)\le C(1+|\xi|^p)$ for some $0<c\le C$, so that the functionals are defined on  $W^{1,p}(\Omega_\e;\mathbb R^m)$. The compactness argument above applies verbatim, with $p$ in place of $2$.

\medskip\noindent
{\em The convex (or scalar) case.} We first analyze the case when $f$ is convex. For ease of notation, we consider the scalar case $m=1$, when the functionals are defined on $W^{1,p}(\Omega_\e)$. We note that the convexity assumption is not necessary if $m=1$, or in the vector case if $d=2$ (when $\omega$ is an interval and the limit object is a curve), up to a relaxation argument: it suffices to replace the function $\overline f$ defined below by its convex envelope.

Using the function $v$ as above we rewrite
$$
F_\e(u)=\int_{\Omega} f\Big(\frac{\partial v}{\partial x_1},\ldots, \frac{\partial v}{\partial x_{d-1}},\frac1\e\frac{\partial v}{\partial x_d}\Big)dx
=\int_{\Omega} f\Big(\nabla' v,\frac1\e\frac{\partial v}{\partial x_d}\Big)dx.
$$
Since we have no control on the last component, we simply minimize it out, and define
$\overline f\colon\mathbb R^{d-1}\to \mathbb R$ by
$$
\overline f(\xi')= \min\{f(\xi',\xi_d): \xi_d\in \mathbb R\}
$$
(the minimum is attained by coercivity).
The function $\overline f$ is convex on $\mathbb R^{d-1}$ and satisfies a $p$-growth condition, so that the corresponding integral functional is lower semicontinuous and, if $u_\e\to u$ in the sense of dimension reduction and  setting $v_\e(x)=u_\e(x',\e x_d)$ as above, we have
$$
\liminf_{\e\to 0} F_\e(u_\e)\ge \liminf_{\e\to 0} \int_{\Omega} \overline f(\nabla' v_\e)dx\ge \int_{\Omega} \overline f(\nabla' v)dx=\int_{\omega} \overline f(\nabla' u)dx'.
$$
To prove the converse inequality, it is enough to consider the case of a linear target function
$u(x')=\langle \xi',x'\rangle$. Correspondingly there exists a minimizer $\xi_d$ such that
$\overline f(\xi')= f(\xi',\xi_d)$, or, equivalently, such that, setting $\xi=(\xi',\xi_d)$ we have $ f(\xi)= \overline f(\xi')$. The recovery sequence is simply given by
$$
u_\e(x) = \langle \xi',x'\rangle + \xi_d x_d= \langle \xi,x\rangle.
$$
Note that correspondingly, we have $v_\e(x) = \langle \xi',x'\rangle + \e \xi_d x_d$, so that the dimension-reduction convergence to $u$ is easily checked. For a general (smooth) target $u$ the same construction is applied pointwise, $u_\e(x)=u(x')+x_d\,\xi_d(x')$ with $\xi_d(x')$ a minimizer for $\overline f(\nabla'u(x'))$, after a mollification of the (measurable) function $\xi_d$; the proof is completed by an approximation argument (see the Appendix). The construction of the recovery sequence shows an {\bf optimization of the derivative in the small direction at scale $\e$}.

\medskip\noindent
{\em The vector (quasiconvex) case.} In the vector case of functions in $W^{1,p}(\Omega_\e;\mathbb R^m)$, with $m>1$, we may assume that $f\colon \mathbb M^{m\times d}\to\mathbb R$ is quasiconvex. We can still minimize out the derivative in the $x_d$-direction and in the same way as above define $\overline f\colon\mathbb M^{m\times (d-1)}\to\mathbb R$.
Unfortunately, $\overline f$ may turn out not to be quasiconvex (easy examples may be constructed
for $f$ polyconvex, with $d=m=2$, such that $\overline f$ is not convex; see the Appendix). Hence, a further relaxation is necessary. If we denote by $Q\overline f$ the ($m\times (d-1)$)-quasiconvex envelope of $\overline f$, then, in the notation already used in the convex case, we have
$$
\liminf_{\e\to 0} F_\e(u_\e)\ge  \liminf_{\e\to 0} \int_{\Omega} Q\overline f(\nabla' v_\e)dx\ge\int_{\Omega} Q\overline f(\nabla' v)dx=\int_{\omega} Q\overline f(\nabla' u)dx'.
$$
Here the second inequality follows from the lower semicontinuity of the functional $v\mapsto\int_\Omega Q\overline f(\nabla'v)dx$, which holds since $Q\overline f$ is also quasiconvex as a function of $\xi\in\mathbb M^{m\times d}$ independent of the last column (it suffices to apply the definition of quasiconvexity on each slice $\{x_d=t\}$ and integrate in $t$).
To check the optimality of the lower bound we have to construct recovery sequences. It suffices to do this for an affine target function $u(x')= \xi' x'$. To that end we use the formula
$$
Q\overline f(\xi')=\lim_{S\to+\infty}
\frac{1}{S^{d-1}}\inf\Big\{   \int_{Q'_{S}}\overline f(\xi'+\nabla' w(y'))\, dy': w\in W^{1,p}_{\rm loc}(\mathbb R^{d-1};\mathbb R^m), w \hbox{ $S$-periodic}  \Big\},
$$
where $Q'_S=(0,S)^{d-1}$,
and for every $y'$ we let $\xi_d(y')\in\mathbb R^m$ be such that $f(\xi'+\nabla' w(y'),\xi_d(y'))=\overline f(\xi'+\nabla' w(y'))$.  Finally, for given $\delta=\delta_\e$, we set
$$
u_\e(x) = \xi' x'+\delta w\Big(\frac{x'}\delta\Big)+ x_d \xi_d\Big(\frac{x'}\delta\Big);
$$
that is,
$$
v_\e(x) = \xi' x'+\delta w\Big(\frac{x'}\delta\Big)+ \e x_d \xi_d\Big(\frac{x'}\delta\Big).
$$
Up to an approximation argument (we may take $w$ smooth and mollify the measurable function $\xi_d$), we may suppose that $y'\mapsto \xi_d(y')$ is a $C^1$-function, so that
$$
\nabla' v_\e(x) = \xi'+\nabla' w\Big(\frac{x'}\delta\Big)+ \frac\e\delta x_d \nabla' \xi_d\Big(\frac{x'}\delta\Big),\qquad \frac1\e\frac{\partial v_\e}{\partial x_d}(x)=\xi_d\Big(\frac{x'}\delta\Big),
$$
and this gives a recovery sequence if $\e<\!<\delta$, because this condition ensures that the last term in the sum expressing $\nabla' v_\e(x)$ is asymptotically negligible. The construction of the recovery sequence shows an {\bf optimization of the derivative in the small direction at scale $\e$ periodically modulated at a scale $\delta>\!>\e$}.  This characterization can also be expressed through a homogenization-type formula
in which oscillations in the $x'$-direction and optimization in the $x_d$-direction compete, showing that
\begin{eqnarray*}
&&Q\overline f(\xi')=\lim_{S\to+\infty}
\frac{1}{S^{d-1}}\inf\Big\{   \int_{Q'_{S}\times (0,1)}f((\xi',0)+\nabla w(y))\, dy
\\
&& \hskip 2.5cm: w\in W^{1,p}_{\rm loc}(\mathbb R^{d-1}\times(0,1);\mathbb R^m), w \hbox{ $S$-periodic in the $x'$-directions}  \Big\}
\end{eqnarray*}
for $\xi'\in \mathbb M^{m\times (d-1)}$.

\begin{figure}[h!]
	\centerline{\includegraphics [width=5.2in]{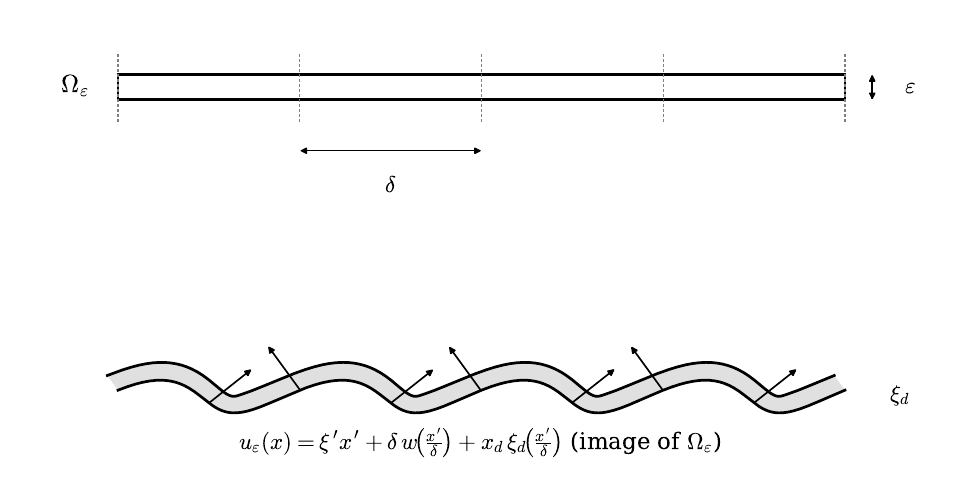}}
	\caption{Recovery sequence in the vector case. The film $\Omega_\e$ (top, thickness $\e$) is divided into periodicity cells of size $\delta\gg\e$ in the $x'$-direction. The image of $\Omega_\e$ under $u_\e$ (bottom) is a thin band, whose midline oscillates at scale $\delta$ around the affine map $\xi'x'$ through $\delta w(x'/\delta)$, and whose thickness vector $\e\,\xi_d(x'/\delta)$ is optimized (the arrows) pointwise at the same scale $\delta$.}\label{Fig31}
\end{figure}

\subsection{Perforated domains and double porosity}
{\em Perforated Domains.} A different type of varying domain, of oscillating type, is given by periodically perforated domains. These are obtained by taking a periodic connected open set $E$ of $\mathbb R^d$ with Lipschitz boundary, and considering energies on functions defined on the domains $\Omega_\e=\Omega\cap  \e E$, with $\Omega$ an open set of $\mathbb R^d$. For simplicity we consider the case of $\Omega=\mathbb R^d$, in order to remove the technical issue of a lack of regularity of $\Omega_\e$ at points in $\partial\Omega$. As a model case we consider functionals of the form
\begin{equation}\label{perfe}
F^E_\e(u)=\int_{\e E} |\nabla u|^pdx +\int_{\e E} |u|^pdx,
\end{equation}
with $p>1$, with domain $W^{1,p}(\e E)$. As for thin films, it is convenient to interpret these energies as defined on a common space. To that end, we suppose that there exists a linear {\em extension operator} $T\colon W^{1,p}(E)\to W^{1,p}(\mathbb R^d)$ such that

(i) (extension) $Tu(x)=u(x)$ if $x\in E$

(ii) (bound on the norm of the gradient) there exists a constant $C$ such that
$$
\int_{\mathbb R^d}|\nabla Tu|^pdx\le C\int_{E} |\nabla u|^pdx
$$

(iii) (bound on the norm of the function) there exists a constant $C$ such that
$$
\int_{\mathbb R^d}|Tu|^pdx\le C\int_{E} |u|^pdx.
$$

\noindent Note that the linear operator $T_\e\colon W^{1,p}(\e E)\to W^{1,p}(\mathbb R^d)$, defined
by setting $T_\e u(x)$ as the value at $x/\e$ of the function obtained by applying $T$ to $y\mapsto u(\e y)$, satisfies the properties above, with the same constant $C$, with $\e E$ in place of $E$. We say that $u_\e\in W^{1,p}(\e E)$ converges to $u\in W^{1,p}(\mathbb R^d)$ if $T_\e u_\e$ converges weakly to $u$ in $W^{1,p}(\mathbb R^d)$. With these extension operators we obtain some equi-coerciveness properties. Indeed, if $F^E_\e(u_\e)\le S<+\infty$ then $T_\e u_\e$ is a sequence bounded in $W^{1,p}(\mathbb R^d)$, from which we can extract a sequence weakly converging to some $u\in W^{1,p}(\mathbb R^d)$ (and strongly in $L^p_{\rm loc}(\mathbb R^d)$). It must be noted that this function $u$ does not depend on the particular choice of the extension operator. Indeed if $T'$ is another extension operator, and $T'_\e u_\e\to u'$ then we deduce that $u'=u$, letting $\e\to0$ in the equality $T_\e u_\e\chi_{\e E}= T'_\e u_\e\chi_{\e E}$, using the strong convergence in $L^p_{\rm loc}$ of $T_\e u_\e$ and $T'_\e u_\e$, and the weak$^*$ convergence in $L^\infty$ of $\chi_{\e E}$ to the strictly positive constant $\theta=|E\cap Q_1|$, where $Q_1=(0,1)^d$.

With this improved coerciveness, the $\Gamma$-limit can be computed as for the homogenization of usual (pointwise) coercive integrals, obtaining
$$
F^E_{\rm hom}(u)=\int_{\mathbb R^d} f^E_{\rm hom}(\nabla u)dx+ \theta\int_{\mathbb R^d}|u|^pdx,
$$
where
$$
 f^E_{\rm hom}(\xi):=\inf\Big\{\int_{Q_{1}\cap E}|\xi+\nabla v(y)|^p dy: v\in W^{1,p}_{\rm loc}(E), v \hbox{ 1-periodic}   \Big\}.
$$
The second integral in $F^E_{\rm hom}$ comes directly from the continuous convergence of the second term in \eqref{perfe}.

\begin{figure}[h!]
	\centerline{\includegraphics [width=5.4in]{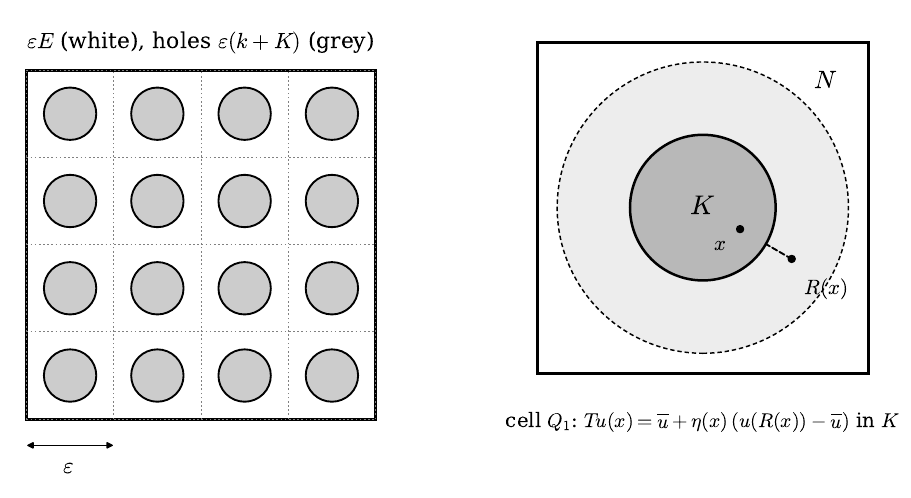}}
	\caption{Left: a periodically perforated domain $\e E$, obtained by removing the sets $\e(k+K)$ (in grey). Right: a cell $Q_1$ with the perforation $K$ and the neighbourhood $N$ of $\partial K$ in $E$ used in the reflection construction of the extension operator: inside $K$, $Tu$ is obtained from the values of $u$ at the reflected point $R(x)\in N$, cut off to the constant $\overline u$ away from $\partial K$.}\label{Fig32}
\end{figure}

\bigskip\noindent
{\em High-contrast media.} The argument above suggests that the usual homogenization process for integral functionals can be carried out also in the case when the energy densities are coercive on an infinite Lipschitz domain. We now show that this fact must be interpreted with care. We will consider a small variation of the energy above, taking
$$
F_\e(u)=\int_{\e E} |\nabla u|^pdx +\e^p \int_{\mathbb R^d\setminus\e E} |\nabla u|^pdx+\int_{\mathbb R^d} |u|^pdx,
$$
now defined on $W^{1,p}(\mathbb R^d)$.
Note that, if we omit the second integral and replace the third term with the integral $\int_{\e E} |u|^pdx$, we obtain \eqref{perfe}.

To understand why the contribution of the second term in $F_\e$ is not asymptotically negligible,
we will consider the simplest case of disjoint bounded perforations; that is,
$$
E=\mathbb R^d\setminus (\mathbb Z^d+ K),
$$
where $d\ge 2$ and $K\subset(0,1)^d$ is the closure of a bounded open set of class $C^2$ with connected boundary (so that the sets $k+K$, $k\in\mathbb Z^d$, are pairwise disjoint and $E$ is connected). In this case, the existence of an extension operator $T$ is obtained by a textbook ``periodic reflection argument'' in a neighbourhood of $\partial K$ (see the Appendix), and $\theta=1-|K|$. Hence, if we take $u_\e$ such that $T_\e(u_\e)\to u$ we can write
$$
F_\e(u_\e)=F^E_\e(u_\e) +\sum_{k\in\mathbb Z^d}\int_{\e(k+K)} (\e^p |\nabla u_\e|^p+|u_\e|^p)dx.
$$
In order to estimate the sum, we note that, up to a slight variation of boundary values (made precise in the Appendix), we can suppose that $u_\e= \overline u^\e_k$ on $\partial (\e k+\e K)$, where $\overline u^\e_k=\e^{-d}\theta^{-1}\int_{\e(k+(Q_1\setminus K))}u_\e\,dx$.
We then have
\begin{eqnarray*}
\int_{\e(k+K)} (\e^p |\nabla u_\e|^p+|u_\e|^p)dx\ge \min\Big\{\int_{\e(k+K)} (\e^p |\nabla v|^p+|v|^p)dx
:v=\overline u^\e_k \hbox{ on } \partial (\e k+\e K)\Big\}\\
=\e^d\min\Big\{\int_{K} (|\nabla w|^p+|w|^p)dx
:w=\overline u^\e_k \hbox{ on } \partial  K\Big\}= C_K\e^d|\overline u^\e_k|^p ,
\end{eqnarray*}
where $
C_K=\min\big\{\int_{K} (|\nabla w|^p+|w|^p)dx
:w=1 \hbox{ on } \partial  K\big\}$. By plugging these terms into the estimate of $F_\e(u_\e)$, we deduce that
$$
\liminf_{\e\to 0} F_\e(u_\e)\ge F^E_{\rm hom} (u)+C_K\int_{\mathbb R^d}|u|^pdx.
$$
This inequality is actually sharp. This can be proved by choosing $u_\e$ as a recovery sequence for $F^E_{\rm hom}(u)$ in $\e E$, and using optimal $w(x/\e)$ inside the sets $\e(k+K)$. Hence the $\Gamma$-limit of $F_\e$ is
$$
\int_{\mathbb R^d} f^E_{\rm hom}(\nabla u)dx+ (1-|K|+C_K)\int_{\mathbb R^d}|u|^pdx.
$$

Hence, the introduction of a critical coefficient $\e^p$ in the perforation produces a non-trivial interaction between the gradient and the lower-order term. The nature of this interaction is clearer if we slightly modify
$F_\e$ as
$$
F^\lambda_\e(u)=\int_{\e E} |\nabla u|^pdx +\e^p \int_{\mathbb R^d\setminus\e E} |\nabla u|^pdx+\int_{\e E} |u|^pdx+\lambda \int_{\mathbb R^d\setminus \e E} |u|^pdx,
$$
in which case, by the same argument, the $\Gamma$-limit is
$$
\int_{\mathbb R^d} f^E_{\rm hom}(\nabla u)dx+ (1-|K|+C_K(\lambda))\int_{\mathbb R^d}|u|^pdx,
$$
where the coefficient of the extra term is
$$
C_K(\lambda)=\min\Big\{\int_{K} (|\nabla w|^p+\lambda |w|^p)dx
:w=1 \hbox{ on } \partial  K\Big\}.
$$
In particular, for $p=2$ we can write
$$
C_K(\lambda)=\sqrt\lambda \min\Big\{\int_{K} \Big(\frac1{\sqrt\lambda}|\nabla w|^2+\sqrt\lambda |w|^2\Big)dx :w=1 \hbox{ on } \partial  K\Big\}.
$$
In terms of $w$ this is a singular-perturbation problem with a single well at $0$, and the boundary condition $w=1$ on $\partial K$; hence, as $\lambda\to+\infty$, the minimum produces a boundary layer of thickness $1/\sqrt\lambda$, with energy tending to a constant times the $\mathcal H^{d-1}$ measure of $\partial K$; namely, $C_K(\lambda)=\sqrt\lambda\,(\mathcal H^{d-1}(\partial K)+o(1))$ (see the Appendix). Summarizing, this {\em double porosity} limit exhibits a {\bf periodic behaviour at scale $\e$ with a concentration effect (at scale $
\e/\sqrt\lambda$) within the perforation}.

\begin{figure}[h!]
	\centerline{\includegraphics [width=5.4in]{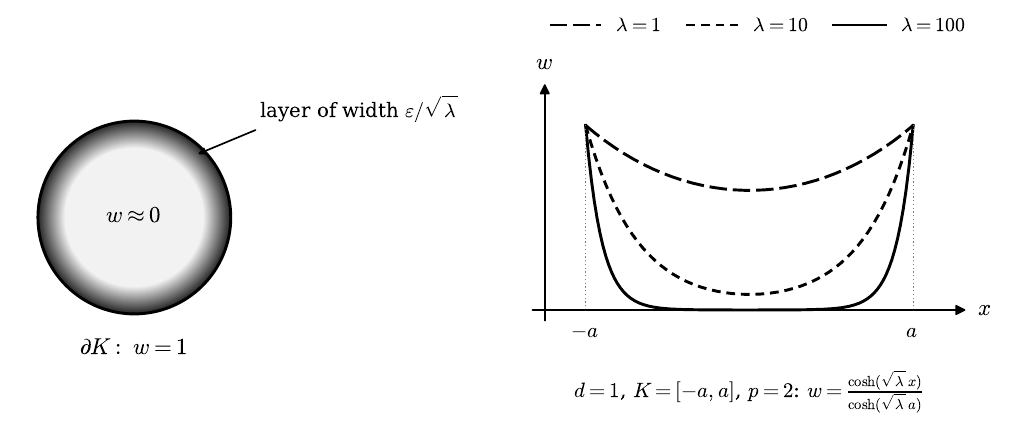}}
	\caption{High-contrast media with the term $\lambda\int|u|^p$ in the perforations. Left: the optimal $w$ for $C_K(\lambda)$ equals $1$ on $\partial K$ and is concentrated in a boundary layer of width $1/\sqrt\lambda$ (width $\e/\sqrt\lambda$ in the original variables). Right: the explicit profile $w(x)=\cosh(\sqrt\lambda x)/\cosh(\sqrt\lambda a)$ for $d=1$, $K=[-a,a]$, $p=2$, for increasing values of $\lambda$.}\label{Fig33}
\end{figure}

\subsection{Appendix: relevant tools of Chapter 8}\label{App-Ch8}

\paragraph{References.} The derivation of lower-dimensional models by $\Gamma$-convergence starts with the elastic string of Acerbi, Buttazzo and Percivale \cite{AcBP91}, and with the nonlinear membrane model of Le Dret and Raoult \cite{LDR95}, in which the density is obtained by minimizing out the normal derivative and then taking the quasiconvex envelope, as in the vector case of Section~8.1. A general treatment, including asymptotic developments, is in Anzellotti, Baldo and Percivale \cite{AnBP94}, and the case of inhomogeneous films, where oscillations of the energy density compete with the thickness, is in Braides, Fonseca and Francfort \cite{BFF00}; see also the chapter on dimension-reduction problems of \cite{Bra02}. Homogenization in perforated domains is due to Cioranescu and Saint Jean Paulin \cite{CSJP79}; the extension operators in general periodic domains are in Acerbi, Chiad\`o Piat, Dal Maso and Percivale \cite{ACDP92}, and a complete treatment of connected media (Chapter~19 and Appendix~B) is in \cite{BD98}; see also the chapter on Dirichlet problems in perforated domains of \cite{Bra02}. The critical size of the perforations, giving rise to an extra zero-order term, is the ``strange term'' of Cioranescu and Murat \cite{CM82}. The high-contrast scaling $\e^p$ and the double-porosity models are discussed by Arbogast, Douglas and Hornung \cite{ADH90}, Zhikov \cite{Zhi00}, and Bouchitt\'e and Felbacq \cite{BoF04}, who also consider the resonant regime.

\paragraph{Dimension reduction: the convex case.}
We prove the upper bound of Section~8.1 in the scalar case with $f$ convex; the lower bound is proved in the text. We assume $c|\xi|^p-1\le f(\xi)\le C(1+|\xi|^p)$. Since $f$ is convex and finite it is continuous, $\overline f$ is convex (partial minimization of a convex function) and $c|\xi'|^p-1\le\overline f(\xi')\le C(1+|\xi'|^p)$.

\begin{proposition}\label{dimred}
Let $F_\e$ be as in Section~8.1 with $f$ convex. Then the $\Gamma$-limit of $F_\e$ with respect to the dimension-reduction convergence is $F(u)=\int_\omega\overline f(\nabla'u)\,dx'$ if $u\in W^{1,p}(\omega)$, and $+\infty$ otherwise.
\end{proposition}

\begin{proof}
The compactness and the lower bound were proved in the text, together with the fact that the limit belongs to $W^{1,p}(\omega)$ if the energies are bounded. As for the upper bound, since the $\Gamma$-limsup $F''$ is lower semicontinuous and $u\mapsto\int_\omega\overline f(\nabla'u)dx'$ is continuous with respect to the strong convergence in $W^{1,p}(\omega)$, and since $\omega$ is Lipschitz, it suffices to prove $F''(u)\le\int_\omega\overline f(\nabla'u)dx'$ for $u\in C^\infty(\overline\omega)$. Let $g\colon\omega\to\mathbb R$ be a Borel function such that $f(\nabla'u(x'),g(x'))=\overline f(\nabla'u(x'))$ for every $x'$ (take, e.g., the minimizer of smallest modulus, which is a Borel function of $\xi'$ since the set of minimizers is a compact interval depending in an upper semicontinuous way on $\xi'$). By coercivity $g$ is bounded, as $\nabla'u$ is bounded. Let $g_\eta=g*\rho_\eta$ be mollifications (extending $g$ by $0$ outside $\omega$), so that $g_\eta$ are smooth, $\|g_\eta\|_\infty\le\|g\|_\infty$ and $g_\eta\to g$ a.e.~as $\eta\to0$. For fixed $\eta$ we set
$$
u_{\e}(x)=u(x')+x_d\,g_\eta(x'),\qquad v_\e(x)=u(x')+\e x_d\,g_\eta(x').
$$
Then $v_\e\to u$ uniformly, and
$$
F_\e(u_\e)=\int_\Omega f\big(\nabla'u+\e x_d\nabla'g_\eta,\ g_\eta\big)dx\ \longrightarrow\ \int_\Omega f(\nabla'u,g_\eta)\,dx=\int_\omega f(\nabla'u,g_\eta)\,dx'
$$
as $\e\to0$, since the integrand is uniformly continuous on a bounded set. Hence $F''(u)\le \int_\omega f(\nabla' u,g_\eta)dx'$ for every $\eta$. As $\eta\to0$, $f(\nabla'u,g_\eta)\to f(\nabla'u,g)=\overline f(\nabla'u)$ a.e., with the uniform bound $C(1+\|\nabla'u\|_\infty^p+\|g\|^p_\infty)$, so by dominated convergence $F''(u)\le\int_\omega\overline f(\nabla'u)dx'$.
\end{proof}

If $f$ is not convex, but $m=1$ or $d=2$, the same holds with $\overline f$ replaced by its convex envelope $C\overline f$ (which coincides with $Q\overline f$ in these cases, since quasiconvexity reduces to convexity on $\mathbb M^{1\times(d-1)}$ and on $\mathbb M^{m\times1}$). The lower bound is the same as above, since $f\ge\overline f\ge C\overline f$ and $C\overline f$ is convex, while the upper bound is obtained with the oscillating recovery sequences described below for the vector case.

\paragraph{Dimension reduction: the vector case.}
\begin{lemma}\label{slices}
Let $g\colon\mathbb M^{m\times(d-1)}\to\mathbb R$ be quasiconvex with $p$-growth and set $\hat g(\xi',\xi_d)=g(\xi')$ for $(\xi',\xi_d)\in\mathbb M^{m\times d}$. Then $\hat g$ is quasiconvex, so that $v\mapsto\int_\Omega g(\nabla'v)dx$ is sequentially lower semicontinuous for the weak convergence in $W^{1,p}(\Omega;\mathbb R^m)$.
\end{lemma}

\begin{proof}
Let $\varphi\in W^{1,\infty}_0(Q_1;\mathbb R^m)$, $Q_1=(0,1)^d$, and $\xi=(\xi',\xi_d)$. For a.e.~$t\in(0,1)$ the function $\varphi(\cdot,t)$ belongs to $W^{1,\infty}_0((0,1)^{d-1};\mathbb R^m)$ and its gradient is $\nabla'\varphi(\cdot,t)$. By Fubini's theorem and the quasiconvexity of $g$ on each slice,
$$
\int_{Q_1}\hat g(\xi+\nabla\varphi)dx=\int_0^1\int_{(0,1)^{d-1}}g(\xi'+\nabla'\varphi(y',t))dy'dt\ge g(\xi').
$$
The lower semicontinuity is then the classical theorem of Morrey and Acerbi--Fusco for quasiconvex integrands with $p$-growth.
\end{proof}

\noindent {\em The upper bound.} Let $u(x')=\xi'x'$, let $\nu>0$, and let $w$ be a smooth $S$-periodic function with $S^{1-d}\int_{Q'_S}\overline f(\xi'+\nabla'w)dy'\le Q\overline f(\xi')+\nu$ (smooth functions are dense for the formula, by the $p$-growth of $\overline f$). As in the scalar case, let $g$ be a bounded Borel selection of minimizers, $f(\xi'+\nabla'w(y'),g(y'))=\overline f(\xi'+\nabla'w(y'))$, and let $g_\eta$ be an $S$-periodic mollification of $g$. With $\delta=\delta_\e$, $\e/\delta\to0$, we put
$$
v_\e(x)=\xi'x'+\delta w\Big(\frac{x'}\delta\Big)+\e x_d\,g_\eta\Big(\frac{x'}\delta\Big).
$$
Then $v_\e\to u$ uniformly and
$$
F_\e(u_\e)=\int_\Omega f\Big(\xi'+\nabla'w\big(\tfrac{x'}\delta\big)+\tfrac\e\delta x_d\nabla'g_\eta\big(\tfrac{x'}\delta\big),\ g_\eta\big(\tfrac{x'}\delta\big)\Big)dx\to|\Omega|\,S^{1-d}\int_{Q'_S}f(\xi'+\nabla'w,g_\eta)\,dy',
$$
by the uniform continuity of $f$ on bounded sets and the Riemann--Lebesgue lemma applied to the periodic function $y'\mapsto f(\xi'+\nabla'w(y'),g_\eta(y'))$. Letting $\eta\to0$ by dominated convergence, and then $\nu\to0$ we get $F''(u)\le|\Omega|\,Q\overline f(\xi')$. The passage from affine to general $u\in W^{1,p}(\omega;\mathbb R^m)$ follows the usual localization argument of $\Gamma$-convergence (the $\Gamma$-limsup is the restriction to open sets of a measure, so that it can be computed on piecewise-affine functions, see \cite{BD98}), and is not repeated here.

\paragraph{An example where $\overline f$ is not convex.}
Let $d=m=2$, $p=4$ and $f(\xi)=|\xi|^4+(\det\xi-1)^2$, where $|\xi|$ is the Euclidean norm of the matrix. Then $f$ is polyconvex, being a convex function of $(\xi,\det\xi)$, and has $4$-growth in the above sense. Write $\xi=(a|b)$ with $a=\xi'$ and $b=\xi_2\in\mathbb R^2$. Then $\overline f(0)=\min_b(|b|^4+1)=1$. If $a\neq0$, $s=|a|$, and $b=t\,a^\perp/s$ with $t>0$, then $\det\xi=st$, $|\xi|^2=s^2+t^2$ and
$$
f(a|b)=(s^2+t^2)^2+(st-1)^2=1-2st+s^2t^2+(s^2+t^2)^2.
$$
For $t=s^{1/3}$ this equals $1-s^{4/3}+3s^{8/3}+s^{4}$, which is strictly less than $1$ for $s$ small and positive. Hence $\overline f(a)<\overline f(0)$ for small $a\ne0$. On the other hand $\overline f$ is even: $f(-a|-b)=f(a|b)$ gives $\overline f(-a)=\overline f(a)$. A convex even function satisfies $\overline f(0)\le\frac12(\overline f(a)+\overline f(-a))=\overline f(a)$, so $\overline f$ cannot be convex.

\paragraph{The extension operator for disjoint perforations.}
Let $d\ge2$ and let $K\subset(0,1)^d$ be the closure of a bounded open set of class $C^2$ with connected boundary, $E=\mathbb R^d\setminus(\mathbb Z^d+K)$, $Q_1=(0,1)^d$.

\begin{lemma}\label{ext}
There exists a linear extension operator $T\colon W^{1,p}(E)\to W^{1,p}(\mathbb R^d)$ satisfying (i)--(iii) of Section~8.2, which is local on cells: $Tu$ on $k+K$ depends only on $u$ in a fixed neighbourhood $k+N$ of $k+\partial K$ in $E$. More precisely, there exist a connected open set $N\subset Q_1\setminus K$ with $\partial K\subset\overline N$ and $C>0$ such that, with $\overline u=\frac1{|N|}\int_Nu\,dx$,
$$
\int_K|\nabla Tu|^pdx\le C\int_N|\nabla u|^pdx,\qquad \int_K|Tu-\overline u|^pdx\le C\int_N|u-\overline u|^pdx\le C\int_N|\nabla u|^pdx.
$$
\end{lemma}

\begin{proof}
Since $\partial K$ is of class $C^2$ and compact, there is $d_0>0$ such that the signed distance is $C^1$ in the tubular neighbourhood $\{{\rm dist}(x,\partial K)<d_0\}$, which is contained in $(0,1)^d$, and every point there can be written uniquely as $\pi-s\nu(\pi)$ with $\pi\in\partial K$, $\nu$ the outer normal, $s\in(-d_0,d_0)$ ($s>0$ in $K$). Let $R(\pi-s\nu)=\pi+s\nu$ be the reflection; it is a bi-Lipschitz map exchanging $K_{d_0}=\{x\in K:s<d_0\}$ and $N=\{x\notin K:{\rm dist}(x,K)<d_0\}$, the identity on $\partial K$. The set $N$ is homeomorphic to $\partial K\times(0,d_0)$, hence connected. Let $\zeta\in C^\infty(\mathbb R)$ with $\zeta=1$ on $(-\infty,d_0/4]$ and $\zeta=0$ on $[d_0/2,+\infty)$, and set $\eta(x)=\zeta(s(x))$ on $K_{d_0}$. We define $Tu=u$ on $E$, and on $K$
$$
Tu(x)=\overline u+\eta(x)\big(u(R(x))-\overline u\big)\ \hbox{ if }x\in K_{d_0},\qquad Tu(x)=\overline u\ \hbox{ otherwise}.
$$
Then $Tu$ is linear in $u$, is the identity on $E$, and has the same trace as $u$ on $\partial K$ since $\eta=1$ and $R={\rm id}$ there; moreover $Tu=\overline u$ near $\partial K_{d_0}$ inside $K$, so that $Tu\in W^{1,p}(Q_1)$. Since $\nabla Tu=\eta\,(\nabla u\circ R)\nabla R+(u\circ R-\overline u)\otimes\nabla\eta$, a change of variables and the Poincar\'e--Wirtinger inequality in the connected Lipschitz set $N$ give
$$
\int_K|\nabla Tu|^pdx\le C\int_N|\nabla u|^pdx+C\int_N|u-\overline u|^pdx\le C'\int_N|\nabla u|^pdx,
$$
and similarly for $\int_K|Tu-\overline u|^p$. Finally $\int_K|Tu|^pdx\le C|\overline u|^p+C\int_N|u-\overline u|^pdx\le C\int_N|u|^pdx$ by Jensen's inequality. Repeating the construction in each cell $k+Q_1$ (where $Tu=u$ near the boundary of the cell since $K\subset(0,1)^d$) gives (i)--(iii).
\end{proof}

For $T_\e$ the estimates are scaled: if $N_k^\e=\e(k+N)$, $K^\e_k=\e(k+K)$ and $\overline u_k$ is the average of $u$ on $N_k^\e$, then
$$
\e^p\int_{K_k^\e}|\nabla T_\e u|^pdx+\int_{K^\e_k}|T_\e u-\overline u_k|^pdx\le C\e^p\int_{N^\e_k}|\nabla u|^pdx.\eqno(8.1)
$$

\paragraph{The high-contrast limit.}
For $g$ in the trace space of $W^{1,p}(K)$ (the interior of $K$ is of class $C^2$) we set
$$
\Phi(g)=\min\Big\{\int_K(|\nabla w|^p+|w|^p)dx:\ w\in W^{1,p}(K),\ w=g\hbox{ on }\partial K\Big\},
$$
so that $C_K=\Phi(1)$. The function $\Phi^{1/p}$ is the norm of the quotient of $W^{1,p}(K)$ by $W^{1,p}_0(K)$: it is $p$-homogeneous, $\Phi(tg)=|t|^p\Phi(g)$, and satisfies the triangle inequality, since if $w_1,w_2$ are optimal for $g_1,g_2$ then $w_1+w_2$ is a competitor for $g_1+g_2$ and Minkowski's inequality applies. Moreover $\Phi({\rm tr}\,W)^{1/p}\le\|W\|_{W^{1,p}(K)}$ for every $W\in W^{1,p}(K)$. Hence
$$
\big|\Phi(g)^{1/p}-C_K^{1/p}|c|\big|\le\Phi(g-c)^{1/p}\le\|W-c\|_{W^{1,p}(K)}\eqno(8.2)
$$
for every constant $c$ and every $W$ with trace $g$.

\begin{proposition}\label{highcontrast}
Under the assumptions of Section~8.2 on $K$, the functionals $F_\e$ $\Gamma$-converge, with respect to the convergence $T_\e u_\e\rightharpoonup u$ in $W^{1,p}(\mathbb R^d)$, to
$$
F(u)=\int_{\mathbb R^d}f^E_{\rm hom}(\nabla u)dx+(1-|K|+C_K)\int_{\mathbb R^d}|u|^pdx.
$$
\end{proposition}

\begin{proof}
We use the $\Gamma$-convergence of $F^E_\e$ to $F^E_{\rm hom}$, stated in the text, as known. Let $u_\e\in W^{1,p}(\mathbb R^d)$, $U_\e=T_\e u_\e$ (where $T_\e$ is applied to $u_\e|_{\e E}$), and let $\overline u_k$ be the average of $u_\e$ on $N^\e_k$. Note that $U_\e$ and $u_\e$ have the same trace on $\partial K^\e_k$, since they coincide in $\e E$, so that, setting $W_k(y)=U_\e(\e k+\e y)$ and $w_k(y)=u_\e(\e k+\e y)$ for $y\in K$, we have $\Phi({\rm tr}\,w_k)=\Phi({\rm tr}\,W_k)$ and
$$
\int_{K^\e_k}(\e^p|\nabla u_\e|^p+|u_\e|^p)dx=\e^d\int_K(|\nabla w_k|^p+|w_k|^p)dy\ge\e^d\Phi({\rm tr}\,w_k).
$$
Set $a_k=\Phi({\rm tr}\,w_k)^{1/p}$, $b_k=C_K^{1/p}|\overline u_k|$ and $\rho_k=\|W_k-\overline u_k\|_{W^{1,p}(K)}$. By (8.2), $|a_k-b_k|\le\rho_k$, and by (8.1) and a change of variables,
$$
\e^d\rho_k^p=\e^p\int_{K^\e_k}|\nabla U_\e|^pdx+\int_{K^\e_k}|U_\e-\overline u_k|^pdx\le C\e^p\int_{N^\e_k}|\nabla u_\e|^pdx.
$$
Since the sets $N^\e_k$ are pairwise disjoint and contained in $\e E$, we get $\sum_k\e^d\rho_k^p\le C\e^p\int_{\e E}|\nabla u_\e|^pdx=C\e^pG_\e$, where $G_\e=\int_{\e E}|\nabla u_\e|^pdx$. By Minkowski's inequality in $\ell^p$ with weights $\e^d$,
$$
\Big|\Big(\sum_k\e^da_k^p\Big)^{1/p}-\Big(\sum_k\e^db_k^p\Big)^{1/p}\Big|\le C\e\,G_\e^{1/p}.\eqno(8.3)
$$
Let $\overline u^\e$ be the piecewise-constant function equal to $\overline u_k$ on $\e(k+Q_1)$. Since $N^\e_k$ is a subset of fixed proportion of the cell $\e(k+Q_1)$ and $U_\e|_{\e(k+Q_1)}$ is controlled by $u_\e$ in $\e(k+Q_1)\cap\e E$ through (ii), the Poincar\'e inequality on the cells gives $\|\overline u^\e-U_\e\|_{L^p(\mathbb R^d)}\le C\e G_\e^{1/p}$, and $\sum_k\e^db_k^p=C_K\|\overline u^\e\|^p_{L^p}$.

{\em Lower bound.} We may assume $F_\e(u_\e)\le S$ along the sequence, so that $U_\e\rightharpoonup u$ weakly in $W^{1,p}$ and strongly in $L^p_{\rm loc}$. Then $\overline u^\e\to u$ in $L^p_{\rm loc}$ and by (8.3), restricting the sums to the cells contained in a bounded set $A$,
$$
\liminf_{\e\to0}\sum_k\int_{K^\e_k}(\e^p|\nabla u_\e|^p+|u_\e|^p)dx\ge C_K\int_A|u|^pdx.
$$
Letting $A\uparrow\mathbb R^d$, and using $\liminf F^E_\e(u_\e)\ge F^E_{\rm hom}(u)$, we get $\liminf F_\e(u_\e)\ge F^E_{\rm hom}(u)+C_K\int|u|^p$.

{\em Upper bound.} Let $u\in W^{1,p}(\mathbb R^d)$ and let $v_\e\in W^{1,p}(\e E)$ be a recovery sequence for $F^E_{\rm hom}(u)$; then $\int_{\e E}|v_\e|^pdx\to\theta\int|u|^pdx$ and $\int_{\e E}|\nabla v_\e|^p dx\to \int f^E_{\rm hom}(\nabla u)dx$ since, the sum converging and each of the two terms having the corresponding liminf inequality, each term converges. We define $u_\e=v_\e$ in $\e E$ and, in each $K^\e_k$, $u_\e(\e k+\e y)=w_k(y)$ where $w_k$ is the minimizer for $\Phi({\rm tr}\,W_k)$, $W_k$ being built as above from $T_\e v_\e$; then $u_\e\in W^{1,p}(\mathbb R^d)$ since the traces on $\partial K^\e_k$ coincide, and the sum over $k$ in the definition of $F_\e(u_\e)$ is exactly $\e^d\sum_ka_k^p$. By (8.3) and the estimate on $\overline u^\e-T_\e v_\e$ (note that now $\overline u_k$ is the average of $v_\e$ on $N^\e_k$), we have
$$
\limsup_{\e\to0}\Big(\sum_k\e^da_k^p\Big)^{1/p}\le C_K^{1/p}\limsup_{\e\to0}\|\overline u^\e\|_{L^p(\mathbb R^d)}.
$$
Finally, by the Poincar\'e inequality on the sets $\e(k+Q_1)\setminus K^\e_k$ (connected Lipschitz sets of size $\e$) we have $\big|(\int_{\e E}|v_\e|^p)^{1/p}-(\theta\sum_k\e^d|\overline u_k|^p)^{1/p}\big|\le C\e G_\e^{1/p}$, from which $\|\overline u^\e\|^p_{L^p}\to\int_{\mathbb R^d}|u|^pdx$ (the factor $\theta$ cancels, since $\int_{\e E}|v_\e|^p\to\theta\int|u|^p$ and $\theta\sum_k\e^d|\overline u_k|^p=\theta\|\overline u^\e\|_{L^p}^p$). Therefore $\limsup F_\e(u_\e)\le F^E_{\rm hom}(u)+C_K\int|u|^pdx$.
\end{proof}

The same proof applies to $F^\lambda_\e$, with $\Phi$ replaced by $\Phi_\lambda(g)=\min\{\int_K(|\nabla w|^p+\lambda|w|^p):w={\rm tr}^{-1}g\}$ and $C_K$ by $C_K(\lambda)=\Phi_\lambda(1)$.

\paragraph{The coefficient $C_K(\lambda)$.}
Let $p=2$ and $d=1$, with $K=[-a,a]$ (this is only a cell computation, since $E$ is not connected if $d=1$). The Euler--Lagrange equation of the problem defining $C_K(\lambda)$ is $w''=\lambda w$ with $w(\pm a)=1$, so that $w(x)=\cosh(\sqrt\lambda x)/\cosh(\sqrt\lambda a)$, and integrating by parts
$$
C_K(\lambda)=\int_{-a}^a(w'^2+\lambda w^2)dx=\big[ww'\big]_{-a}^a=2\sqrt\lambda\tanh(\sqrt\lambda a)\ \sim\ 2\sqrt\lambda=\sqrt\lambda\,\mathcal H^0(\partial K)\quad(\lambda\to+\infty).
$$
The profile $w$ is concentrated in a layer of thickness $1/\sqrt\lambda$ around $\pm a$, which, in the variable $x$ of the functional $F^\lambda_\e$, has width $\e/\sqrt\lambda$.

In general, for $p>1$ and $K$ as in Section~8.2, one has
$$
\lim_{\lambda\to+\infty}\frac{C_K(\lambda)}{\lambda^{1/p'}}=c_p\,\mathcal H^{d-1}(\partial K),\qquad c_p=(p-1)^{-1/p'},\quad p'=\tfrac p{p-1},
$$
where $c_p$ is the value of the one-dimensional problem $\min\{\int_0^\infty(|w'|^p+|w|^p)ds:w(0)=1\}$ attained at $w(s)=e^{-s/(p-1)^{1/p}}$ (indeed, the first integral of the Euler--Lagrange equation is $(p-1)|w'|^p=|w|^p$, and then the energy is $\frac p{p-1}\int_0^\infty|w|^pds=(p-1)^{-1/p'}$). Moreover $C_K(\lambda)\le\lambda|K|$ (take $w\equiv1$), so $C_K(\lambda)\to0$ as $\lambda\to0$. We sketch the proof of the limit. Let $s$ denote the distance from $\partial K$ in $K$. {\em Upper bound}: take $w=\zeta(s)\,e^{-\kappa s}$ with $\kappa=(\lambda/(p-1))^{1/p}$ and a cut-off $\zeta$ equal to $1$ for $s<d_0/2$; for $s<d_0/2$ one has $|w'|^p+\lambda|w|^p=\frac{p\lambda}{p-1}e^{-p\kappa s}$, and the volume element is $(1+O(s))\,d\mathcal H^{d-1}ds$, so that the energy is $\lambda^{1/p'}c_p\,\mathcal H^{d-1}(\partial K)(1+O(\kappa^{-1}))+O(e^{-\kappa d_0})$. {\em Lower bound}: for every $\pi\in\partial K$, along the normal segment $s\in(0,d_0)$ the one-dimensional problem on $(0,d_0)$ with $w(0)=1$ and free end has minimum at least $c_p\lambda^{1/p'}(1-O(e^{-\kappa d_0}))$, and integrating in $\pi$ with the Jacobian $\ge1-Cd_0$ we get $\liminf_\lambda C_K(\lambda)\lambda^{-1/p'}\ge c_p(1-Cd_0)\mathcal H^{d-1}(\partial K)$ for every small $d_0$. For $p=2$ this reads $C_K(\lambda)=\sqrt\lambda\,(\mathcal H^{d-1}(\partial K)+o(1))$ and the optimal $w$ decays as $e^{-\sqrt\lambda s}$, which gives the concentration scale $\e/\sqrt\lambda$.

\section{Non-local energies}
In this section we consider functionals defined by double integrals, first remarking that lower-semicontinuity and convergence conditions may lead to ``non-representability''. This fact is made precise in the first section by a counterexample.

\subsection{Non-locality brings non-representability}\label{sec-nonrep}

\paragraph{Relaxation.} We first note that a sufficient condition for functionals of the form
\begin{equation}
F(u)=\int_{\Omega\times\Omega} f(u(x)-u(y))\,dxdy
\end{equation}
to be weakly lower semicontinuous with respect to the convergence in $L^1(\Omega)$ is that $f$ is (lower semicontinuous and) convex. This follows from the remark that if $u_j\wto u$ then the functions $U_j(x,y)= u_j(x)-u_j(y)$ converge weakly in $L^1(\Omega\times\Omega)$ to $U(x,y)= u(x)-u(y)$, and using the known weak lower semicontinuity theorem for convex energies. Under some mild assumptions, convexity turns out also to be necessary. Conversely, if $f$ is not convex then the lower-semicontinuous and convex envelope of $f$ is far from giving the integrand of the  
lower-semicontinuous envelope $\overline F$ of $F$. Indeed, $\overline F$ may not even be representable as a double integral.

To make a simple example, we consider an ``extreme'' $f$ given by
$$
f(z)=
\begin{cases} 0 &\hbox{ if $z\in\{-1,1\}$}\cr
1  &\hbox{ if $z=0$}\cr
+\infty &\hbox{ if $z\not\in\{-1,0,1\}$}.
\end{cases}
$$
Note that if $f$ is the increasing limit of a sequence $f_n$ of continuous functions, then by monotone convergence the functionals with integrand $f_n$ $\Gamma$-converge to $\overline F$ (the lower-semicontinuous envelope of $F$ computed below), so that the counterexample holds also for $f_n$, with $n$ large (see the Appendix).

Let $u_n\wto u$ in $L^1(\Omega)$ and $F(u_n)\le S<+\infty$. Then there exist
$c_n$ such that $u_n(x)\in\{c_n-\frac12,c_n+\frac12\}$ almost everywhere.
This implies that the oscillation of $u$ is at most $1$, and we can write 
$$F(u_n)= |\{u_n(x)=c_n-\tfrac12\}|^2+  |\{u_n(x)=c_n+\tfrac12\}|^2.$$ Note that such a quantity is minimal if $|\{u_n(x)=c_n-\tfrac12\}|=|\{u_n(x)=c_n+\tfrac12\}|$, so that it is always at least $\frac12|\Omega|^2$.

We examine two cases, which are the only ones needed in the following (the intermediate values $\frac12<{\rm osc}\,u<1$ will not be used): first that of shallow oscillation; that is, osc$\,u\le 1/2$. Up to considering $u-|\Omega|^{-1}\int_\Omega u\,dy$ in place of $u$, we can suppose that $\int_\Omega u\,dy=0$ and $-\frac12\le u\le \frac12$. If $u_n\wto u$ with $u_n(x)\in\{-\frac12,\frac12\}$ then from $$\tfrac12|\{u_n(x)=\tfrac12\}|-\tfrac12|\{u_n(x)=-\tfrac12\}|=\int u_n
dy\to\int udy=0, $$ we deduce that $F(u_n)\to\frac12|\Omega|^2$, which implies that $
\overline F(u)\le \frac12|\Omega|^2$. Since the converse inequality is always valid, we deduce that
$\overline F(u)=\frac12|\Omega|^2$.
 
Second, we examine the case of maximal oscillation, osc$\,u=1$. Up to a translation it is sufficient to examine the case ess-inf$\,u=0$ and ess-sup$\,u=1$. 
In this case, any approximating sequence $u_n$ is such that the corresponding $c_n$ tend to $1/2$, and we have $|\{u_n=c_n+\tfrac12\}|\to \int_\Omega udy$, $|\{u_n=c_n-\tfrac12\}|\to |\Omega|-\int_\Omega udy$, so that 
$$
\lim_{n\to+\infty}F(u_n)=\lim_{n\to+\infty}
\big( |\{u_n=c_n+\tfrac12\}|^2+|\{u_n=c_n-\tfrac12\}|^2\big)= \Big(\int_\Omega u\,dy\Big)^2+\Big(|\Omega|-\int_\Omega u\,dy\Big)^2.
$$
Since this value is independent of the sequence $u_n$, it equals $\overline F(u)$.

We now suppose that there exists $g$ such that
$$
\overline F(u)=\int_{\Omega\times\Omega} g(u(x)-u(y))\,dxdy,
$$
where we may suppose that $g$ is even, since the double integral is invariant under the exchange of $x$ and $y$.
The function $u=0$ is of the first type considered, so that $\overline F(0)=\frac12|\Omega|^2$ would give $g(0)=\frac12$. For $t\in[0,1]$ we consider $A_t\subset \Omega$ with $|A_t|=t|\Omega|$ and $u_t=\chi_{A_t}$, so that $\int_\Omega u_tdy=t|\Omega|$. This is a function with maximal oscillation so that we would obtain
$$
|\Omega|^2(t^2+(1-t)^2)=\overline F(u)=|\Omega|^2\big(g(0)(t^2+(1-t)^2)+g(1)(2t-2t^2)\big).
$$
Recalling that $g(0)=\frac12$, we deduce that
$$
g(1)=\frac{t^2+(1-t)^2}{4(t-t^2)}.
$$
This is a contradiction since the right-hand side is not a constant.

\begin{figure}[t]
	\centerline{\includegraphics [width=5.8in]{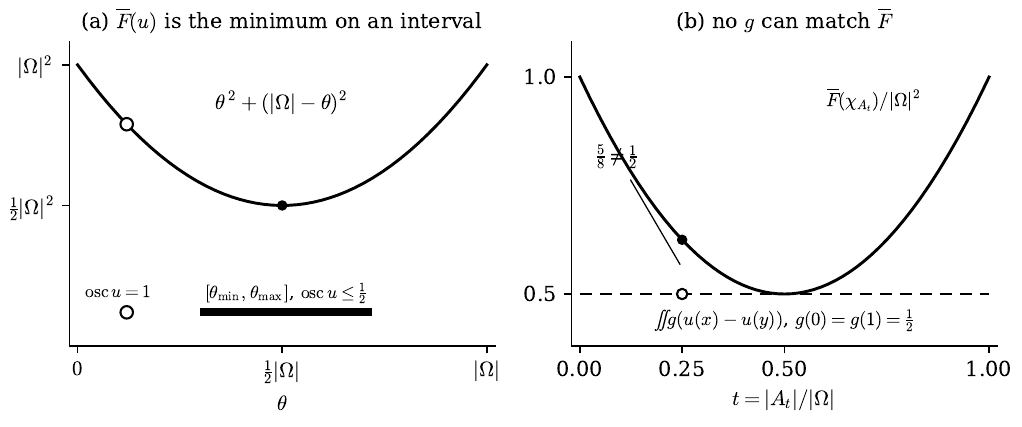}}
	\caption{Non-representability. (a) For $|\Omega|$ fixed, $\overline F(u)$ is the minimum of $\theta^2+(|\Omega|-\theta)^2$ over $[\theta_{\min},\theta_{\max}]$: when ${\rm osc}\,u\le\frac12$ the interval contains $\frac12|\Omega|$ and $\overline F(u)=\frac12|\Omega|^2$, while for ${\rm osc}\,u=1$ it reduces to a point and the value is read on the parabola. (b) The values $\overline F(\chi_{A_t})/|\Omega|^2=t^2+(1-t)^2$ (solid) and the only double integral compatible with the values at $u=0$ and $u=\chi_{A_{1/2}}$, which is constant equal to $\frac12$ (dashed): they differ at $t=\frac14$, where $\frac58\neq\frac12$.}\label{Fig39}
\end{figure}

\smallskip
The constructions in this counterexample show the {\em non-locality} of recovery sequences.
Indeed for $u$ constant, the recovery sequences oscillate between $u-\frac12$ and $u+\frac12$, while if $u$ is a non-constant characteristic function, then the only recovery sequence is essentially the function $u$ itself.

\goodbreak
\subsection{Fractional problems and their singular limits}\label{sec-frac}
Fractional Sobolev spaces may be defined in several ways. A variational definition is through the so-called Gagliardo seminorms
$$
[u]_{s,p}(\Omega)=\Big(\iint_{\Omega\times\Omega} \frac{|u(x)-u(y)|^p}{|x-y|^{d+sp}}dxdy
\Big)^{\frac1p},
$$
defined for $s\in(0,1)$ and $p\ge 1$. The space of functions $u\in L^p(\Omega)$ such that $[u]_{s,p}(\Omega)<+\infty$ is denoted by $W^{s,p}(\Omega)$. Hence, a ``natural'' class of functionals on $W^{s,p}(\Omega)$ consists of double integrals dominated by the corresponding Gagliardo seminorm.

If $\Omega$ is a bounded Lipschitz open set, the space $W^{s,p}(\Omega)$ is compactly embedded in $L^p(\Omega)$. As a consequence the limits of some classes of integrals are easily described. For example, if $a_h=a_h(x,y)$ is a family of coefficients with $0<\alpha\le a_h\le\beta<+\infty$ for all $h$ and weakly* converging to $a_\infty=a_\infty(x,y)$, then we have
\begin{equation}\label{5b}
\Gamma\hbox{-}\lim_{h\to+\infty}\iint_{\Omega\times\Omega}a_h(x,y) \frac{|u(x)-u(y)|^p}{|x-y|^{d+sp}}dxdy= \iint_{\Omega\times\Omega}a_\infty(x,y) \frac{|u(x)-u(y)|^p}{|x-y|^{d+sp}}dxdy.
\end{equation}
This holds in particular for oscillating coefficients $a_h(x,y)=a (hx,hy)$ with $a$ $1$-periodic, for which $a_\infty$ is the average
\begin{equation}\label{5c}
\overline a=\iint_{(0,1)^d\times(0,1)^d}a(x,y)dxdy;
\end{equation} 
that is, {\bf oscillations are averaged out}. This is easily seen by fixing $\delta>0$ and, given $u_h$ converging to $u$, estimating 
\begin{eqnarray*}
&&\liminf_{h\to+\infty}\iint_{\Omega^2}a_h(x,y) \frac{|u_h(x)-u_h(y)|^p}{|x-y|^{d+sp}}dxdy\\
&\ge&\liminf_{h\to+\infty} \iint_{(\Omega\cap B_{1/\delta})^2\cap\{\delta<|x-y|<1/\delta\}}a_h(x,y) \frac{|u_h(x)-u_h(y)|^p}{|x-y|^{d+sp}}dxdy\\
&
=&\iint_{(\Omega\cap B_{1/\delta})^2\cap\{\delta<|x-y|<1/\delta\}}a(x,y) \frac{|u(x)-u(y)|^p}{|x-y|^{d+sp}}dxdy,
\end{eqnarray*}
by weak/strong convergence. Then the lower bound is achieved by letting $\delta\to 0$.
The upper bound is trivial; that is, we can take $u_h=u$.

\bigskip
The limit as $s\to1$ of the $(s,p)$-Gagliardo seminorm is singular, and its behaviour has been highlighted in a seminal paper by Bourgain, Brezis and Mironescu and is driven by {\bf concentration on the diagonal}. We briefly explain this limit process by examining the pointwise behaviour for a function in $C^\infty_c(\Omega)$ and for $p>1$. The first remark is that,
thanks to the decay of the kernel in the integral, and simply estimating $|u(x)-u(y)|^p\le 2^{p-1}(|u(x)|^p+|u(y)|^p)$, we can restrict the interactions to $|x-y|<r$ for any $r>0$ (actually, even $r=r_s\to0$, as long as $r_s^p\gg 1-s$ and $(1-s)|\log r_s|\to0$). Then, using the Taylor expansion, the symmetry of the kernel, and integrating in polar coordinates,  we have
$$
\int_\Omega \Big(\int_{B_r} \frac{|u(x+\xi)-u(x)|^p}{|\xi|^{d+sp}}d\xi\Big)dx\sim
\int_\Omega \Big(\int_{B_r} \frac{|\langle\nabla u(x),\xi\rangle|^p}{|\xi|^{d+sp}}d\xi\Big)dx
$$
$$
=\int_\Omega |\nabla u(x)|^p\Big(\int_{B_r} \frac{|\langle e_1,\xi\rangle|^p}{|\xi|^{d+sp}}d\xi\Big)dx=\int_\Omega |\nabla u(x)|^pdx\int_{S^{d-1}} |\langle e_1,\nu\rangle|^p d\mathcal H^{d-1}\int_0^r\rho^{p-sp-1}d\rho
$$
$$
\sim  \frac1{1-s}C_{d,p}\int_\Omega |\nabla u(x)|^pdx, \qquad C_{d,p}= \frac1{p}\int_{S^{d-1}} |\langle e_1,\nu\rangle|^pd\mathcal H^{d-1}
$$
as $s\to 1$, so that we have the {\em Bourgain--Brezis--Mironescu formula} \cite{BBM01,Pon04}
$$
\lim_{s\to 1} (1-s)\iint_{\Omega\times\Omega} \frac{|u(x)-u(y)|^p}{|x-y|^{d+sp}}dxdy= C_{d,p}\int_\Omega |\nabla u(x)|^pdx
$$

The effect of the limit as $s\to 0$ has been studied by Maz'ya and Shaposhnikova, highlighting {\bf concentration at infinity}. To that end it is convenient to consider $\Omega=\mathbb R^d$, and scale by $s$. Then,  for $p>1$ and $u\in C^\infty_c(\mathbb R^d)$ one notes that the double integral is asymptotically negligible on every bounded set. Since $u=0$ outside any ball $B_R$ for $R$ large, we have 
$$
s\iint_{\mathbb R^{2d}} \frac{|u(x)-u(y)|^p}{|x-y|^{d+sp}}dxdy\sim
2s\int_{\mathbb R^{d}}|u(x)|^pdx \int_{\mathbb R^{d}\setminus B_R} \frac{1}{|\xi|^{d+sp}}d\xi 
$$
$$
=2s\mathcal H^{d-1}(S^{d-1})\int_{\mathbb R^{d}}|u(x)|^pdx \int_R^{+\infty}\rho^{-sp-1}d\rho\sim C'_{d,p}\int_{\mathbb R^{d}}|u(x)|^pdx,
$$
with $C'=\frac{2}{p}\mathcal H^{d-1}(S^{d-1})$, the factor $2$ coming from the equality of the integrals on $\mathbb R^d\times (\mathbb R^{d}\setminus B_R)$ and on $(\mathbb R^{d}\setminus B_R)\times \mathbb R^d$.
The limit
$$
\lim_{s\to 0^+}s\iint_{\mathbb R^{2d}} \frac{|u(x)-u(y)|^p}{|x-y|^{d+sp}}dxdy= C'_{d,p}\int_{\mathbb R^{d}}|u(x)|^pdx
$$
is called the {\em Maz'ya--Shaposhnikova formula} \cite{MS02}.

\bigskip
Both formulas hold as $\Gamma$-limits in the appropriate topologies, so that variational problems with Gagliardo seminorms can be approximated by usual local problems.
In both extreme regimes we have concentration to a local functional. Note however that in the second case, the limit exists also for bounded $\Omega$, but in this case we simply have
$$
\lim_{s\to 0^+}\iint_{\Omega\times\Omega} \frac{|u(x)-u(y)|^p}{|x-y|^{d+sp}}dxdy= \iint_{\Omega\times\Omega} \frac{|u(x)-u(y)|^p}{|x-y|^{d}}dxdy.
$$
The right-hand double integral is finite for example on Lipschitz functions.

\subsection{Scale interactions for double integrals: homogenization}\label{sec-hom-frac}
We now examine the homogenization of Gagliardo seminorms when an oscillating coefficient is coupled with the singular limits $s\to 0$ and $s\to 1$ already examined in the previous section. Let $a\colon\mathbb R^d\times\mathbb R^d\to\mathbb R$ be $1$-periodic in both variables, with $0<\alpha\le a\le\beta$, and consider the oscillating Gagliardo energies
\begin{equation}\label{osc-gagl}
F_{\e,s}(u)=\iint_{\Omega\times\Omega} a\Big(\frac x\e,\frac y\e\Big)\frac{|u(x)-u(y)|^2}{|x-y|^{d+2s}}\,dxdy,
\end{equation}
with $\Omega$ either $\mathbb R^d$ or a bounded Lipschitz subset of $\mathbb R^d$. Two parameters are now present: the period $\e$ of the oscillation and the exponent $s=s(\e)$ governing the singular limit. We only treat the quadratic case $p=2$.

\paragraph{The case of $s$ fixed.} If $s\in(0,1)$ is kept fixed as $\e\to 0$, or $s\to s_o\in (0,1)$, the $\Gamma$-limit with respect to the strong $L^2$-topology is simply given by \eqref{5b} and \eqref{5c}, with no genuine homogenization effect: the contribution of pairs close to the diagonal is negligible for the strong topology, and away from the diagonal the oscillating coefficient is simply replaced by its average $\overline a$ through a Riemann--Lebesgue argument, exactly as in Section~\ref{hom-sec}. The interesting regime, where oscillations and the singular limits $s\to0,1$ genuinely interact, appears only when $s=s(\e)$ tends to an endpoint simultaneously with $\e\to0$.

\paragraph{Three scales of interaction.} As for the homogeneous case of Section~\ref{sec-frac}, write $u(x+\xi)-u(x)$ for the relevant difference quotients and distinguish three ranges for $|\xi|=|x-y|$:

\smallskip
1) {\bf microscopic interactions}, where $|\xi|\ll\e$;

\smallskip
2) {\bf mesoscopic interactions}, where $\e\lesssim|\xi|\lesssim 1$;

\smallskip
3) {\bf macroscopic interactions}, where $|\xi|\gg1$.

\smallskip\noindent
The relative weight of these three contributions, and hence the nature of the limit, is governed by a single scalar parameter comparing $\e$ with $s$ (or $1-s$), exactly as the comparison between $\e$ and $\delta$ governed the periodic arrays of Section~\ref{strangeterms}.1.

\medskip\noindent
{\em The parameter as $s\to0$.} The mass of the (homogeneous) kernel beyond scale $\e$ is
\begin{equation}\label{kmass0}
s\int_{\mathbb R^d\setminus B_\e}\frac{1}{|\xi|^{d+2s}}\,d\xi=\frac{\omega_d}2\,\e^{-2s},
\end{equation}
where $\omega_d=\mathcal H^{d-1}(S^{d-1})$, so the relevant parameter is
$$
\lambda=\lim_{\e\to0}\e^{-2s}\in[1,+\infty].
$$
More precisely, in units of $\frac{\omega_d}2$ the kernel has mass $1$ on the macroscopic range $|\xi|>1$ and mass $\e^{-2s}-1\to\lambda-1$ on the mesoscopic range $\e<|\xi|<1$, while the microscopic range $|\xi|<\e$ gives a vanishing contribution. Hence, when $\lambda<+\infty$, mesoscopic and macroscopic interactions both contribute at the same order; for $\lambda=1$ (that is, $s\ll|\log\e|^{-1}$) the mesoscopic interactions are negligible, while for $\lambda=+\infty$ (that is, $s\gg|\log\e|^{-1}$) they dominate over the macroscopic ones and force the oscillations to average out.

\smallskip\noindent
{\em The parameter as $s\to1$.} Symmetrically, approximating $u(x+\xi)-u(x)\sim\langle\nabla u(x),\xi\rangle$ for $|\xi|$ small, the microscopic kernel mass is
\begin{equation}\label{kmass1}
(1-s)\int_{B_\e}\frac{|\langle\nabla u(x),\xi\rangle|^2}{|\xi|^{d+2s}}\,d\xi=\frac{\omega_d}{2d}|\nabla u(x)|^2\,\e^{2(1-s)},
\end{equation}
so that the relevant parameter is now
$$
\lambda=\lim_{\e\to0}\e^{2(1-s)}\in[0,1],
$$
where, in units of $\frac{\omega_d}{2d}|\nabla u(x)|^2$, the microscopic range $|\xi|<\e$ has mass $\lambda$ and the mesoscopic range $\e<|\xi|<1$ has mass $1-\lambda$, while the macroscopic range gives a vanishing contribution. Hence $\lambda=0$ corresponds to $1-s\gg|\log\e|^{-1}$ (microscopic interactions are negligible), and $\lambda=1$ to $1-s\ll|\log\e|^{-1}$ (they are the only ones that matter).

\paragraph{Three coefficients associated to $a$.} \begin{figure}[h!]
	\centerline{\includegraphics [width=5.8in]{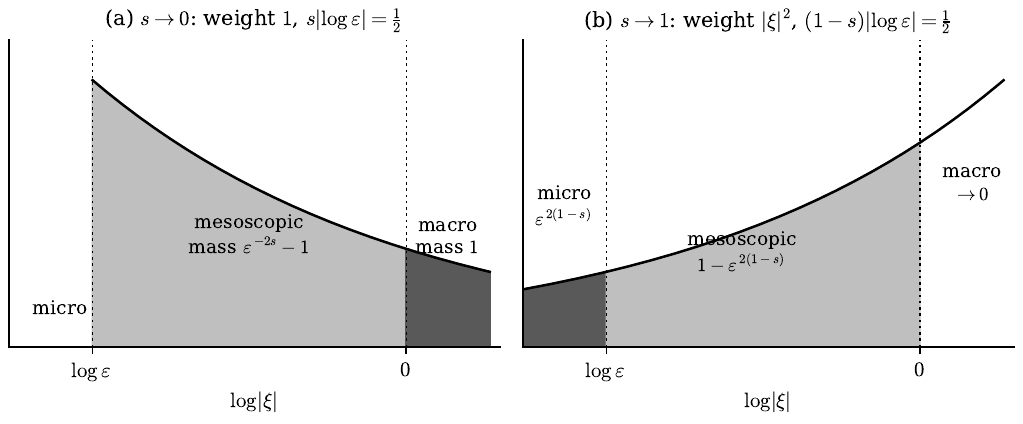}}
	\caption{The mass of the kernel $|\xi|^{-d-2s}$ distributed over the scales $\log|\xi|$, for $\e=e^{-30}$ and $s|\log\e|=\frac12$, resp.~$(1-s)|\log\e|=\frac12$. (a) $s\to0$: the weight is $1$; the macroscopic range (dark) has mass $1$, the mesoscopic one (light) $\e^{-2s}-1$, and the microscopic one is negligible. (b) $s\to1$: with the weight $|\xi|^2$ the microscopic range (dark) has mass $\e^{2(1-s)}$, the mesoscopic one $1-\e^{2(1-s)}$ and the macroscopic one is negligible. Masses are in units of $\frac{\omega_d}{2}$, resp.~$\frac{\omega_d}{2d}|\nabla u|^2$.}\label{Fig40}
\end{figure}

To state the limits we need, besides the average $\overline a$, two further quantities built from $a$ (here $Q=(0,1)^d$):

\smallskip
$\bullet$ the {\bf contracted average} $a^c(x)=\int_Q(a(x,y)+a(y,x))\,dy$, which describes the averaged interaction of the point $x$ with ``points at infinity'' at the macroscopic level, and its {\bf harmonic mean} $\underline a$, defined by $\underline a^{-1}=\int_Q a^c(x)^{-1}dx$ and characterized, as in Section~\ref{sec-L2hom}, by
$$
\Gamma\hbox{-}\lim_{\e\to0}\int_\Omega a^c\Big(\frac x\e\Big)|u|^2dx=\underline a\int_\Omega|u|^2dx\qquad\hbox{(weak $L^2(\Omega)$-convergence);}
$$

\smallskip
$\bullet$ the {\bf diagonal value} $a(x,x)$ (for which we assume $a$ continuous), governing microscopic interactions, together with the {\bf homogenized matrix} $A_{\rm hom}$ obtained, as in Section~\ref{sec-W1phom}, from the coefficient $a(x,x)$:
$$
\Gamma\hbox{-}\lim_{\e\to0}\int_\Omega a\Big(\frac x\e,\frac x\e\Big)|\nabla u|^2dx=\int_\Omega\langle A_{\rm hom}\nabla u,\nabla u\rangle\,dx.
$$

\paragraph{A Maz'ya--Shaposhnikova-type homogenization formula.} As $s\to0$ with $\e^{-2s}\to\lambda\in[1,+\infty]$, define, for $\lambda<+\infty$,
\begin{equation}\label{alam}
a_\lambda=\inf_{\substack{\phi\in L^2(Q)\\ \int_Q\phi=0}}\Big\{\int_Q a^c(x)|1+\phi(x)|^2dx+(\lambda-1)\iint_{Q\times Q}a(x,y)|\phi(x)-\phi(y)|^2dxdy\Big\},
\end{equation}
and $a_{+\infty}=\int_Q a^c(x)\,dx=2\overline a$, with $\overline a$ as in \eqref{5c}. Then
\begin{equation}\label{MSH}
\Gamma\hbox{-}\lim_{\e\to0}s\iint_{\mathbb R^d\times\mathbb R^d}a\Big(\frac x\e,\frac y\e\Big)\frac{|u(x)-u(y)|^2}{|x-y|^{d+2s}}\,dxdy=\frac{\omega_d}2\,a_\lambda\,\|u\|_{L^2(\mathbb R^d)}^2
\end{equation}
with respect to the weak $L^2(\mathbb R^d)$-convergence; this homogenized Maz'ya--Shaposhnikova formula is proved in \cite{BPSP26}. The optimal $\phi$ in \eqref{alam} is a {\bf corrector}: if $A\subset\mathbb R^d$ is a bounded Lipschitz open set and we extend the optimal $\phi$ periodically, the functions $u_\e(x)=\chi_A(x)\big(1+\phi(\frac x\e)\big)$ realize the limit on $\chi_A$. The formula shows that oscillations are felt \emph{both} at the mesoscopic and the macroscopic level, while the microscopic ones are negligible; the term $(\lambda-1)$ in \eqref{alam}, consistent with \eqref{kmass0}, tunes the relative weight of the mesoscopic optimization. For $a\equiv1$ we recover $a_\lambda=2$ and the classical Maz'ya--Shaposhnikova formula $\omega_d\|u\|_{L^2}^2$.

The map $\lambda\mapsto a_\lambda$ is non-decreasing, continuous and concave on $[1,+\infty]$. Its two extreme values correspond to a {\bf separation of scales}: for $\lambda=1$ (that is, $|\log\e|\ll s^{-1}$) mesoscopic optimization is negligible and $a_1=\underline a$, the harmonic mean of $a^c$, so that the limit is obtained by first averaging the macroscopic interactions and only then homogenizing with respect to the weak $L^2$-topology; for $\lambda=+\infty$ (that is, $|\log\e|\gg s^{-1}$) one may instead first let $\e\to0$ at $s>0$ fixed, replacing $a$ by $\overline a$ as in the case of fixed $s$ above, and only then apply the (homogeneous) Maz'ya--Shaposhnikova formula of Section~\ref{sec-frac}, giving $2\overline a$.

\paragraph{A Bourgain--Brezis--Mironescu-type homogenization formula.} As $s\to1$ with $\e^{2(1-s)}\to\lambda\in[0,1]$, define the matrix
\begin{equation}\label{Alam}
A_\lambda=(1-\lambda)\,\overline a\,{\rm Id}+\lambda\,A_{\rm hom}.
\end{equation}
Then, for $\Omega$ Lipschitz and $a$ continuous,
\begin{equation}\label{BBMH}
\Gamma\hbox{-}\lim_{\e\to0}(1-s)\iint_{\Omega\times\Omega}a\Big(\frac x\e,\frac y\e\Big)\frac{|u(x)-u(y)|^2}{|x-y|^{d+2s}}\,dxdy=\frac{\omega_d}{2d}\int_\Omega\langle A_\lambda\nabla u,\nabla u\rangle\,dx,
\end{equation}
equivalently with respect to the weak or the strong $L^2(\Omega)$-convergence. In sharp contrast with the previous case, here the effect of the oscillations is felt \emph{only} at the microscopic level: formula \eqref{Alam} is an exact convex combination of the ``averaged'' energy of Section~\ref{sec-frac} (coefficient $\overline a$, corresponding to no homogenization at all) and the classical homogenized energy with coefficient $a(x,x)$, since mesoscopic interactions see $u$ through its gradient alone and are indifferent to the fine structure of the corrector, by the same Riemann--Lebesgue mechanism of \eqref{5b}. The correctors are precisely those of the classical homogenization problem for $A_{\rm hom}$: if $\phi_\xi$ is $1$-periodic and minimizes $\int_Q a(y,y)|\xi+\nabla\phi|^2dy$, then $u_\e(x)=\langle\xi,x\rangle+\e\,\phi_\xi(\frac x\e)$ realizes the limit on the affine function $u^\xi(x)=\langle\xi,x\rangle$. For $a\equiv1$, \eqref{BBMH} reduces to the classical Bourgain--Brezis--Mironescu formula.

As before, the extreme values of $\lambda$ correspond to separation of scales: $\lambda=0$ ($|\log\e|\gg(1-s)^{-1}$) gives $A_0=\overline a\,{\rm Id}$, obtainable by first letting $\e\to0$ at fixed $s$ (replacing $a$ by $\overline a$) and then applying the homogeneous Bourgain--Brezis--Mironescu formula; $\lambda=1$ ($|\log\e|\ll(1-s)^{-1}$) gives $A_1=A_{\rm hom}$, obtainable by first applying the (homogeneous-in-$\e$) Bourgain--Brezis--Mironescu formula at fixed $\e$, which concentrates on the diagonal and produces the coefficient $a(\frac x\e,\frac x\e)$, and only then homogenizing the resulting local functional as in Section~\ref{sec-W1phom}.

\begin{figure}[h!]
	\centerline{\includegraphics [width=5.8in]{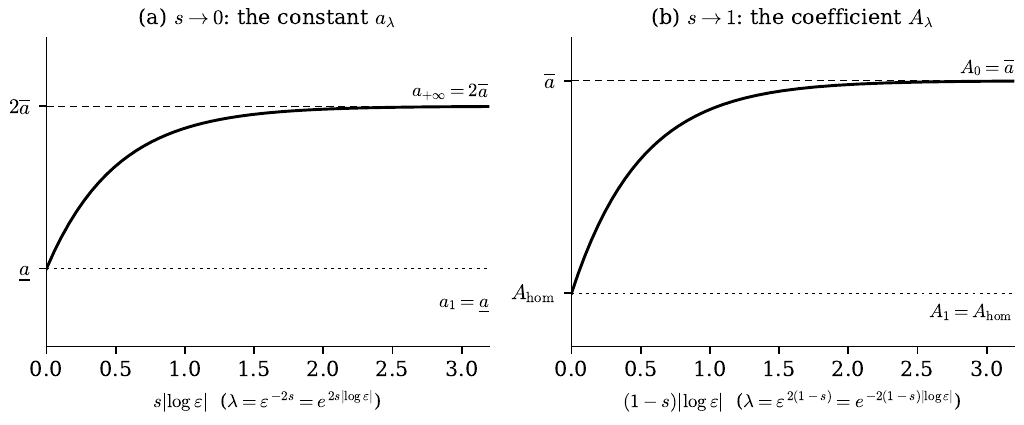}}
	\caption{A numerical example in dimension $d=1$ with $a(x,y)=1+0.4\cos 2\pi x+0.4\cos 2\pi y$ (so that $\overline a=1$, $a^c(x)=2(1+0.4\cos2\pi x)$, $\underline a\approx1.83$ and $A_{\rm hom}=0.6$ from $a(x,x)=1+0.8\cos 2\pi x$). (a) The constant $a_\lambda$ as a function of $s|\log\e|$ ($\lambda=\e^{-2s}$), computed by solving the minimum problem \eqref{alam} on a grid: it is increasing and concave in $\lambda$, with $a_1=\underline a$ and $a_\lambda\to2\overline a$. (b) The coefficient $A_\lambda=(1-\lambda)\overline a+\lambda A_{\rm hom}$ as a function of $(1-s)|\log\e|$ ($\lambda=\e^{2(1-s)}$).}\label{Fig41}
\end{figure}

\paragraph{A unified picture.} The formulas above describe different sides of a single phenomenon, and can be organized into a joint asymptotic description valid for $\e\to0$ and $s\to s_0\in[0,1]$ simultaneously, in the spirit of an \emph{expansion by $\Gamma$-convergence}: the functionals $s(1-s)F_{\e,s}$ (with $F_{\e,s}$ as in \eqref{osc-gagl} and $\Omega=\mathbb R^d$) admit, for every joint regime of $\e\to0$ and $s\to s_0$, a $\Gamma$-limit coinciding with one of
$$
s(1-s)\,\overline a\iint_{\mathbb R^d\times\mathbb R^d}\frac{|u(x)-u(y)|^2}{|x-y|^{d+2s}}dxdy \quad(\hbox{fixed } 0<s_0<1),
$$
$$
a_{\e^{-2s}}\frac{\omega_d}2\int_{\mathbb R^d}|u|^2dx\quad(s\to0),\qquad\qquad \frac{\omega_d}{2d}\int_{\mathbb R^d}\langle A_{\e^{2(1-s)}}\nabla u,\nabla u\rangle dx\quad(s\to1),
$$
the three regimes matching continuously at the boundary values $\lambda\in\{1,+\infty\}$ and $\lambda\in\{0,1\}$ respectively, where $a_\lambda$ and $A_\lambda$ are given by \eqref{alam} and \eqref{Alam}. In this sense the {\bf scale at which oscillations are felt migrates continuously}, as $s$ moves from $0$ to $1$, from the macroscopic and mesoscopic range (Maz'ya--Shaposhnikova regime) to the purely microscopic one (Bourgain--Brezis--Mironescu regime), with the mesoscopic range always present as a transitional mechanism except in the two separated-scales limits.

\subsection{Scale interactions for double integrals: dimension reduction}\label{sec-thin-frac}
We close this chapter by examining the interaction between the non-local scale $s$ and a genuinely geometric scale, that of dimension reduction already considered in Section~\ref{sec-thin} for local (Sobolev) energies. We follow here the analysis of Braides, Pinamonti and Solci \cite{BPS26}.

Let $\omega\subset\mathbb R^{d-1}$ be a bounded Lipschitz open set, and for $\e>0$ consider the thin film $\Omega_\e=\omega\times(0,\e)$ already introduced in Section~\ref{sec-thin}. For fixed $s\in(0,1)$ we consider the squared Gagliardo seminorm on $\Omega_\e$
$$
F_{\e,s}(u)=\iint_{\Omega_\e\times\Omega_\e}\frac{|u(x)-u(y)|^2}{|x-y|^{d+2s}}\,dxdy,
$$
and ask for its asymptotic behaviour, suitably scaled, as $\e\to 0$. As in Section~\ref{sec-thin}, we rescale the thin direction by setting $v_\e(x',x_d)=u_\e(x',\e x_d)$ for $x=(x',x_d)\in\Omega=\omega\times(0,1)$, so that $u_\e$ is identified with a function on the fixed domain $\Omega$. Unlike the local case, however, no Fubini-type decomposition of $F_{\e,s}$ into a horizontal and a vertical part is available, since the Gagliardo seminorm couples points at all mutual distances, including those connecting different values of $x_d$: this is the reason why the analysis of the non-local case is considerably more delicate, and why the scale $s$ interacts non-trivially with the thin-film scale $\e$.

\paragraph{Compactness and the scaling of the energy.}
In terms of $v_\e$ the energy reads
$$
F_{\e,s}(u_\e)=\e^2\iint_{\Omega\times\Omega}\frac{|v_\e(x)-v_\e(y)|^2}{\big(|x'-y'|^2+\e^2(x_d-y_d)^2\big)^{\frac{d+2s}2}}\,dxdy .
$$
A function depending on $x_d$ only pays, for each pair of values $x_d\neq y_d$, a kernel integrated over the planar variable $x'-y'$ equal to a constant times $(\e|x_d-y_d|)^{-1-2s}$, so that its energy is of order $\e^{1-2s}$ (the factor $\e^2$ is compensated by $\e^{-1-2s}$). Since $1-2s$ is smaller than both $2$ and $3-2s$, such oscillations in the thin direction are much more expensive than the scalings $\e^2$, $\e^2|\log\e|$, $\e^{3-2s}$ found below, which are those of functions depending on $x'$ only. This is the mechanism that makes the dimension-reduction scaling coercive: in each of the three regimes below, if $v_\e$ is bounded in $L^2(\Omega)$ and the scaled energies are bounded, then, up to subsequences, $v_\e$ converges strongly in $L^2(\Omega)$ to a function $v=v(x')$ independent of $x_d$, which is the dimension-reduction convergence of Section~\ref{sec-thin} (in its $L^2$ version). The scalings, the form of the limit and the domain of the limit functional depend qualitatively on whether $s$ is smaller, equal to, or larger than $\frac12$.

\paragraph{The subcritical case $s<\frac12$: gain of differentiability.}
For $s<\frac12$ the relevant interactions occur at planar distance of order one: for $|x'-y'|\gg\e$ the integral of the kernel over the thickness is $\e^2|x'-y'|^{-d-2s}$ up to lower-order terms, while pairs with $|x'-y'|\lesssim\e$ give a contribution which is $o(\e^2)$ on functions in $H^{s+\frac12}(\omega)$, since the Gagliardo energy of such functions is concentrated at no particular small scale. In this regime the correct dimension-reduction scaling is $\e^2$: if we set
$$
F^s_\e(u)=\frac1{\e^2}\,F_{\e,s}(u),
$$
then $F^s_\e$ is equi-coercive with respect to dimension-reduction convergence, and if $v_\e\to v=v(x')$ in $L^2(\Omega)$, then
$$
\Gamma\hbox{-}\lim_{\e\to 0}F^s_\e(u_\e)=\iint_{\omega\times\omega}\frac{|v(x')-v(y')|^2}{|x'-y'|^{(d-1)+2(s+\frac12)}}\,dx'dy'=:F^s(v),
$$
with domain $H^{s+\frac12}(\omega)$. The striking feature of this formula is that the exponent in the kernel can be rewritten as $(d-1)+2(s+\frac12)$, since $d+2s=(d-1)+2(s+\frac12)$: hence $F^s$ is the squared Gagliardo seminorm defining $H^{s+\frac12}(\omega)$. In other words, {\bf the dimension-reduction process produces a gain of $\frac12$ in the differentiability exponent}: a thin film of $H^s$-functions converges, after rescaling, to the fractional space $H^{s+\frac12}(\omega)$ on the limit cross-section. This is the non-local counterpart of the classical fact that dimension reduction of $H^1(\Omega_\e)$ produces $H^1(\omega)$ with no change of exponent; here, instead, the vertical fractional oscillations are converted into extra planar regularity.

As $s\to\frac12^-$ the limit functionals $F^s$ themselves undergo a singular-limit transition, since $s+\frac12\to1^-$ and $F^s$ is a Gagliardo seminorm of order $s+\frac12$ on $\omega\subset\mathbb R^{d-1}$: by the Bourgain--Brezis--Mironescu formula in dimension $d-1$, after a further scaling by $(\frac12-s)$ we have
$$
\lim_{s\to\frac12^-}(\tfrac12-s)\,F^s(v)=\frac{\sigma_{d-2}}{2(d-1)}\int_\omega |\nabla' v|^2\,dx',
$$
where $\sigma_{k}=\mathcal H^{k}(S^{k})$, with domain $H^1(\omega)$: the gain in differentiability saturates exactly at the classical exponent $1$ as $s\to\frac12$, anticipating the local behaviour that dominates the critical and supercritical regimes.

\paragraph{The critical case $s=\frac12$: a logarithmic scale.}
At $s=\frac12$ the dimension-reduction scaling itself acquires a logarithmic correction. Setting
$$
F^{1/2}_\e(u)=\frac1{\e^2|\log\e|}\,F_{\e,1/2}(u),
$$
we have dimension-reduction compactness in $L^2(\Omega)$, with limit directly in $H^1(\omega)$ (no intermediate fractional space appears), and
$$
\Gamma\hbox{-}\lim_{\e\to 0}F^{1/2}_\e(u_\e)=\frac{\sigma_{d-2}}{d-1}\int_\omega |\nabla' v|^2\,dx'.
$$
The heuristic explanation of the logarithmic correction is geometric: the interactions that drive the limit are now {\bf mesoscopic}, that is, at every scale $r$ between $\e$ and $1$, and contribute to the energy in comparable amounts at each dyadic scale; integrating over all such scales --- for instance by decomposing the relevant interactions on cylinders of the form $(B'_r\setminus B'_\e)\times(0,\e)$, on which the kernel can be treated as if it only depended on the planar variable $x'$ --- produces the factor $|\log\e|$. This is the same mechanism, read in the non-local setting, that produces logarithmic scales in Section~\ref{strangeterms} for perforated domains with $p=d$ (compare formula~\eqref{mindd}): there too, scale-invariance of the relevant interactions is responsible for the logarithmic correction.

\paragraph{The supercritical case $\frac12<s<1$: competition of microscopic and mesoscopic interactions.}
For $s>\frac12$, a function of $x'$ in $H^1(\omega)$ has a finite energy at the thin-film scale, and the interactions between points at planar distance of order $\e$ (which we call {\bf microscopic}, as they take place within a few multiples of the film's own thickness) are no longer negligible with respect to the mesoscopic ones. Indeed, after Taylor expansion of $v$ the energy is governed by the integral of $|h'|^2$ against the kernel obtained by integrating $|x-y|^{-d-2s}$ over the thickness. This kernel is of order $\e^2|h'|^{-d-2s}$ for $|h'|\gg\e$, and of order $\e|h'|^{-(d-1)-2s}$ for $|h'|\lesssim\e$: both the microscopic range and the mesoscopic range --- whose contribution, since the resulting integral in $|h'|$ is $\int_\e^1r^{-2s}dr\approx\e^{1-2s}/(2s-1)$, is concentrated at distances comparable with $\e$ --- give contributions of order $\e^{3-2s}$, which is therefore the correct dimension-reduction scale in this regime. The precise result is the following: for $\frac12<s<1$ the functionals $\e^{-(3-2s)}F_{\e,s}$ are equi-coercive with respect to dimension-reduction convergence in $L^2(\Omega)$, and their $\Gamma$-limit has domain $H^1(\omega)$, on which
$$
\Gamma\hbox{-}\lim_{\e\to 0}\frac1{\e^{3-2s}}F_{\e,s}(u_\e)=\kappa_{s,d}\int_\omega |\nabla' v|^2\,dx',
$$
where
$$
\kappa_{s,d}=\frac{J_s}{(2s-1)(3-2s)(1-s)},\qquad J_s=\int_{\mathbb R^{d-1}}\frac1{(1+|\xi'|^2)^{\frac d2+s}}\,d\xi'=\frac{\pi^{\frac{d-1}2}\,\Gamma(s+\frac12)}{\Gamma(s+\frac d2)}.
$$
The constant $\kappa_{s,d}$ reflects a kernel in which both the microscopic and the mesoscopic range of interactions contribute, each weighted according to $s$. It correctly interpolates between the two limiting cases already discussed: at the lower end,
$$
\lim_{s\to\frac12^+}(2s-1)\,\kappa_{s,d}=\frac{\sigma_{d-2}}{d-1},
$$
matching the critical formula above, while at the upper end,
$$
\lim_{s\to1^-}(1-s)\,\kappa_{s,d}=\frac{\omega_d}{2d},
$$
matching the Bourgain--Brezis--Mironescu constant on $\Omega\subset\mathbb R^d$ (Section~\ref{sec-frac}) composed with the classical ($p=2$) dimension reduction of Section~\ref{sec-thin}. In this regime the two scales --- the non-local scale $s\to1$ and the geometric scale $\e\to0$ --- merge into a single local limit in $H^1(\omega)$, mediated by a kernel that remembers both.

\begin{figure}[h!]
	\centerline{\includegraphics [width=5.8in]{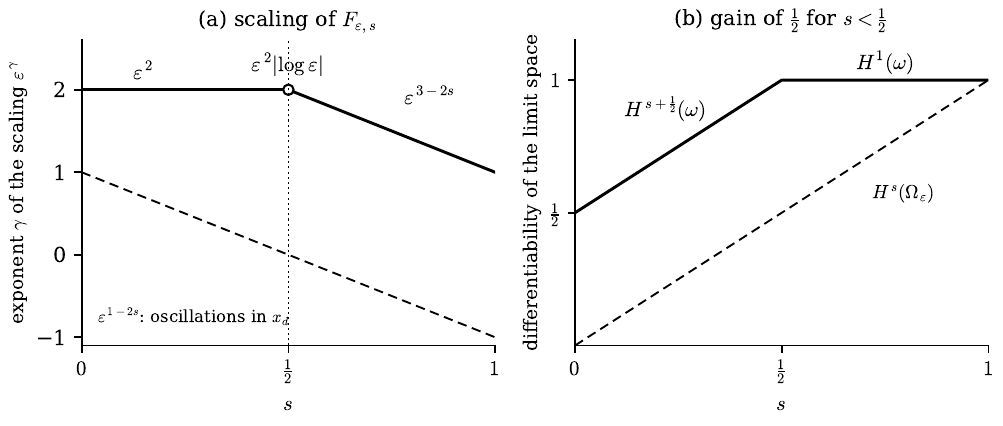}}
	\caption{Dimension reduction of Gagliardo seminorms on thin films. (a) The exponent $\gamma$ of the scaling $\e^{\gamma}$ of $F_{\e,s}$: $\gamma=2$ for $s<\frac12$, $\e^2|\log\e|$ at $s=\frac12$, $\gamma=3-2s$ for $s>\frac12$; the dashed line is the order $\e^{1-2s}$ of the energy of functions oscillating in the thin direction, which is always larger. (b) Differentiability of the limit space on $\omega$: $H^{s+\frac12}(\omega)$ for $s<\frac12$ (a gain of $\frac12$), $H^1(\omega)$ for $s\ge\frac12$; the dashed line is the original exponent $s$.}\label{Fig42}
\end{figure}

\paragraph{The joint limit $\e\to0$, $s\to1$, and separation of scales.}
The results above concern $\e\to0$ at $s$ fixed. It is natural to ask what happens if $s\to1$ and $\e\to0$ simultaneously, so as to connect directly with the Bourgain--Brezis--Mironescu formula of Section~\ref{sec-frac}. Letting $s\to1$ in the supercritical scaling, the correct combined scaling is $(1-s)\e^{2s-3}$: the functionals $(1-s)\e^{-(3-2s)}F_{\e,s}$ are equi-coercive with respect to dimension-reduction convergence in $L^2(\Omega)$, with $\Gamma$-limit of domain $H^1(\omega)$ given by
$$
\Gamma\hbox{-}\lim_{\e\to 0}(1-s)\,\frac{F_{\e,s}(u_\e)}{\e^{3-2s}}=\frac{\omega_d}{2d}\int_\omega |\nabla' v|^2\,dx',
$$
independently of the relative speed at which $s\to1$ and $\e\to0$. The limit functional coincides with the one obtained by first letting $s\to1$ (recovering, by Bourgain--Brezis--Mironescu, the Dirichlet integral $\frac{\omega_d}{2d}\int_{\Omega_\e}|\nabla u|^2dx$, divided by $(1-s)$) and then performing the classical (local) dimension reduction of Section~\ref{sec-thin} on the resulting Dirichlet integral, with the usual membrane scaling $\frac1\e$: the two singular limits ``commute'', in the sense that the joint limit has the same form as the iterated one. The scalings, however, are not the same: since
$$
\frac{1-s}{\e^{3-2s}}=\e^{-2(1-s)}\,\frac{1-s}{\e},
$$
the membrane scaling $(1-s)/\e$ differs from the joint one by the factor $\e^{2(1-s)}$, and so it gives the same limit only in the regime $1-s\ll\frac1{|\log\e|}$, where $\e^{2(1-s)}\to1$. In the critical regime $1-s\sim\frac1{|\log\e|}$, with $\e^{2(1-s)}\to\lambda\in(0,1)$, the membrane-scaled energies converge to $\lambda$ times the same limit, and if $1-s\gg\frac1{|\log\e|}$ they are of lower order than the energy of the joint limit and tend to $0$ on functions of $x'$ only. This is a form of {\bf separation of scales}, already encountered for finite-perimeter concentration in Chapter~\ref{conc-sec} and for the competing effects of Section~\ref{sec-competing}: the two mechanisms --- concentration of non-local interactions as $s\to1$, and geometric collapse of the thin direction as $\e\to0$ --- do not interact at leading order, provided the joint scaling correctly accounts for both, and the parameter $\e^{2(1-s)}$ measures how far the iterated scaling is from it.

\medskip
Summarizing, the analysis of Gagliardo seminorms on thin films exhibits a genuinely three-scale phenomenology: {\bf finite-distance} interactions govern the subcritical regime $s<\frac12$ and produce a gain in differentiability; {\bf mesoscopic} interactions, spread logarithmically over all scales between $\e$ and $1$, govern the critical case $s=\frac12$; and {\bf microscopic} interactions, at the scale $\e$ itself, come into play for $s>\frac12$ (at the same order $\e^{3-2s}$ as the mesoscopic ones) and dominate only in the joint limit $s\to1$, $\e\to0$, always producing in the end a local, classical Dirichlet-type limit on the reduced domain $\omega$.

\subsection{Appendix: relevant tools of Chapter 9}\label{App-Ch9}

\paragraph{References.} The lower semicontinuity and relaxation of double integrals, and the fact that the relaxation of a double integral may fail to be a double integral, are studied by Mora-Corral and Tellini \cite{MCT20} and by Kreisbeck and Zappale \cite{KZ20}, using Young measures; the simple example of Section~\ref{sec-nonrep} is taken from \cite{Bra24}. A comprehensive reference on fractional Sobolev spaces and Gagliardo seminorms is \cite{DNPV12}. The limit as $s\to1$ is due to Bourgain, Brezis and Mironescu \cite{BBM01} (see also Ponce \cite{Pon04} for its formulation as a $\Gamma$-limit), and the limit as $s\to0$ to Maz'ya and Shaposhnikova \cite{MS02}. The interaction of these limits with oscillating coefficients is studied in \cite{BBD23,BMS26} for $s\to1$, where the formula for $A_\lambda$ is proved in the regime $\e^{2(1-s)}\to\lambda$, and in \cite{BPSP26} for $s\to0$, where the formula for $a_\lambda$ is proved. The dimension reduction of Gagliardo seminorms on thin films is treated in \cite{BPS26}.

\paragraph{Relaxation of double integrals.}
\begin{proposition}\label{dbl-lsc}
Let $f\colon\mathbb R\to[0,+\infty]$ be convex and lower semicontinuous. Then $F(u)=\iint_{\Omega\times\Omega}f(u(x)-u(y))\,dxdy$ is lower semicontinuous with respect to the weak convergence in $L^1(\Omega)$.
\end{proposition}

\begin{proof}
Let $u_j\rightharpoonup u$ in $L^1(\Omega)$ and $U_j(x,y)=u_j(x)-u_j(y)$. For every $\varphi\in L^\infty(\Omega\times\Omega)$ we have
$$
\iint U_j\varphi\,dxdy=\int_\Omega u_j(x)\Big(\int_\Omega\varphi(x,y)dy\Big)dx-\int_\Omega u_j(y)\Big(\int_\Omega\varphi(x,y)dx\Big)dy\ \longrightarrow\ \iint U\varphi\,dxdy,
$$
since the two marginals of $\varphi$ belong to $L^\infty(\Omega)$. Hence $U_j\rightharpoonup U$ in $L^1(\Omega\times\Omega)$. Since $U\mapsto\iint f(U)dxdy$ is convex and lower semicontinuous with respect to the strong convergence in $L^1$, it is lower semicontinuous with respect to the weak convergence (Mazur's lemma).
\end{proof}

We now compute the relaxation $\overline F$ in the example of Section~\ref{sec-nonrep}, where $f=0$ on $\{-1,1\}$, $f(0)=1$ and $f=+\infty$ elsewhere, in full generality.

\begin{proposition}\label{dbl-relax}
Let $u\in L^1(\Omega)$ and let $\theta_{\min}=\int_\Omega(u-{\rm ess\,inf}\,u)dx$, $\theta_{\max}=|\Omega|+\int_\Omega(u-{\rm ess\,sup}\,u)dx$. Then $\overline F(u)=+\infty$ if ${\rm osc}\,u>1$, while otherwise $\theta_{\min}\le\theta_{\max}$ and
$$
\overline F(u)=\min\big\{\theta^2+(|\Omega|-\theta)^2:\ \theta\in[\theta_{\min},\theta_{\max}]\big\}.
$$
In particular $\overline F(u)=\frac12|\Omega|^2$ if ${\rm osc}\,u\le\frac12$, and $\overline F(u)=\theta_{\min}^2+(|\Omega|-\theta_{\min})^2$ if ${\rm osc}\,u=1$.
\end{proposition}

\begin{proof}
If $F(u_n)<+\infty$ then $u_n$ takes only two values $c_n\pm\frac12$ almost everywhere, and, setting $\theta_n=|\{u_n=c_n+\frac12\}|$, $F(u_n)=\theta_n^2+(|\Omega|-\theta_n)^2$ and $\int_\Omega u_n=(c_n-\frac12)|\Omega|+\theta_n$. If $u_n\rightharpoonup u$ then $(c_n)$ is bounded; up to subsequences $c_n\to c$ and, since the weak limit of functions with values in $[c_n-\frac12,c_n+\frac12]$ takes values in $[c-\frac12,c+\frac12]$, we have ${\rm ess\,sup}\,u-\frac12\le c\le{\rm ess\,inf}\,u+\frac12$, which is possible only if ${\rm osc}\,u\le1$. Moreover $\theta_n\to\theta=\int_\Omega u-(c-\frac12)|\Omega|$, and when $c$ varies in the interval above, $\theta$ varies in $[\theta_{\min},\theta_{\max}]$. Hence $\liminf F(u_n)\ge\min\{\theta^2+(|\Omega|-\theta)^2:\theta\in[\theta_{\min},\theta_{\max}]\}$. Conversely, given $c$ in the interval, we write $u=c-\frac12+\lambda$ with $\lambda\in[0,1]$, and consider $u_n=c-\frac12+\chi_{E_n}$, where $\chi_{E_n}\rightharpoonup\lambda$ (e.g.~by oscillating sets); then $|E_n|\to\int\lambda=\theta$ and $F(u_n)\to\theta^2+(|\Omega|-\theta)^2$. The last claims follow since $\theta^2+(|\Omega|-\theta)^2$ is minimal for $\theta=\frac12|\Omega|$, and $\theta_{\min}\le\frac12|\Omega|\le\theta_{\max}$ if ${\rm osc}\,u\le\frac12$, while $\theta_{\min}=\theta_{\max}$ if ${\rm osc}\,u=1$.
\end{proof}

\noindent The proof that $\overline F$ is not a double integral given in Section~\ref{sec-nonrep} uses only the values of $\overline F$ at $u=0$ and at the characteristic functions $u_t=\chi_{A_t}$, which are given by this proposition: $\overline F(0)=\frac12|\Omega|^2$ and $\overline F(u_t)=|\Omega|^2(t^2+(1-t)^2)$. Note that the representation $\overline F(u)=\iint g(u(x)-u(y))$ can always be assumed with $g$ even, by symmetrization, so that $\overline F(u_t)=|\Omega|^2\big(g(0)(t^2+(1-t)^2)+2t(1-t)g(1)\big)$.

\begin{proposition}\label{dbl-cont}
Let $f_n$ be continuous (or just Borel) functions with $f_n\uparrow f$, with $f$ as above, and let $F_n(u)=\iint f_n(u(x)-u(y))dxdy$. Then there exists $n_0$ such that, for $n\ge n_0$, the lower-semicontinuous envelope $\overline{F_n}$ is not a double integral $\iint g_n(u(x)-u(y))$.
\end{proposition}

\begin{proof}
Since $F_n\uparrow F$, by the monotone convergence of $\Gamma$-limits we have $\overline{F_n}\uparrow\overline F$ (this is the property used in Section~\ref{sec-nonrep} to transfer the counterexample). Suppose by contradiction that $\overline{F_n}=\iint g_n(u(x)-u(y))$ for infinitely many $n$, with $g_n$ even, and consider only these $n$. Setting $a_n=g_n(0)$ and $b_n=g_n(1)$, evaluation at $u=0$ gives $|\Omega|^2a_n=\overline{F_n}(0)\to\overline F(0)=\frac12|\Omega|^2$, so $a_n\to\frac12$. Evaluation at $u_{1/2}$ gives $\frac12(a_n+b_n)\to\frac12$, so $b_n\to\frac12$. Finally, at $u_{1/4}$ we get
$$
\frac{\overline{F_n}(u_{1/4})}{|\Omega|^2}=\frac58a_n+\frac38b_n\to\frac12,\qquad\hbox{while}\qquad\frac{\overline F(u_{1/4})}{|\Omega|^2}=\frac58,
$$
a contradiction with $\overline{F_n}(u_{1/4})\to\overline F(u_{1/4})$.
\end{proof}

\paragraph{Kernel masses at different scales.}
For $0<r_1<r_2$ and $s\in(0,1)$, integrating in polar coordinates, with $\omega_d=\mathcal H^{d-1}(S^{d-1})$,
$$
\int_{r_1<|\xi|<r_2}\frac{d\xi}{|\xi|^{d+2s}}=\frac{\omega_d}{2s}\big(r_1^{-2s}-r_2^{-2s}\big),\qquad
\int_{r_1<|\xi|<r_2}\frac{|\langle e,\xi\rangle|^2}{|\xi|^{d+2s}}d\xi=\frac{\omega_d}{2d(1-s)}\big(r_2^{2(1-s)}-r_1^{2(1-s)}\big)
$$
for every unit vector $e$ (we have used $\int_{S^{d-1}}\nu_1^2\,d\mathcal H^{d-1}=\omega_d/d$). These are the identities behind the three ranges of interactions of Section~\ref{sec-hom-frac}. As $s\to0$, with the scaling $s$, the masses of the macroscopic range $|\xi|>1$, of the mesoscopic range $\e<|\xi|<1$ and of the microscopic range $|\xi|<\e$ (where we use the weight $|\xi|^2/\e^2$ coming from the regularity of $u$ at scale $\e$) are, in units of $\omega_d/2$,
$$
1,\qquad \e^{-2s}-1\to\lambda-1,\qquad O(s\,\e^{-2s})\to0.
$$
As $s\to1$, with the scaling $1-s$ and weight $\langle\nabla u,\xi\rangle^2$, the masses of the microscopic range, of the mesoscopic range and of the macroscopic range $|\xi|>1$ are, in units of $\omega_d|\nabla u|^2/(2d)$,
$$
\e^{2(1-s)}\to\lambda,\qquad 1-\e^{2(1-s)}\to1-\lambda,\qquad O(1-s)\to0.
$$
This explains why in the first case the oscillations of the coefficient are felt at the mesoscopic and macroscopic level, and in the second one only at the microscopic level, and the role of the quantity $s|\log\e|$ (resp.~$(1-s)|\log\e|$) in distinguishing the regimes: $\e^{-2s}\to1$ if $s|\log\e|\to0$, and $\to+\infty$ if $s|\log\e|\to+\infty$; $\e^{2(1-s)}\to1$ if $(1-s)|\log\e|\to0$, and $\to0$ if $(1-s)|\log\e|\to+\infty$.

\paragraph{Formal derivation of the homogenization formulas.} Let $A$ be a bounded Lipschitz set and $u_\e(x)=\chi_A(x)(1+\phi(x/\e))$, with $\phi$ $1$-periodic with zero mean. In the limit $s\to0$, the double integral $s\iint a(\frac x\e,\frac y\e)|u_\e(x)-u_\e(y)|^2|x-y|^{-d-2s}dxdy$ is the sum of the contribution of the pairs with $x\in A$ and $y\notin A$ (and symmetrically $x\notin A$, $y\in A$), which are macroscopic interactions, and that of the pairs with both $x$ and $y$ in $A$. The first one, since the mass of the kernel is concentrated at infinity and the coefficient there is equidistributed in $y/\e$, gives (in the units above)
$$
\frac{\omega_d}2\int_Aa^c\Big(\frac x\e\Big)\Big|1+\phi\Big(\frac x\e\Big)\Big|^2dx\ \longrightarrow\ \frac{\omega_d}2|A|\int_Qa^c(x)|1+\phi(x)|^2dx,
$$
where we have used that the weight of $x\in A$ is $\int(a(\frac x\e,\frac y\e)+a(\frac y\e,\frac x\e))dy$ averaged over $y$, that is $a^c(x/\e)$. The second one, using the mesoscopic range $\e<|x-y|<1$ where $(x/\e,y/\e)$ is equidistributed in $Q\times Q$, gives $\frac{\omega_d}2(\lambda-1)|A|\iint_{Q\times Q}a(x,y)|\phi(x)-\phi(y)|^2dxdy$. Summing the two terms and optimizing in $\phi$ we get $\frac{\omega_d}2a_\lambda|A|$, with $a_\lambda$ as in \eqref{alam}. The formula for $A_\lambda$ is obtained in the same way: in the limit $s\to1$ for $u_\e(x)=\langle\xi,x\rangle+\e\phi(x/\e)$, the microscopic interactions $|x-y|\ll\e$, of mass $\lambda$, see $a(\frac x\e,\frac x\e)$ and the gradient $\xi+\nabla\phi(x/\e)$, giving $\lambda\langle A_{\rm hom}\xi,\xi\rangle$ after optimization in $\phi$, while the mesoscopic interactions, of mass $1-\lambda$, see only the average $\overline a$ of the coefficient and the gradient $\xi$ of the macroscopic part of $u_\e$, giving $(1-\lambda)\overline a|\xi|^2$.

\begin{proposition}\label{alam-prop}
The function $\lambda\mapsto a_\lambda$ is non-decreasing and concave on $[1,+\infty)$, $a_1=\underline a$ and $\lim_{\lambda\to+\infty}a_\lambda=\int_Qa^c=2\overline a$.
\end{proposition}

\begin{proof}
For each $\phi$ the function inside the infimum in \eqref{alam} is affine and non-decreasing in $\lambda$, so $a_\lambda$ is concave and non-decreasing. For $\lambda=1$, setting $v=1+\phi$, we have to minimize $\int_Qa^cv^2$ under the constraint $\int_Qv=1$. By the Cauchy--Schwarz inequality
$1=(\int_Qv)^2\le\int_Qa^cv^2\int_Q(a^c)^{-1}$, so that $a_1\ge\underline a$, with equality for $v=\underline a/a^c$ (which has $\int_Qv=1$). Finally, $a_\lambda\le\int_Qa^c$ for all $\lambda$ (take $\phi=0$) and, if $\phi_\lambda$ is a competitor with energy at most $\int_Qa^c+1$, since $a\ge\alpha$ and since $\iint_{Q\times Q}|\phi(x)-\phi(y)|^2\,dxdy=2\|\phi\|^2_{L^2(Q)}$ for $\phi$ of zero mean, we get $2\alpha(\lambda-1)\|\phi_\lambda\|_{L^2}^2\le C$, so $\phi_\lambda\to0$ in $L^2(Q)$ and $\liminf a_\lambda\ge\int_Qa^c$.
\end{proof}

\paragraph{The dimension reduction of Gagliardo seminorms: formal computation.}
Let $v\in C^2(\overline\omega)$ and $u(x)=v(x')$. Then, integrating in the thin direction,
$$
F_{\e,s}(u)=\iint_{\omega\times\omega}|v(x')-v(y')|^2k_\e(x'-y')\,dx'dy',\qquad k_\e(z')=\int_0^\e\!\!\int_0^\e\frac{da\,db}{(|z'|^2+(a-b)^2)^{\frac{d+2s}2}}.
$$
If $|z'|\gg\e$ then $k_\e(z')\sim\e^2|z'|^{-d-2s}$, and the pairs with $|z'|\gg\e$ give, for $s<\frac12$, $\e^2\iint_{\omega\times\omega}|v(x')-v(y')|^2|x'-y'|^{-d-2s}$, since the contribution of the pairs with $|z'|\lesssim\e$, which is of order $\e^{3-2s}$ by the computation below, is of lower order. This gives the scale $\e^2$ and the limit $F^s(v)$ in Section~\ref{sec-thin-frac}, whose kernel exponent is $d+2s=(d-1)+2(s+\frac12)$. For the pairs with $|z'|\lesssim\e$ we use $|v(x')-v(y')|^2\sim\langle\nabla'v(x'),z'\rangle^2$ and integrate first in $z'\in\mathbb R^{d-1}$ at a fixed value of $h=|a-b|$, with the change of variable $z'=h\zeta$:
$$
\int_{\mathbb R^{d-1}}\frac{\langle\nabla'v,z'\rangle^2}{(|z'|^2+h^2)^{\frac{d+2s}2}}dz'=\frac{|\nabla'v|^2}{d-1}\,h^{1-2s}\,I_s,\qquad I_s=\int_{\mathbb R^{d-1}}\frac{|\zeta|^2}{(1+|\zeta|^2)^{\frac d2+s}}d\zeta,
$$
where $I_s<+\infty$ if and only if $s>\frac12$. Since $\int_0^\e\!\int_0^\e|a-b|^{1-2s}da\,db=\e^{3-2s}/((1-s)(3-2s))$, we obtain the scale $\e^{3-2s}$ and the limit
$$
\frac{F_{\e,s}(u)}{\e^{3-2s}}\ \to\ \frac{I_s}{(d-1)(1-s)(3-2s)}\int_\omega|\nabla'v|^2dx'.
$$
This coincides with the constant $\kappa_{s,d}$ of Section~\ref{sec-thin-frac}. Indeed, set $b=\frac d2+s$ and $J_s=\int_{\mathbb R^{d-1}}(1+|\zeta|^2)^{-b}d\zeta$. Integrating by parts, using $\zeta(1+|\zeta|^2)^{-b}=-\frac1{2(b-1)}\nabla(1+|\zeta|^2)^{-(b-1)}$ (the boundary terms vanish since $s>\frac12$) we get
$$
I_s=\frac{d-1}{d+2s-2}\int_{\mathbb R^{d-1}}\frac{d\zeta}{(1+|\zeta|^2)^{b-1}}=\frac{d-1}{d+2s-2}\,(J_s+I_s),
$$
since $(1+|\zeta|^2)^{-(b-1)}=(1+|\zeta|^2)^{-b}+|\zeta|^2(1+|\zeta|^2)^{-b}$. Solving for $I_s$ gives $I_s=\frac{d-1}{2s-1}J_s$, and, by the standard formula $\int_{\mathbb R^n}(1+|\zeta|^2)^{-b}d\zeta=\pi^{n/2}\Gamma(b-\frac n2)/\Gamma(b)$,
$$
J_s=\pi^{\frac{d-1}2}\frac{\Gamma(s+\frac12)}{\Gamma(s+\frac d2)}.
$$
Hence
$$
\kappa_{s,d}=\frac{J_s}{(2s-1)(3-2s)(1-s)},
$$
and
$$
\lim_{s\to1^-}(1-s)\kappa_{s,d}=J_1=\frac{\pi^{d/2}}{2\Gamma(1+\frac d2)}=\frac{\omega_d}{2d},\qquad \lim_{s\to\frac12^+}(2s-1)\kappa_{s,d}=J_{\frac12}=\frac{\pi^{\frac{d-1}2}}{\Gamma(\frac{d+1}2)}=\frac{\sigma_{d-2}}{d-1}.
$$
At the critical exponent the integral in $z'$ diverges logarithmically: the same computation, for $\e<|z'|<1$ and with $k_\e(z')\sim\e^2|z'|^{-d-1}$, gives $\e^2\int_\omega\int_{\e<|z'|<1}\frac{\langle\nabla'v,z'\rangle^2}{|z'|^{d+1}}dz'dx'=\e^2|\log\e|\,\frac{\sigma_{d-2}}{d-1}\int_\omega|\nabla'v|^2dx'$, which is the critical scale $\e^2|\log\e|$ of Section~\ref{sec-thin-frac}. Finally, for $s<\frac12$ and $\sigma=s+\frac12\to1^-$, the Bourgain--Brezis--Mironescu formula in dimension $d-1$ applied to $F^s$ (with the exponent $(d-1)+2\sigma$ and $1-\sigma=\frac12-s$) gives the limit $\frac{\sigma_{d-2}}{2(d-1)}\int_\omega|\nabla'v|^2dx'$ stated there.

We note that the optimal recovery sequence is $u(x)=v(x')$ without a corrector in the thin direction. Indeed, for the quadratic form obtained by adding to $\langle\nabla'v(x'),y'-x'\rangle$ a perturbation $w(y)-w(x)$ with $w$ defined on $\mathbb R^{d-1}\times(0,\e)$, the Euler--Lagrange equation at $w=0$ is
$$
\int\langle\nabla'v(x'),y'-x'\rangle\,\frac{dy}{|x-y|^{d+2s}}=0,
$$
which is satisfied since the kernel is even in $y'-x'$; by the convexity of the quadratic form, $w=0$ is the minimizer. The compactness of sequences with bounded scaled energy in the sense of dimension reduction and the lower bound are proved in \cite{BPS26}.

\paragraph{The joint limit $s\to1$, $\e\to0$ and the membrane scaling.}
By the formal computation above, for smooth $v$, the energies $(1-s)\e^{-(3-2s)}F_{\e,s}(u)$ converge to $(1-s)\kappa_{s,d}\int_\omega|\nabla'v|^2$, and $(1-s)\kappa_{s,d}\to\omega_d/(2d)$ as $s\to1$: the limit does not depend on the relative speed of $s\to1$ and $\e\to0$ because the integral $I_s$ is uniformly bounded as $s\to1$ and the contribution of the pairs at distance of order $1$ is of order $\e^2\ll\e^{3-2s}/(1-s)$. The rigorous proof of this statement, including the lower bound, is in \cite{BPS26}. The membrane scaling $(1-s)/\e$, obtained by first sending $s\to1$ at fixed $\e$ (so that $(1-s)F_{\e,s}\to\frac{\omega_d}{2d}\int_{\Omega_\e}|\nabla u|^2$) and then $\e\to0$ in the local dimension reduction of Section~\ref{sec-thin}, differs from the joint scaling by the factor $\e^{3-2s}/\e=\e^{2(1-s)}$:
$$
\frac{(1-s)}{\e}F_{\e,s}(u)=\e^{2(1-s)}\cdot\frac{(1-s)}{\e^{3-2s}}F_{\e,s}(u).
$$
Hence the two scalings give the same limit if $(1-s)|\log\e|\to0$ (the factor tends to $1$); if $\e^{2(1-s)}\to\lambda\in(0,1)$, that is, if $1-s\sim|\log\e|^{-1}$, the membrane-scaled energies converge to $\lambda$ times the limit of the joint scaling, and, if $(1-s)|\log\e|\to+\infty$, to $0$. This is the critical regime, in which the mesoscopic interactions between the scales $\e$ and $1$, which are not seen when $s\to1$ first, become relevant.

\vfill\eject
%
%
%
%
%
%
%
%
%
%
%
%

\section{Discrete-to-continuum homogenization}\label{disc-sec}
We now consider functionals whose domain consists of functions defined on lattices depending on a small scaling parameter, as for instance functions defined on $\e\mathbb Z^d$. In this sense, we can say that the scale is {\em a priori} given by a geometrical constraint. The analysis as $\e$ goes to $0$ of these energies is a particular kind of homogenization, 
giving a suitable continuum energy which is referred to as the {\bf discrete-to-continuum limit}. 

The first problem to face is that of giving a meaning to the convergence of discrete functions to a continuous limit. This can be done by identifying discrete functions with particular interpolations, depending on the structure and the scaling of the energies.  

\subsection{Bulk scaling}
Since our focus is on scaling effects, and not on the structure of the limit energies depending on the structure of the original lattice, we only consider functions $u$ defined on the lattice $\e \mathbb Z^d$. The discrete analog of integral functionals in the simplest case and in dimension $1$ is then an energy of the form 
\begin{equation}\label{bulkdisc1}
F_\e(u)=\e \sum_{i\in I_\e(\Omega)} f\Big(\frac{u_i-u_{i-1}}{\e}\Big), 
\end{equation}
where $f\colon \mathbb R\to [0,+\infty)$ is continuous with $|z|^p\le f(z)\le c(1+|z|^p)$ for some $p>1$ and $c>0$, 
$\Omega$ is an open subset of $\mathbb R$, and  
$u\colon\Omega\cap\e\mathbb Z\to \mathbb R$. The sum is taken over the set 
$I_\e(\Omega)=\{i\in \mathbb Z: \e i, \e (i-1)\in \Omega\}$, 
and the notation $u_i$ denotes the value $u(\e i)$. 
This can be viewed as a finite-difference counterpart of the functional 
$\int_{\Omega}f(u^\prime)\, dt$. Indeed, 
the scalings (of the variable and of the energy) in \eqref{bulkdisc1} are designed in such a way that,  
identifying $u$ with its piecewise-affine interpolation, we have  
$$\int_{\Omega}f(u^\prime)\, dt=\sum_{i\in I_\e(\Omega)}\e \,f\Big(\frac{u_i-u_{i-1}}{\e}\Big)+o(1)=F_\e(u)+o(1)$$ 
as $\e$ goes to $0$. The error term can also be set to $0$ by suitably extending the interpolation at the endpoints of the intervals of $\Omega$. Note however that, even though we identify lattice functions with a subset of $L^1(\Omega)$, for each $\e>0$ the functional $F_\e$ is defined on a finite-dimensional space, and hence it is a continuous functional if $f$ is continuous, regardless of any convexity condition. 

The lattice size $\e$ introduces a {\bf constraint on the scale of oscillations} that is overcome only in the limit as $\e$ tends to $0$. This is clear if we remark that
$$
\Gamma\hbox{-}\lim_{\e\to 0} F_\e(u)= \int_\Omega f^{**}(u')\,dt
$$
for $u\in W^{1,p}(\Omega)$ (and $+\infty$ otherwise), with respect to the convergence in $L^1(\Omega)$ of the piecewise-affine interpolations; that is, convexification appears only as a limit process. The proof of this fact is easily obtained by observing that 
$F_\e(u)\ge \int_\Omega f^{**}(u')\,dt$ for all $\e$ and $u$ (by Jensen's inequality, since $f\ge f^{**}$ and the interpolation of $u$ has derivative $\frac{u_i-u_{i-1}}\e$ on each interval). Since $f(z)\ge|z|^p$, a sequence with bounded energy is bounded in $W^{1,p}(\Omega)$, and $u\mapsto\int_\Omega f^{**}(u')\,dt$ is lower semicontinuous along such sequences, as $f^{**}$ is convex; hence the inequality holds also for the $\Gamma$-liminf.
Conversely, if $u\in C^2(\overline\Omega)$ then $\lim_{\e\to 0} F_\e(u)=
\int_\Omega f(u')\,dt$, so that $\Gamma$-$\limsup_{\e\to 0} F_\e(u)\le \int_\Omega f(u')\,dt$, from which the inequality 
$\Gamma$-$\limsup_{\e\to 0} F_\e(u)\le \int_\Omega f^{**}(u')\,dt$ follows by the relaxation argument of Section~\ref{Pre-Relax} (and the growth condition on $f$). 
As a consequence, we have an {\bf asymptotic formula for the convex envelope} of a continuous $f$ as
\begin{eqnarray*}
f^{**}(z)&=&\lim_{N\to+\infty}{1\over N}\min\bigg\{\sum_{i=1}^N f(u_i-u_{i-1}): u_0=0, u_N= Nz\bigg\}\\
&=&\lim_{N\to+\infty}{1\over N}\min\bigg\{\sum_{i=1}^N f(z_i): {1\over N}\sum_{i=1}^N z_i= z\bigg\},
\end{eqnarray*}
where the limit exists, and is equal to the infimum over $N$, by the subadditivity of $N\mapsto\min\{\sum_{i=1}^Nf(z_i):\sum_{i=1}^Nz_i=Nz\}$. In general it cannot be replaced by a minimum problem with a fixed number of variables: the convex envelope is obtained by convex combinations $\sum_j\lambda_jf(w_j)$ of at most two values, and these can be realized by averages of $N$ values only in the limit when the weights $\lambda_j$ are irrational. For instance, for $f(z)=(z^2-1)^2$ and $z$ irrational in $(-1,1)$ we have $f^{**}(z)=0$, but the minimum is strictly positive for every $N$, since the average of $N$ values in $\{-1,1\}$ is rational.

\begin{figure}[h!]
	\centerline{\includegraphics [width=5.8in]{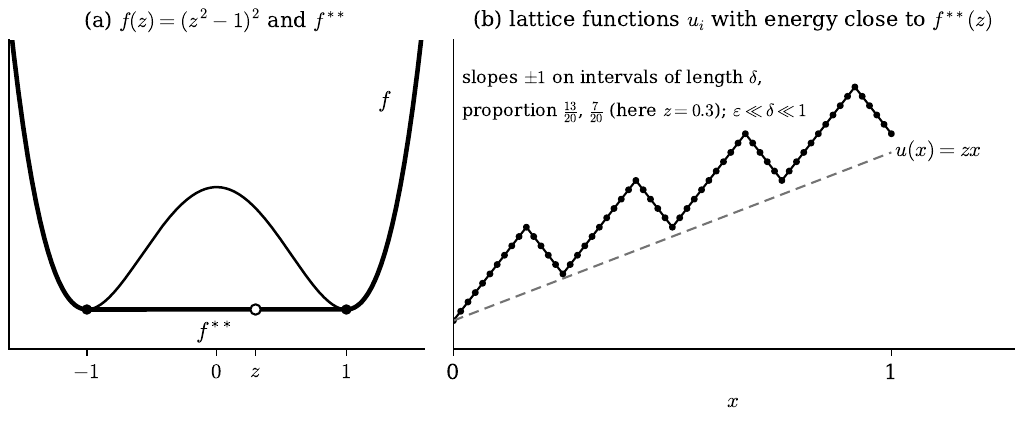}}
	\caption{Convexification at the lattice level. (a) A double-well potential $f(z)=(z^2-1)^2$ and its convex envelope $f^{**}$, which vanishes on $[-1,1]$. (b) A lattice function $u_i$ (dots) whose increments are $\pm\e$, that is, with slopes $\pm1$ at the lattice scale $\e$, in proportion $\frac{13}{20}$ and $\frac7{20}$ on intervals of length $\delta$, with $\e\ll\delta\ll1$: it converges to the affine function $u(x)=zx$ with $z=0.3$, which is not a well of $f$, and its energy tends to $f^{**}(z)=0<f(z)$.}\label{Fig44}
\end{figure}

\bigskip
The discrete setting allows one to consider functionals that have a structure more general than simple discretizations of continuum functionals, by taking into account {\em long-range interactions}; that is, in the notation introduced, taking into account terms depending on differences $u_i-u_j$ also when $|i-j|>1$. 
Again, in the one-dimensional setting we can consider
$$
F_\e(u)=\sum_{k\ge 1}\sum_{i\in I^k_\e(\Omega)}\e\, f_k\Big(\frac{u_i-u_{i-k}}{\e}\Big),
$$
with each $f_k$ continuous and non-negative, and $f_1(z)\ge c|z|^p-C$ for some $p>1$, $c>0$, $C\ge0$,
where now $I^k_\e(\Omega)=\{i\in \mathbb Z: \e i, \e (i-k)\in \Omega\}$. 
Formally, this looks like a discretization of $\int_\Omega \sum_k f_k(k u')dt$.
Indeed, this is the pointwise limit of $F_\e$, after identifying a function $u\in W^{1,p}(\Omega)$ with the function $u\colon\Omega\cap \e\mathbb Z\to\mathbb R$ given by the continuous representative of $u$, but in general 
the corresponding relaxation, with integrand the convex envelope of the function $f_0(z):= \sum_k f_k(k z)$, is only an upper bound. 

If we suppose that $\sum_k k\, f_k(k z)\le c(1+|z|^p)$ for all $z$ (the factor $k$ takes into account that $k$ pairs $(i,i-k)$ interact across any given point; this condition controls the effect of the boundary, and of the interactions at long range), we can repeat the arguments of the homogenization of Chapter~\ref{hom-sec} (see \cite{BGe02,AC04}), and deduce that, for $u\in W^{1,p}(\Omega)$,
$$
\Gamma\hbox{-}\lim_{\e\to 0} F_\e(u)= \int_\Omega f_{\rm hom}(u')\,dt,
$$
where
\begin{eqnarray*}
f_{\rm hom}(z)=\lim_{N\to+\infty}{1\over N}\min\bigg\{\sum_{k=1}^N\sum_{i=k}^N f_k(u_i-u_{i-k}): u_0=0, u_N= Nz\bigg\}.
\end{eqnarray*}
For this function we have the inequalities
$$
\sum_{k=1}^\infty f_k^{**}(kz)\le f_{\rm hom}(z)\le f^{**}_0(z),
$$
with $f_0$ giving the pointwise limit, while the first inequality is obtained by separately optimizing oscillations for each $f_k$. 

\bigskip
We now show that indeed we cannot reduce to a single cell-problem formula in a simple but meaningful example, where
$f_k=0$ except for $k=1$ and $k=M>1$, generalizing the case when only $f_1=f\neq 0$. We take $\Omega=[0,1]$; this is not an open set, but the only change is that we also consider $i=0$ and possibly $i=1/\e$ if it is an integer, so that our functionals are
$$
F_\e(u)=\sum_{i=1}^{N_\e}\e f_1\Bigl({u_i-u_{i-1}\over\e}\Big)+ \sum_{i=M}^{N_\e}\e f_M\Bigl({u_i-u_{i-M}\over\e}\Big),\qquad\hbox{ where $N_\e=\lfloor {1\over\e}\rfloor$.}
$$
In this case we optimize nearest-neighbour interactions at fixed $M$-neighbour interactions. To check a lower bound we start by writing
\begin{eqnarray*}
F_\e(u)&\ge&{1\over M}\sum_{k=1}^M
\sum_{j=1}^{\lfloor{N_\e\over M}\rfloor-1}\e M\biggl(f_M\Bigl({u_{k+Mj}-u_{k+M(j-1)}\over\e}\Big)\\
&&\hskip3cm+\sum_{i=k+M(j-1)+1}^{k+Mj}{1\over M} f_1\Bigl({u_i-u_{i-1}\over\e}\Big)\bigg);
\end{eqnarray*}
that is, we subdivide the $M$-neighbour interactions into $M$ nonintersecting sublattices $\e k+\e M\mathbb Z$ with $k\in\{1,\ldots, M\}$, and consider each nearest-neighbour interaction as an intermediate interaction between two points in one of these sublattices. Since in this way each nearest-neighbour interaction is accounted for $M$ times, we have to divide it by $M$. In order to avoid the boundary, the sum is taken slightly inside $[0,1]$.
At this point we define the {\em effective $M$-neighbour potential}
$$
f_{\rm eff} (z)=f_M ( M z)+{1\over M}\min\bigg\{\sum_{i=1}^M f_1(z_i): \sum_{i=1}^M z_i= Mz
\bigg\},
$$
and we have the estimate
\begin{eqnarray*}
F_\e(u)\ge{1\over M}\sum_{k=1}^M
\sum_{j=1}^{\lfloor{N_\e\over M}\rfloor-1}\e M\, f_{\rm eff}\Bigl({u_{k+Mj}-u_{k+M(j-1)}\over M\e}\Big).\end{eqnarray*}
The right-hand side of this inequality is a sum of nearest-neighbour energies on disjoint lattices.
Now, if we have a sequence $u_\e$ with bounded energy converging to $u$ we note that the interpolations $u^k_\e$ on each sublattice $\e k+\e M\mathbb Z$ also converge to the same $u$ (indeed, by Hölder's inequality and $f_1(z)\ge c|z|^p-C$, $\|u^k_\e-u_\e\|_{L^p}\le CM\e\|u_\e'\|_{L^p}\le C'\e$). By the result on nearest-neighbour energies, we then have
$$
\liminf_{\e\to0} F_\e(u_\e)\ge {1\over M}\sum_{k=1}^M \int_{ (\eta, 1-\eta)} f_{\rm eff}^{**} (u')dt=\int_{ (\eta, 1-\eta)} f_{\rm eff}^{**} (u')dt
$$
for all $\eta>0$, and we can eventually let $\eta\to 0$ to obtain the lower bound.
The lower bound looks optimistic, since it relies on separate estimates for each $k$ that might need incompatible constructions for the recovery sequence, but this is not the case, as minimizers for $f_{\rm eff} (z)$ can be regarded as periodic functions. 

To check the upper bound it suffices to check that $\Gamma$-$\limsup_{\e\to 0} F_\e(u)\le \int_{(0,1)} f_{\rm eff} (u')dt$; 
we do it for the case $u(x)=zx$, which exemplifies the construction at the discrete level. With fixed $z$ let $(z_1,\ldots, z_M)$ be a minimizer for $f_{\rm eff}(z)$, and consider its $M$-periodic extension to $\mathbb Z$. Then we consider the functions $\overline u_\e\colon\mathbb Z\to \mathbb R$ defined by $\overline u_\e(0)=0$ and $\overline u_\e(\e i)=\overline u_\e(\e (i-1))+\e z_i$. We have $\overline u_\e\to u$ and 
$$
\limsup_{\e\to 0} F_\e(\overline u_\e)= f_{\rm eff }(z)=\int_{(0,1)} f_{\rm eff} (u')dt,
$$
as desired. The convexification is obtained in the same way: if $f^{**}_{\rm eff}(z)=\theta f_{\rm eff}(z')+(1-\theta)f_{\rm eff}(z'')$, for $u(x)=zx$ we alternate on mesoscopic intervals the periodic constructions with slopes $z'$ and $z''$ in proportions $\theta$ and $1-\theta$, and a general $u$ is treated by density of piecewise-affine functions. In conclusion 
$$
\Gamma\hbox{-}\lim_{\e\to 0} F_\e(u)= \int_{(0,1)} f^{**}_{\rm eff} (u')dt.
$$
The form of the integrand highlights the superposition of two phenomena: {\bf optimization of lattice interactions on the period $M$ at scale $\e$ (microscopic scale)}, summarized in the formula for $f_{\rm eff}$, and {\bf optimization of oscillations at scale much larger than $\e$ (mesoscopic scale)}, producing the convexification.

\begin{figure}[H]
	\centerline{\includegraphics [width=4.9in]{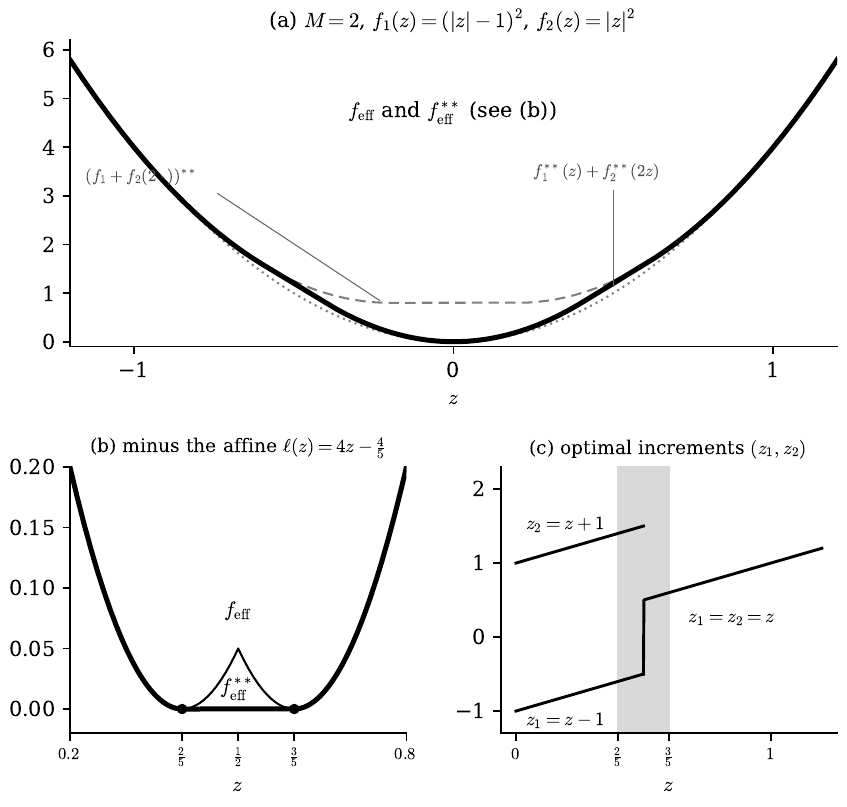}}
	\caption{The example with $M=2$, $f_1(z)=(|z|-1)^2$, $f_2(z)=|z|^2$. (a) The effective potential $f_{\rm eff}$ and its convex envelope $f_{\rm eff}^{**}$. (b) The same functions after subtraction of the affine function $4z-\frac45$. (c) The optimal increments $(z_1,z_2)$ of the microscopic pattern as functions of $z$.}\label{Fig43}
\end{figure}

\medskip\noindent
{\bf An explicit example.} We take $M=2$, $f_1(z)=(|z|-1)^2$ and $f_2(z)=|z|^2$. The nearest-neighbour energy has two wells at $\pm1$, and
$$
f_{\rm eff}(z)=4z^2+\frac12\min\{f_1(a)+f_1(b):a+b=2z\}.
$$
Let $z\ge0$. If $a,b\ge0$ the minimum is attained at $a=b=z$ by the convexity of $f_1$ on $[0,+\infty)$, with value $2(1-z)^2$; if $a<0\le b$ the best choice for $z\le1$ is $(a,b)=(z-1,z+1)$, where both increments are as close as possible to the wells, with value $2z^2$. For $|z|\le\frac12$ the second ``dimerized'' pattern is better and for $|z|\ge\frac12$ the first, ``uniform'', one is, so that
$$
f_{\rm eff}(z)=\begin{cases}5z^2&\hbox{if }|z|\le\frac12\\ 5z^2-2|z|+1&\hbox{if }|z|\ge\frac12.\end{cases}
$$
This function is not convex, since the microscopic pattern changes abruptly at $|z|=\frac12$, producing a concave corner. Its convex envelope is obtained by replacing $f_{\rm eff}$ in the interval $[\frac25,\frac35]$ by the common tangent at the two points $\frac25$ and $\frac35$,
$$
f_{\rm eff}^{**}(z)=\begin{cases}5z^2&\hbox{if }|z|\le\frac25\\ 4|z|-\frac45&\hbox{if }\frac25\le|z|\le\frac35\\ 5z^2-2|z|+1&\hbox{if }|z|\ge\frac35,\end{cases}
$$
see Fig.~\ref{Fig43}: for $z\in(\frac25,\frac35)$ the optimal sequences oscillate, at a mesoscopic scale, between the macroscopic slopes $\frac25$ and $\frac35$, and at the microscopic scale each of these slopes is realized by the optimal $2$-periodic pattern. The functions $f_1^{**}(z)+f_2^{**}(2z)$ and $(f_1+f_2(2\,\cdot\,))^{**}$ bound $f_{\rm eff}^{**}$ from below and above, and both inequalities are strict (for instance at $z=\frac12$, where $f_{\rm eff}^{**}=\frac65$ is between $1$ and $\frac54$).

In Fig.~\ref{Fig43}(a) we picture the effective potential $f_{\rm eff}$ and its convex envelope $f_{\rm eff}^{**}$, which are hardly distinguishable at this scale, together with the lower bound $f_1^{**}(z)+f_2^{**}(2z)=4z^2$ (dotted, for $|z|\le1$) and the upper bound $(f_1+f_2(2\,\cdot\,))^{**}$ (dashed), which is flat on $[-\frac15,\frac15]$. In Fig.~\ref{Fig43}(b), the same functions after subtraction of the affine function $4z-\frac45$ tangent to $f_{\rm eff}$ at $\frac25$ and $\frac35$: $f_{\rm eff}$ has two wells, at $\frac25$ and $\frac35$, separated by a bump of height $\frac1{20}$ at $\frac12$, and $f_{\rm eff}^{**}$ (thick) vanishes on the whole interval $[\frac25,\frac35]$. In Fig.~\ref{Fig43}(c) we picture the optimal increments $(z_1,z_2)$ of the microscopic pattern as functions of $z$: for $z<\frac12$ they alternate between $z-1$ and $z+1$, for $z>\frac12$ they are both equal to $z$; the shaded interval is the range of macroscopic slopes where mesoscopic oscillations between the two patterns are energetically favourable.

\bigskip
In the $d$-dimensional case, the bulk scaling is $\e^d$, and the energies take the form
$$
F_\e(u)=\sum_{k\in \mathbb Z^d\setminus\{0\}}\sum_{i\in I^k_\e(\Omega)}\e^d\, f_k\Big(\frac{u_i-u_{i-k}}{\e}\Big),
$$
where $u\colon\Omega\cap\e\mathbb Z^d\to\mathbb R$ and $I^k_\e(\Omega)=\{i\in\mathbb Z^d:\e i,\e(i-k)\in\Omega\}$ (the interaction between $i$ and $i-k$ is counted both for $k$ and for $-k$, which can be accounted for by a redefinition of $f_k$). Under suitable growth conditions, analogous to those above, the continuum limit is of the form $\int_\Omega f_{\rm hom}(\nabla u)dx$, where $f_{\rm hom}$ is characterized by the $d$-dimensional analog of the asymptotic formula in the $1$-dimensional case,
$$
f_{\rm hom}(z)=\lim_{N\to+\infty}\frac1{N^d}\min\Big\{\sum_{k}\ \sum_{i,i-k\in Q_N}f_k(u_i-u_{i-k}):\ u_i=\langle z,i\rangle\ \hbox{for $i\in\partial Q_N$}\Big\},
$$
where $Q_N=\{0,\dots,N\}^d$ and $\partial Q_N$ is the set of its points with some coordinate equal to $0$ or $N$.


\subsection{Spin systems: surface scaling and the geometry of Wulff shapes}\label{ssec-spin}

All the functionals considered so far in this section are defined on lattice functions with values in $\mathbb R$ (or $\mathbb R^m$), for which the natural scaling of the energy is the {\bf bulk scaling} $\e\sum_i(\cdots)$, designed precisely so that nontrivial limits are integrals over $\Omega$. A genuinely different kind of discrete system is obtained by constraining the lattice function to take only two values, say $u_i\in\{-1,1\}$: these are {\bf spin systems}, the prototype being the Ising model of ferromagnetism. A $\{-1,1\}$-valued function carries no information beyond the location of the interface between the two phases, so the bulk scaling above is now the wrong one: with that scaling, configurations of finite energy would be forced to be eventually constant, and any genuine interface would cost an infinite amount as $\e\to0$: a jump between $-1$ and $1$ at neighbouring sites contributes at least $\e f(2/\e)\ge2^p\e^{1-p}\to+\infty$. The correct scaling is instead a {\bf surface scaling}.

\paragraph{Ising systems and their $\Gamma$-limit.} Let $\alpha_k\ge0$, $k\in\mathbb Z^d\setminus\{0\}$, be a {\em ferromagnetic} interaction (favouring equal neighbouring spins, $\alpha_{-k}=\alpha_k$) satisfying the decay condition $\sum_k\alpha_k\|k\|<+\infty$ and a coercivity condition, which holds for instance if $\alpha_{e_j}>0$ for $j=1,\dots,d$ (for the sharp condition, a connectedness property of the set of $k$ with $\alpha_k>0$, see \cite{BC24}), and for $u\colon\e\mathbb Z^d\to\{-1,1\}$ consider
$$
E_\e(u)=\sum_{i,j\in\mathbb Z^d}\e^{d-1}\,\alpha_{i-j}(u_i-u_j)^2.
$$
Identifying $u$ with the piecewise-constant function equal to $u_i$ on the cell $\e i+[0,\e)^d$ and writing $u=2\chi_A-1$, a family $u^\e$ of equibounded energy converges, up to subsequences, to (the characteristic function of) a set of finite perimeter $A$, and
$$
\Gamma\hbox{-}\lim_{\e\to0}E_\e(u^\e)=F(A)=\int_{\Omega\cap\partial^\ast A}\varphi(\nu)\,d\mathcal H^{d-1},\qquad \varphi(\nu)=4\sum_{k\in\mathbb Z^d}\alpha_k\,|\langle k,\nu\rangle|.
$$
The formula for $\varphi$ (and the correct power $d-1$ of $\e$) is due to Braides and Chambolle \cite{BC24} and has a simple heuristic justification (made rigorous in the Appendix): along every lattice line in direction $k$ crossing a macroscopic interface with normal $\nu$, a spin configuration with finite energy must flip at least once, at a cost $4\alpha_k$; the number of such lines crossing a unit $(d-1)$-area of interface is proportional to $|\langle k,\nu\rangle|$, and summing over $k$ gives $\varphi(\nu)$ (both $k$ and $-k$ are counted, since the sum is over ordered pairs). In short, {\bf a purely discrete, possibly long-range, two-valued system homogenizes into a genuinely anisotropic continuum perimeter}: a surface-tension counterpart to the bulk homogenization formulas of Chapter~\ref{hom-sec}.

\begin{figure}[h!]
	\centerline{\includegraphics [width=5.8in]{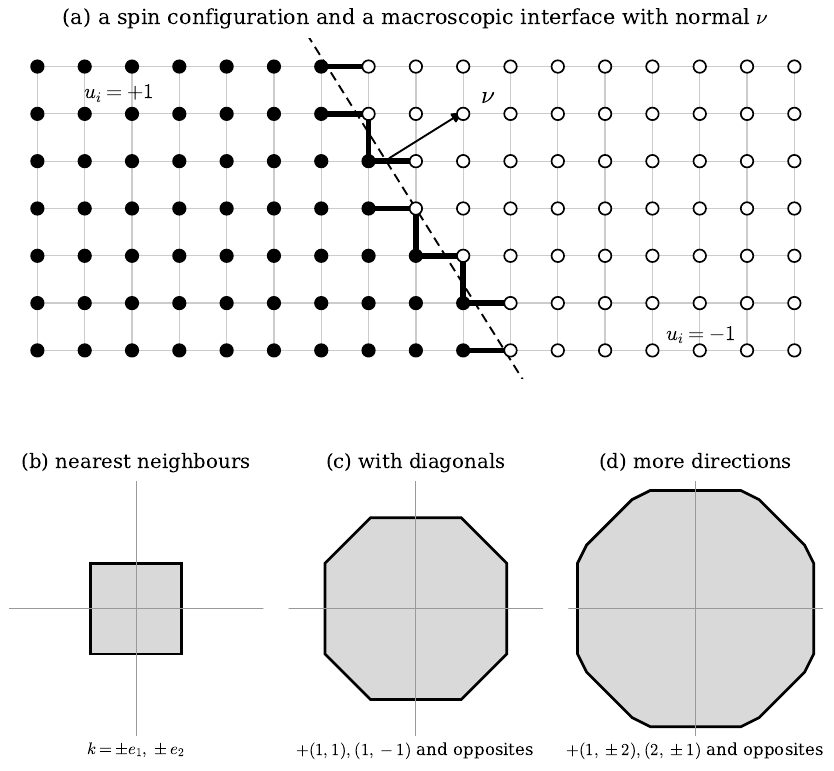}}
	\caption{Ising systems. (a) A spin configuration on $\e\mathbb Z^2$ ($+1$ full dots, $-1$ empty dots) and a macroscopic interface with normal $\nu$; the thick segments are the broken nearest-neighbour bonds, in the directions $e_1$ and $e_2$. Their number per unit length of interface is proportional to $|\nu_1|$ and $|\nu_2|$, and for a general interaction direction $k$ to $|\langle k,\nu\rangle|$, which gives $\varphi(\nu)=4\sum_k\alpha_k|\langle k,\nu\rangle|$. (b)--(d) The corresponding Wulff shapes, Minkowski sums of the segments $[-8\alpha_kk,8\alpha_kk]$ over the directions $k$ with $\alpha_k\neq0$ (one for each pair $\pm k$): the square for nearest-neighbour interactions, the octagon if interactions along the diagonals are added, and polygons with more sides if further rational directions are added; all sides are parallel to the interaction directions.}\label{Fig45}
\end{figure}

\paragraph{Wulff shapes.} Being even, convex and positively $1$-homogeneous, $\varphi$ is itself the support function of a convex body $W_\varphi\subset\mathbb R^d$, the {\bf Wulff shape}: the (unique, up to translation and dilation to fix the volume) minimizer of $F$ at fixed volume, i.e.\ the equilibrium shape of a macroscopic droplet of the ferromagnetic phase. For nearest-neighbour interactions alone ($\alpha_k=\alpha_1>0$ if $\|k\|=1$, $0$ otherwise) one finds $\varphi(\nu)=8\alpha_1(|\nu_1|+\cdots+|\nu_d|)=8\alpha_1\|\nu\|_1$ (the factor $8$ is $4$, coming from $(u_i-u_j)^2=4$, times $2$, since both $k=e_j$ and $k=-e_j$ contribute), whose Wulff shape is the axis-aligned cube $[-8\alpha_1,8\alpha_1]^d$: the classical fact that the zero-temperature nearest-neighbour Ising droplet is a cube, not a ball --- the lattice anisotropy survives, undiluted, in the continuum limit. Switching on further, weaker, interactions along more directions $k$ adds facets to $W_\varphi$ and lets it approximate a ball ever more closely; but, as we are about to see, no \emph{finite} or even countable combination of discrete directions can ever round off every corner exactly.

\paragraph{The zonotope/zonoid classification.} The formula for $\varphi$ has a transparent geometric meaning: writing $\nu_k=k/\|k\|$,
$$
\varphi(z)=4\sum_k\alpha_k\|k\|\,|\langle\nu_k,z\rangle|=\int_{S^{d-1}}|\langle\nu,z\rangle|\,d\mu(\nu),$$
with $\mu=4\sum_k\alpha_k\|k\|\,\delta_{\nu_k}$,
which is exactly the support function of the Minkowski sum of the segments $[-4\alpha_k\|k\|\,\nu_k,4\alpha_k\|k\|\,\nu_k]$. Following \cite{BC24}, a {\bf rational zonotope} is a Minkowski sum of finitely many segments in rational directions, that is, a convex body whose support function is $\int_{S^{d-1}}|\langle\nu,z\rangle|\,d\mu$ with $\mu=\sum_{i=1}^n\lambda_i(\delta_{\nu_i}+\delta_{-\nu_i})$, $\lambda_i>0$ and $\nu_i=w_i/\|w_i\|$ with $w_i\in\mathbb Z^d\setminus\{0\}$ ({\em rational directions}), while a {\bf rational zonoid} is the same with $\mu=\sum_{i=1}^\infty\lambda_i(\delta_{\nu_i}+\delta_{-\nu_i})$, with $\sum_i\lambda_i<+\infty$. Then the Wulff shapes of finite-range homogeneous Ising systems are the rational zonotopes, and those of general homogeneous systems are the rational zonoids (with a connectedness condition on the support of $\mu$ if coercivity is required). The generating measures are atomic and concentrated on rational directions; general {\bf zonoids}, with arbitrary generating measure on the sphere, are Hausdorff limits of rational zonotopes. Hence {\em every} zonoid is approximated arbitrarily well by Wulff shapes of finite-range Ising systems, but only the rational ones are exactly attained: for instance, in dimension $d\ge2$ the ball, whose generating measure is uniform on the sphere, is not the Wulff shape of any homogeneous Ising system. In dimension $2$ every centered convex symmetric set is a zonoid, so that homogeneous Ising systems can approximate {\em any} planar anisotropy; in dimension $d\ge3$, however, zonotopes are characterized among polytopes by the property that all their $2$-dimensional faces are centrally symmetric (in particular, all their faces are themselves zonotopes), and a polytope that is a zonoid is a zonotope. This excludes elementary shapes such as the regular octahedron, whose faces are triangles: {\em no} homogeneous Ising system, however many interaction ranges it uses, can produce an octahedral Wulff shape.

\begin{figure}[h!]
	\centerline{\includegraphics [width=5.8in]{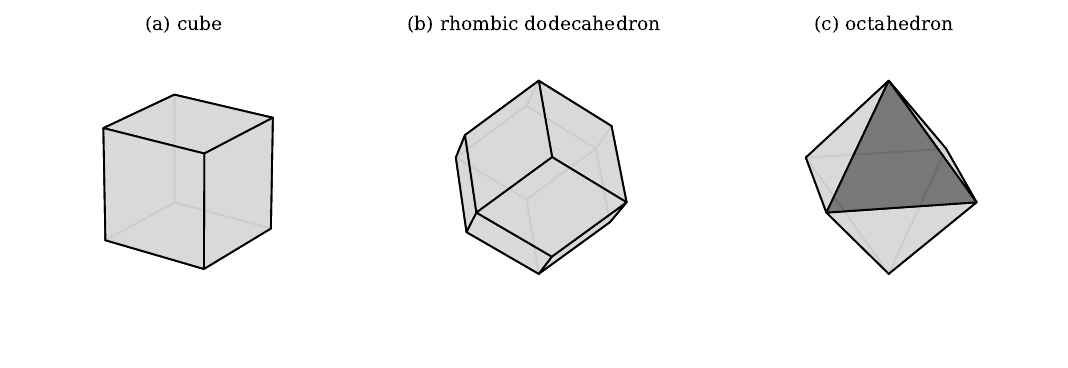}}
	\caption{Wulff shapes in dimension $3$. The cube (a zonotope generated by the three nearest-neighbour directions) and the rhombic dodecahedron (a zonotope generated by the four directions $(\pm1,\pm1,1)$, all faces being rhombi) can be Wulff shapes of homogeneous Ising systems. The regular octahedron cannot: its faces are triangles (one is shaded), which are not centrally symmetric, so that it is not a zonotope.}\label{Fig46}
\end{figure}

\paragraph{Recovering all shapes: periodicity.} This obstruction is removed exactly by the mechanism of Chapter~\ref{hom-sec}: replacing a {\em homogeneous} system $\alpha_{i-j}$ by a {\em periodic} one $a_{ij}$ of period $N$ and homogenizing, the Wulff shape of the $\Gamma$-limit is, for finite-range interactions, a polytope \cite{BC24}, which we call a {\bf zonotope of order $N$}; for $N=1$ these are the rational zonotopes, but for $N>1$ they are no longer constrained to have zonotopic faces. Letting $N\to\infty$, zonotopes of order $N$ approximate {\em all} bounded centered convex bodies: every macroscopic crystalline anisotropy, octahedron included, is the limit, in the Hausdorff sense, of Wulff shapes of periodic spin systems of sufficiently large period, even though it is unreachable from any homogeneous one.

\subsection{Appendix: relevant tools of Chapter \ref{disc-sec}}\label{App-Ch10}

\paragraph{References.} The $\Gamma$-limit of discrete energies with long-range interactions is treated by Braides and Gelli \cite{BGe02}, and the integral representation of continuum limits of discrete energies with superlinear growth by Alicandro and Cicalese \cite{AC04}; a textbook presentation of the convexification effect for lattice energies is in \cite{Bra02}. Discrete systems with values in $\{-1,1\}$ and their surface scaling are studied in \cite{ABCS,BC07}, and the identification of the Wulff shapes of Ising systems with zonoids, the notion of rational zonotope and zonoid, and the zonotopes of order $N$ for periodic systems are due to Braides and Chambolle \cite{BC24}. For the theory of zonotopes and zonoids (in particular the characterization of zonotopes by the central symmetry of their $2$-dimensional faces) see Bolker \cite{Bol69} and Schneider \cite[Section~3.5]{Sch14}.

\paragraph{The asymptotic formula for the convex envelope.}
\begin{proposition}\label{fss-prop}
Let $f\colon\mathbb R\to[0,+\infty)$ be continuous with $f(z)\ge|z|^p$, $p>1$, and let $m_N(z)=\min\{\sum_{i=1}^Nf(z_i):\sum_{i=1}^Nz_i=Nz\}$. Then
$$
f^{**}(z)=\lim_{N\to+\infty}\frac{m_N(z)}N=\inf_N\frac{m_N(z)}N.
$$
\end{proposition}

\begin{proof}
By Jensen's inequality, $\sum_if(z_i)\ge\sum_if^{**}(z_i)\ge Nf^{**}(z)$ when $\sum_iz_i=Nz$, so $m_N(z)\ge Nf^{**}(z)$. Moreover $m_{N+M}(z)\le m_N(z)+m_M(z)$ (concatenate two competitors), so that by Fekete's lemma $m_N(z)/N$ converges to $\inf_Nm_N(z)/N\ge f^{**}(z)$. Conversely, since $f$ is continuous and coercive, $f^{**}(z)=\lambda f(w_1)+(1-\lambda)f(w_2)$ with $z=\lambda w_1+(1-\lambda)w_2$ and $\lambda\in[0,1]$ (Carathéodory's theorem in dimension $1$). Given $N$, let $n=\lfloor\lambda N\rfloor$ and take $z_i=w_1$ for $i\le n$, $z_i=w_2$ for $n<i<N$ and $z_N=Nz-nw_1-(N-n-1)w_2=(\lambda N-n)(w_1-w_2)+w_2$, so that $|z_N|\le C=|w_1-w_2|+|w_2|$. Then $m_N(z)\le nf(w_1)+(N-n-1)f(w_2)+\max_{|w|\le C}f\le Nf^{**}(z)+C'$, with $C'$ independent of $N$, and the conclusion follows by dividing by $N$.
\end{proof}

\noindent If $z$ is such that the weights $\lambda,1-\lambda$ above are irrational (as for the double-well potential $f(z)=(z^2-1)^2$ and $z$ irrational in $(-1,1)$, where $w_{1,2}=\pm1$), the minimum $m_N(z)$ is strictly larger than $Nf^{**}(z)$ for every $N$, since the sets of values taken by the minimizers cannot be reduced to $\{w_1,w_2\}$: this is why the limit cannot be replaced by a minimum on a bounded number of variables.

\paragraph{The example with $M=2$.}
\begin{proposition}\label{feff-prop}
Let $f_1(z)=(|z|-1)^2$, $f_2(z)=|z|^2$ and $f_{\rm eff}(z)=f_2(2z)+\frac12\min\{f_1(a)+f_1(b):a+b=2z\}$. Then
$$
f_{\rm eff}(z)=\begin{cases}5z^2&\hbox{if }|z|\le\frac12\\ 5z^2-2|z|+1&\hbox{if }|z|\ge\frac12,\end{cases}
\qquad
f_{\rm eff}^{**}(z)=\begin{cases}5z^2&\hbox{if }|z|\le\frac25\\ 4|z|-\frac45&\hbox{if }\frac25\le|z|\le\frac35\\ 5z^2-2|z|+1&\hbox{if }|z|\ge\frac35,\end{cases}
$$
and, with $f_0(z)=f_1(z)+f_2(2z)$,
$$
f_1^{**}(z)+f_2^{**}(2z)\le f_{\rm eff}^{**}(z)\le f_0^{**}(z)\quad\hbox{for all }z,
$$
where $f_1^{**}(z)=((|z|-1)^+)^2$, $f_2^{**}=f_2$ and $f_0^{**}(z)=\frac45$ if $|z|\le\frac15$, $f_0^{**}(z)=f_0(z)=5z^2-2|z|+1$ otherwise; the inequalities are strict, for instance, at $z=\frac12$, where the three terms are $1$, $\frac65$ and $\frac54$.
\end{proposition}

\begin{proof}
All functions are even, so let $z\ge0$ and set $h(a)=f_1(a)+f_1(2z-a)$. If $a,b=2z-a$ are both non-negative, then $h(a)\ge2f_1(z)=2(1-z)^2$ by the convexity of $f_1$ on $[0,+\infty)$, with equality for $a=b=z$. Both negative is impossible, as $a+b=2z\ge0$. If $a=-c<0\le b=2z+c$, then $h=(c-1)^2+(2z+c-1)^2$, which for $z\le1$ is minimized on $c\ge0$ at $c=1-z$ with value $2z^2$ (while for $z>1$ the minimum on $c\ge0$ is at $c=0$ and falls in the previous case). Hence $\min h=2\min\{z,1-z\}^2$ for $z\le1$ and $\min h=2(z-1)^2$ for $z\ge1$, and the formula for $f_{\rm eff}$ follows from $f_{\rm eff}=4z^2+\frac12\min h$. The affine function $\ell(z)=4z-\frac45$ satisfies $f_{\rm eff}-\ell=5(z-\frac25)^2\ge0$ on $[0,\frac12]$ and $f_{\rm eff}-\ell=5(z-\frac35)^2\ge0$ on $[\frac12,+\infty)$, so that $\ell\le f_{\rm eff}$ with equality exactly at $\frac25$ and $\frac35$. Hence the function $g$ defined by the three-line formula is convex (it is of class $C^1$ with non-decreasing derivative: $10z$ up to $\frac25$, then $4$, then $10z-2\ge4$ for $z\ge\frac35$, and even), $g\le f_{\rm eff}$ and $g=f_{\rm eff}$ outside $(\frac25,\frac35)$ and at the endpoints; therefore $g=f_{\rm eff}^{**}$: indeed $g\le f_{\rm eff}^{**}$ since $g$ is a convex function below $f_{\rm eff}$, while $f_{\rm eff}^{**}\le f_{\rm eff}=g$ outside $(\frac25,\frac35)$, and in $(\frac25,\frac35)$ we have $f_{\rm eff}^{**}\le\ell=g$, since $\ell$ is the chord between the values of $f_{\rm eff}$ at $\frac25$ and $\frac35$. For the inequalities, $f_{\rm eff}\ge f_2(2z)+\frac12\cdot2f_1^{**}(z)=f_2(2z)+f_1^{**}(z)$ (since $\frac12\min\{f_1(a)+f_1(b):a+b=2z\}\ge f_1^{**}(z)$ by convexity of $f_1^{**}$) and $f_{\rm eff}\le f_0$ (take $a=b=z$); taking convex envelopes gives the two inequalities, and the expressions of $f_1^{**},f_0^{**}$ are computed in the same way as that of $f_{\rm eff}^{**}$ ($f_0(z)=5z^2-2|z|+1$ has a concave corner at $0$, and its tangent line from the minimum point $|z|=\frac15$ is horizontal).
\end{proof}

\paragraph{The $\Gamma$-limit of Ising systems.}
\begin{proposition}\label{ising-prop}
Let $E_\e$ and $\varphi$ be as in Section~\ref{ssec-spin}, $\alpha_k\ge0$ with $\sum_k\alpha_k\|k\|<+\infty$, and let $u^\e\colon\e\mathbb Z^d\to\{-1,1\}$ converge in $L^1_{\rm loc}(\mathbb R^d)$ to $u=2\chi_A-1$. Then $\liminf_{\e\to0}E_\e(u^\e)\ge\int_{\partial^*A}\varphi(\nu)\,d\mathcal H^{d-1}$. If $A$ has Lipschitz boundary and $u^\e=u|_{\e\mathbb Z^d}$, then $\lim_{\e\to0}E_\e(u^\e)=\int_{\partial A}\varphi(\nu)\,d\mathcal H^{d-1}$.
\end{proposition}

\begin{proof}
Write $E_\e=\sum_k\alpha_kG^k_\e$, with $G^k_\e(u)=\sum_i\e^{d-1}(u_i-u_{i-k})^2$. Since $|u_i-u_{i-k}|\in\{0,2\}$ we have $(u_i-u_{i-k})^2=2|u_i-u_{i-k}|$, so that, for the piecewise-constant interpolation $u_\e$,
\begin{eqnarray*}
G^k_\e(u^\e)&=&\frac2\e\int|u_\e(x)-u_\e(x-\e k)|\,dx\ \ge\ 2\int\psi(x)\,\frac{u_\e(x)-u_\e(x-\e k)}\e\,dx\\
&=&2\int u_\e(x)\,\frac{\psi(x)-\psi(x+\e k)}\e\,dx
\end{eqnarray*}
for every $\psi\in C^1_c(\mathbb R^d)$ with $|\psi|\le1$. As $\e\to0$ the right-hand side tends to $-2\int u\,\langle\nabla\psi,k\rangle dx$, and taking the supremum in $\psi$ we get
$$
\liminf_\e G^k_\e(u^\e)\ge2|\langle Du,k\rangle|(\mathbb R^d)=4\int_{\partial^*A}|\langle k,\nu\rangle|\,d\mathcal H^{d-1}.
$$ Since all terms are non-negative, the $\liminf$ of the sum is not less than the sum of the $\liminf$, and the lower bound follows. For the upper bound, $\e^{d-1}\#\{i:u_i\neq u_{i-k}\}\to\int_{\partial A}|\langle k,\nu\rangle|\,d\mathcal H^{d-1}$ for a set with Lipschitz boundary (count the lattice segments of direction $k$ crossing $\partial A$) and this quantity is bounded by $C\|k\|$ uniformly in $\e$; by the summability of $\alpha_k\|k\|$ and the dominated convergence theorem we may exchange the limit and the sum over $k$.
\end{proof}

\noindent The compactness of sequences with bounded energy, which requires the coercivity (connectedness) of the system, and the optimality of the lower bound for general sets of finite perimeter are proved in \cite{BC24}.

\paragraph{The faces of a zonotope.}
Let $Z=\sum_{i=1}^n[-s_i,s_i]$ be a zonotope, with $s_i\in\mathbb R^d\setminus\{0\}$. Its support function is $h_Z(\nu)=\sum_i|\langle s_i,\nu\rangle|$, and the face of $Z$ with outer normal $\nu$ is the Minkowski sum of the faces of the segments $[-s_i,s_i]$ with the same normal: this is the whole segment if $\langle s_i,\nu\rangle=0$, and one of its endpoints otherwise. Hence every face of $Z$ is, up to translation, the sum of the segments $[-s_i,s_i]$ parallel to it, that is, it is itself a zonotope; in particular every $2$-dimensional face is a centrally symmetric polygon. The regular octahedron has triangular faces and is therefore not a zonotope; by the characterization of zonotopes among polytopes through the central symmetry of $2$-faces \cite{Sch14}, and the fact that a polytope which is a zonoid is a zonotope, it is not a zonoid either, and then (Section~\ref{ssec-spin}) not the Wulff shape of any homogeneous Ising system.

\goodbreak

\vfill\eject

%
%

%
%
%
%

\printindex

\end{document}